\documentclass{article}
\usepackage[utf8]{inputenc}
\usepackage{aaPreamble}

\title{Explicit Birational Geometry of Fano threefold complete intersections}
\author{Tiago Duarte Guerreiro}
\date{April 2022}

\begin{document}

\maketitle

\begin{abstract}
We complete the analysis on the birational rigidity of quasismooth Fano 3-fold deformation families appearing in the Graded Ring Database as a complete intersection. When such a deformation family $X$ has Fano index at least 2 and is minimally embedded in a weighted projective space in codimension 2, we determine which cyclic quotient singularity is a maximal centre. If a cyclic quotient singularity is a maximal centre, we construct a Sarkisov link to a non-isomorphic Mori fibre space or a birational involution. This allows, in particular, the construction of new examples of Fano 3-folds of codimension 6 which are realised as complete intersections in fake weighted projective spaces. We define linear cyclic quotient singularities on $X$ and prove that these are maximal centres by explicitly computing Sarkisov links centred at them. It turns out that each $X$ has a linear cyclic quotient singularity leading to a new birational model. As a consequence, we show that if $X$ is birationally rigid then its Fano index is 1. If the new birational model is a strict Mori fibre space, we determine its fibration type explicitly. In this case, a general member of $X$ is birational to a del Pezzo fibration of degrees 1, 2 or 3 or to a conic bundle $Y/S$ where $S$ is a weighted projective plane with at most $A_2$ singularities.  
\end{abstract}

\tableofcontents

\section{Introduction}  \label{sect:intro}

The \emph{Minimal Model Program} aims at finding \emph{good} representatives of a given birational equivalence class. 
It is a largely successful attempt to extend the birational classification of algebraic surfaces carried out by the Italian school of geometers, notably with the efforts of Castelnuovo, Enriques and Severi at the beginning of the last century, to higher dimensions.

In today's language, we can consider smooth projective surfaces divided into two large classes. Those that admit local pluri-canonical sections and those that do not. In the first case, there is a unique "good" representative called a \emph{minimal model} and is characterised by the fact that it is the only surface in its birational equivalence class for which its canonical divisor is nef. In the second case, the surfaces are \emph{uniruled}, that is, covered by rational curves, and there is no unique good representative. 

With the appearance of the Minimal Model Program new objects of foundational importance come to light. For instance, if $W$ is a uniruled variety, then there is a sequence of birational transformations from $W$ that result in a \emph{Mori fibre space} $Y/S$ \cite[Corollary~1.3.3]{BCHM}. That is, $W$ admits a birational model $Y$ with mild singularities and a fibre structure $Y \rightarrow S$ where the relative Picard rank $\rho_{Y/S}$ is 1 and its fibres have ample anticanonical divisor. Fano varieties of Picard rank 1  are important examples of Mori fibre spaces.


There is no guarantee, however, that, given $W$ as above, the Minimal Model Program outputs a unique birational model and so it becomes natural to study the relations between these. The uniqueness of the birational model is encoded in the notion of birational rigidity. Loosely, a Fano variety $X$ of Picard rank 1 is \emph{birationally rigid} if it possesses no birational maps $\sigma \colon X \rat Y/S$ to different Mori fibre spaces. A notable and very recent result in this direction that took over one hundred years to establish and many researchers is that any smooth hypersurface $X \subset \mathbb{P}^{n+1}$ of dimension at least $3$ and degree $n+1$ is birationally (super)-rigid, a stronger form of birational rigidity. See \cite{BGH} for a complete proof and historical account of this problem.

Understanding the birational geometry of a Mori fibre space is understanding the birational maps between itself and its models and this is exactly the content of the Sarkisov Program. The fundamental result here says that a birational map between Mori fibre spaces factors into a finite sequence of \emph{Sarkisov links} \cite{cortiSP, genSP}. 

We are specially interested in the case when $X$ is a mildly singular 3-dimensional Fano complete intersection. If $\iota_X$ is the maximum positive integer for which $-K_X=\iota_XA$, where $A \in \Cl(X)$ is an ample generator, we consider $X$ to be polarised by $A$, that is, we consider $X$ together with an embedding into a weighted projective space via the ring,
\[
R(X,A):=\bigoplus_{n\geq 0}H^0(X,\mathcal{O}_X(nA)).
\] 
In the case of a Fano 3-fold, it is known that $R(X,A)$ is a finitely generated $\mathbb{C}$-algebra. This approach produces lists of Hilbert series which are collectively known as the Graded Ring Database. See \cite{GRDB} and \cite{reidgrdb}. Kawamata established in \cite{Kawamatabound} the boundedness of Fano 3-folds implying, in particular, that there is a finite number of such deformation families. Notice that a much more general statement on the boundedness of $d$-dimensional Fano varieties as been recently achieved in the seminal work of Birkar, see \cite{BirkarboundI, BirkarboundII}. 

Given the finiteness of these deformation families, one can hope to understand their birational geometry. This is a very challenging and active area of research even today. 

The following result is the starting point of our paper. It is the build up of work by many authors such as Iskovskikh and Manin \cite{iskmanin} - where it was proved that the groups of birational and biregular automorphisms of a smooth quartic threefold coincide - Reid, Pukhlikov, Corti \cite{CPR} and others. We state it in its final form:  
\begin{Thm}[{\cite{okadaI,chel,hamidhyp}}] \label{thm:start}
Let $X$ be a quasismooth member of a deformation Fano 3-fold appearing as a complete intersection in the Graded Ring Database.
\begin{enumerate}
	\item If $X$ is a hypersurface, then it is birationally rigid if and only if it has Fano index 1.
	\item If $X$ has codimension 2, then there are exactly 19 deformation families of Fano index 1 which are birationally rigid. 
	\item If $X$ has codimension 3, then it is the intersection of three quadrics in $\mathbb{P}^6$ and it is not birationally rigid.
	\end{enumerate}
\end{Thm}




\section*{Main results}  \label{sect:mainthm}

Throughout this paper, we consider the case when $X$ is a 3-fold with Fano index $\iota_X \geq 2$ and is embedded as a codimension 2 in a weighted projective space. In particular, we prove that in the conditions of Theorem \ref{thm:start}, if $X$ has Fano index at least 2, then it is not birationally rigid. We construct birational maps between $X$ and other Mori fibre spaces in the framework of the Sarkisov Program and present these constructions explicitly. As an immediate corollary we obtain that, in this context, birational rigidity is an exclusive property of Fano 3-folds with Fano index 1. 

Let $X \hookrightarrow \mathbb{P}(a_0,\ldots,a_5)$ be a Fano 3-fold deformation family embedded with codimension 2 and Fano index at least 2. There are 40 such deformation families in the Graded Ring Database which we label by $I=\{86, \ldots,125\}$, see \cite{brownsuI,brownsuII}. We sometimes abuse notation and write $X\in I$ to mean that $X$ is a member of a deformation family indexed by $I$. Moreover, we always assume $X$ to be quasismooth, unless stated otherwise. See \cite[Definition~3.1.5]{DolgachevWPS} for a precise definition. Excluded from the analysis is Family 86: the smooth complete intersection of two quadrics in $\mathbb{P}^5$ which is known to be rational. See for instance \cite{2quadr} for a modern treatment and generalisations. Let
\begin{align*}
I_{Cb}&=\{87,112,113,118,119\}\\
I_{dP}&=\{88,89,90,91,103,114,116,120,121,122,123,124,125\}
\end{align*}
and
\begin{align*}
I_{nS}&=I_{Cb} \cup I_{dP}\\
I_{S}&=I \setminus I_{nS}.
\end{align*}

\begin{Def}[{\cite[Definition~1.4]{pfaff}}]
 A Fano variety $X$ is \textbf{birationally solid} or \textbf{solid} if there is no birational map to a Mori fibre space $Y/S$ where $\dim S > 0$.
\end{Def}
 The main theorem of this paper is the following:

\begin{Thm} \label{thm:main}
Suppose $X$ is a quasismooth Fano 3-fold deformation family index by $I$. 
\begin{enumerate}
	\item If $X \in I_{nS}$ then $X$ is not solid. In particular,
	\begin{enumerate}
		\item if $X \in I_{Cb}$ then there is a Sarkisov link from $X$ to a conic bundle $Y/S$, where $S$ is smooth or has $A_1$ or $A_2$ singularities.
		\item if $X \in I_{dP}$ then there is a Sarkisov link from $X$ to a del Pezzo fibration whose generic fibre has degree 1, 2 or 3.
	\end{enumerate}
		\item If $X \in I_{S}$ then there is a Sarkisov link from $X$  to a terminal non-quasismooth weighted hypersurface or codimension 2 Fano 3-fold.
\end{enumerate}
\end{Thm}


 Notice that the only generality assumption for $X$ in Theorem \ref{thm:main} is quasismoothness. The following are immediate consequences of Theorem \ref{thm:main}.
\begin{Cor}
Let $X \in I$. Then $X$ is birationally non-rigid. 
\end{Cor}

\begin{Cor}
There are exactly 19 quasismooth Fano 3-fold weighted complete intersections of codimension 2 which are birationally rigid and they all have Fano index 1. 
\end{Cor}

Putting together Theorem \ref{thm:main} with \cite{chel,okadaI,hamidhyp} and \cite[Example~1.4.4]{3quad} we have a complete picture on the birational rigidity of Fano 3-folds which appear as a complete intersection in the Graded Ring Database.

\begin{Thm}
Let $X$ be a quasismooth member of a deformation Fano 3-fold which appears as a complete intersection in the Graded Ring Database. Then $X$ is birationally rigid if and only if $X$ is
\begin{enumerate}
	\item A hypersurface of Fano index 1 or
	\item One of 19 deformation families of Codimension 2 and Fano index 1. 
\end{enumerate}
In particular, if $X$ is birationally rigid it has Fano index 1.
\end{Thm}

A Sarkisov link $\sigma \colon X \rat Y/S$ to a Mori fibre space $Y/S$ is initiated by a divisorial extraction from a maximal centre $\Gamma \subset X$, see \cite{cortiSP}. We also have the following classification of maximal centres:

\begin{Thm} \label{thm:mainEx}
Let $X$ be a general member of a deformation family in $I_{S}$ and $\Gamma \subset X$ a cycle. 
\begin{itemize}
	\item If $\Gamma$ is an irreducible curve, then it is not a maximal centre of $X \in I_{S} \setminus \{92\}$.
	\item If $\Gamma$ is a smooth point, then it is not a maximal centre of $X \in I_{S} \setminus \{92, 94, 115\}$ provided the generality assumptions on table \ref{tab:gen} hold.
	\item If $\Gamma$ is a singular point, then it is not a maximal centre if and only if the last column of table \ref{tab:big} is empty. 
	\end{itemize}
\end{Thm}

Finally, we make the following prediction supported by Theorem \ref{thm:mainEx} which will be subject to investigation in the future:

\begin{Conj} \label{conj:solid}
Let $X$ be a general member in $I_S$. Then $X$ is solid.
\end{Conj}
A particular consequence of the above conjecture is that $X$ is not rational.

\paragraph{Structure of the paper.} This paper is roughly divided in two parts: First, the construction of Sarkisov links spanning Chapters 4, 5 and 6. Secondly, the exclusion of subvarieties as maximal centres in Chapter 7. It starts with Chapters 1, 2 and 3 which are introductory and this is where notation and most of the theorems used throughout are stated and proved. It ends with the Big Table in chapter 8 which summarises the results obtained.

In Chapter 3 we introduce the objects we are interested in and the main techniques. In particular, we explain how to use the natural polarisation of $X$ to obtain a Sarkisov link from the blow up of a terminal cyclic quotient singularity in $X$ in Subsection \ref{sub:toricemb}. This is done by understanding how sections lift under the Kawamata blowup which is explained in Lemma \ref{lem:utbl} and Corollary  
\ref{cor:genlift}. We then explain how to continue the Sarkisov link under certain useful assumptions in Lemmas \ref{lem:anticanonicalwalliso} and \ref{lem:discr}. We define \emph{linear} cyclic quotient singularities w.r.t. $X$ in Definition \ref{def:linear}. These are relevant since there is always a Sarkisov link centred on a linear cyclic quotient singularity. 

In Chapter 4 we construct Sarkisov links to Mori fibrations with positive dimensional base, that is, either to conic bundles or del Pezzo fibrations. See Theorem \ref{thm:nS}. In Chapters 5 and 6 we construct Sarkisov links to other Fano 3-folds. The results are achieved by first playing a 2-ray game on a rank two toric variety and then restricting it to our 3-fold, as explained in Subsection \ref{sub:toricemb}. In Chapter 6, in particular, we construct new examples of Fano 3-folds of codimension 6 which are realised as complete intersections embedded in a fake weighted projective space. See, for instance, Example \ref{ex:cod6}.

In Chapter 7, certain subvarieties are excluded as maximal centres using different methods. See Theorem \ref{thm:mainEx}. This paves the way to Conjecture \ref{conj:solid}.

The paper ends with a summary of all the results obtained in The Big Table.

\paragraph{Acknowledgements} The author would like to thank Hamid Abban for the continuous support and encouragement. We would also like to thank Gavin Brown, Livia Campo, Ivan Cheltsov, Takuzo Okada, Erik Paemurru and Alan Thompson  for many conversations and interest. This work has been partially supported by the EPSRC grant EP/V048619/1 and EPSRC grant EP/V055399/1.





\section{Notation}  \label{sect:notation}

We denote by $\mathbb{P}$ the $5$-dimensional weighted projective space $\mathbb{P}(a_0,\ldots,a_5)=\Proj \mathbb{C}[x_0,\ldots,x_5]$ with unordered weights.
Let $\Phi \colon T \rightarrow \mathbb{P}$ be the toric $(b_0, \ldots, b_4)$-weighted blowup of the coordinate point $\mathbf{p}_{\xi}$. If this is the $i$-th coordinate point, then we substitute $a_i$ by $a_{\xi}$ and $x_i$ by $\xi$. We specify $T$ by its weight matrix,
\[
\begin{array}{cccc|cccccc}
             &       & u     & \xi & x_0     &   & x_4 & \\
\actL{T}     &  \lBr &  0    & a_{\xi} & a_0     & \ldots           &  a_4 & \actR{}\\
             &       & -a_{\xi}  & 0   & b_0 & \ldots & b_4 &  
\end{array}
\]
where the vertical bar represents the irrelevant ideal of $T$ which in this case is $(u,\xi) \cap (x_0,\ldots,\widehat{\xi},\ldots,x_4)$ and each row is a representation of a $\mathbb{C}^*$-action on $\mathbb{C}^{7}$. In other words, the variety $T$ is defined by the geometric quotient
\[
T=\frac{\mathbb{C}^{7}\setminus \mathcal{Z}((u,\xi)\cap(x_0,\ldots,\widehat{\xi},\ldots,x_4))}{\mathbb{C}^*\times \mathbb{C}^*}.
\]
We always consider a well-formed isomorphic model of $T$ (See \cite[Definition~2.11]{hamidplia} . The Cox ring of $T$ is $\Cox(T)=\mathbb{C}[u,\xi,x_0,x_1,x_2,x_3,x_4]$ with the grading of each variable in the corresponding column of the matrix above. We denote the exceptional divisor of $\Phi$ by $E$. This is the effective divisor given by $(u=0)$ and is isomorphic to $\mathbb{P}(b_0,\ldots,b_4)$. Let $A$ be the $\mathbb{Q}$-Cartier divisor $\frac{1}{a_{\xi}} \cdot (\xi=0)$ in $\mathbb{P}$. Then, $\Pic(T)_{\mathbb{Q}}$ is generated by (the strict transform of) $A$ and $E$.

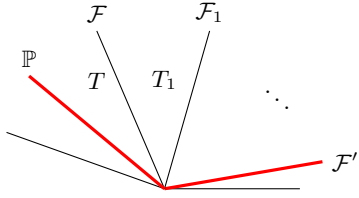
\begin{figure}%
\centering
\begin{tikzpicture}[scale=3,font=\small]
  \coordinate (A) at (0, 0);
  \coordinate [label={left:}] (E) at (-0.7, 0.25);
  \coordinate [label={above:$\mathbb{P}$}] (K) at (-0.6,0.5);
	\coordinate [label={above:$\mathcal{F}$}] (5) at (-0.3, 0.7);
	\coordinate [label={above:$\mathcal{F}_1$}] (2) at (0.2,0.7);
	\coordinate [label={right:$\mathcal{F}'$}] (6) at (0.7,0.12);
	\coordinate [label={right:}] (4) at (0.6,0);
	
		\coordinate [label={$T$}] (T) at (-0.3,0.4);
		\coordinate [label={$T_1$}] (T1) at (0,0.4);
		\coordinate [label={$\ddots$}] (T2) at (0.5,0.3);
  \draw  (A) -- (E);
  \draw  (A) -- (4);
	\draw [very thick,color=red] (A) -- (6);
	\draw (A) -- (5);
	\draw (A) -- (2);
	\draw [very thick,color=red](A) -- (K);
\end{tikzpicture}
\caption{A representation of the chamber decomposition of the cone of effective divisors of $T$.}%
\label{fig:genNot}%
\end{figure}

We make repeated use of several cones of $T$ as depicted in Figure \ref{fig:genNot}. Let $N^1(T)$ be the finite dimensional vector space over $\mathbb{R}$ of $\mathbb{Q}$-divisors modulo numerical equivalence.  The effective cone of $T$, $\Eff(T) \subset N^1(T)$, is generated by the outermost rays in Figure \ref{fig:genNot} and it has a decomposition into smaller chambers. 
Recall that the cone of Movable divisors of $T$, $\Mov(T) \subset \Eff(T)$, that we depict in red, is the cone of divisors $D$ for which the base locus of $|D|$ has codimension at least 2 in $T$. There is a finite number of small birational maps $f_j \colon T \dashrightarrow T_j$ for which $\Mov(T)=\bigcup_{j}f_j^*(\Nef(T_j))$ .A small $\mathbb{Q}$-factorial modification $\tau_i : T_i \rightarrow T_{i+1} $ is described by a $\mathbb{C}^*$-action and follows the notation in \cite{brownflips}. A ray in $\Eff(T)$ given by the vanishing of the section $s$ is written as $\mathbb{R}_+[s]$.

Let $X$ be a $\mathbb{Q}$-factorial terminal Fano 3-fold deformation family of Fano index $\iota_X$ at least two embedded in $\mathbb{P}$ with codimension two. We write it as 
\[
X \colon (f=g=0) \subset \mathbb{P}
\]
where $f$ and $g$ are quasi-homogeneous equations of degrees $d_1$ and $d_2$, respectively. For a polynomial $h=h(x_0,\ldots,x_n)$ and a monomial $x_0^{c_0} \cdots x_n^{c_n}$ of degree $\deg h$ we write $x_0^{c_0} \cdots x_n^{c_n} \in h$ (respectively $x_0^{c_0} \cdots x_n^{c_n} \not \in h$) if the coefficient of the monomial is non-zero in $h$ (respectively zero). For a polynomial $h$, if $D$ is a $\mathbb{Q}$-Cartier divisor in the algebraic variety $W$ we write $h \in |D|$ to say that $(h=0)$ is a member of the linear system given by $D$.

 The weights $a_i$ of the weighted projective space inside which $X$ sits have very particular features. In fact, if $\iota_X\not =4$, one of the weights is equal to $\iota_X$, and exactly another one a multiple of $\iota_X$. We always set $x_0$ and $x_1$ to be the variables with such weights, respectively.

We call $Y$ the proper transform of $X$ via $\Phi$ and $\varphi$ the restriction $\Phi\vert_Y$. We write it as
\[
Y \colon (\widetilde{f}=\widetilde{g}=0) \subset T.
\]
\section{Preliminaries}  \label{sect:prelim}

\subsection{Terminal Singularities} \label{sect:termsing}

It turns out that in order to generalise the results on the birational classification of surfaces to higher dimensions, one needs to allow some mild singularities. 

\begin{Def}
We say that a normal projective variety $X$ is \textbf{$\mathbb{Q}$-factorial} if for any Weil divisor $D$ on $X$ there is a positive integer $m$ for which $mD$ is Cartier.
\end{Def}

In particular, we can take the pull-back of $D$ with respect to any morphism $W \rightarrow X$ and define the intersection number $D \cdot C$ with any curve $C \subset X$. Next we define the smallest class of singularities needed to run the Minimal Model Program for non-singular varieties.

\begin{Def}[Terminal Singularities]
Suppose $X$ is a normal projective variety. Then, $X$ has \textbf{terminal singularities} if 
\begin{itemize}
	\item $X$ is $\mathbb{Q}$-Gorenstein, that is, there is a positive integer $m$ for which the Weil divisor $mK_X$ is Cartier;
	\item if $f\colon Y \rightarrow X$ is a resolution of $X$ and $\{ E_i\}$ the collection of all exceptional prime divisors of $f$ then 
	\[
	K_Y = f^*(K_X) + \sum a_iE_i
	\]
	with strictly positive discrepancies $a_i >0$. 
\end{itemize}
We say that $X$ is in the \textbf{Mori Category} if it is a normal projective variety with $\mathbb{Q}$-factorial and terminal singularities. We say $X$ has \textbf{canonical singularities} if $a_i \geq 0$. 
\end{Def}
The complete list of terminal 3-fold singularities is known and is due to the works of Reid \cite{ReidI}, Danilov \cite{DanilovI}, Morrison and Stevens \cite{MorrisonStevens}, Mori \cite{MoriTerminal} and others. In dimension 2 isolated terminal singularities are smooth points.  See \cite{reidyoung} for a friendly introduction to terminal and, more generally, canonical singularities.

\begin{Def}
Let $F \in \mathbb{C}\{x_1,x_2,x_3,x_4 \}$ be a convergent power series around $0$. Then $F=0$ is a \textbf{compound du Val singularity} (or cDV) if $F$ is of the form
\[
h(x,y,z)+tg(x,y,z,t)=0
\] 
where $h(x,y,z)=0$ defines a canonical surface singularity (also known as du Val singularity). 
\end{Def}



Recall that $\bm{\mu}_r$ denotes the cyclic group of $r$th roots of unity. Define the action of $\bm{\mu}_r$ on $\mathbb{C}^4$ with coordinates $x_1,\,x_2,\,x_3,\,x_4$ by
\begin{align*}
\bm{\mu}_r \times \mathbb{C}^4 &\longrightarrow \mathbb{C}^4 \\
(\epsilon, (x_1,x_2,x_3,x_4)) &\longmapsto (\epsilon^{\alpha_1}x_1,\epsilon^{\alpha_2}x_2,\epsilon^{\alpha_3}x_3,\epsilon^{\alpha_4}x_4)
\end{align*} 
where $\epsilon$ is a primitive $r$ root of unity and $\alpha_i$ are integers. We denote such an action by
\[
\frac{1}{r}(\alpha_1,\alpha_2,\alpha_3,\alpha_4)
\] 
and its orbit space by
\[
\mathbb{C}^4/\bm{\mu}_r(\alpha_1,\alpha_2,\alpha_3,\alpha_4).
\]
Let $F \in \mathbb{C}\{x_1,x_2,x_3,x_4 \}$ be a convergent power series around $0$ and suppose that $F$ is equivariant with respect to this action. Then $\bm{\mu}_r$ acts on the germ of the hypersurface $(F(x_1,x_2,x_3,x_4)=0) \subset \mathbb{C}^4$ and we can take the quotient
\[
(F(x_1,x_2,x_3,x_4)=0)/\bm{\mu}_r(\alpha_1,\alpha_2,\alpha_3,\alpha_4).
\]

The main theorem of \cite{reidterminal} is

\begin{Thm}[{\cite[Main Theorem~I]{reidterminal}}]
Every terminal 3-fold singularity over $\mathbb{C}$ is analytically isomorphic to
\[
(F(x_1,x_2,x_3,x_4)=0)/\bm{\mu}_r(\alpha_1,\alpha_2,\alpha_3,\alpha_4)
\]
where $F$ defines a cDV singularity.
\end{Thm}

It turns out that only very few actions produce terminal singularities. A complete list can be found in \cite[Theorem~5.43]{kollarmori}. Of special importance for us is the case of \emph{cyclic quotient singularities}.

\begin{Def}
Let $F \in \mathbb{C}\{x_1,x_2,x_3,x_4 \}$ be a convergent power series around $0$ and suppose $F=0$ is the germ of a smooth point. Then 
\[
(F(x_1,x_2,x_3,x_4)=0)\big/\bm{\mu}_r(\alpha_1,\alpha_2,\alpha_3,\alpha_4)
\]
is called a \textbf{cyclic quotient singularity} and we write it as
\[
\frac{1}{r}(a_1,a_2,a_3).
\]
\end{Def}

The following is adapted from \cite[Section~5]{reidyoung}.
\begin{Thm} \label{thm:cqs}
Let $r\geq 2$. A cyclic quotient singularity is terminal if and only if it is of type
\[
\frac{1}{r}(1,a,r-a) 
\]
for some $0<a<r$ where $\gcd(a,r)=1$.
\end{Thm}

\subsection{Fano 3-folds in the Graded Ring Database} \label{sect:term}

\begin{Def}
A normal projective variety over $\mathbb{C}$ is called a $\mathbb{Q}$-\textbf{Fano $3$-fold} or \textbf{Fano $3$-fold} if $X$ is a 3-dimensional variety with $\mathbb{Q}$-factorial terminal singularities and such that $-K_X$ is ample. The \textbf{Fano index} of $X$ is 
\[
\iota_X := \max \{q \in \mathbb{Z}_{>0} \, | \, -K_X \sim_{\mathbb{Q}} qA,\, \text{$A \in \Cl(X)$ is a $\mathbb{Q}$-Cartier Divisor}\}.
\]  
If $\iota_X \geq 2$ we say that $-K_X$ is \textbf{divisible} in $\Cl(X)$ and that $X$ has \textbf{higher index}.
\end{Def}

Let $A \in \Cl(X)$ be an ample divisor for which $-K_X \sim_{\mathbb{Q}} \iota_XA$. We consider $X$ to be polarised by $A$, that is, we consider $X$ together with an embedding into weighted projective space given by the ring of sections,
\[
R(X,A):=\bigoplus_{m \geq 0} H^0(X,\mathcal{O}_X(mA)).
\] 
Notice that $R(X,-K_X) \subset R(X,A)$. These are finitely generated over $\mathbb{C}$ and we realise $X$ as
\[
X:=\Proj R(X,A).
\]
A choice of generators $x_0,\ldots,x_N$ where $x_i \in H^0(X,\mathcal{O}_X(a_iA))$ determines a map 
\[
X \hookrightarrow \Proj\mathbb{C}[x_0, \ldots,x_N] = \mathbb{P}(a_0,\ldots,a_N)=\mathbb{P}.
\]  We define the codimension of the 3-fold $X$ (with respect to this embedding) to be $N-3$. A subvariety $X \subset \mathbb{P}$ of codimension $N-3$ is a \textbf{weighted complete intersection (WCI) of multidegree $(d_1,\ldots, d_{N-3})$} if its weighted homogeneous ideal $I_X \subset \mathbb{C}[x_0,\ldots,x_N]$ is generated by a regular sequence of homogeneous elements of degrees $d_i$ in  $\mathbb{C}[x_0,\ldots,x_N]$.  Unless otherwise stated, we always assume that $\mathbb{P}(a_0,\ldots,a_N)$ is \textbf{well formed}, that is, we have $\gcd(a_0,\ldots, \widehat{a_i},\ldots, a_N) = 1$ for each $i$. Each weighted projective space is isomorphic to a well formed weighted projective space, see \cite[1.3.1]{DolgachevWPS}. A subvariety $X \subset \mathbb{P}$ is \textbf{well formed} if $\mathbb{P}$ is well formed and $\codim_X(X \cap \mathbb{P}_{\text{sing}}) \geq 2$.

\begin{Def}[{\cite[Definition~3.1.5]{DolgachevWPS}}]
Consider the weighted projective space $\mathbb{P}$ and the natural quotient map $\pi \colon \mathbb{A}^{n+1} \setminus 0 \rightarrow \mathbb{P}$. A subvariety $X \subset \mathbb{P}$ is said to be \textbf{quasismooth} if its affine cone $C_X:=\pi^{-1}X$ is smooth outside the origin.
\end{Def}

\begin{Thm}[{\cite[Theorem~3.2.4]{DolgachevWPS}}] \label{thm:pic1}
Assume $X \subset \mathbb{P}$ is a quasismooth and well formed WCI of dimension at least 3. Then $\Pic(X)=\mathbb{Z}$.
\end{Thm}

When $X$ is quasismooth its singularities are exclusively due to the $\mathbb{C}^{*}$-action defining the weighted projective space $\mathbb{P}$. If $X$ is additionally terminal, its singularities are terminal cyclic quotient singularities as in Theorem \ref{thm:cqs}.

We have a description of the deformation families $X$ which are quasismooth Fano 3-folds whose general member is not a linear cone and can be realised as a complete intersection in some weighted projective space.  Iano-Fletcher has constructed lists of the possible Hilbert series of anticanonically embedded quasismooth Fano 3-fold complete intersections, see \cite[Tables~5 and 6]{fletcher}. It was later proved that the list of corresponding quasismooth Fano 3-folds derived from those tables are complete, see \cite{ccc}.  On the other hand, with a similar approach in \cite{brownsuI, brownsuII}, the authors look at the possible Hilbert series associated to $R(X,A)$ where $X$ is a Fano 3-fold with Fano index at least 2. They found that there is a total of 1964 Hilbert series and obtained the following important result as a consequence:
\begin{Thm}[{\cite[Theorem~1]{brownsuII}}]
Let $X$ be a Fano 3-fold. If $-K_X$ is divisible then $H^0(X,-K_X) \not = 0$.
\end{Thm}

The codimension $c$ of a quasismooth Fano $d$-fold in weighted projective space which is not a linear cone is bounded above by $d$, that is, we have $c \leq d$, see \cite{ccc}. In particular, with the assumption of quasismoothness and the exception of degenerate cases such as linear cones, there are no Fano $3$-folds of codimension four or higher which can be realised as a complete intersection in weighted projective space.

%



The following result is crucial in the construction of Sarkisov links from $X$, see Lemma \ref{lem:iso}, as well as for the exclusion of smooth points as maximal centres, see Subsection \ref{subsec:smooth}.

\begin{Lem} \label{lem:awaybase}
Let 
\[
X \subset \Proj \mathbb{C}[x_0,\ldots,x_5] \simeq \mathbb{P}(a_0,\ldots,a_5)
\]
be a quasismooth WCI Fano 3-fold deformation family for which $\iota_X \geq 2$. Then, the rational curve $\Gamma : (x_2=x_3=x_4=x_5=0) \simeq \mathbb{P}^1$ does not intersect $X$.
\end{Lem} 

\begin{proof}
This is a consequence of quasismoothness of $X$. Using \cite[Proposition~4.2]{ccc}, one can see that $d_1$ and $d_2$ need to be multiples of $\iota_X$ if $\iota_X \not =4$. On the other hand, if $\iota_X=4$, then exactly two weights are multiples of 2 and $d_1$ and $d_2$ are even. Suppose $x_0^{\alpha_0}x_1^{\alpha_1}x_{\mu}$ is a monomial in $f$. Then,
\[
\wt(x_{\mu}) = d_1 - \alpha_0\wt(x_0)- \alpha_1\wt(x_1) \equiv 
\begin{cases*}
        0 \pmod{\iota_X}, & \text{for } $\iota_X \not = 4$\\
        0 \pmod{2}, & \text{for } $\iota_X  = 4$
\end{cases*} 
\] 
 Hence, $x_{\mu} \in \{x_0,x_1 \}$ and 
\begin{equation*}
\frac{\partial f}{\partial x_{\mu'}}\bigg\rvert_{\Gamma},\, \frac{\partial g}{\partial x_{\mu'}}\bigg\rvert_{\Gamma}
\end{equation*}
vanish identically if $x_{\mu'} \not \in  \{x_0,x_1 \}$. In other words, 
\[
J(X)|_{\Gamma} = 
\begin{pmatrix}
\partial_{x_0}f(x_0,x_1) & \partial_{x_1}f(x_0,x_1) & 0& 0& 0& 0\\
\partial_{x_0}g(x_0,x_1) & \partial_{x_1}g(x_0,x_1) & 0& 0& 0& 0
\end{pmatrix}
\]
where the column $0 \leq i \leq 5$ is the partial derivative with respect to $x_i$. Let $h \in \mathbb{C}[x_0,x_1]$ be a common factor of $f|_{\Gamma}$ and $g|_{\Gamma}$. We claim that if $X$ is quasismooth, then $h$ is the constant non-zero polynomial. Write $f|_{\Gamma}=hf'$ and $g|_{\Gamma}=hg'$. Then, the first minor $M$ of $J(X)|_{\Gamma}$ is
\begin{align*}
\partial_{x_0}f|_{\Gamma}\partial_{x_1}g|_{\Gamma}- \partial_{x_1}f|_{\Gamma}\partial_{x_0}g|_{\Gamma} &= f'g'(\partial_{x_0}h\partial_{x_1}h-\partial_{x_1}h\partial_{x_0}h)\\
&+h^2(\partial_{x_0}f'\partial_{x_1}g'-\partial_{x_1}f'\partial_{x_0}g')  \\ 
&+hf'(\partial_{x_0}h\partial_{x_1}g'-\partial_{x_1}h\partial_{x_0}g') \\
&+hg'(\partial_{x_0}f'\partial_{x_1}h-\partial_{x_1}f'\partial_{x_0}h).
\end{align*}
Clearly, the first line, $\partial_{x_0}h\partial_{x_1}h-\partial_{x_1}h\partial_{x_0}h$, vanishes and $M \in (h)$. Since $X$ is quasismooth, it follows that $h \in \mathbb{C}^*$.
\end{proof}

\subsection{The Sarkisov Program} \label{sect:The Sarkisov Program}

Let $W$ be a smooth projective variety. If $K_W$ is not pseudo-effective, then by \cite[Corollary~1.3.3]{BCHM}, the Minimal Model Program produces a birational model $V$ of $W$ which admits a very special fibre structure $V \rightarrow S$ called a Mori fibre space.

\begin{Def}
A \textbf{Mori fibre space} is an extremal contraction $f \colon V \rightarrow S$ between normal varieties for which $f_*\mathcal{O}_V=\mathcal{O}_S$ and
\begin{itemize}
	\item $V$ has at worst $\mathbb{Q}$-factorial, terminal singularities;
	\item $-K_V$ is relatively ample for $f$;
	\item The relative Picard rank satisfies $\rho (V/S)=1$ and $\dim S < \dim V$.
\end{itemize}
We sometimes denote the Mori fibre space $f \colon V \rightarrow S$ by $V/S$. A \textbf{$\mathbb{Q}$-Fano variety} is a Mori fibre space $f \colon V \rightarrow S$ where the base is a point. We denote it simply by $V$.
\end{Def}

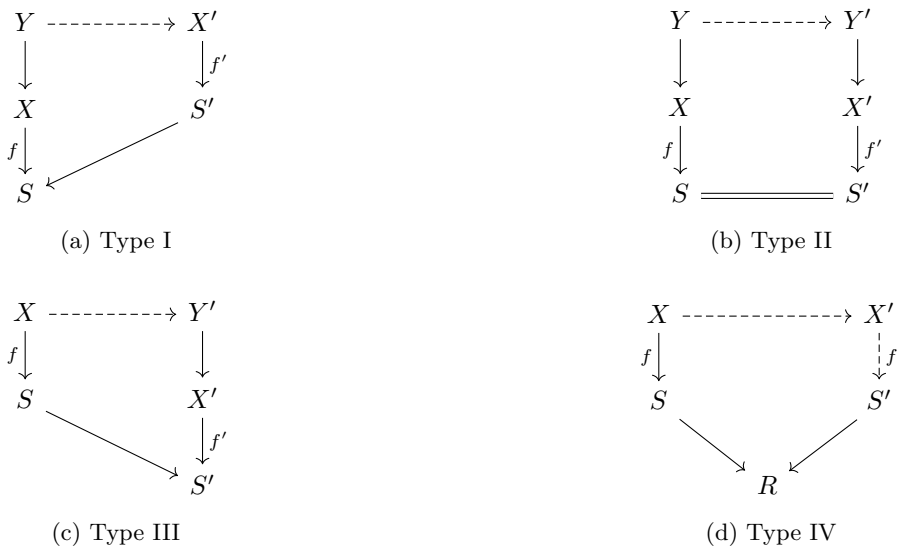
\begin{figure*}[h!]
        \centering
        \begin{subfigure}[b]{0.475\textwidth}
            \centering
\begin{tikzcd}
Y\ar[d] \ar[rr, dashed]&  &  X' \ar[d, "f'"]    \\
X \ar[d, swap,"f"]  & & S'\ar[dll]\\
S &  &    
\end{tikzcd}

            \caption[Link1]%
            {{\small  Type I}}    
            \label{fig:1}
        \end{subfigure}
        \hfill
        \begin{subfigure}[b]{0.475\textwidth}  
            \centering 
\begin{tikzcd}
Y\ar[d] \ar[rr, dashed]&  &  Y' \ar[d]    \\
X \ar[d, swap,"f"]  & & X'\ar[d,"f'"]\\
S \ar[rr,equal] &  & S'    
\end{tikzcd}
            \caption[]%
            {{\small Type II}}    
            \label{fig:2}
        \end{subfigure}
        \vskip\baselineskip
        \begin{subfigure}[b]{0.475\textwidth}   
            \centering 
\begin{tikzcd}
X\ar[d, "f",swap] \ar[rr, dashed]&  &  Y' \ar[d]    \\
S \ar[drr]   & & X'\ar[d, "f'"]\\
 &  & S'   
\end{tikzcd}
            \caption[]%
            {{\small Type III}}    
            \label{fig:3}
        \end{subfigure}
        \hfill
        \begin{subfigure}[b]{0.475\textwidth}   
            \centering 
\begin{tikzcd}
X\ar[d,swap, "f"] \ar[rr, dashed]&  &  X' \ar[d, "f'", dashed]    \\
S \ar[dr]  & & S'\ar[dl]\\
 & R &    
\end{tikzcd}
            \caption[]%
            {{\small Type IV}}    
            \label{fig:4}
        \end{subfigure}
        \caption{\small The four types of Sarkisov links between Mori fibre spaces.} 
        \label{fig:4types}
    \end{figure*}

 The Sarkisov Program is the study of the relations between birational Mori fibre spaces. The following is the fundamental result of the Sarkisov Program: 

\begin{Thm}[{\cite[Theorem~1.1]{genSP}, \cite{cortiSP} for $3$-folds}] \label{thm:sp}
Two Mori fibre spaces are birational if and only if they are connected by a finite sequence of elementary Sarkisov links.
\end{Thm}
Recall that an elementary Sarkisov link $\sigma \colon X \dashrightarrow X'$ between two Mori fibre spaces $f \colon X \rightarrow S$ and $f' \colon X' \rightarrow S'$ is one of four types as in figure \ref{fig:4types}. 

Every non-horizontal arrow is an extremal contraction. The dashed arrows are (finite) sequences of anti-flips, flops or flips, in that order. Notice that a Sarkisov link from a $\mathbb{Q}$-Fano variety can only be of Type I or Type II.  The following definition captures the uniqueness of $\mathbb{Q}$-Fano varieties and is essentially due to Corti. See \cite[Definition~1.3]{cortising} and \cite[Definition~1.1.2]{chel}.
\begin{Def}
The $\mathbb{Q}$-Fano variety $X$ is \textbf{birationally rigid} if whenever there is birational map $\sigma \colon X \rat X'$ to a Mori fibre space $X'/S$, then $X$ and $X'$ are biregular (and in particular $X'/S$ is a $\mathbb{Q}$-Fano variety). If, in addition, the birational automorphism group of $X$ coincides with its biregular automorphism group we say that $X$ is \textbf{birationally super-rigid}.  
\end{Def} 

\subsection{Maximal extractions and Mori Dream Spaces}

Let $X$ be a $\mathbb{Q}$-Fano variety and $\mathcal{H}$ a movable linear system on $X$. Let $n \in \mathbb{Q}_{>0}$ be such that $\mathcal{H}+nK_X \sim_{\mathbb{Q}}0$.
\begin{Def}
The \textbf{canonical threshold} of the pair $(X,\mathcal{H})$ is 
\[
c(X,\mathcal{H}):= \max\{ \lambda\,\,|\,\, K_X+\lambda\mathcal{H}\, \text{is canonical} \} \in \mathbb{Q}_{>0}.
\]
A \textbf{maximal singularity} of $(X,\mathcal{H})$ is an extremal divisorial contraction $Y\rightarrow X$ with exceptional divisor $E$ such that 
\[
\frac{1}{n} > c(X,\mathcal{H}).
\]
Notice that $ c(X,\mathcal{H})= \frac{a_E(K_X)}{m_E(\mathcal{H})}$ where $m_E(\mathcal{H})$ is the multiplicity of $\mathcal{H}$ along $E$ and $a_E(K_X)$ is the discrepancy of $Y \rightarrow X$. A subvariety $Z \subset X$ is called a \textbf{maximal centre} if there is a maximal extraction $Y \rightarrow X$ centred at $Z$, that is, a maximal singularity of $\mathcal{H}$ on $X$, $h \colon Y\rightarrow X$, for which $h(E)=Z$.
\end{Def} 

\begin{Rem} \label{rem:maxcentre}
Suppose that $Z \subset X$ is a maximal centre. By definition there is a maximal singularity $h \colon Y \rightarrow X$ centred at $Z$ with respect to a movable linear system $\mathcal{H}$ as above. In particular, the pair $(X, \frac{1}{n}\mathcal{H})$ is not canonical where $n$ is such that $\mathcal{H}+nK_X \sim_{\mathbb{Q}}0$. Hence, according to the Noether-Fano-Iskovskihk inequality, see \cite[Theorem~1.26]{chelshra}, $X$ is not birationally superrigid. It follows that there must be a Mori fibre space $X' \rightarrow S'$ and a birational map $\sigma \colon X \dashrightarrow X'$ which is not an isomorphism. The Sarkisov Program says that $\sigma$ can be factored as a finite sequence of elementary Sarkisov links starting with the maximal singularity $h \colon Y \rightarrow X$, see \cite[Proposition~2.10]{cortiSP}.
\end{Rem}


Recall the definition of a Mori Dream Space

\begin{Def}[{\cite[Definition~1.10]{mdsGIT}}]
A normal projective variety $X$ is a \textbf{Mori Dream Space} if the following hold
\begin{itemize}
	\item $X$ is $\mathbb{Q}$-factorial and $\Pic(X)_{\mathbb{Q}} = N^1(X)$;
	\item $\Nef(X)$ is the affine hull of finitely many semi-ample line bundles;
	\item There exists a finite collection of small $\mathbb{Q}$-factorial modifications $\tau_i \colon X \dashrightarrow X_i$ such that each $X_i$ satisfies the previous point and $\overline{\Mov}(X)= \bigcup \tau_i^*(\Nef(X_i))$.
\end{itemize}
\end{Def}


The following result says, in particular, that the only way for a divisorial contraction $f \colon Y \rightarrow X$ from a $\mathbb{Q}$-Fano variety to initiate a Sarkisov link is for $Y$ to be a Mori Dream Space.

\begin{Lem}[{\cite[Lemma~2.9]{hamidquartic}}] \label{lem:linkmds}
Let $X$ be a $\mathbb{Q}$-Fano 3-fold and $\varphi \colon Y \rightarrow X$ be a divisorial extraction. Then $\varphi$ initiates a Sarkisov link if and only if the following hold:
\begin{enumerate}
    \item $Y$ is a Mori Dream Space;
    \item If $\tau \colon Y \dashrightarrow Y'$ is a small birational map and $Y'$ is $\mathbb{Q}$-factorial, then $Y'$ is terminal;
    \item $[-K_Y] \in \Int(\overline{\Mov}(Y))$
\end{enumerate}
\end{Lem}

Putting together remark \ref{rem:maxcentre} and lemma \ref{lem:linkmds}, we have the following corollary:

\begin{Cor}\label{cor:exmds}
Suppose $\mathcal{H}$ is a movable linear system in the $\mathbb{Q}$-Fano 3-fold $X$. Let $\varphi \colon Y \rightarrow X$ be a maximal singularity of the pair $(X,\mathcal{H})$. Then $Y$ is a Mori Dream Space.
\end{Cor}


Suppose $X$ and $\varphi \colon Y \rightarrow X$ are as in corollary \ref{cor:exmds}. Then $\rho(Y)=2$ and we can play the 2-ray game on $Y$. The following explanation is based on the proof of \cite[Lemma~2.9]{hamidquartic}. Indeed, let
\[
\overline{\Mov}(Y) = \mathbb{R}_+[M_1]+\mathbb{R}_+[M_1] \subset \overline{\Eff}(Y) = \mathbb{R}_+[D_1]+\mathbb{R}_+[D_2].
\]
where $M_i$ and $D_i$ are effective divisors. Assume that $\varphi \colon Y \rightarrow X$ is a maximal extraction from a maximal centre $Z \subset X$. Since $X$ is a $\mathbb{Q}$-Fano 3-fold, $\varphi$ is the ample model (see \cite[Def~3.6.5]{BCHM}  for a precise definition) of the big divisor $M_1$. Hence $M_1 \not = D_1$. Let $\epsilon$ be a small positive rational number and $Y'$ be the ample model of $M_2+\epsilon M_1$. Then, there is a birational map $\tau \colon Y \dashrightarrow Y'$ which is a finite sequence of anti-flips, flops or flips in one-to-one correspondence to the Mori Chambers of $\overline{\Mov}(Y)$, see \cite[Proposition~1.11]{mdsGIT}. Moreover, since $\tau$ is small, the movable cones of $Y$ and $Y'$ can be identified and, in particular, by \cite[Theorem~3.2]{hamidmds}, $-K_{Y'} \in \Mov(Y')$. If $Y$ has terminal singularities, then $Y'$ has terminal singularities if and only if the anti-flips are terminal. There are two cases to consider for the end of a link: 
\begin{enumerate}
\item Suppose $M_2=D_2$. Let $S'$ be the ample model of $M_2$. Since $M_2 \in \partial \overline{\Eff}$, $M_2$ is not big and $\dim S' < \dim Y'$. Moreover, since $-K_{Y'} \in \Mov(Y')$ the fibration $\varphi' \colon Y' \rightarrow S'$ is $K_{Y'}$-negative. We conclude that $\sigma \colon X \dashrightarrow Y'/S'$ is a Sarkisov link of Type I.
\item Suppose $M_2\not=D_2$. Let $X_2$ be the ample model of $M_2$. Then, the morphism $\varphi' \colon Y' \rightarrow X'$ contracts $D_2$ to a centre in $X'$. Moreover, $X'$ has terminal singularities since $\varphi'$ is $K_{Y'}$-negative. Hence $X'$ is a $\mathbb{Q}$-Fano 3-fold and $\sigma \colon X \dashrightarrow X'$ is a Sarkisov link of type II.  
\end{enumerate}

\paragraph{The Kawamata blowup.} A maximal centre $\Gamma$ on a $\mathbb{Q}$-Fano 3-fold $X$ can be one of three types: A curve, a smooth point or a singular point. See Chapter \ref{sect:excl} for when $\Gamma$ is either of the first two cases. The following result of Kawamata classifies extractions in the Mori category from germs of 3-fold terminal cyclic quotient singularities.  

\begin{Thm}[{\cite{kawamata}}] \label{thm:kwbl}
Let $(\mathbf{p}\in X) \sim \frac{1}{r}(1,a,r-a)$ be the germ of a terminal cyclic quotient singularity. If $\varphi \colon Y \rightarrow X$ is a divisorial extraction centred along $\Gamma$ where $\mathbf{p} \in \Gamma$, then $\Gamma= \mathbf{p}$ and $\varphi$ is the weighted blowup with weights $\frac{1}{r}(1,a,r-a)$ and discrepancy $\frac{1}{r}$. Moreover, the exceptional divisor is $E \simeq \mathbb{P}(1,a,r-a) \subset Y$.
\end{Thm}
 
The weighted blowup of the previous theorem is called the \textbf{Kawamata blowup centred at $\mathbf{p}$}. An immediate corollary is the following

\begin{Cor}[{\cite[Corollary~3.4.3]{CPR}}] \label{cor:kwblcor}
Suppose $X$ is a 3-fold with only terminal quotient singularities. If $\Gamma \subset X$ is the centre of a divisorial extraction in the Mori category then $\Gamma \subset X \setminus X_{\text{sing}}$.
\end{Cor}

\subsection{Toric Embedding and the 2-ray game} \label{sub:toricemb}

By \cite[Proposition~2.11]{mdsGIT}, the birational contractions of a Mori Dream Space $Y$ are induced from toric geometry. In particular, there is an embedding $Y \subset T$ into a quasismooth projective toric variety for which every Mori chamber of $Y$ is a finite union of Mori chambers of $T$. We briefly explain the idea of the embedding $Y \subset T$ we use. Let $X \subset \mathbb{P}$ be a quasismooth $\mathbb{Q}$-Fano 3-fold. If $\mathbf{p} \in X \subset \mathbb{P}$ is a germ of a cyclic quotient singularity point, we perform a toric blowup  $\Phi \colon T \rightarrow \mathbb{P}$ centred at $\mathbf{p}$ with the constraint that it restricts locally around $\mathbf{p}$ to the unique Kawamata blowup, see Theorem \ref{thm:kwbl}, 
\[
\varphi \colon E \subset Y \rightarrow \mathbf{p} \in X, 
\]
where $Y:= \Phi_*^{-1}X$ and $\varphi$ is the restriction $\Phi|_{\Phi_*^{-1}X}$. This is summarised in the following diagram,
\[
\begin{tikzcd}[column sep = 4em]
E \arrow[swap]{d}{\varphi|_E} \arrow[swap,hookrightarrow]{r}{} & Y \arrow[swap]{d}{\varphi} \arrow[swap,hookrightarrow]{r}{} \arrow[dashed]{r}{} & T \arrow[swap]{d}{\Phi}    \\
 \mathbf{p} \arrow[swap,hookrightarrow]{r}{} & X \arrow[hookrightarrow]{r}{} & \mathbb{P}   
\end{tikzcd}
\]

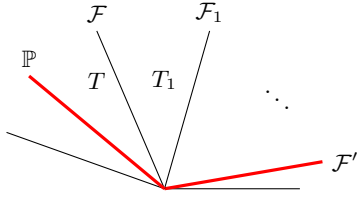
\begin{figure}%
\centering
\begin{tikzpicture}[scale=3,font=\small]
  \coordinate (A) at (0, 0);
  \coordinate [label={left:}] (E) at (-0.7, 0.25);
  \coordinate [label={above:$\mathbb{P}$}] (K) at (-0.6,0.5);
	\coordinate [label={above:$\mathcal{F}$}] (5) at (-0.3, 0.7);
	\coordinate [label={above:$\mathcal{F}_1$}] (2) at (0.2,0.7);
	\coordinate [label={right:$\mathcal{F}'$}] (6) at (0.7,0.12);
	\coordinate [label={right:}] (4) at (0.6,0);
	
		\coordinate [label={$T$}] (T) at (-0.3,0.4);
		\coordinate [label={$T_1$}] (T1) at (0,0.4);
		\coordinate [label={$\ddots$}] (T2) at (0.5,0.3);
  \draw  (A) -- (E);
  \draw  (A) -- (4);
	\draw [very thick,color=red] (A) -- (6);
	\draw (A) -- (5);
	\draw (A) -- (2);
	\draw [very thick,color=red](A) -- (K);
\end{tikzpicture}
\caption{A representation of the chamber decomposition of the cone of effective divisors of $T$. The red subcone is the cone of movable divisors of $T$.}%
\label{fig:gen}%
\end{figure}

The constraint that $\Phi$ should restrict to a Kawamata blowup allows us to determine  the rank two toric variety $T$ up to isomorphism. This is done in detail in Lemma \ref{lem:utbl}. Once $T$ is determined, we look at the closure of the movable and effective cones of $T$, denoted $\overline{\Mov}(T)$ and $\overline{\Eff}(T)$, respectively and how these are related. Notice that since $\Phi$ is a divisorial extraction from $\mathbb{P}$, the movable cone of $T$ is strictly contained in the effective cone of $T$. The chambers $f_i^*(\Nef(T_i))$ and their faces form a fan supported at $\overline{\Mov}(T)$. The number of cones in this fan is finite by \cite[Proposition~1.11~(2)]{mdsGIT} and are in one-to-one correspondence with contracting rational maps $g \colon T \rat Q$, with $Q$ normal and projective via
\[
g \colon T \rat Q  \quad \leftrightsquigarrow \quad g^*(\Nef(Q)) \subset \overline{\Mov}(T)
\]
by \cite[Proposition~1.11~(3)]{mdsGIT}. In particular, if $T_i$ and $T_{i+1}$ are the ample models (see \cite[Definition~3.6.5]{BCHM}) of adjacent chambers, then they are related by a small $\mathbb{Q}$-factorial modification that we denote by $\tau_i$. Since in our case $\Pic(T)=2$, we can think of these cones in $\mathbb{R}^2$ and hence there is a completely determined sequence of birational transformations 
\[
\begin{tikzcd} 
& T \arrow[swap]{dl}{\Phi} \arrow[swap]{dr}{f} \arrow[dashed]{rr}{\tau} & &  T_1\arrow[swap]{dr}{f_1} \arrow{dl}{g} \arrow[dashed]{r}{\tau_1} & \cdots \arrow[dashed]{r}{\tau_n}   & T'\arrow{dr}{\Phi'}\arrow{dl}{g_{n-1}}  & \\
 \mathbb{P} &  & \mathcal{F} &  & \cdots  &  & \mathcal{F}' 
\end{tikzcd}
\]
starting from $\Phi \colon T \rightarrow \mathbb{P}$ called the \textbf{2-ray game}. As explained, these maps are induced from a chamber decomposition of $\overline{\Eff}(T)$ such as in figure \ref{fig:gen}. 
There are two possible outcomes of running the 2-ray game on $T$. 
\begin{enumerate}
	\item The closure of the movable cone of $T$ is stricly contained in the interior of the cone of pseudo-effective divisors of $T$. In this case, the map $\Phi'$ is birational and is, in fact, a divisorial contraction to a closed subvariety. This is the case depicted in figure \ref{fig:gen}
	\item The closure of the movable cone of $T$ is \emph{not} stricly contained in the interior of the cone of pseudo-effective divisors of $T$. In this case, there is a divisor class $D$ in the boundary of $\overline{\Mov}(T)$ which is not big. Hence the map $\Phi'$ given by the positive multiples $|mD|$ is not birational and $\Phi'$ is a fibration.
\end{enumerate}

\begin{Ex} \label{ex:bundle}
Let $a\geq -2$ be an integer and consider the Toric variety $T= \Proj_{\mathbb{P}^1} \mathcal{E}$, where 
\[
\mathcal{E}=\mathcal{O}_{\mathbb{P}^1}\oplus\mathcal{O}_{\mathbb{P}^1}(a)\oplus\mathcal{O}_{\mathbb{P}^1}(a).
\]
Let $M$ be the tautological line bundle on $T$ and $L$ a fibre of the natural projection to $\mathbb{P}^1$.  Then, $\Pic(T)=\mathbb{Z}[M]+ \mathbb{Z}[L]$ and
\[
\overline{\Mov}(T)=  \mathbb{R}_{+}[L]+\mathbb{R}_{+}[M-aL]. 
\]
Hence, $-K_T \in |3M+(2-2a)L|$ is in the interior of the movable cone when $a>-2$ and in the boundary otherwise. We play the 2-ray game on $T$ and show a plethora of different behaviors depending on $a$.

Suppose $\overline{\Mov}(T)=\overline{\Eff}(T)$. Then $a\geq 0$. 

\begin{itemize}
	\item If $a=0$, then $T \simeq \mathbb{P}^1 \times \mathbb{P}^2$ and it comes with the two natural projections.
	\item If $a=1$, there is a birational involution on $\mathbb{P}^1$ obtained in the following way. $T$ is the ample model for $\epsilon L+ (M-aL) $ and let $T'$ be the ample model for $L+ \epsilon (M-aL) $, where $\epsilon$ is a small positive rational number. Since $-K_T$ is in the interior of the mobile cone of $T$, there is a small birational modification $\tau \colon T \dashrightarrow T'$. The map $\tau$ contracts the base locus of the linear system $|M-aL|$ and extracts the base locus of the linear system $|L|$ which are both small and trivial against the anticanonical divisor. Hence $\tau$ is a flop, in fact the Atyiah flop.
	\item If $a\geq 2$, with the same notation for the ample models as before, there is a non-terminal anti-flip $\tau \colon T \dashrightarrow T'$, where $T'$ is a $\mathbb{P}(1,1,a)$ fibration over $\mathbb{P}^1$. In this case, $\tau$ extracts a line of singularities to $T'$ and the game goes out of the Mori category. 
	\end{itemize}
	Suppose $\overline{\Mov}(T)\not =\overline{\Eff}(T)$. Then $a\leq  -1$. 
 
	\begin{itemize}
		\item  If $a=-1$, let $X'$ be the ample model for $M-aL=M+L$. Then $X' \simeq \mathbb{P}^3$ and $T$ admits a divisorial contraction from $M\simeq \mathbb{P}^1 \times \mathbb{P}^1$ to a smooth rational curve $T \rightarrow \mathbb{P}^3$.
		\item  If $a=-2$, let $X'$ be the ample model for $M-aL=M+2L$. Then $X' \simeq \mathbb{P}^3(1,1,-a,-a)$ and $T$ admits a divisorial contraction to a rational curve $\varphi \colon T \rightarrow \mathbb{P}^3(1,1,-a,-a)$. Notice that $-K_T \sim 3(M-aL)$ and therefore $-K_T$ is in the boundary of the mobile cone of $T$. It is easily checked that $K_T-\varphi^*(K_{\mathbb{P}(1,1,-a,-a)}) \sim 0$, i.e., it is a crepant blowup, hence, the game goes out of the Mori category. Indeed $\mathbb{P}(1,1,-a,-a)$ contains a line of $A_1$ singularities so it is not terminal.  
	\end{itemize}
\end{Ex}

\paragraph{Variation of GIT.} For a Mori Dream Space, the GIT chambers and Mori chambers coincide, see \cite[Theorem~2.3]{mdsGIT}. Hence, the 2-ray game on a rank 2 toric variety $T$ can be obtained by variation of GIT or vGIT as explained in \cite[Section~4]{2raybrown}. In order to illustrate this we use the case of $a=1$ in example \ref{ex:bundle}. In this case, the Cox ring of $T$ is $\mathbb{C}[y_1,y_2,t,x_1,x_2]$ with weight system
\[
\begin{array}{cccc|cccc}
             &       & y_1  & y_2 &   t & x_1 & x_2  & \\
\actL{T}   &  \lBr &  0 & 0 &   1 & 1 & 1  &   \actR{.}\\
             &       & 1 & 1 &   0 & -1 & -1 &  
\end{array}
\] 
and irrelevant ideal $(y_1,y_2) \cap (t,x_1,x_2)$ as shown by the vertical bar. That is, we have an action of $(\mathbb{C}^*)^2$ on $\mathbb{C}^5$ given by
\begin{align*}
\lambda \cdot (y_1,y_2,t,x_1,x_2) &= (y_1,y_2, \lambda t, \lambda x_1,\lambda x_2) \\
\mu \cdot (y_1,y_2,t,x_1,x_2) &= (\mu y_1,\mu y_2,  t, \mu^{-1} x_1,\mu^{-1} x_2).  
\end{align*}
The GIT chambers of $T$ are 

\begin{figure}[h]%
\centering
\begin{tikzpicture}[scale=3,font=\small]
  \coordinate (A) at (0, 0);
  \coordinate [label={left:$(0,1)$}] (E) at (0,0.5);
  \coordinate [label={above:$(1,0)$}] (K) at (0.5,0);
	\coordinate [label={right:$(1,-1)$}] (5) at (0.5,-0.5);
	\draw  (A) -- (E);
  \draw  (A) -- (K);
	\draw (A) -- (5);
	\end{tikzpicture}
\caption{A representation of the GIT chamber decomposition of $T$.}%
\label{fig:git}%
\end{figure}
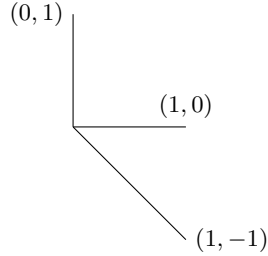
Choosing a character in the ray generated by $(0,1)$ constructs a variety isomorphic to $\mathbb{P}^1$. This is given by 
\[
\Proj \bigoplus_{m\geq 1} H^0(T,\mathcal{O}(0,m)) = \Proj \mathbb{C}[y_1,y_2] \simeq \mathbb{P}^1
\]
and similarly with the ray generated by $(1,-1)$. The ray $(1,0)$ constructs the non-$\mathbb{Q}$-factorial variety
\[
\Proj \bigoplus_{m\geq 1} H^0(T,\mathcal{O}(m,0)) = \Proj \mathbb{C}[t,y_1x_1,y_1x_2,y_2x_1,y_2x_2] \simeq (u_{11}u_{22}-u_{12}u_{21}=0) \subset \mathbb{P}^4
\]
where $u_{ij} := y_ix_j$. This is a quadric cone $Q$ whose only singular point is $p=(1:0:0:0:0)$. On the other hand, choosing characters in the interior of the two GIT chambers constructs varieties related by a small $\mathbb{Q}$-factorial modification. Choosing a character in the interior of the cone generated by $\left(\begin{smallmatrix} 0 \\ 1 \end{smallmatrix} \right)$ and $\left(\begin{smallmatrix} 1 \\ 0 \end{smallmatrix} \right)$ constructs a variety isomorphic to $T$ and choosing a character in the interior of the cone generated by $\left(\begin{smallmatrix} 1 \\ 0 \end{smallmatrix} \right)$ and $\left(\begin{smallmatrix} 1 \\ -1 \end{smallmatrix} \right)$ constructs a variety isomorphic to $T'$, where $T'$ has the same Cox ring of $T$ but its irrelevant ideal is $(y_1,y_2,t) \cap (x_1,x_2)$. Consider the map $f \colon T \rightarrow Q$ given by
\[
(y_1,y_2,t,x_1,x_2) \mapsto (t,y_1x_1,y_1x_2,y_2x_1,y_2x_2).
\] 
and similarly $g\colon T' \rightarrow Q$. Then $f$ contracts the locus $C^{-} \colon \mathbb{P}^1 \simeq (x_1=x_2=0) \subset T$ to $p \in Q$ and $g$ contracts $C^{+} \colon \mathbb{P}^1 \simeq (y_1=y_2=0) \subset T'$. Notice that $f$ and $g$ are isomorphisms away from the contracted loci. Hence, the map $\tau \colon T \rat T'$ replacing $C^-$ by $C^+$ is a small $\mathbb{Q}$-factorial modification that we denote by $(-1,-1,1,1)$. In fact, it is the Atyiah flop since $K_T \cdot C^- = K_{T'} \cdot C^+ = 0$. Moreover, the map 
\[
T' \rightarrow \mathbb{P}^1, \quad (y_1,y_2,t,x_1,x_2) \mapsto (x_1,x_2)
\]
is a fibration whose fibres are isomorphic to $\mathbb{P}^2$.

\paragraph{Restriction to $Y$.}

Our goal in playing the 2-ray game on $T$ is to restrict it to $Y$ and obtain a Sarkisov link.

\begin{Def}[{\cite[Definition~3.5]{hamidmds}}]
Let $Y \subset T$ be $\mathbb{Q}$-factorial projective varieties and suppose $\Pic(T)=2$. We say that \textbf{$Y$ 2-ray follows $T$} or that \textbf{$Y$ follows $T$} if the 2-ray game on $T$ restricts to a 2-ray game on $Y$ where $Y_i \subset T_i$ is given by the ideal $I_Y$ and the small $\mathbb{Q}$-factorial modification $T_i \rat T_{i+1}$ over $\mathcal{F}_i$ restricts to a small $\mathbb{Q}$-factorial modification $Y_i \rat Y_{i+1}$ over $\mathcal{F}_i|_{Y_i}$. 
\end{Def}


The following lemma is analogous to \cite[Lemma~3.6]{hamidmds}.

\begin{Lem} \label{lem:utbl}
Let $X\in I$ be a quasismooth $\mathbb{Q}$-Fano 3-fold where 
\[
X:=X_{d_1,d_2} \colon (f=g=0) \subset \Proj \mathbb{C}[x_0,\ldots,x_4,\xi] \simeq \mathbb{P}(a_0, \ldots,a_4, a_{\xi}) :=\mathbb{P}
\]
with no specified order on the weights and $\mathbf{p_{\xi}} \in X$ be a coordinate point whose germ is a terminal cyclic quotient singularity. Then, there is a unique toric blowup $\Phi: T \rightarrow \mathbb{P}$ centred at $\mathbf{p_{\xi}} \in \mathbb{P}$ that restricts to the Kawamata blowup of $X$ centred at $\mathbf{p_{\xi}} \in X$. 
\end{Lem}

\begin{proof} 
The orbifold point $\mathbf{p_{\xi}} \in \mathbb{P}$ is the cyclic quotient singularity $\mathbb{C}^4/\bm{\mu}_{a_{\xi}}:=\frac{1}{a_{\xi}}(a_0,\ldots,a_4)$. We consider the weighted blowup of the affine space $\mathbb{C}^4/\bm{\mu}_{a_{\xi}}$ at the origin with weights $\sigma:=\frac{1}{a_{\xi}}(b_0, \ldots,b_4)$. The cone corresponding to $\mathbb{C}^4/\bm{\mu}_{a_{\xi}}$ is 
\[
C=\mathbb{R}_+[e_0]+ \cdots + \mathbb{R}_+[e_4] 
\]
where the $e_i$ are the canonical coordinates of $\mathbb{Z}^5$. The weighted blowup of $\mathbb{C}^4/\bm{\mu}_{a_{\xi}}$ at the origin corresponds to a subdivision of $C$ into the five cones
\[
C_k=\mathbb{R}_+[e_0]+ \cdots+\mathbb{R}^k_+[\sigma] + \cdots+\mathbb{R}_+[e_4], \quad k \in \{0, \ldots, 4 \}.
\]
That is, for each $k$ we have $\mathbb{R}_+[\sigma]$ instead of $\mathbb{R}_+[e_k]$. The blowup map is the natural map between the refined fan and the original one and the exceptional divisor is isomorphic to $\mathbb{P}(a_0, \ldots, a_4)$.
We now determine the $b_i$.

Since $X$ is quasismooth and $\mathbf{p_{\xi}}\in  X$, the equations defining $X$ are of the form 
\begin{align*}
f\colon \xi^{r_1}x_3+f(x_0,x_1,x_2,x_3,x_4,\xi)&=0\\
g\colon \xi^{r_2}x_4+g(x_0,x_1,x_2,x_3,x_4,\xi)&=0
\end{align*}
where $r_1$ and $r_2$ are the highest exponents of $\xi$ in $f$ and $g$, respectively. Then, $\mathbf{p_{\xi}} \sim \frac{1}{a_{\xi}}(a_0,a_1,a_2)$, where the first component is $a_0=\iota_X$ and is due to the local polarisation of the singularity by the anticanonical class.
From terminality of $\mathbf{p_{\xi}}$, we have that $\gcd(\iota_X,a_{\xi})=1$ and so there is a (unique) inverse $0 < k < a_{\xi}$ of $\iota_X \pmod {a_{\xi}}$. Using the automorphism of $\bm{\mu}_{a_{\xi}}$ given by $\mu \mapsto \mu^{k}$ and changing coordinates in $\mathbb{C}^3$ the cyclic quotient singularity becomes
\[
\mathbf{p_{\xi}} \sim \frac{1}{a_{\xi}}(\iota_X,a_1,a_2) \sim \frac{1}{a_{\xi}}(1,\overline{ka_1},a_{\xi}-\overline{ka_1}). 
\]
Hence by Theorem \ref{thm:kwbl}, the Kawamata blowup $\varphi \colon E \subset Y \rightarrow X$ of $\mathbf{p_{\xi}}$ is the weighted blowup with weights $\frac{1}{a_{\xi}}(1,\overline{ka_1},a_{\xi}-\overline{ka_1})$ and the exceptional divisor is $E \simeq \mathbb{P}(1,\overline{ka_1},a_{\xi}-\overline{ka_1})$. We conclude that $b_0=1,\, b_1=\overline{ka_1},\, b_2=a_{\xi}-\overline{ka_1}$.


Let $T_0$ be the (non well-formed, see \cite[Definition~2.11]{hamidplia}) rank two toric variety
\[
\begin{array}{cccc|ccccccc}
             &       & u  & \xi &   x_0 & x_1 & x_2 & x_3 & x_4 & \\
\actL{T_0}   &  \lBr &  0 & a_{\xi} &   \iota_X & a_1 & a_2 & a_3 & a_4 &   \actR{}\\
             &       & -a_{\xi} & 0 &   1 & \overline{ka_1} & a_{\xi}-\overline{ka_1} & b_3 & b_4 &  
\end{array}
\]
and consider the blowup map.
\begin{align*}
 \Phi \colon T_0 &\rightarrow \mathbb{P}(\iota_X, a_1, a_2, a_3, a_4, \xi )\\
(u,\xi, x_0, x_1, x_2, x_3, x_4) &\mapsto (x_0u^{1/{a_{\xi}}}: x_1u^{\overline{ka_1}/{a_{\xi}}} : x_2u^{(a_{\xi}-\overline{ka_1})/{a_{\xi}}}: x_3u^{b_3/{a_{\xi}}}:x_4u^{b_4/{a_{\xi}}}:{\xi} ).
\end{align*}
The map $\Phi$ realises the $\frac{1}{a_{\xi}}(1,\overline{ka_1}, a_{\xi}-\overline{ka_1}, b_3, b_4)$-weighted blowup explained in the beginning of the proof. We now determine $b_3$ and $b_4$ using the fact that the restriction of $\Phi$ to the Kawamata blowup of $X$ centred at $\mathbf{p_{\xi}}$ has discrepancy $\frac{1}{a_{\xi}}$.

Define $Y_0= \Phi^{-1}_{|_X}$. Then $Y_0$ is given by 
\begin{align*}
\Phi^*f &\colon \xi^{r_1}x_3u^{b_3/a_{\xi}}+f'(u,x_0,x_1,x_2,x_3,x_4,\xi)=0\\
\Phi^*g&\colon \xi^{r_2}x_4u^{b_4/a_{\xi}}+g'(u,x_0,x_1,x_2,x_3,x_4,\xi)=0.
\end{align*}

Let $m_f$ and $m_g$ be the highest rational numbers for which $\Phi^*f=u^{m_f}\widetilde{f}$ and $\Phi^*f=u^{m_g}\widetilde{g}$, respectively, where $\widetilde{f},\, \widetilde{g} \in \mathbb{C}[u,x_0,x_1,x_2,x_3,x_4,\xi]$. Hence,  the inequalities $b_3/a_{\xi} \geq m_f$ and $b_4/a_{\xi} \geq m_g$ hold. We can well-form $T_0$ (and henceforth write it as $T$) and define $Y$ as 
\[
\begin{array}{cccc|ccccccc}
             &       & u  & \xi &   x_0 & x_1 & x_2 & x_3 & x_4 & \\
\actL{Y \colon (\widetilde{f}=\widetilde{g}=0)\subset T }   &  \lBr &  0 & a_{\xi} &   \iota_X & a_1 & a_2 & a_3 & a_4 &   \actR{.}\\
             &       & 1 & k &   \frac{k\iota_X-1}{a_{\xi}} & \frac{ka_1-\overline{ka_1}}{a_{\xi}} & \frac{ka_2-\overline{ka_2}}{a_{\xi}} & \frac{ka_3-b_3}{a_{\xi}} & \frac{ka_4-b_4}{a_{\xi}} &  
\end{array}
\] 

In particular, the bi-degrees of $\widetilde{f}$ and $\widetilde{g}$ are $\left(\begin{smallmatrix} d_1 \\ \frac{kd_1-m_fa_{\xi}}{a_{\xi}} \end{smallmatrix} \right)$ and $\left(\begin{smallmatrix} d_2 \\ \frac{kd_2-m_ga_{\xi}}{a_{\xi}} \end{smallmatrix} \right)$, respectively. It is easy to compute $-K_{T}$ and $-K_{\mathbb{P}(a_0,\ldots,a_{\xi})}$. Using the adjunction formula,
\begin{align*}
-K_Y &\sim \mathcal O \bigg(\sum{a_i} - (d_1+d_2), \frac{k(\sum{a_i} - (d_1+d_2))-1 - (b_3+b_4)+a_{\xi}(m_f+m_g)}{a_{\xi}}\bigg) \\ &\sim \mathcal O \bigg(\iota_X, \frac{k\iota_X-1 - (b_3+b_4)+a_{\xi}(m_f+m_g)}{a_{\xi}}\bigg), 
\end{align*}
and 
\[
\Phi^*(-K_{X}) \sim \mathcal O \bigg(\iota_X,\frac{k\iota_X}{a_{\xi}}\bigg).
\]
Therefore,
\[
-K_Y-\Phi^*(-K_{X})=\bigg(\frac{-(b_3+b_4)+a_{\xi}(m_f+m_g)}{a_{\xi}}-\frac{1}{a_{\xi}}\bigg)E
\]
where $E$ is the irreducible exception divisor of $\Phi$. This map restricts to the Kawamata blowup centred at $\mathbf{p_{\xi}}$ if and only if 
\[
(-b_2/a_{\xi}+m_f)+(-b_4/a_{\xi}+m_g)=0 
\]
Since both these terms are non-negative, it follows that $b_3=a_{\xi}m_f$ and $b_4=a_{\xi}m_g$ and, moreover, $Y$ is the proper transform of $X$ under $\Phi$. 
\end{proof}

\begin{Cor} \label{cor:genlift}
Let $\varphi \colon Y \rightarrow X$ be the Kawamata blowup of $X$ centred at $\mathbf{p_{\xi}} \in X$. Then $\Pic(Y) = \frac{\mathbb{Z}}{\iota_X}[-K_Y]+\frac{\mathbb{Z}}{\iota_X}[E]$. Moreover, $x_0$ lifts to an anticanonical section of $Y$ and 
\begin{align*}
x_i &\in H^0\bigg(Y,-\frac{a_i}{\iota_X}K_Y-\frac{\alpha_i}{\iota_X} E\bigg),\,\,  1\leq i\leq 5 
\end{align*} 
for some $\alpha_i \in \mathbb{Z}$ depending on the order of vanishing of $x_i$ at $E$.
\end{Cor}

\begin{proof}
Let $\mathbf{p_{\xi}}\sim \frac{1}{a_{\xi}}(a_0,a_1,a_2)\sim \frac{1}{a_{\xi}}(1,\overline{ka_1}, a_{\xi}-\overline{ka_1})$ be the cyclic quotient singularity we blowup. Notice that we are not assuming $a_{\xi} \geq a_i$. Recall that the discrepancy of the Kawamata blowup $\varphi \colon Y \rightarrow X$ is $\frac{1}{r}$, where $r$ is the index of the cyclic quotient singularity, i.e.,
\[
K_Y=\varphi^*K_X+\frac{1}{r}E.
\]

Let $\nu_i$ be the vanishing order of $x_i$ at $E$. Then,
\begin{align*}
x_i &\in H^0(Y,a_i\varphi^*(A)-\nu_iE) \\
		&= H^0\bigg(Y,\frac{a_i}{\iota_X}\varphi^*(-K_X)-\nu_iE\bigg) \\ 
    &= H^0\bigg(Y,\frac{a_i}{\iota_X}(-K_Y+\frac{1}{a_{\xi}}E)-\nu_iE\bigg)\\
    &= H^0\bigg(Y,-\frac{a_i}{\iota_X}K_Y+\frac{a_i-\iota_Xa_{\xi}\nu_i}{\iota_Xa_{\xi}}E\bigg).
\end{align*}

 Locally at $\mathbf{p_{\xi}}$, the local coordinates of the blowup $\varphi \colon Y \rightarrow X$ are given by
\[
\bigg(\frac{x_0^k}{{\xi}^{k_0}}:\frac{x_1^k}{{\xi}^{k_1}}:\frac{x_2^k}{{\xi}^{k_2}} \bigg) \in \mathbb{P}(1,\overline{ka_1}, a_{\xi}-\overline{ka_1})
\] 
where $k,\,k_i$ are uniquely defined positive integers by terminality of $\mathbf{p_{\xi}}$ and the exceptional divisor is $E \simeq\mathbb{P}(1,\overline{ka_1}, a_{\xi}-\overline{ka_1})$.
Hence, $x_0,\,x_1,\,x_2$ vanish at $E$ with orders $\frac{1}{a_{\xi}},\,\frac{\overline{ka_1}}{a_{\xi}},\,\frac{a_{\xi}-\overline{ka_1}}{a_{\xi}}$, respectively. Hence we have,
\begin{align*}
x_0 &\in H^0(Y,-K_Y)\\ 
x_1 &\in H^0\bigg(Y,-\frac{a_1}{\iota_X}K_Y+\frac{a_1-\iota_X\overline{ka_1}}{\iota_Xa_{\xi}}E\bigg)\\ 
x_2 &\in H^0\bigg(Y,-\frac{a_2}{\iota_X}K_Y+\frac{a_2-\iota_Xa_{\xi}+\iota_X\overline{ka_1}}{\iota_Xa_{\xi}}E\bigg).
\end{align*}
By the definition of $k$, see the proof of Lemma \ref{lem:utbl}, it follows that $\iota_X\overline{ka_1} \equiv a_1 \pmod{a_{\xi}}$. In other words, $a_1-\iota_X\overline{ka_1}$ is a multiple of $a_{\xi}$ and therefore there is $\alpha_1 \in \mathbb{Z}$ such that $\frac{a_1-\iota_X\overline{ka_1}}{\iota_Xa_{\xi}} = -\frac{\alpha_1}{\iota_X}$. Similarly, using the fact that $a_1+a_2 \equiv 0 \pmod{a_{\xi}}$, see Theorem \ref{thm:cqs}, it follows that $a_2-\iota_Xa_{\xi}+\iota_X\overline{ka_1} \equiv 0 \pmod{a_{\xi}}$ since $\iota_X\overline{ka_1} \equiv a_1 \pmod{a_{\xi}}$.

From Lemma \ref{lem:utbl}, it follows that $x_3$ and $x_4$ vanish at $E$ with order $m_f$ and $m_g$, respectively. Hence,
\begin{align*}
x_3 &\in H^0\bigg(Y,-\frac{a_3}{\iota_X}K_Y+\frac{a_3-a_{\xi}\iota_Xm_f}{\iota_Xa_{\xi}}E\bigg)\\
x_4 &\in H^0\bigg(Y,-\frac{a_4}{\iota_X}K_Y+\frac{a_4-a_{\xi}\iota_Xm_g}{\iota_Xa_{\xi}}E\bigg).
\end{align*}
We show that $a_{\xi}\iota_Xm_f \equiv a_3 \pmod{a_{\xi}}$ and therefore that there is $\alpha_3 \in \mathbb{Z}$ for which $x_3 \in H^0\Big(Y,-\frac{a_3}{\iota_X}K_Y-\frac{\alpha_3}{\iota_X}E\Big)$. By Lemma \ref{lem:awaybase}, $f$ has monomials of the form $x_0^{\beta_0}x_1^{\beta_1}$ where $\beta_0 a_0 + \beta_1 a_1 = d_1$ and $\beta_i \geq 0 $ are integers not both zero. Its pull-back via $\varphi$ vanishes at $E$ with order $\frac{\beta_0\overline{ka_0}+\beta_1\overline{ka_1}}{a_{\xi}}$. Notice that $\beta_0\overline{ka_0}+\beta_1\overline{ka_1} \equiv kd_1 \pmod{a_{\xi}}$. By definition of $m_f$, we have that $a_{\xi}m_f-\beta_0\overline{ka_0}+\beta_1\overline{ka_1} \equiv 0 \pmod{a_{\xi}}$. In particular, it follows that 
\[
\iota_X a_{\xi}m_f \equiv k\iota_Xd_1 \equiv d_1=ra_{\xi}+a_3 \equiv a_3 \pmod{a_{\xi}}.
\]
The reasoning for $x_4$ is completely analogous.

On the other hand, ${\xi}$ does not vanish at $E$, that is, $\nu_5 =0$ and 
\begin{align*}
{\xi} &\in H^0\bigg(Y,-\frac{a_{\xi}}{\iota_X}K_Y+\frac{1}{\iota_X}E\bigg). 
\end{align*}

It follows that $\widetilde{f} \ni {\xi}^{r_1}x_3$ and $\widetilde{g}\ni {\xi}^{r_2}x_4$ as in lemma \ref{lem:utbl} are divisors in the linear systems,
\[
\widetilde{f} \in \bigg|-\frac{d_1}{\iota_X}K_Y + \Big(\frac{d_1}{a_{\xi}\iota_X} -m_f \Big)E \bigg|, \qquad \widetilde{g} \in \bigg|-\frac{d_2}{\iota_X}K_Y + \Big(\frac{d_2}{a_{\xi}\iota_X} -m_g \Big)E \bigg|. 
\]
\end{proof}

\begin{Lem} \label{lem:anticanonicalwalliso}
Let $\varphi \colon Y \rightarrow X$ be the Kawamata blowup of $X$ centred at $\mathbf{p_{\xi}} \in X$. Suppose $Y\colon(\widetilde{f}=\widetilde{g}=0) \subset T$ is a Mori Dream space with $-K_Y \in \Int\overline{\Mov}(Y)$. If $x_1 \in H^0(X,-r_1K_X)$ lifts to a pluri-anticanonical section on $Y$ and 
\[
\widetilde{f} \in |-m_1K_Y|,\, \widetilde{g} \in |-m_2K_Y|
\]
for some positive integers $m_1$ and $m_2$, then there is a rational number $\epsilon > 0$ for which the ample models of $-K_Y+\epsilon E$ and $-K_Y-\epsilon E$ are isomorphic.
\end{Lem}
\begin{proof}
Let $X \subset \Proj \mathbb{C}[x_0,\ldots,x_4,\xi]$. We consider the Kawamata blowup $\varphi \colon E \subset Y \rightarrow X$ centred at $\mathbf{p_{\xi}}$. By Lemma \ref{cor:genlift} we know how the sections lift which depends on $X$. The exceptional divisor is the effective divisor $E \colon (u=0)$ in $Y$.  Let $y_i \in \{u,\xi,x_0,x_1,x_2,x_3,x_4 \}$ be the sections in $H^0(Y,-m_1K_Y+m_2E)$ for some $m_i>0$ and similarly, $z_i \in H^0(Y,-n_1K_Y-n_2E)$, $n_i >0$. The 3-fold $Y$ is contained in a toric variety $T$ from Lemma \ref{lem:utbl} whose Cox ring is isomorphic to $\mathbb{C}[u,x_0,\ldots,x_4,\xi]$. The bi-grading is known and given from the bi-grading of the sections $x_i,\,\xi$ of $Y$. Consider the ray $\mathbb{R}_+[x_0]$ in $\mathbb{R}^2$. Since $x_1$ lifts to a pluri-anticanonical section on $Y$ the ray above coincides with $\mathbb{R}_+[x_1]$.
 We play the 2-ray game on $T$. A small $\mathbb{Q}$-factorial modification 
 \[
        \begin{tikzcd}[ampersand replacement=\&, column sep = 2em]
             T_k   \ar[rr, dashed ] \ar[dr, swap, "\displaystyle{f_k}" ] \& {} \& T_{k+1} \ar[ld, "\displaystyle{g_k}" ] \\
             {} \& \mathcal{F}_k \& {}
        \end{tikzcd}
    \]

between ample models in the movable cone of $T$ is given by the the composition of small contractions $f_k$ and $g_k$. The irrelevant ideals of $T_k$ and $T_{k+1}$ are $(u,\xi,\ldots, y_i,\ldots) \cap (x_0,x_1, \ldots, z_i, \ldots)$ and $(u,\xi,\ldots, y_i,\ldots,x_0,x_1) \cap (\ldots, z_i, \ldots)$, respectively. In particular, 
\[
f_k \colon T_k \rightarrow \mathcal{F}_k,\, \mathcal{F}_k = \Proj \bigoplus_{m\geq 1} H^0(T, -mK_Y) \simeq \Proj \mathbb{C}[x_0,x_1, \ldots, y_i^{\alpha_i}z_i^{\beta_i},\ldots],
\]
where $\alpha_i,\beta_i >0$. Then, it is clear that $f_k$ contracts the exceptional set $E_k \colon (z_i=0) \subset T_k$ to $\Proj \mathbb{C}[x_0,x_1]$. Let $Y_k \subset T_k$ be the 3-fold obtained by restricting the 2-ray game on $T$ to $Y$. In particular it is defined by the same equations as $Y$. We claim that the locus contracted by $f_k$ is disjoint from $Y_k$.

By assumption $Y$ (and hence $Y_k$) is defined by the intersection of effective divisors $\widetilde{f}=0,\, \widetilde{g}=0$ in the pluri-anticanonical linear systems, $|-m_1K_Y|$ and $|-m_2K_Y|$, respectively. In particular, a monomial is in $f^*$ only if it is purely in the ideal $(x_0,x_1)$ or contains $y_i$ \emph{and} $z_i$. Hence,
\[
Y_k|_{E_k} \colon (\widetilde{f}|_{E_k}=\widetilde{g}|_{E_k} = 0) \subset T_k \quad \iff \quad Y_k|_{E_k} \colon (\widetilde{f}(x_0,x_1)=\widetilde{g}(x_0,x_1) = 0) \subset T_k.
\]  
Suppose $Y_k|_{E_k}$ is not empty. Then $ (f(x_0,x_1)=g(x_0,x_1)=0) \subset X$ is not empty which is a contradiction by Lemma \ref{lem:awaybase}. Hence, $Y_k$ and $Y_{k+1}$ are ample models in the same ample cone. It follows that they are isomorphic.
\end{proof}

\begin{Lem} \label{lem:discr}
Let $\sigma \colon X \dashrightarrow X'$ be an elementary Sarkisov link between $\mathbb{Q}$-Fano 3-folds initiated by a divisorial extraction $\varphi \colon E \subset Y \rightarrow X$. Then, there is a birational map $g \colon Y \dashrightarrow X'$ which is the composition of small $\mathbb{Q}$-factorial modifications followed by a divisorial contraction $\varphi' \colon E' \subset Y' \rightarrow X'$. If, in addition, $\varphi'$ is a divisorial contraction to a point it has discrepancy
\[
a= \frac{m_2}{n_2m_1-n_1m_2}
\]
where $m_i$ and $n_i$ are positive rational numbers such that $\varphi'^*(-K_{X'}) \sim -m_1K_Y-m_2E$ and $E' \sim -n_1K_Y-n_2E$.
\end{Lem}

\begin{proof}
By assumption, the Sarkisov link $\sigma$ is initiated by $\varphi \colon Y \rightarrow X$. By Lemma \ref{lem:linkmds} we have that $Y$ is a Mori Dream space and the chambers of $\overline{\Eff}(Y)$ correspond to the contraction $\varphi$ and the elementary contractions which define $g$. By the assumption that $X'$ is a $\mathbb{Q}$-Fano 3-fold it follows that $\overline{\Mov}(Y)$ is strictly contained in $\overline{\Eff}(Y)$. In particular, every movable divisor is big. Hence $g$ is a finite composition of small $\mathbb{Q}$-factorial modifications $\tau \colon Y \rat Y'$ in the Mori Category, followed by a $K_{Y'}$-negative contraction $\varphi \colon E' \subset Y' \rightarrow X'$. Moreover, the class of $-K_Y$ is in the interior of $\overline{\Mov}(Y)$. See also the paragraph after Corollary \ref{cor:exmds}.  


 Write
\[
-K_{Y'} \sim \varphi'^*(-K_{X'})-aE'
\]
where $E'$ is the exceptional divisor. By Lemma \ref{cor:genlift}, there are $m_i,\,n_i \in \frac{1}{\iota_X}\mathbb{Z}$ such that 
\[
\varphi'^*(-K_{X'}) \sim -m_1K_Y-m_2E \quad \text{and} \quad E' \sim -n_1K_Y-n_2E.
\] 
Hence,
\[
-K_{Y'} \sim (an_1-\lambda m_1)K_Y+(an_2-\lambda m_2)E.
\]
Since $Y$ and $Y'$ are connected by small birational maps the closure of their movable cones are equal and $-K_Y$ can be identified with $-K_{Y'}$. Under this identification and using the above,
\[
0 \sim (an_1-\lambda m_1+1)K_Y+(an_2-\lambda m_2)E.
\]
which is possible if and only if its coefficients are zero, that is, $an_1-\lambda m_1+1=an_2-\lambda m_2=0$. Solving for $a$ yields the discrepancy
\[
a= \frac{m_2}{n_2m_1-n_1m_2}.
\]
Notice that $n_2m_1-n_1m_2>0$ since $\varphi'^*(-K_{X'})$ is an effective divisor i.e., it is in the strict interior of $\overline{\Eff}(Y)$.

If $-K_Y \in \partial \overline{\Mov}(Y)$, then $-K_Y \sim_{\mathbb{Q}} M_2$ and $m_2=0$. It follows that we have a crepant blowup 
\[
-K_{Y'} \sim \varphi'^*(-K_{X'})
\]
to a strictly canonical singularity by definition. Otherwise $-K_Y \in \Int\overline{\Mov}(Y)$ and $m_2>0$, which in turn implies that the discrepancy of $\varphi'$ is strictly positive.
\end{proof}

\subsubsection{Linear Cyclic Quotient Singularities} \label{subsec:lin}

\begin{Def} \label{def:linear}
Let $(f,g)$ be the defining ideal of 
\[
X \subset \Proj \mathbb{C}[x_0,\ldots,x_4,\xi] \simeq \mathbb{P}(a_0,\ldots, a_4,a_{\xi})
\]
 (with no specified order of the weights) and suppose that there exists a change of variables for which $f$ and $g$ can be written as 
\begin{align*}
\xi x_3 + f_{d_1} &=0\\
\xi x_4 + g_{d_2} &=0
\end{align*}
where $f_{d_1},\,g_{d_2} \in \mathbb{C}[x_0,\ldots,x_4]$. If $\wt(\xi ) = a_{\xi} > 0$ the cyclic quotient singularity $\mathbf{p_{\xi}} = (0:\cdots : 1) \in X$ is said to be \textbf{linear} with respect to $X$.
\end{Def}

Interestingly, for each $X \in I$ there is at least one linear cyclic quotient singularity, see Lemma \ref{lem:polyform}. Moreover, blowing up a linear cyclic quotient singularity always initiates a Sarkisov link as we shall see. Notice that in \cite[Definition~2.1]{okadaI}, the author defines a \emph{distinguished singular point} to be a point for which there is a Sarkisov link (which is not an involution) to a Fano hypersurface or to a del Pezzo fibration. These distinguished singular points are precisely the linear cyclic quotient singularities in our notation. It turns out, however, that, in our case, not all distinguished points will be linear.

The following notation will be used throughout this paper. Suppose $\mathbf{p_{\xi}}\in X$ is \emph{linear} with respect to $X$. Then we can write the defining equations $(f,g) \in \mathbb{C}[x_0,x_1,x_2,x_3,x_4,\xi]$ of $X$ as 
\begin{align*}
\xi x_3 + f_{d_1} &=0\\
\xi x_4 + g_{d_2} &=0
\end{align*}
inside $\Proj \mathbb{C}[x_0,\ldots,x_4,\xi] \simeq \mathbb{P}(a_0,\ldots, a_4,a_{\xi})$ with no specified order of the weights. We define the $x_i$ with respect to $\xi$ in the following way.
\begin{itemize}
	\item By \cite{brownsuI} and \cite{brownsuII}, we have that $H^0(X,-K_X) \not =0$. In fact, it follows that in codimension 2 there are exactly two sections $x_0 \in H^0(X,-K_X)$ and $x_1 \in H^0(X,-mK_X)$ for some $m \in \mathbb{Z}_{>0}$ if $\iota_X \not = 4$ and $m \in \frac{1}{2}\mathbb{Z}_{>0}$ if $\iota_X = 4$.
	\item The point $\mathbf{p_{\xi}}$ is a cyclic quotient singularity of type $\frac{1}{a_{\xi}}(a_0,a_1,a_2)$ with local coordinates $x_0,\,x_1,\,x_2$. Since $x_0$ and $x_1$ have been assigned in the previous bullet point, the variable $x_2$ is defined to be the only local coordinate of $\mathbf{p_{\xi}}$ which is not $x_0$ or $x_1$.
	\item The variables $x_3$ and $x_4$ are the only variables for which $\xi x_3 \in f$ and $\xi x_4 \in g$. 
\end{itemize}
     
The following observation will be useful later. If $\mathbf{p_{\xi}} \sim \frac{1}{a_{\xi}}(1,a,b)$ is a terminal cyclic quotient singularity we have $a+b \equiv 0 \pmod{a_{\xi}}$. If, in addition, it is also linear we have $a_1+a_2 = a_{\xi}$ as can be seen from $\sum a_i-d_1-d_2= \iota_X$.



\begin{Lem} \label{lem:polyform}
Each $\mathbb{Q}$-Fano 3-fold $X \in I$ has a linear cyclic quotient singularity. Moreover, at least one of these has the highest possible index singularity in $X$.
\end{Lem}

\begin{proof}
Let $\mathbf{p_{\xi}}$ be the coordinate point $(0:\ldots : 0:1)$ and suppose that $a_{\xi} \geq a_i$. Then $2a_{\xi} \geq d_1$ and $2a_{\xi} \geq d_2$. By quasismoothness at $\mathbf{p_{\xi}}$ and since $X$ is not a linear cone, it follows that either $\xi^2 \in f$ (or $g$) or  ${\xi}x_3\in f$ and $\xi x_4\in g$. If $2a_{\xi} > d_1$ and $2a_{\xi} > d_2$, then     
\begin{align*}
f\colon \xi(x_3+\widetilde{f}(x_0,x_1,x_2,x_4))+f(x_0,x_1,x_2,x_3,x_4)&=0\\
g\colon \xi(x_4+\widetilde{g}(x_0,x_1,x_2,x_3))+g(x_0,x_1,x_2,x_3,x_4)&=0.
\end{align*}
The change of coordinates are therefore $ x_{3} \mapsto x_{3} - \widetilde{f}(x_0,x_1,x_2,x_4)$ followed by  $ x_4 \mapsto x_{4} - \widetilde{g}(x_0,x_1,x_2,x_3)$.

On the other hand, suppose that $2a_{\xi} = d_2$. If $d_1<d_2$, one can see by inspection that there is another monomial, say $x_4$, such that $2a_4=d_2$ and no others. By quasismoothness then either $\xi x_4 \in g$ and $(\xi+x_4)x_3 \in f$ and the conclusion follows or  $\xi^2+x_4^2 \in g$ and $(\xi+x_4)x_3 \in f$. The complex quadric $\xi^2+x_4^2$ has rank 2 and is therefore equivalent to $\xi x_4$.

If $d_1=d_2$ then we are in the case of family 91. Then $X:=X_{6,6} \subset \Proj \mathbb{C}[x,y,z,t,v,w] \simeq \mathbb{P}(1,2,2,3,3,3)$ is 
\begin{align*}
f \colon f_2(t,v,w)+l(t,v,w)f_3(x,y,z)+f_6(x,y,z)&=0\\
g \colon g_2(t,v,w)+l'(t,v,w)g_3(x,y,z)+g_6(x,y,z)&=0
\end{align*}
where $f_2,\,g_2$ are quadratic and $l,\,l'$ are linear forms. Consider the plane $\Pi \colon (x=y=z=0) \simeq \mathbb{P}^2$. If $X$ is quasismooth its singularities are all in $X \cap \Pi$ and are $4 \times \frac{1}{3}(1,1,2)$. Since $f_3$ and $g_3$ are not linear in any $x,\,y,\,z$ we have 
\[
J(X)|_{\Pi} = 
\begin{pmatrix}
0 & 0 & 0 & \partial_tf_2 & \partial_vf_2 & \partial_wf_2\\
0 & 0 & 0 & \partial_tg_2 & \partial_vg_2 & \partial_wg_2
\end{pmatrix}
\]
We write it in the reduced form,
\[
J(X)|_{\Pi} = 
\begin{pmatrix}
 \partial_tf_2 & \partial_vf_2 & \partial_wf_2\\
 \partial_tg_2 & \partial_vg_2 & \partial_wg_2
\end{pmatrix}.
\]
It is clear that $f_2$ and $g_2$ can have no non-trivial common factor $h \in \mathbb{C}[t,v,w]$. The idea is the same as in the proof of Lemma \ref{lem:awaybase}.


We claim that $f|_{\Pi}=0$ and $g|_{\Pi}=0$ are quadrics in $\mathbb{P}^2$ of rank at least two. Indeed, suppose that $f|_{\Pi}=0$ is a quadric of rank $0$ or $1$. Clearly if the rank is $0$, $X$ is not quasismooth. Then, $f|_{\Pi}=f_2(t,v,w)=w^2$ and $w=0$ and $X$ is not quasismooth at the two points $(g_2(t,v,w)=w=0) \subset \mathbb{P}^2$. The same conclusion follows from $g$. Hence we can write $f_2(t,v,w) = wv+\alpha t^2$ for some $\alpha \in \mathbb{C}^*$ and 
\[
J(X)|_{\Pi} = 
\begin{pmatrix}
 2\alpha t & w & v\\
 \partial_tg_2 & \partial_vg_2 & \partial_wg_2
\end{pmatrix}
\]
Suppose $g_2(t,v,w)=t^2+g_2'(t,v,w)$. Then $g_2 \not =0$ and at least one of its monomials contain $v$ or $w$ since $g_2=0$ defines a quadric of rank at leat 2. Also, since $g_2$ and $f_2$ have no non-trivial common factor, we can change variables to write 
\begin{align*}
f_2&=wv+\alpha t^2+f_2'(t,v,w) \\
g_2&=wt+g_2''(t,v,w) 
\end{align*}
and $\mathbf{p_{w}}=(0:0:0:0:0:1) \in X$ is a linear cyclic quotient singularity.

Suppose on the other hand that $g_2(t,v,w)=g_2(v,w)$. Then $g_2= vw$ and we claim that $\alpha \not = 0$. This is clear since in this case, after a row operation,
\[
J(X)|_{\Pi} = 
\begin{pmatrix}
 2\alpha t & 0 & 0\\
 0 &  w & v
\end{pmatrix}
\]
We can change variables to write 
\begin{align*}
f_2&=wt+f_2'(t,v,w) \\
g_2&=wv+g_2''(t,v,w) 
\end{align*}
and $\mathbf{p_{w}}=(0:0:0:0:0:1) \in X$ is a linear cyclic quotient singularity.
\end{proof}


As we mentioned before, of the whole basket of singularities of $X$ the linear ones are particularly nice from the point of view of birational geometry. The fact that all the codimension 2 Fano 3-folds admit such singularities is, in fact, only true for higher index Fano 3-folds. Moreover, for index 1 it is also true that the families with such singularities are not birationally rigid. See \cite{okadaI}.

\begin{Lem} \label{lem:lift}
Let $\varphi : Y \rightarrow X$ be the Kawamata blowup of $X$ centred at the linear cyclic quotient singularity $\mathbf{p_{\xi}} \in X$. Suppose $\iota_X \not =4$. Then, $x_0$ and $x_1$ lift to plurianticanonical sections of $Y$ and 
\[
x_i \in H^0\Big(Y,-\frac{a_i}{\iota_X}K_Y+\frac{\alpha_i}{\iota_X} E\Big),\,\,  2\leq i\leq 5 
\]
for some $\alpha_i \in \mathbb{Z}_{<0},\,\, 2 \leq i\leq 4$ and $\alpha_5 := a_{\xi}>0$. Moreover $Y$ is the intersection of two irreducible hypersurfaces $\widetilde{f}=0$ and $\widetilde{g}=0$ in the pluri-anticanonical linear systems 
\[
\widetilde{f} \in \Big|-\frac{d_1}{\iota_X}K_Y \Big|, \quad \widetilde{g}\in \Big|-\frac{d_2}{\iota_X}K_Y \Big|.
\]
\end{Lem}

\begin{proof}
Let $\mathbf{p_{\xi}}\sim \frac{1}{a_{\xi}}(a_0,a_1,a_2)\sim \frac{1}{a_{\xi}}(1,\overline{ka_1}, a_{\xi}-\overline{ka_1})$ be the linear cyclic quotient singularity we blowup. Notice that we are \emph{not} assuming $a_{\xi} \geq a_i$. Since $\mathbf{p_{\xi}}$ is terminal we have $a_1+a_2 \equiv 0 \pmod{a_{\xi}}$ but by linearity in fact $a_1+a_2=a_{\xi}$. Moreover, by definition of $k$, $k\iota_X \equiv 1 \pmod{a_{\xi}} $. Hence, $\overline{ka_1} = \frac{a_1}{\iota_X}$ since $\frac{a_1}{\iota_X} < a_{\xi}$. Notice that is an integer since $\iota_X \not = 4$. By Corollary \ref{cor:genlift}, $x_0$ is an anticanonical section of $Y$ and 
\begin{align*}
x_1 & \in H^0\bigg(Y,-\frac{a_1}{\iota_X}K_Y+\frac{a_1-\iota_X\overline{ka_1}}{\iota_Xa_{\xi}}E\bigg)    = H^0\bigg(Y,-\frac{a_1}{\iota_X}K_Y\bigg) \\
x_2 & \in H^0\bigg(Y,-\frac{a_2}{\iota_X}K_Y+\frac{a_2-\iota_Xa_{\xi}+\iota_X\overline{ka_1}}{\iota_Xa_{\xi}}E\bigg) = H^0\bigg(Y,-\frac{a_2}{\iota_X}K_Y-\frac{\iota_X-1}{\iota_X}E\bigg)
\end{align*} 
Recall that $x_3$ and $x_4$ vanish at $E$ with order $m_f$ and $m_g$, respectively. We focus on $m_f$ since the computation of $m_g$ is completely analogous. By quasismoothness of $X$, there are monomials in the ideal $(x_0,x_1)$ in $f$ (and $g$). When pulled back by the Kawamata blowup, the monomials $f(x_0,x_1)$ in $\overline{f}$ become $u^{d_1/\iota_Xa_{\xi}}f(x_0,x_1)$, where $u=0$ is the exceptional divisor of $\varphi$. By the proof of Lemma \ref{lem:polyform} (i.e., by the fact that $\mathbf{p_{\xi}}$ is linear) we can assume that $d_1<\iota_Xa_{\xi}$. Then, by the definition of $m_f$ in the proof of Lemma \ref{lem:utbl}, it follows automatically that $m_f=\frac{d_1}{\iota_Xa_{\xi}}$. Hence,
\[
x_3 \in H^0\bigg(Y,-\frac{a_3}{\iota_X}K_Y+\frac{a_3-a_{\xi}\iota_Xm_f}{\iota_Xa_{\xi}}E\bigg)   = H^0\bigg(Y,-\frac{a_3}{\iota_X}K_Y-\frac{1}{\iota_X}E\bigg) \\
\]
In the same way
\[
x_4 \in H^0\bigg(Y,-\frac{a_4}{\iota_X}K_Y-\frac{1}{\iota_X}E\bigg).
\]
Moreover,
\[
{\xi} \in H^0\bigg(Y, -\frac{a_{\xi}}{\iota_X}K_Y+\frac{1}{\iota_X}E\bigg).
\]
In particular, it follows from the computation of $m_f$ and $m_g$ that $\widetilde{f} \in \Big|-\frac{d_1}{\iota_X}K_Y\Big|$ and $\widetilde{g} \in \Big|-\frac{d_2}{\iota_X}K_Y\Big|$. 
\end{proof}

\begin{figure}%
\centering
\begin{tikzpicture}[scale=3]
  \coordinate (A) at (0, 0);
  \coordinate [label={left:$E$}] (E) at (0, 1);
  \coordinate [label={left:$-K_Y$}] (K) at (1, 1);

\coordinate [label={left:${\xi}$}] (5) at (0.5,1);
\coordinate [label={right:$x_2$}] (2) at (1.2,0.8);
\coordinate [label={right:$x_4$}] (4) at (1.2,0.5);
\coordinate [label={right:$x_3$}] (3) at (1.2,0.2);

  \draw  (A) -- (E);
  \draw  (A) -- (K);
	\draw [very thick,color=red](A) -- (5);
	\draw (A) -- (2);
	\draw (A) -- (3);
	\draw [very thick,color=red](A) -- (4);
  
\end{tikzpicture}
\caption{A representation of the chamber decomposition of the cone of effective divisors of $Y$. Notice that the rays corresponding to the characters $x_2,\,x_3,\,x_4$ may coincide.}%
\label{fig:chambdecomplin}%
\end{figure}
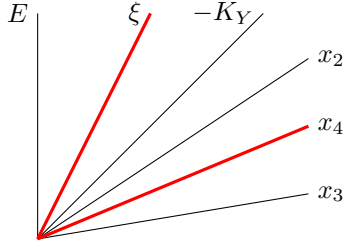


Suppose $\mathbf{p_{\xi}} \in X$ is a linear cyclic quotient singularity and let $Y\subset T$ be the proper transform of $X$ via the Kawamata blowup centred at $\mathbf{p}$, as in Lemma \ref{lem:utbl}. We now prove that the variation of GIT in $T$ with respect to the first wall crossing restricts to an isomorphism on $Y$.

\begin{Lem} \label{lem:iso}
Let $X \in I$ and suppose $\mathbf{p_{\xi}}$ is the germ of a linear cyclic quotient singularity. Let $ \Phi \colon T \rightarrow \mathbb{P}$ be the unique toric blowup centred at $\mathbf{p_{\xi}}$ and $Y \subset T$ the proper transform of $X$ via $\Phi$. Then the map $T \rat T_1$ arising from crossing the first wall in the movable cone of $T$ restricts to an isomorphism $Y \xrightarrow{\sim} Y_1$.
\end{Lem}

\begin{proof}
Suppose $\iota_X \not = 4$ or $X$ is a quasismooth member of family $X_{8,12} \subset \mathbb{P}(1,3,4,4,5,7)$. Then this is a direct consequence of Lemma \ref{lem:lift} and  Lemma \ref{lem:anticanonicalwalliso}.

Suppose on the other hand that $\iota_X=4$ except for the one family mentioned before. At the toric level we have the following diagram
 \[
        \begin{tikzcd}[ampersand replacement=\&, column sep = 2em]
             T_1   \ar[rr, dashed] \ar[dr, swap, "\displaystyle{f_1}" ] \& {} \& T_2 \ar[ld, "\displaystyle{g_1}" ] \\
             {} \& \mathcal{F} \& {}
        \end{tikzcd}
    \]
where $f_1$ contracts the locus $(x_1=x_2=x_3=x_4=0)$ and $g_1$ contracts the locus $(u=\xi=0)$. We have
\[
 \mathcal{F} = \Proj \bigoplus_{m\geq 0} H^0(T, mD) \simeq \Proj \mathbb{C}[x_0, \ldots, u_i,\ldots],
\]
where $D$ is the divisor given by $(x_0=0)$ and $\{u_i\}$ are monomials in $(u,{\xi}) \cap (x_1,x_2,x_3,x_4)$. It follows that $f_1$ contracts $(x_1=x_2=x_3=x_4={\xi}=0)$ to the point $p_0=(1:0:\cdots :0)$ in $\mathcal{F}$ and $g_1$ extracts from it $(u={\xi}=0)$. From quasismoothness of $X$ we have $x_0^{\alpha} \in f$. Hence, 
\[
x_0^{\alpha}u^{\frac{\iota_X\alpha}{\iota_X a_{\xi}}}= x_0^{\alpha}u^{\frac{d_1}{\iota_X a_{\xi}}}
\]
and $x_0^{\alpha} \in \widetilde{f}$ since $d_1 < \iota_X a_{\xi}$. Hence, the preimages of the restriction maps are empty in $Y_1$ and $Y_2$ and the conclusion follows.
\end{proof}

\section{Non Solid Families}  \label{sect:nonsolid}

In this section we prove that every quasismooth member $X$ of a deformation family in $I_{nS}$ admits a structure of a strict Mori fibre space, that is, $X$ is birational to a conic bundle over a normal surface or to a del Pezzo fibration $Y/\mathbb{P}^1$, establishing the first part of Theorem \ref{thm:main}. We use the following result in order to establish non-solidity of most deformation families in $I_{nS}$. 
\begin{Thm}{{\cite[Corollary~2.3]{liviatiago}}} \label{thm:nonsolid}
Let $X$ be a  $\mathbb{Q}$-factorial terminal Fano $d$-fold with $d\geq 3$ and $A \in \Cl(X)$ a generator of the Class group of $X$. Consider the embedding given by the ring of sections of $A$
\begin{equation*}
    X \hookrightarrow \mathbb{P}(a_0,\ldots , a_N)
\end{equation*}
 Suppose that there are $a_i,\, a_j$ such that  $l=\lcm(a_i,a_j) < \iota_X$.\ Then, $X$ is non-solid.\ 
\end{Thm}

\begin{Cor} \label{cor:nonsolid 2}
Let $X$ be a family in Table \ref{tab:dimlinsys}.\ Then $X$ is not solid.\
\end{Cor}

\begin{table}[ht]
    \centering
    \begin{tabular}[t]{c c c}
        \toprule
        Fam.                     & $\iota_X$  & $(a_i,a_j)$     \\ \midrule
        87   & 2 & $(1,1)$ \\
        88  & 2 & $(1,1)$\\
        89  & 2 & $(1,1)$\\
        90   & 2& $(1,1)$ \\
        \bottomrule
    \end{tabular}\hfill%
    \begin{tabular}[t]{c c c}
        \toprule
        Fam.                    & $\iota_X$  & $(a_i,a_j)$     \\ \midrule
        112   & 3 &  $(1,1)$ \\
        113   & 3& $(1,2)$ \\
        114  & 3 & $(1,2)$\\
        118  & 4 & $(1,2)$\\
        \bottomrule
    \end{tabular}\hfill%
    \begin{tabular}[t]{c c c}
        \toprule
        Fam.                     & $\iota_X$  & $(a_i,a_j)$  \\ \midrule
        119   & 4 & $(1,2)$ \\
        120   & 4 & $(1,3)$\\
        121  & 4 & $(1,3)$\\
        \bottomrule
    \end{tabular}
\caption{Families $X$ for which there are $(a_i,a_j)$ for which  $\lcm(a_i,a_j)<\iota_X$.}
\label{tab:dimlinsys}
\end{table}

We go further and construct for each $X$ in Table \ref{tab:dimlinsys} an explicit birational map to a strict Mori fibre space. 

\subsection{Conic Bundles}

In this section we prove that each quasismooth member of a family in the next table is birational to a conic bundle. This is done by constructing a Sarkisov link which is initiated by the Kawamata blowup of the respective cyclic quotient singularity.
\begin{center}
\begin{tabular}{lccccc} \toprule
    $X$ & 87 & 112 & 113& 118 & 119  \\
    $\mathbf{p}$  & $\frac{1}{3}(1,1,2)$  &$\frac{1}{5}(1,1,4)$ & $\frac{1}{4}(1,1,3)$ & $\frac{1}{3}(1,1,2)$ & $\frac{1}{7}(1,3,4)$ \\
		\bottomrule
		\end{tabular} 
\end{center}




\begin{Prop} \label{prop:elemcb}
Let $X$ be a quasismooth member of families 87 or 112. Then, there is an elementary Sarkisov link 
\[
\sigma \colon X \dashrightarrow Y/S
\]
to a conic bundle over $S=\mathbb{P}(\iota_X-1,1,1)$.

\end{Prop}

\begin{proof}
Any quasismooth member in each of these families can be written as 
\begin{align*}
\xi x_3 + f_{d_1}(x_0,\ldots, x_4) &=0\\
\xi x_4 + g_{d_2}(x_0,\ldots, x_4) &=0
\end{align*}
inside $\mathbb{P}:=\mathbb{P}(a_0,\ldots,a_4,a_{\xi})$ with homogeneous variables $x_0,\ldots,x_4,\xi$. Notice that $\iota_X$ is either 2 or 3, depending on whether $X$ is a member of family 87 or 112, respectively. Then $\mathbf{p_{\xi}} \in X$ is a linear cyclic quotient singularity of type $\frac{1}{a_{\xi}}(1,1,a_{\xi}-1)$ and local coordinates $x_0,\,x_1,\,x_2$. We consider the toric blowup $\Phi \colon T \rightarrow \mathbb{P}$ that restricts to the unique Kawamata blowup $\varphi \colon E \subset Y \rightarrow X$ centred at $\mathbf{p_{\xi}}$. Locally around $\mathbf{p_{\xi}}$ the map $\varphi$ is the graph of the rational map to $\mathbb{P}(1,1,a_{\xi}-1)$ given by
\[
\bigg(\frac{x_0^2}{\xi}:\frac{x_1^2}{\xi}:x_2^2 \bigg).
\]
By Lemma \ref{lem:lift}, we know exactly how the sections $x_i,\,\xi \in H^0(X,m_iA)$ of $X$ lift. For instance, the sections $x_0,\,x_1 \in H^0(X,-K_X)$ lift to anti-canonical sections of $Y$. Moreover, Lemma \ref{lem:lift} implies that, for $i \in \{2,3,4\}$, the sections  $x_i$ are all in $H^0(Y,mD)$ for some $m>0$ and effective divisor $D$ if and only if 
\[
a_4(\iota_X-1)=a_3(\iota_X-1)=a_2
\] 
and $D \sim -\frac{a_4}{\iota_X}K_Y-\frac{1}{\iota_X}E$. This is satisfied exactly by families 87 and 112. Hence, we can write the weight system of $T \supset Y$ as
\[
\begin{array}{cccc|ccccccc}
             &       & u  & \xi &   x_0 & x_1 & x_2 & x_3 & x_4 & \\
\actL{ T }   &  \lBr &  0 & a_{\xi} &   \iota_X & a_1 & a_2 & a_3 & a_4 &   \actR{}\\
             &       & 1 & 2 &   1 & 1 & 0 & 0 & 0 &  
\end{array}
\] 
where $E \colon (u=0)$ is the exceptional divisor of $\varphi$. The defining equations of $Y$ are $(\widetilde{f}=\widetilde{g}=0) \subset T$ and $\widetilde{f} \in |\frac{d_1}{\iota_X}K_Y|$ and $\widetilde{g} \in |\frac{d_2}{\iota_X}K_Y|$. In particular, every monomial purely in $(x_2,x_3,x_4)$, say $x_2^{\alpha_2}x_3^{\alpha_3}x_4^{\alpha_4} \in f$, lifts to $x_2^{\alpha_2}x_3^{\alpha_3}x_4^{\alpha_4}u^2 \in \widetilde{f}$. Hence we can write $\widetilde{f}$ and $\widetilde{g}$ as
\begin{align*}
\xi x_3 +\widetilde{f_{d}}(x_0,x_1)+ u^2\widetilde{f_{d}}(x_2,x_3,x_4)+ u(x_0+x_1)\widetilde{f_{d-\iota_X}}(x_2,x_3,x_4) &=0 \\
\xi x_4 +\widetilde{g_{d}}(x_0,x_1)+ u^2\widetilde{g_{d}}(x_2,x_3,x_4)+ u(x_0+x_1)\widetilde{g_{d-\iota_X}}(x_2,x_3,x_4) &=0.
\end{align*}
By Lemma \ref{lem:iso}, the first wall crossing $T \rat T'$, where $T'$ has the same Cox ring of $T$ but its irrelevant ideal is $(u,\xi,x_0,x_1) \cap (x_2,x_3,x_4)$, restricts to an isomorphism at the level of 3-folds. Moreover there is a map 
\[
\Phi' \colon T' \rightarrow \Proj\bigoplus_{l\geq 1} H^0(T',lD) \simeq \mathbb{P}(a_2,a_3,a_4) = \mathbb{P}(\iota_X-1,1,1)
\]
given by $(u,\xi,x_0,x_1,x_2,x_3,x_3) \mapsto (x_2,x_3,x_3)$. Each fibre is isomorphic to the 3-dimensional weighted projective space $\mathbb{P}(1,2,1,1)$ with homogeneous variables $u, \, \xi, \, x_0, \,x_1$. When we restrict $\Phi'$ to $Y \subset T'$ we have a fibration to $\mathbb{P}(\iota_X-1,1,1)$. Over a generic point it is possible to eliminate $\xi$ globally, and we end up with a conic given by 
\[
C\colon (\widetilde{f}(x_0,x_1)+\widetilde{g}(x_0,x_1)+u^2+u(x_0+x_1)=0) \subset \mathbb{P}^2.
\]
The claim follows.
\end{proof}

\begin{Prop} \label{prop:113cb}
Let $X$ be a quasismooth member of family 113. Then, there is a composition of elementary Sarkisov links to a conic bundle
\[
\sigma \colon X \rat Y'/\mathbb{P}^2.
\] 
\end{Prop}

\begin{proof}
By Theorem \ref{thm:cod2curve}, $X$ is birational to a special (singular) member of family 87, that we call $Z_{4,4}$, which has a $\frac{1}{3}(1,1,2)$ cyclic quotient singularity. Following the proof of Proposition \ref{prop:elemcb}, we see that it applies to $Z_{4,4}$ without any changes. This is mainly because all the key monomials are still in $Z_{4,4}$. We conclude that $X$ admits a conic bundle structure. 
\end{proof}

There are two families of index $4$ whose quasismooth members admit a conic bundle structure. These are families 118 and 119. 
\begin{Prop} \label{prop:118_119CB}
Let $X$ be a quasismooth member of families 118 or 119. Then, there is an elementary Sarkisov link 
\[
\sigma \colon X \dashrightarrow Y/S
\]
to a conic bundle over $S=\mathbb{P}(1,2,3)$.
\end{Prop}

\begin{proof}
Contrary to Proposition \ref{prop:elemcb} we do each family separately. This is not strictly necessary for the proof but the computations become clearer.
\paragraph{Family 118.} Let $X \subset \mathbb{P}(1,2,3,3,4,5):=\mathbb{P}$ with homogeneous variables $x,\,y,\,z,\,t,\,v,\,w$ be a quasismooth member of family 118. After a change of variables we can write $X$ as in Table \ref{tab:big} since, by quasismoothness, $f$ contains a rank 2 quadric in $t$ and $z$. Hence, the two $\frac{1}{3}(1,1,2)$ cyclic quotient singularities are the coordinate points $\mathbf{p_z}$ and $\mathbf{p_t}$ and at least one of them is linear. We construct a Sarkisov link initiated by blowing up the linear one which we assume w.l.o.g to be the $\mathbf{p_t}$ point. Let $\Phi \colon T \rightarrow \mathbb{P}$ be the toric blowup which restricts to the unique Kawamata blowup $\varphi \colon E \subset Y \rightarrow X$ centred at $\mathbf{p_t}$. Locally around $\mathbf{p_t}$ the Kawamata blowup is
\[
\bigg(\frac{v}{t},x,y\bigg) \in E \simeq \mathbb{P}(1,1,2).
\]
Then, by Corollary \ref{cor:genlift}, 
\[
v \in H^0(Y,-K_Y),\,\, x \in H^0\bigg(Y,-\frac{1}{4}K_Y-\frac{1}{4}E \bigg),\,\, y \in H^0\bigg(Y,-\frac{1}{2}K_Y-\frac{1}{2}E \bigg) 
\]
since these sections vanish at $E$ with orders $\frac{1}{3}$, $\frac{1}{3}$ and $\frac{2}{3}$ respectively. Hence, $m_f=1$ since by quasismoothness of $X$, $f$ contains at least one of the monomials $vy$ or $y^3$ and $m_g=\frac{2}{3}$ since, again by quasismoothness of $X$, $g$ contains the monomial $v^2$, where $m_f$ and $m_g$ are as in Lemma \ref{lem:utbl}. We can now also compute 
\[
w \in H^0\bigg(Y,-\frac{5}{4}K_Y-\frac{1}{4}E \bigg),\,\, z \in H^0\bigg(Y,-\frac{3}{4}K_Y-\frac{3}{4}E \bigg). 
\]
Hence, we can write $T$ as 
\[
\begin{array}{cccc|ccccccc}
             &       & u  & t &   v & w & x & y & z & \\
\actL{T}   &  \lBr &  0 & 3 &   4 & 5 & 1 & 2 & 3 &   \actR{.}\\
             &       & 1 & 1 &   1 & 1 & 0 & 0 & 0 &  
\end{array}
\] 
The movable cone of $T$ is generated by $\mathbb{R}_+[t=0]+\mathbb{R}_+[x=0]$, that is, by the rays $(1,3)$ and $(0,1)$ in $\mathbb{R}^2$. Adjacent Mori chambers in $\Mov(T)$ are related by small $\mathbb{Q}$-factorial modifications by \cite[Proposition~1.11]{mdsGIT}. We multiply the weight system of $T$ by
	\[
\begin{pmatrix}
    1        & -4 \\
    0      & 1  
\end{pmatrix} \in \SL_2(\mathbb{Z}).
\]
to get the isomorphic toric variety
\[
\begin{array}{cccc|ccccccc}
             &       & u  & t &   v & w & x & y & z & \\
\actL{T}   &  \lBr &  -4 & -1 &   0 & 1 & 1 & 2 & 3 &   \actR{.}\\
             &       & 1 & 1 &   1 & 1 & 0 & 0 & 0 &  
\end{array}
\] 
The irrelavant ideal of $T$ is $I=(u,t)\cap(v,w,x,y,z)$ as shown by the vertical bar in the matrix.

Choose a character in the interior of the GIT chamber generated by the rays $\mathbb{R}_+[t=0]$ and $\mathbb{R}_+[v=0]$, for instance, $(-1,2)$. Then we have 
\[
T=\Proj \bigoplus_{m\geq 1} H^0(T,m\mathcal{O}_T(-1,2))
\]
In the same way we define $T_1$
\[
T_1=\Proj \bigoplus_{m\geq 1} H^0(T,m\mathcal{O}_T(1,2)).
\]
It has the same Cox ring as $T$ but its irrelevant ideal is $I_1=(u,t,v)\cap(w,x,y,z)$. Hence there is a diagram
\[
        \begin{tikzcd}[ampersand replacement=\&, column sep = 2em]
             T   \ar[rr, dashed, "\displaystyle{\tau}" ] \ar[rd,swap, "\displaystyle{f}" ]   \&{} \&T_1 \ar[ld, "\displaystyle{g}" ]\\
						{}   \& \mathcal{F} \& {} 
        \end{tikzcd}
    \]
where, similarly,
\[
\mathcal{F}=\Proj \bigoplus_{m\geq 1} H^0(T,m\mathcal{O}_T(0,1))=\Proj \mathbb{C}[v,tw,tx,\ldots,u^3z^4].
\]
Notice that the ray $\mathbb{R}_+[v=0]$ separates the GIT chambers of $T$ and $T_1$. The map $f \colon T \rightarrow \mathcal{F}$ is then given in coordinates by 
\[
(u,t,v,w,x,y,z) \mapsto (v,tw,tx,t^2y,t^3z,\ldots,u^3z^4).
\]
This is a small contraction over $p_v=(1,0,\cdots,0) \in \mathcal{F}$. Indeed, 
\[
f^{-1}p_v=\{tw=tx=t^2y=t^3z=\cdots=u^3z^4 \} = \{u=t=0 \} \cup \{w=x=y=z=0 \}
\]
However any point of the form $(0,0,1,w_0,x_0,y_0,z_0)$ is contained in the zero set of the irrelavant ideal of $T$. Hence $f^{-1}p_v \simeq \mathbb{P}^1 \subset T$. Away from $p_v$ the map $f$ is one to one. We conclude that $f$ is a small contraction as we claimed. Similarly, $g^{-1}p_v \simeq \mathbb{P}(1,1,2,3)$ and is one to one otherwise. The map $\tau \colon T \rat T_1$ is the composition $\tau = g^{-1} \circ f$ and we denote it by $(-4,-1,1,1,2,3)$ following \cite{brown}.  As before, we multiply  the matrix $T_1$ by
	\[
\begin{pmatrix}
    1        & -1 \\
    0      & 1  
\end{pmatrix} \in \SL_2(\mathbb{Z})
\]
and get the isomorphic toric variety 
\[
\begin{array}{ccccc|cccccc}
             &       & u  & t &   v & w & x & y & z & \\
\actL{T_1}   &  \lBr &  -5 & -2 &  -1 & 0 & 1 & 2 & 3 &   \actR{.}\\
             &       & 1 & 1 &   1 & 1 & 0 & 0 & 0 &  
\end{array}
\] 
Similarly, we have a diagram
\[
        \begin{tikzcd}[ampersand replacement=\&, column sep = 2em]
             T_1   \ar[rr, dashed, "\displaystyle{\tau_1}" ] \ar[rd,swap, "\displaystyle{f_1}" ]   \&{} \&T_2 \ar[ld, "\displaystyle{g_1}" ]\\
						{}   \& \mathcal{F}_1 \& {} 
        \end{tikzcd}
    \]
where 
\[
\mathcal{F}_1=\Proj \bigoplus_{m\geq 1} H^0(T_1,m\mathcal{O}_{T_1}(w=0))=\Proj \mathbb{C}[w,vx,ty,\ldots,u^3z^5].
\]
In this case $\tau_1=g_1^{-1} \circ f_1$ is the small $\mathbb{Q}$-factorial modification $(-5,-2,-1,1,2,3)$. As can be seen in Figure \ref{fig:cb118}, the cones $\Eff(T)$ and $\Mov(T)$ share the boundary ray generated by $(1,0)$. Hence we have a fibration, which, in this case, is 
\[
\Phi' \colon T_2 \rightarrow \mathcal{F}_3,\quad \mathcal{F}_3 = \Proj \bigoplus_{m\geq 1} H^0(Y,m\mathcal{O}(1,0))=\mathbb{P}(1,2,3)
\]
written in coordinates as 
\[
(u,t,v,w,x,y,z) \mapsto (x,y,z).
\]

\begin{figure}%
\centering
\begin{tikzpicture}[scale=3,font=\small]
  \coordinate (A) at (0, 0);
  \coordinate [label={left:$u$}] (E) at (-0.7, 0.25);
  \coordinate [label={above:$\mathbb{P}$}] (K) at (-0.5,0.5);
	\coordinate [label={above:$\mathcal{F}$}] (5) at (0, 0.6);
	\coordinate [label={right:$\mathcal{F}_1$}] (2) at (0.5,0.5);
	\coordinate [label={right:$\mathbb{P}(1,2,3)$}] (4) at (0.5,0);
	
		\coordinate [label={$T$}] (T) at (-0.25,0.5);
		\coordinate [label={$T_1$}] (T1) at (0.25,0.5);
		\coordinate [label={$T_2$}] (T2) at (0.5,0.2);
  \draw  (A) -- (E);
  \draw [very thick,color=red] (A) -- (4);
	\draw (A) -- (5);
	\draw (A) -- (2);
	\draw [very thick,color=red](A) -- (K);
\end{tikzpicture}
\caption{A representation of the chamber decomposition of the cone of effective divisors of $T$. The red subcone is the cone of movable divisors of $T$.}%
\label{fig:cb118}%
\end{figure}
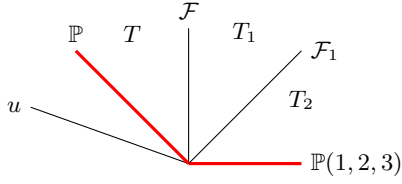

By adjunction $-K_T \sim 2\mathcal{O}(1,2)$ and therefore, $-K_T \in \Int\Mov(T)$. Hence, $\Phi'$ is a $K_T$-negative contraction. Its fibres are each isomorphic to $\mathbb{P}^3$. We have therefore the following diagram,
\[
        \begin{tikzcd}[ampersand replacement=\&, column sep = 2em]
      \& T \ar[rr, dashed, "\displaystyle{\tau}" ]\ar[rd,swap, "\displaystyle{f}" ] \ar[dl,swap, "\displaystyle{\Phi}" ]  \& \& T_1 \ar[dl, "\displaystyle{g}" ] \ar[rd,swap, "\displaystyle{f_1}" ] \ar[rr, dashed, "\displaystyle{\tau_1}" ]\& \& T_2 \ar[ld, "\displaystyle{g_1}" ] \ar[rd, "\displaystyle{\Phi'}" ] \&\\
\mathbb{P} \&   \& \mathcal{F}  \&  \& \mathcal{F}_1 \&  \& \mathbb{P}(1,2,3)
        \end{tikzcd}
    \] 
We want to restrict the above diagram to $X \in \mathbb{P}$. Recall that, by construction of $\Phi$, this restricts to the unique Kawamata blowup $\varphi \colon Y \rightarrow X$. The defining equations of $Y \colon (\widetilde{f}=\widetilde{g}=0) \subset T$ are such that $\widetilde{f} \in |-\frac{3}{2}K_Y-\frac{1}{2}E|$ and $\widetilde{g} \in |-2K_Y|$. By Lemma \ref{lem:iso}, the map $T \rat T_1$ restricts to an isomorphism on $Y$. On the other hand, since $wx \in f$ and $wt \in g$ by quasismoothness of $X$ at $\mathbf{p_w}$ and $wx \in |-\frac{3}{2}K_Y-\frac{1}{2}E|$ and $wt \in |-2K_Y|$ it follows that $wx \in \widetilde{f}$ and $wt \in \widetilde{g}$. Hence, locally around this point, we can eliminate $x$ and $t$. The restriction is a flip of type $(-5,-1,2,3)$. Let $Y_2$ be the resulting 3-fold. The fibration $\Phi'\colon T_2 \rightarrow \mathbb{P}(1,2,3)$ restricts to a fibration $\varphi'$ on $Y_2$ whose generic fibres are smooth conics in $\mathbb{P}^2$. That is, we have the following Sarkisov link,
\[
        \begin{tikzcd}[ampersand replacement=\&, column sep = 3em]
      \& Y \ar[rr, "\displaystyle{\simeq}" ] \ar[dl,swap, "\displaystyle{\Phi}" ]  \& \& Y_1   \ar[rr, dashed, "\displaystyle{(-5,-1,2,3)}" ]\& \& Y_2  \ar[rd, "\displaystyle{\Phi'}" ] \&\\
X \&   \&   \&  \&  \&  \& \mathbb{P}(1,2,3)
        \end{tikzcd}
    \]

Notice that the base of the fibration is a surface with two du Val singularities of types $A_1$ and $A_2$. Indeed by \cite[Theorem~1.2.7]{MoriProkConicsI}, the base of a $\mathbb{Q}$-conic bundle can have at most type $A$ du Val singularities.

\paragraph{Family 119.} We proceed analogously as in the case of family 118. Let $X \subset \mathbb{P}(1,2,3,4,5,7):=\mathbb{P}$ with homogeneous variables $x,\,y,\,z,\,t,\,v,\,w$ be a quasismooth member of family 119. Any such member can be written as in Table \ref{tab:big}. Let $\Phi \colon T \rightarrow \mathbb{P}$ be the toric blowup which restricts to the unique Kawamata blowup $\varphi \colon E \subset Y \rightarrow X$ centred at $\mathbf{p_w}$. Locally around $\mathbf{p_w}$ the Kawamata blowup is the graph of the rational map to $\mathbb{P}(1,3,4)$ given by
\[
\bigg(\frac{t^2}{w},\frac{v^2}{w},y^2\bigg) \in E \simeq \mathbb{P}(1,3,4).
\]
Then, by Corollary \ref{cor:genlift}, 
\[
t \in H^0(Y,-K_Y),\,\, v \in H^0\bigg(Y,-\frac{5}{4}K_Y-\frac{1}{4}E \bigg),\,\, y \in H^0\bigg(Y,-\frac{1}{2}K_Y-\frac{1}{2}E \bigg). 
\]
since these sections vanish at $E$ with orders $\frac{1}{7}$, $\frac{3}{7}$ and $\frac{4}{7}$ respectively. Hence, $m_f=\frac{2}{7}$ since by quasismoothness of $X$, $f$ contains the monomial $t^2$ and $m_g=\frac{6}{7}$ since, again by quasismoothness of $X$, $g$ contains the monomial $v^2$, where $m_f$ and $m_g$ are as in Lemma \ref{lem:utbl}. We can now also compute 
\[
x \in H^0\bigg(Y,-\frac{1}{4}K_Y-\frac{1}{4}E \bigg),\,\, z \in H^0\bigg(Y,-\frac{3}{4}K_Y-\frac{3}{4}E \bigg). 
\]
Hence, we can write $T$ as 
\[
\begin{array}{cccc|ccccccc}
             &       & u  & w &   t & v & x & y & z & \\
\actL{T}   &  \lBr &  0 & 7 &   4 & 5 & 1 & 2 & 3 &   \actR{.}\\
             &       & 1 & 2 &   1 & 1 & 0 & 0 & 0 &  
\end{array}
\] 
We play a 2-ray game on $T$ and get the following diagram
\[
        \begin{tikzcd}[ampersand replacement=\&, column sep = 2em]
      \& T \ar[rr, dashed, "\displaystyle{\tau}" ]\ar[rd,swap, "\displaystyle{f}" ] \ar[dl,swap, "\displaystyle{\Phi}" ]  \& \& T_1 \ar[dl, "\displaystyle{g}" ] \ar[rd,swap, "\displaystyle{f_1}" ] \ar[rr, dashed, "\displaystyle{\tau_1}" ]\& \& T_2 \ar[ld, "\displaystyle{g_1}" ] \ar[rd, "\displaystyle{\Phi'}" ] \&\\
\mathbb{P} \&   \& \mathcal{F}  \&  \& \mathcal{F}_1 \&  \& \mathbb{P}(1,2,3)
        \end{tikzcd}
    \] 
by multiplying the matrix of $T$ by the same two matrices we used for family 118. In this case the small $\mathbb{Q}$-factorial modifications $\tau$ and $\tau_1$ are $(-4,-1,1,1,2,3)$ and $(-5,-3,-1,1,2,3)$, respectively. The map $\Phi' \colon T_1 \rightarrow \mathbb{P}(1,2,3)$ is a fibration given in coordinates by 
\[
(u,w,t,v,x,y,z) \mapsto (x,y,z).
\] 
Each fibre is isomorphic to $\mathbb{P}(1,1,1,2)$. We restrict this diagram to $X$. The map $\Phi$ restricts to the Kawamata blowup $\varphi \colon E \subset Y \rightarrow X$ centred at $\mathbf{p_w}$ by construction. Moreover, by the proof of Corollary \ref{cor:genlift}, we have $\widetilde{f} \in |-2K_Y|$ and $\widetilde{g} \in |-\frac{5}{2}K_Y-\frac{1}{2}E|$. Therefore $t^2 \in \widetilde{f}$ and $v^2 \in \widetilde{g}$. Hence, both $\tau$ and $\tau_1$ restrict to isomorphisms on $Y$. On the other hand, the restriction of $\Phi'$ to $Y \subset T_2$ is a map $\varphi'$ whose generic fibre is a conic in $\mathbb{P}^2$ with variables $v,\,u,\,t$. Hence we have the Sarkisov link,
\[
        \begin{tikzcd}[ampersand replacement=\&, column sep = 3em]
      \& Y \ar[rr, "\displaystyle{\simeq}" ] \ar[dl,swap, "\displaystyle{\varphi}" ]  \& \& Y_1   \ar[rr, "\displaystyle{\simeq}" ]\& \& Y_2  \ar[rd, "\displaystyle{\varphi'}" ] \&\\
X \&   \&   \&  \&  \&  \& \mathbb{P}(1,2,3)
        \end{tikzcd}
    \]  
\end{proof}

The conclusion is the following:

\begin{Thm} \label{thm:cb}
Let $X$ be a quasismooth member of a family in $I_{Cb}$. Then, there is a Sarkisov link 
\[
\sigma \colon X \dashrightarrow Y/S
\]
to a conic bundle over $S$ where 
\[
S = \begin{cases*}
                    \mathbb{P}(\iota_X-1,1,1) & if  $X \in \{87,112,113 \}$  \\
                     \mathbb{P}(1,2,3) & if $X \in \{118,119 \}$
                 \end{cases*}
\]
In particular, the base of the conic bundle has only isolated $A_1$ or $A_2$ du Val singularities.
\end{Thm}

\subsection{del Pezzo Fibrations}
Let $X$ be a quasismooth member of a family in $I_{dP}$. We prove that blowing up the cyclic quotient singularity in $X$ as in the next table initiates a Sarkisov link to a del Pezzo fibration.

\vspace{-0.5cm}

\begin{center}
\resizebox{\textwidth}{!}{\begin{tabular}{lccccccccccccc} \toprule
    $X$ & 88 & 89  & 90 & 91 & 103 & 114 & 116 &120 & 121 & 122 & 123 & 124 & 125\\
    $\mathbf{p}$  & $\frac{1}{3}(1,1,2)$  &$\frac{1}{5}(1,1,4)$ & $\frac{1}{5}(1,2,3)$ & $\frac{1}{3}(1,1,2)$  & $\frac{1}{7}(1,2,5)$ & $\frac{1}{5}(1,1,4)$ &$\frac{1}{7}(1,1,6)$ &$\frac{1}{7}(1,1,6)$ &$\frac{1}{7}(1,2,5)$   & $\frac{1}{7}(1,3,4)$ & $\frac{1}{9}(1,4,5)$ &  $\frac{1}{7}(1,3,4)$ & $\frac{1}{8}(1,1,7)$ \\ 
		\bottomrule
		\end{tabular}}
\end{center}



\paragraph{Families 88, 90, 91, 103, 114, 116 and 120.} Let $X$ be a quasismooth member in any one of these families. Then, we can write $X$ as
\begin{align*}
\xi x_3 + f(x_0,x_1) + x_4^{\alpha} l_{\mu} + f(x_3,x_2)+f_{d_1}&=0\\  
\xi x_4 + g(x_0,x_1) + g(x_3,x_2)+ g_{d_2}&=0
\end{align*}
where
\begin{itemize}
  \item $X\subset \mathbb{P}(a_0,a_1,\ldots,a_4, a_{\xi}) =: \mathbb{P}$, with no order on the weights and homogeneous variables $x_0,\,x_1,\ldots, x_4,\,\xi$;
	\item $x_{\mu}$ is either $x_4$ or $x_2$; In the second case we have $x_4^2x_3^{\beta} \in g_{d_2}$ for some $\beta$.
	\item $f_{d_1}, \, g_{d_2} \in \mathbb{C}[x_0, \ldots x_4]$ and contain no monomials purely in $(x_0,x_1)$ nor in $(x_3,x_2)$; 
\end{itemize}
Then $\mathbf{p_{\xi}} \in X$ is a linear cyclic quotient singularity of type $\frac{1}{a_{\xi}}(1,\overline{ka_1},\overline{ka_2})$ w.r.t. $X$. Notice that for each of these families we have 
\[
a_3(\iota_X-1) = a_2, \quad a_4 > a_3.
\]
As usual, we consider the toric blowup $\Phi \colon T \rightarrow \mathbb{P}$ whose restriction to $X$ is the unique Kawamata blowup $\varphi \colon E \subset Y \rightarrow X$ centred at $\mathbf{p_{\xi}}$. We know how the sections $H^0(X,m_iA)$ of $X$ lift to $Y$ by Lemma \ref{lem:lift}. Notice that in particular,
\[
x_3 \in H^0\bigg(Y,-\frac{a_3}{\iota_X}K_Y-\frac{1}{\iota_X}E \bigg),\, x_2 \in H^0\bigg(Y,-\frac{a_2}{\iota_X}K_Y-\frac{\iota_X-1}{\iota_X}E \bigg)
\]
span the same ray in $\mathbb{R}^2$ if and only if $a_3(\iota_X-1) = a_2$ which is the case for $X$ as we noticed previously. Hence we can write the weight system of $T$ in clockwise order as
		\[
\begin{array}{cccc|ccccccc}
             &       & u  & \xi &   x_0 & x_1 & x_4 & x_3 & x_2 & \\
\actL{T}   &  \lBr &  0 & a_{\xi} &   \iota_X & a_1 & a_4 & a_3 & a_3(\iota_X-1) &   \actR{}\\
             &       & 1 & k &   \frac{k\iota_X-1}{a_{\xi}} & \frac{k\iota_X-1}{a_{\xi}}\cdot \frac{a_1}{\iota_X} & \frac{a_4(k\iota_X-1)-a_{\xi}}{a_{\xi}\iota_X} & \frac{a_3(k\iota_X-1)-a_{\xi}}{a_{\xi}\iota_X} & \frac{a_3(k\iota_X-1)-a_{\xi}}{a_{\xi}\iota_X}\cdot (\iota_X-1) &  
\end{array}
\] 
as in the proof of Lemma \ref{lem:utbl}. Hence,
\[
\mathbb{R}_+[\xi]+\mathbb{R}_+[x_3]=\overline{\Mov}(T) \subset \overline{\Eff}(T) = \mathbb{R}_+[u]+\mathbb{R}_+[x_3].
\]
The movable cone of $T$ is subdivided into three chambers, the first of which is separated by the ray generated by $\mathbb{R}_+[x_0]=\mathbb{R}_+[x_1]$. Lemma \ref{lem:iso} implies that crossing this wall restricts to an isomorphism on $Y$, so we can assume $Y \subset T_1$ where $T_1$ has the same Cox ring of $T$ but its irrelevant ideal is $(u,\xi,x_0,x_1) \cap (x_4,x_3,x_2)$.

\begin{Lem} \label{lem:flip1}
There is a small modification $ T_1 \rat T_2$ over the point $\mathbf{p_4}$. It replaces the $3$-fold 
\[
\mathbb{P}\bigg(a_4,\frac{d_2}{\iota_X},1,\frac{a_1}{\iota_X}\bigg) \subset T_1
\]
 with the projective line  
\[
\mathbb{P}\bigg(\frac{a_4-a_3}{\iota_X},\frac{a_4-a_3}{\iota_X}(\iota_X-1)\bigg) \subset T_2. 
\]
\end{Lem}

\begin{proof}
Let $a_4':=\frac{ka_4-m_2}{a_{\xi}} = \frac{a_4(k\iota_X-1)-a_{\xi}}{a_{\xi}\iota_X} \in \mathbb{Z}$. Since $\gcd(a_4,a_4') = 1$ (for instance, from direct verification), there are integers $r,\,s$ for which $ra_4+sa_4'=1$. Consider the matrix
	\[
A=\begin{pmatrix}
    r        & s \\
    -a_4'       & a_4  
\end{pmatrix} \in \SL_2(\mathbb{Z}).
\]
Then,
\[
A\cdot \begin{pmatrix}
    a_4          \\
    a_4'        
\end{pmatrix} = \begin{pmatrix}
    1         \\
    0     
\end{pmatrix}
\]
This means that, by restricting to the open set $(x_4 \not = 0)$, we can eliminate the action given by the first row of $A\cdot T$ and that the small contraction $\tau_1 \colon T_1 \rat T_2$ is given by its second row. We have the following diagram
\[
        \begin{tikzcd}[ampersand replacement=\&, column sep = 2em]
             T_1   \ar[rr, dashed ] \ar[rd,swap, "f_1" ]  \&  \&T_2 \ar[ld, "g_1" ]\\
						{} \& p_{x_4} \in \mathcal{F}_1  \& {} 
        \end{tikzcd}
    \]
where 
\begin{align*}
\mathcal{F}_1 &= \Proj\bigoplus_{m\geq 1} H^0(T_1,m\mathcal{O}(x_4))\\
&=\Proj \mathbb{C}[x_4, \ldots, u_iv_i, \ldots]
\end{align*}
where each $u_i$ is a monomial in $(u,\xi,x_0,x_1)$ and each $v_i$ is a monomial in $(x_3,x_2)$. The contractions $f_1$ and $g_1$ are given in coordinates by
\[
(u,\xi,x_0,x_1,x_4,x_3,x_1) \mapsto (x_4, \ldots, u_iv_i, \ldots).
\]
Thus, $f_1$ contracts the locus $(x_3=x_2=0) \subset T_1$ to $p_{x_4} = (1:0\cdots:0) \in \mathcal{F}_1$ and $g_2$ contracts the locus $(u=\xi=x_0=x_1=0) \subset T_2$ to $p_{x_4}$. Easy but tedious computations show that the contracted loci are isomorphic to
\[
\mathbb{P}\bigg(a_4,\frac{d_2}{\iota_X},1,\frac{a_1}{\iota_X}\bigg) \subset T_1 \quad \text{and} \quad \mathbb{P}\bigg(\frac{a_4-a_3}{\iota_X},\frac{a_4-a_3}{\iota_X}(\iota_X-1)\bigg) \subset T_2,
\]
respectively.  Hence, the small modification $\tau_1$ is 
\[
\bigg(a_4,\frac{d_2}{\iota_X},1,\frac{a_1}{\iota_X}, \frac{a_3-a_4}{\iota_X},\frac{a_3-a_4}{\iota_X}(\iota_X-1)\bigg).
\] 
\end{proof}

\begin{Cor} \label{cor:dpflip}
There is a flip $Y\simeq Y_1 \rat Y_2$ over the point $\mathbf{p_4}$. Moreover, this is a toric flip if and only if there is a monomial $x_4^{\alpha} \in f$ where $\alpha=\iota_X$ and, otherwise, it is a hypersurface flip.
\end{Cor}

\begin{proof} 
We restrict the construction made in Lemma \ref{lem:flip1} to the $3$-fold $Y$. Recall that $\widetilde{f} \in |-\frac{d_1}{\iota_X}K_Y|$ and $\widetilde{g} \in |-\frac{d_2}{\iota_X}K_Y|$. Hence, in a neighborhood of $(x_4\not =0)$, where the small modification $\tau_1$ happens, one can eliminate the variable $\xi$ because $\xi x_4 \in \widetilde{g}$. Hence, $\mathbb{P}(a_4,1,\frac{a_1}{\iota_X}) \subset T_1$ with homogeneous variables $u,\,x_0,\,x_1$ is cut down by 
\begin{equation} \label{eq:cut1}
f(x_0,x_1)+x_{\mu} u^{\beta} =0
\end{equation}
where $x_{\mu}$ is either $0$ or $1$.  In order to have a toric flip, we need to be able to eliminate another variable in $Y$. The only way to do so is if there is $\alpha$ such that $x_4^{\alpha} \in f$ which lifts to $x_4^{\alpha}u \in \widetilde{f}$. Indeed, if $x_4^{\alpha} \in f$, it lifts to
\[
	x_4^{\alpha} u^{\alpha m_g-m_f} = x_4^{\alpha} u^{\frac{\alpha d_2-d_1}{\iota_X a_{\xi}}}. 
	\]
	That is, we have that $x_4^{\alpha}u \in \widetilde{f}$ if and only if $\alpha d_2 - d_1 =  \iota_X a_{\xi}$. Using the fact that $\alpha a_4 = d_1 = a_{\xi}+a_3$ we get that
	\[
	x_4^{\alpha}u \in \widetilde{f} \quad \iff \quad \iota_X a_4 = a_{\xi}+a_3 \quad \iff \quad \alpha=\iota_X.
	\]
	Hence, equation \ref{eq:cut1} becomes
	\[
	f(x_0,x_1)+x_{4}^{\iota_X} u =0
	\]
Using the inverse function theorem, we can eliminate $u$ in an analytic neighborhood of $p_{x_4}$ and the flip becomes
\[
        \begin{tikzcd}[ampersand replacement=\&, column sep = 2em]
            T_1 \supset Y   \ar[rrrrrr, dashed, "\displaystyle{\bigg(1,\frac{a_1}{\iota_X},\frac{a_3-a_4}{\iota_X}, \frac{a_3-a_4}{\iota_X}(\iota_X-1)\bigg)}" ] \ar[rrrd ]  \& {} \&{} \&{} \&{}  \&{} \&Y_2 \subset T_2 \ar[llld ]\\
						{} \&{} \&{}  \& p_{x_4} {} \& {} \& {} 
        \end{tikzcd}
    \]

We have proven that $x_4^{\iota_X} \in f$ is a necessary condition in order to have a flip in this setting. This is also sufficient by \cite[Theorem~7]{brown}: Indeed, this small $\mathbb{Q}$-factorial modification is terminal if and only if $a_4-a_3=\iota_X$, which is always the case and we can conclude that it is a flip by the classification of toric flips in \cite[Theorem~7]{brown}.

If $x_4^{\alpha} \not \in f$, then we prove that we have a hypersurface flip. In this case $x_{\mu} = x_2$ and, since $\widetilde{f} \in |-\frac{d_1}{\iota_X}K_Y|$, we have 
\[
x_4^{\alpha}x_2u^{\beta} \in \bigg|-\frac{d_1}{\iota_X}K_Y-\bigg(\frac{\alpha+1}{\iota_X}- \beta \bigg)E\bigg|.
\]
Since we know the degrees of the equations defining $X$, we can see that $\alpha=1$. Therefore $x_4^{\alpha}x_2u^{\beta} \in \widetilde{f}$ if and only if $\beta = 1$. As we have seen above, the fibre inside $Y_1$ above $p_{x_4}$ is given by $\widetilde{f}(x_0,x_1)=0$ which has degree $d_1/\iota_X$ and contains monomials $x_0^{\alpha_1}x_1^{\alpha_2}$ where not both $\alpha_i$ are zero. In this case we can not eliminate any other variables from $\widetilde{f}$ in an analytic neighborhood of $(x_4 \not = 0)$. Indeed we remain with the equation
\[
(x_2u+f(x_0,x_1)+g(x_0,x_1)+ u^{\beta}f(x_3,x_2)+\cdots = 0 ) \subset \mathbb{C}^5_{ux_0x_1x_3x_2}.
\]
These have weights $\wt(u,x_0,x_1,x_3,x_2) = \Big(a_4,1,\frac{a_1}{\iota_X},\frac{a_3-a_4}{\iota_X}, \frac{a_3-a_4}{\iota_X}(\iota_X-1)\Big)$ . Hence, we a have a hypersurface flip as in the following diagram 
\[
        \begin{tikzcd}[ampersand replacement=\&, column sep = 2em]
             Y   \ar[rrrrrrrr, dashed, "\displaystyle{\bigg(a_4,1,\frac{a_1}{\iota_X},\frac{a_3-a_4}{\iota_X}, \frac{a_3-a_4}{\iota_X}(\iota_X-1);\frac{d_1}{\iota_X}\bigg)}" ] \ar[rrrrd ]  \& {} \&{} \&{} \&{} \&{} \&{} \&{} \&Y_2 \ar[lllld ]\\
						{} \&{} \&{} \& {} \& p_{x_4} {} \& {} \& {} \& {} 
        \end{tikzcd}
    \]
This is indeed a flip of type 1, according to \cite[Theorem~8]{brown}.		
\end{proof}

\begin{Rem} \label{rem:120}
Although a quasismooth member $X$ of family 120 has index 4, since $a_1$ is actually a multiple of 4 it follows that $X$ behaves exactly as a quasismooth Fano of prime index. In particular, Corollary \ref{cor:dpflip} applies and there is a hypersurface flip over $p_{x_4}$, $Y_1\rat Y_2$ of type $(-5,-1,-1,1,3;-2)$.  That is, locally analytically around $p_{x_4}$the flip arises as the quotient of a hypersurface by the $\mathbb{C}^*$-action described. 
\end{Rem}

We now prove that there is a fibration $Y_3 \rightarrow \mathbb{P}^1$ whose fibres are del Pezzo surfaces. 

\begin{Lem} \label{lem:toricfib}
There is a toric fibration $T_2 \rightarrow \mathbb{P}^1(1,\iota_X-1)$ whose fibres are 4-dimensional weighted projective spaces.
\end{Lem}
\begin{proof}
Exactly as in the proof of Lemma \ref{lem:flip1} we find a matrix $A$ such that the last two columns of $T$ are multiples of the vector $\left(\begin{smallmatrix} 1 \\ 0 \end{smallmatrix} \right)$ provided that this is not the case already. If it is not, let $a_3':=\frac{ka_3-m_1}{\iota_X} = \frac{a_4(k\iota_X-1)-a_{\xi}}{a_{\xi}\iota_X}$. Since $\gcd(a_3,a_3') =1$, there are integers $r,\,s$ for which $ra_3+sa_3'=1$. Let $A$ be the matrix
\[
\begin{pmatrix}
    r        & s \\
    -a_3'       & a_3  
\end{pmatrix} \in \SL_2(\mathbb{Z}).
\]
As we already observed, the closures of the movable and effective cones of $T$ share a wall, which after multiplication by $A$ is the ray $\mathbb{R}_+[\left(\begin{smallmatrix} 1 \\ 0 \end{smallmatrix} \right)]$. Then we necessarily have a fibration. This is given by 
\[
\Phi' \colon T_2 \rightarrow \Proj \bigoplus_{m\geq 1} H^0(T_2,m \mathcal{O}\left(\begin{smallmatrix} 1 \\ 0 \end{smallmatrix} \right)) \simeq \Proj \mathbb{C}[x_3,x_2] \simeq \mathbb{P}(1,\iota_X-1). 
\]
The fibres of this projection are copies of a four dimensional weighted projective space whose weights are given by $a_3,-a_3'a_{\xi}+a_3k, \ldots, -a_3'a_4+a_3a_4'$. It is a straightforward computation to see that these are 
\[
\wt(u,\xi,x_0,x_1,x_4)=\Big(a_3,\frac{d_1}{\iota_X},1,\frac{a_1}{\iota_X},\,\frac{a_4-a_3}{\iota_X}\Big).
\]
Hence the fibres of the projection are copies of 
\[
\mathbb{P}\Big(a_3,\frac{d_1}{\iota_X},1,\frac{a_1}{\iota_X},\frac{a_4-a_3}{\iota_X}\Big).
\]
\end{proof}

\begin{Cor} \label{cor:dP}
There is a del Pezzo fibration $Y_2 \rightarrow \mathbb{P}^1(1,\iota_X-1)$.
\end{Cor}
\begin{proof}
We restrict to $Y_2$ the construction made in Lemma \ref{lem:toricfib}. Let $\lambda \in \mathbb{C}$ and 
\[
\mathbf{p}_{\lambda} \colon (\lambda x_3^{\iota_X-1}- x_2 = 0) \subset \mathbb{P}^1(1,\iota_X-1).
\]
be a point in the base. Hence $\mathbf{p}_{\lambda} = (1:\lambda)$ The fibre $F_{\lambda}$ over the generic point $\mathbf{p}_{\lambda}$ is given by 
\[
-x_4(f(x_0,x_1)+x_4^{\alpha}x_{\mu}u+\lambda^{\alpha_1}u^{\beta_1}+u\lambda f_{d_1/\iota_X})+g(x_0,x_1)+u^{\beta_2}+ug_{d_2/\iota_X}=0
\]
inside $\mathbb{P}\Big(a_3,1,\frac{a_1}{\iota_X},\frac{a_4-a_3}{\iota_X}\Big)$ with homogeneous variables $u,\,x_0,x_1,\,x_4$ since $\xi$ can be globally eliminated. We want to notice the following features of this equation. If $x_{\mu} = x_2$, then we have $\lambda^{\beta'}x_4^2 \in g_{d_2/\iota_X}$. Otherwise $-x_4^{\alpha+1}u$ is a monomial in the defining equation of the fibre. Hence, $F_{\lambda}$ is smooth for a generic $\lambda \in \mathbb{C}$.

Moreover, $F_{\lambda}$ is given by an equation of degree $d_2/\iota_X$. Hence, 
\begin{align*}
-K_{F_{\lambda}} &\sim \mathcal{O}\Big(a_3+1+\frac{a_1}{\iota_X}+\frac{a_4-a_3}{\iota_X}-\frac{d_2}{\iota_X}\Big)\\
& =\mathcal{O}\Big(\frac{a_3(\iota_X-1)-a_2+\iota_X}{\iota_X}\Big) \\
& = \mathcal{O}(1) && a_3(\iota_X-1)-a_2=0. 
\end{align*}
The degree of $F_{\lambda}$ can be readily computed and it is the integer $K_{F_{\lambda}}^2 = \frac{d_2}{a_3a_1}$. The special fibre over the point at infinity $\mathbf{p_{\infty}} = (0:1)$ is $F_0$ given by
\begin{align*}
f(x_0,x_1)+x_4^{\alpha}x_{\mu}u+x_3^{\iota_X-1}u^{\beta}+u\lambda x_3^{\alpha_1}f_{d_1} &= 0\\
\xi x_4 + g(x_0,x_1)+u^{\beta_2}+g_{d_2} &=0
\end{align*}
inside $\mathbb{P}\Big(a_3,\frac{d_1}{\iota_X},1,\frac{a_1}{\iota_X},\frac{a_4-a_3}{\iota_X})$ with homogeneous variables $u,\,\xi,\,x_0,x_1,\,x_4$. These are two equations of degrees $d_1/\iota_X$ and $d_2/\iota_X$ having a singularity at $\mathbf{p_{\xi}}$. Similarly we can see that $-K_{F_0} \sim \mathcal{O}(1)$ and every fibre is a del Pezzo surface.
\end{proof}

We have the following conclusion
\begin{Thm} \label{thm:delPezzoMain}
Let $X \in \{88, 90, 103, 114, 116, 120 \}$ be a quasismooth member. Then, there is an elementary Sarkisov link $\sigma \colon X \rat Y/\mathbb{P}^1$ to a fibration into del Pezzo surfaces of integer degree $\frac{d_2}{a_3a_1}$ given by
\[
        \begin{tikzcd}[ampersand replacement=\&, column sep = 3em]
      \& Y \ar[rr, "\displaystyle{\simeq}" ] \ar[dl,swap, "\displaystyle{\varphi}" ]  \& \& Y_1   \ar[rr,dashed, "\displaystyle{\tau_1}" ]\& \& Y_2  \ar[rd, "\displaystyle{\varphi'}" ] \&\\
X \&   \&   \&  \&  \&  \& \mathbb{P}(1,\iota_X-1)
        \end{tikzcd}
    \]  
		where
		\begin{itemize}
			\item $\varphi$ and $\varphi'$ are the restrictions of $\Phi$ and $\Phi'$ to $Y$ and $Y_2$, respectively.
			\item $\tau_1 \colon Y_1 \rat Y_2$ is a toric flip on $\mathbb{C}^4$ if and only if $x_4^{\iota_X} \in f$ and is hypersurface flip on $(h=0) \subset \mathbb{C}^{5}$ otherwise.
		\end{itemize}
\end{Thm}
\begin{proof}
The statement is a build up from the previous results. We only justify that the factorisation of the birational map $\sigma$ is indeed a Sarkisov link. Since $Y$ follows the 2-ray game on $T$ we have that the the cones of pseudo-effective divisors and movable divisors of $Y$ coincide. Moreover, the class of the anticanonical divisor of $Y$, $-K_{Y}$ is in the interior of $\overline{\Mov}(Y)$. Hence, by Lemma \ref{lem:linkmds}, the Kawamata blowup $\varphi \colon Y \rightarrow X$ centred at $\mathbf{p_{\xi}}$ initiates a Sarkisov link.
\end{proof}

For $\iota_X\not =4$, a quasismooth member in $I_{dP}$ satisfying $a_4 = a_3$ and $a_3(\iota_X-1)<a_2$ also admits a del Pezzo fibration. This is only satisfied by family 89, $X_{6,6} \subset \mathbb{P}(1,1,2,2,3,5)$:

\begin{Prop} \label{prop:89}
Let $X$ be a member of family 89. Then it admits an elementary Sarkisov link  to a del Pezzo fibration, $X \rat X'/\mathbb{P}^1$ of degree 3.
\end{Prop}

\begin{proof}
The computations are very similar to the previous with the difference that the variables change slightly, so we are brief. By quasismoothness we can assume $X$ has equations as in Table \ref{tab:big}. We can also assume from quasismoothness that $v^2 \in f$. The proper transform of $X$ via the Kawamata blowup of $\mathbf{p_w}$ is isomorphic to 
		\[
\begin{array}{cccccc|ccccc}
             &       & u  & w &   z & t & v & x & y & \\
\actL{(\widetilde{f}=\widetilde{g}=0) \colon Y_2 \subset T_2 }   &  \lBr &  -3 & -4 &   -1 & -1 & 0 & 1 & 1 &   \actR{.}\\
             &       & 1 & 3 &   1 & 1 & 1 & 0 & 0 &  
\end{array}
\] 
where $\widetilde{f}=wx+v^2u+f'(u,x,y,z,t,v)$ and $\widetilde{g}=wy+g'(z,t)+g''(u,x,y,z,t,v)$. Around a small neighborhood of $\mathbf{p_v}$ we can eliminate $u$ and hence we are contracting the curve $C:(x=y=0) \cong (wy+g'(z,t)+g''(u(z,t,w),0,0,z,t,1)=0)  \subset \mathbb{P}(1_z,1_t,4_w)$ and extracting a $\mathbb{P}^1$ which is a  hypersurface flip of the form $(4,1,1,-1,-1;-3)$ and flipping equation $wx+g'(z,t)+...=0$ which is a type 1 terminal hypersurface flip by \cite[Theorem~8]{brown}, $(Y_2 \subset T_2 )\rat (Y_3 \subset T_3)$. Moreover, the projection $T_3 \rightarrow \mathbb{P}^1$ restricts to a del Pezzo fibration where the fibre above each point $(x_0,y_0)$ is a del Pezzo surface of degree 3 in $\mathbb{P}^3$ given by $\widetilde{f}(u,0,z,t,v,x_0,y_0)+\widetilde{g}(u,0,z,t,v,x_0,y_0)=0$.   
\end{proof}

\subsubsection{The exceptional cases.} So far we have dealt with most cases of quasismooth families $X \in I_{dP}$. We now complete the analysis considering the exceptional cases. These correspond to three subcases: 
\begin{itemize}
\item $X \in I_{dP}\setminus \{120 \}$ and $\iota_X=4$.
\item $X \in I_{dP}$ and the Sarkisov link is initiated by a non-linear cyclic quotient singularity.
\item $X \in I_{dP}$ and there is a composition of more than one elementary Sarkisov link to a del Pezzo fibration.
\end{itemize}
Although by Theorem \ref{thm:cb} familes 113 and 119 are in $I_{cB}$, we include these here and find del Pezzo fibration models. The cases above are realised by the following families:
\begin{center}
\begin{tabular}{lccccccc} \toprule
    $X$ & 113  & 119 &121&122& 123  & 124\\
		$\iota_X$ & 3  & 4 &4&4& 4 & 4\\
    $\mathbf{p_{\xi}}$  & $\frac{1}{2}(1,1,1)$  & $\frac{1}{3}(1,1,2)$ &$\frac{1}{7}(1,2,5)$ &$\frac{1}{7}(1,3,4)$& $\frac{1}{9}(1,4,5)$ & $\frac{1}{11}(1,4,7)$    \\
		 		\bottomrule
		\end{tabular}
\end{center}

\begin{Prop} \label{prop:113_119}
Let $X$ be the general member of families 113 or 119, respectively. Then there is a Sarkisov link
\[
\sigma \colon X \rat Y'/\mathbb{P}^1
\]
to a degree 4 del Pezzo fibration.
\end{Prop}

\begin{proof} We analyse each family separately. 
\paragraph{Family 113.}  Let $X$ be a general quasismooth member of family 113. In general $X$ has four cyclic quotient singularities of type $\frac{1}{2}(1,1,1)$. After a change of coordinates we can assume that at least one of these is a coordinate point, that is, we can write $X$ as 
\begin{align*}
\xi w + tv  + z^3 + f_6(x,z,t,v,w) &=0\\
\xi^2z+\xi g_{4}(x,z,t,v) + t^2+tv+v^2  + g_{6}(x,z,t,v,w) &=0.
\end{align*} 
inside $\mathbb{P}(1,2,2,3,3,4)$ with homogeneous variables $x,\,\xi,\,z,\,t,\,v,\,w$. Hence $\mathbf{p_{\xi}} \in X$ is a cyclic quotient singularity of type $\frac{1}{2}(1,1,1)$. Locally around $\mathbf{p_{\xi}}$, the Kawamata blowup $\varphi \colon E \subset Y \rightarrow X$ centred at this point is the graph of the rational map  to $\mathbb{P}(1,1,1)$ given by
\[
\bigg(x: \frac{t}{\xi} :\frac{v}{\xi} \bigg).
\]
Let $\Phi \colon T \rightarrow \mathbb{P}$ be the toric blowup of $\mathbb{P}$ centred at $\mathbf{p_{\xi}}$ whose restriction to $X$ is $\varphi \colon Y \rightarrow X$. By Corollary \ref{cor:genlift}, we have that the sections $t,\,v \in H^0(X,-K_X)$ lifts to anti-canonical sections in $Y$ and
\[
x \in H^0\bigg(Y,-\frac{1}{3}K_Y-\frac{1}{3}E\bigg),\, \xi \in H^0\bigg(Y,-\frac{2}{3}K_Y+\frac{1}{3}E\bigg). 
\]
Hence $m_f = m_g = 1$ since both $f$ and $g$ contain monomials purely in the ideal $(t,v)$ by quasismoothness of $X$. Hence,
\[
w \in H^0\bigg(Y,-\frac{4}{3}K_Y-\frac{1}{3}E\bigg),\, z \in H^0\bigg(Y,-\frac{2}{3}K_Y-\frac{2}{3}E\bigg).
\]
Following the proof of Lemma \ref{lem:utbl}, the weight system of $T \supset Y$ is 
		\[
\begin{array}{cccc|ccccccc}
             &       & u  & \xi &   t & v & w & x & z & \\
\actL{T}   &  \lBr &  0 & 2 &  3 & 3 & 4 & 1 & 2 &   \actR{.}\\
             &       & 1 & 1 &   1 & 1 & 1 & 0 & 0 &  
\end{array}
\] 
The closure of the cones of movable and effective divisors of $T$ can therefore be described as
\[
\langle \left(\begin{smallmatrix} 2 \\ 1 \end{smallmatrix} \right),\left(\begin{smallmatrix} 1 \\ 0 \end{smallmatrix} \right) \rangle = \overline{\Mov}(T) \subset \overline{\Eff}(T) =  \langle \left(\begin{smallmatrix} 0 \\ 1 \end{smallmatrix} \right),\left(\begin{smallmatrix} 1 \\ 0 \end{smallmatrix} \right) \rangle.
\]
The cone of movable divisors on $T$ is subdivided into three chambers. Hence we have small modifications $T \rat T_1 \rat T_2$. In particular we describe $T_1 \rat T_2$. This map corresponds to crossing the wall given by the ray $\mathbb{R}_+[\left(\begin{smallmatrix} 4 \\ 1 \end{smallmatrix} \right)]$. Let 
\[
\mathcal{F}_1 := \Proj \bigoplus_{m\geq 0} H^0(T_1, m\mathcal{O}_{T_1}\left(\begin{smallmatrix} 4 \\ 1 \end{smallmatrix} \right)) = \Proj  \bigoplus_{m\geq 0} H^0(T_2, m\mathcal{O}_{T_2}\left(\begin{smallmatrix} 4 \\ 1 \end{smallmatrix} \right))
\]
Then, we have the following two small contractions
\[
f_1 \colon T_1 \rightarrow  \mathcal{F}_1 \leftarrow T_2 \colon g_1 
\]
and $\tau_1 \colon T_1 \rightarrow T_2$ is  the composition $ g_1^{-1}\circ f_1$. Notice that the $T_i$ have the same Cox ring of $T$ but differ in their irrelevant ideals. For example, the irrelevant ideal of $T_1$ is $(u,\xi,t,v) \cap (w,x,z)$ and the irrelevant ideal of $T_2$ is $(u,\xi,t,v,w) \cap (x,z)$. The maps $f_1$ and $g_1$ are given in coordinates by
\[
(u,\xi,t,v,w,x,y) \mapsto (w,vx,tx,\ldots,yz,uz^2) \in \mathcal{F}_1.
\]
In particular, $f_1$ contracts the threefold $\mathbb{P}(4_u,2_{\xi},1_t,1_v) \simeq (x=y=0) \subset T_1$ to $p_w \in \mathcal{F}_1$ and $g_1$ extracts from it the locus $\mathbb{P}(1_x,2_z) \simeq (u=\xi=t=v=0) \subset T_2$. Away from $p_w$ both $f_1$ and $g_1$ are isomorphisms. We denote this small contraction by $(-4,-2,-1,-1,1,2)$.

Notice that the closures of the movable and effective cones of $T$ share the wall given by the ray $\mathbb{R}_+[\left(\begin{smallmatrix} 1 \\ 0 \end{smallmatrix} \right)]$. Hence we have the fibration
\[
\Phi' \colon T_2 \rightarrow \mathcal{F}_2
\]
where 
\[
\mathcal{F}_2 = \Proj \bigoplus_{m\geq 1} H^0(T_2, m\mathcal{O}_{T_2}\left(\begin{smallmatrix} 1 \\ 0 \end{smallmatrix} \right)) = \Proj \mathbb{C}[x,z] \simeq \mathbb{P}(1,2).
\]
Each fibre is isomorphic to $\mathbb{P}^4$. We restrict this construction to $Y \subset T$. The defining equations of $Y$ are $\widetilde{f}=0$ and $\widetilde{g}=0$ where $(\widetilde{f}=0),\, (\widetilde{g}=0) \in |-2K_Y|$. By Lemma \ref{lem:anticanonicalwalliso}, the small modification $T \rat T_1$ restricts to an isomorphism on $Y$. On the other hand, $Y|_{x=z=0} \subset T_1$ is isomorphic to the curve $(t^2+tv+v^2 = 0) \subset \mathbb{P}(4_u,1_t,1_v)$ and $Y|_{u=\xi=t=v=0} \subset T_2$ is isomorphic to $\mathbb{P}(1_x,2_z)$. Hence, around the point $p_w \in \mathcal{F}_1$ we have the hypersurface flip $(-4,-1,-1,1,2;-2)$ with flipping equation $zu+t^2+vt+v^2+\cdots = 0$. This is a type 1 hypersurface flip as in \cite[Theorem~8]{brown}. Let $Y'$ be the flipped 3-fold. The map $\Phi'$ restricts to a del Pezzo fibration $Y'/\mathbb{P}(1,2)$. In fact, the fibre $S$ over a generic point is
\begin{align*}
\xi w + tv + f_2(u,t,v ) &=0\\
\xi^2+\xi l(u,t,v,w) + t^2+tv+v^2 + wu + g_{2}(u,t,v,w) &=0
\end{align*} 
in $\mathbb{P}^4$ where $u^2 \in f_2$ or $u^2 \in g_2$ and $l$ is a linear form. Hence, the generic fibre is a degree four del Pezzo surface. Over the point $\mathbf{p_x}$ the fibre, $S_0$, is singular and isomorphic to 
\begin{align*}
\xi w + tv + f_2(u,t,v ) &=0\\
\xi l(u,t,v,w) + t^2+tv+v^2 + g_{2}(u,t,v,w) &=0
\end{align*} 
and $\mathbf{p_w} \in S_0$ is the only singular point which is a du Val $A_1$ singularity if the polynomials $f_2$ and $g_2$ are general.

\paragraph{Family 119.} Let $X$ be a quasismooth member of family 119. Then, $X \subset \mathbb{P}:=\mathbb{P}(1,2,3,4,5,7)$ with homogeneous variables $x,\,y,\,\xi,\,t,\,v,\,w$ can be written as 
\begin{align*}
\xi v + wx  + t^2 + f_8(x,y,t,v ) &=0\\
\xi w + v^2 + y^3t + g_{10}(x,y,t,v,w) &=0.
\end{align*} 
Then $\mathbf{p_{\xi}} \in X$ is a cyclic quotient singularity of type $\frac{1}{3}(1,1,2)$. Locally around $\mathbf{p_{\xi}}$, the Kawamata blowup $\varphi \colon E \subset Y \rightarrow X$ centred at this point is the graph of the rational map  to $\mathbb{P}(1,1,2)$ given by
\[
\bigg(\frac{t}{z}: x :y \bigg).
\]
Let $\Phi \colon T \rightarrow \mathbb{P}$ be the toric blowup of $\mathbb{P}$ centred at $\mathbf{p_{\xi}}$ whose restriction to $X$ is $\varphi \colon Y \rightarrow X$. By Corollary \ref{cor:genlift}, we have that the section $t \in H^0(X,-K_X)$ lifts to an anti-canonical section in $Y$ and
\[
x \in H^0\bigg(Y,-\frac{1}{4}K_Y-\frac{1}{4}E\bigg),\, y \in H^0\bigg(Y,-\frac{1}{2}K_Y-\frac{1}{2}E\bigg),\, y \in H^0\bigg(Y,-\frac{3}{4}K_Y+\frac{1}{4}E\bigg). 
\]
Hence $m_f = \frac{2}{3}$ since $f$ contains the monomial $t^2$ and $m_g = \frac{4}{3}$ since $g$ contains the monomial $v^2$ and no monomial vanishing less than 1 at $E$. Hence,
\[
v \in H^0\bigg(Y,-\frac{5}{4}K_Y-\frac{1}{4}E\bigg),\, w \in H^0\bigg(Y,-\frac{7}{4}K_Y-\frac{3}{4}E\bigg).
\]
Following the proof of Lemma \ref{lem:utbl}, the weight system of $T \supset Y$ is 
		\[
\begin{array}{cccc|ccccccc}
             &       & u  & \xi &   t & v & w & x & y & \\
\actL{T}   &  \lBr &  0 & 3 &  4 & 5 & 7 & 1 & 2 &   \actR{.}\\
             &       & 1 & 1 &   1 & 1 & 1 & 0 & 0 &  
\end{array}
\] 
The closure of the cones of movable and effective divisors of $T$ can therefore be described as
\[
\langle \left(\begin{smallmatrix} 3 \\ 1 \end{smallmatrix} \right),\left(\begin{smallmatrix} 1 \\ 0 \end{smallmatrix} \right) \rangle = \overline{\Mov}(T) \subset \overline{\Eff}(T) =  \langle \left(\begin{smallmatrix} 0 \\ 1 \end{smallmatrix} \right),\left(\begin{smallmatrix} 1 \\ 0 \end{smallmatrix} \right) \rangle.
\]
The cone of movable divisors on $T$ is subdivided into four chambers. Hence we have small modifications $T \rat T_1 \rat T_2 \rat T_3$. In particular we describe $T_2 \rat T_3$. This map corresponds to crossing the wall given by the ray $\mathbb{R}_+[\left(\begin{smallmatrix} 7 \\ 1 \end{smallmatrix} \right)]$. Let 
\[
\mathcal{F}_2 := \Proj \bigoplus_{m\geq 1} H^0(T_2, m\mathcal{O}_{T_2}\left(\begin{smallmatrix} 7 \\ 1 \end{smallmatrix} \right)) = \Proj  \bigoplus_{m\geq 1} H^0(T_3, m\mathcal{O}_{T_3}\left(\begin{smallmatrix} 7 \\ 1 \end{smallmatrix} \right))
\]
Then, we have the following two small contractions
\[
f_2 \colon T_2 \rightarrow  \mathcal{F} \leftarrow T_3 \colon g_2 
\]
and $\tau_2 \colon T_2 \rightarrow T_3$ is  the composition $ g_2^{-1}\circ f_2$. Notice that the $T_i$ have the same Cox ring of $T$ but differ in their irrelevant ideals. For example, the irrelevant ideal of $T_2$ is $(u,\xi,t,v) \cap (w,x,z)$ and the irrelevant ideal of $T_3$ is $(u,\xi,t,v,w) \cap (x,y)$. The maps $f_2$ and $g_2$ are given in coordinates by
\[
(u,\xi,t,v,w,x,y) \mapsto (w,vx^2,tx^3,\ldots,zy^2,uxy^3) \in \mathcal{F}.
\]
In particular, $f_2$ contracts the threefold $\mathbb{P}(7_u,4_{\xi},3_t,2_v) \simeq (x=y=0) \subset T_2$ to $p_w \in \mathcal{F}$ and $g_2$ extracts from it the locus $\mathbb{P}(1_x,2_y) \simeq (u=\xi=t=v=0) \subset T_3$. Away from $p_w$ both $f_2$ and $g_2$ are isomorphisms. We denote this small contraction by $(-7,-4,-3,-2,1,2)$.

Notice that the closures of the movable and effective cones of $T$ share the wall given by the ray $\mathbb{R}_+[\left(\begin{smallmatrix} 1 \\ 0 \end{smallmatrix} \right)]$. Hence we have the fibration
\[
\Phi' \colon T_3 \rightarrow \mathcal{F}_3
\]
where 
\[
\mathcal{F}_3 = \Proj \bigoplus_{m\geq 0} H^0(T_3, m\mathcal{O}_{T_3}\left(\begin{smallmatrix} 1 \\ 0 \end{smallmatrix} \right)) = \Proj \mathbb{C}[x,y] \simeq \mathbb{P}(1,2).
\]
Each fibre is isomorphic to $\mathbb{P}^4$. We restrict this construction to $Y \subset T$. The defining equations of $Y$ are $\widetilde{f}=0$ and $\widetilde{g}=0$ where $(\widetilde{f}=0) \in |-2K_Y|$ and $(\widetilde{g}=0) \in |-\frac{5}{2}K_Y-\frac{1}{2}E|$. Hence $t^2 \in \widetilde{f}$ and $v^2 \in \widetilde{g}$ and the small modifications $T \rat T_1 \rat T_2$ therefore restrict to isomorphisms on $Y$. On the other hand, $Y|_{x=y=0} \subset T_2$ is isomorphic to the curve $(v^3+t^2 = 0) \subset \mathbb{P}(7_u,3_t,2_v)$ and $Y|_{x=y=0} \subset T_2$ is isomorphic to $\mathbb{P}(1_x,2_y)$. Hence, around the point $p_w \in \mathcal{F}$ we have the hypersurface flip $(-7,-3,-2,1,2;-6)$ with flipping equation $xu+t^2+v^3+\cdots = 0$. This is a type 2 hypersurface flip as in \cite[Theorem~8]{brown}. Let $Y'$ be the flipped 3-fold. The map $\Phi'$ restricts to a del Pezzo fibration $Y'/\mathbb{P}(1,2)$. In fact, the fibre $S$ over a generic point is
\begin{align*}
\xi v + wu  + t^2 + f_2(u,t,v ) &=0\\
\xi w + v^2 + tu + g_{2}(t,v,w) &=0
\end{align*} 
in $\mathbb{P}^4$ where $u^2 \in f_2$ or $u^2 \in g_2$. That is, the generic fibre is a degree four del Pezzo surface. Over the point $\mathbf{p_y}$ the fibre, $S_0$, is singular and isomorphic to 
\begin{align*}
\xi v + t^2 + f_2(u,t,v ) &=0\\
\xi w + v^2 + tu + g_{2}(t,v,w) &=0
\end{align*} 
and $\mathbf{p_w} \in S_0$ is the only singular point which is a du Val $A_1$ singularity if the polynomials $f_2$ and $g_2$ are general.
\end{proof}


\begin{Prop} \label{prop:121_122}
Let $X$ be a quasismooth member of families 121 or 122, respectively. Then there is a Sarkisov link
\[
\sigma \colon X \rat Y'/\mathbb{P}^1
\]
to a degree 3 del Pezzo fibration.
\end{Prop}

\begin{proof}
\begin{description}
	\item[Family 121] By Lemma \ref{lem:utbl} the proper transform of $X$ via the Kawamata blowup of $\mathbf{p_w}$ is given by 
			\[
\begin{array}{cccccc|ccccc}
             &       & u  & w &   z & t & v & x & y & \\
\actL{(\widetilde{f}=\widetilde{g}=0) \colon Y_3 \subset T_3 }   &  \lBr &  -6 & -5 &   -2 & -1 & 0 & 1 & 3 &   \actR{.}\\
             &       & 1 & 2 &   1 & 1 & 1 & 0 & 0 &  
\end{array}
\] 
	
	Moreover, by quasismoothness, $t^2 \in f$ and it lifts to a plurianticanonical section $t^2 \in \widetilde{f}$. This induces an isomorphism $Y_2 \rightarrow Y_3$ and by Lemma \ref{lem:iso} $(Y_1\subset T_1) \cong (Y_2\subset T_2)$. Hence $(Y_1\subset T_1) \cong (Y_3\subset T_3)$. Over a small neighborhood of $\mathbf{p_v}$ it is possible to eliminate the variable $z$ since $vz \in f$ by quasismoothness and it lifts to a pluri-anticanonical section $vz \in \widetilde{f}$. The variable $u$ can be eliminated for the same reasons since $v^2u \in \widetilde{g}$. Hence we have a toric flip of type $(-5,-1,1,3)$ over $\mathbf{p_v}$. 
	\item[Family 122]By Lemmas \ref{lem:utbl} and \ref{lem:iso} the proper transform of $X_{122}$ via the Kawamata blowup of $\mathbf{p_w}$ can be written as 
			\[
\begin{array}{ccccc|cccccc}
             &       & u  & w &   z & t & v & x & y & \\
\actL{(\widetilde{f}=\widetilde{g}=0) \colon Y_2 \subset T_2 }   &  \lBr &  -5 & -3 &   -1 & 0 & 0 & 2 & 3 &   \actR{.}\\
             &       & 1 & 2 &   1 & 1 & 1 & 0 & 0 &  
\end{array}
\] 
where $\widetilde{f}=wy+vt+f'$ and $\widetilde{g}=w(v+\alpha t)+z^3+t^2xu+g'$. From $T_2$ there is toric flip based on $\mathbb{P}^1$ which contracts a codimension 2 weighted plane $\mathbb{P}(5,3,1)$, and extracts a projective line. We prove that this restricts to 2 flips over 2 different points. Indeed $Y_2{_{|_{u=w=z=x=y=0}}} : (vt=0) \subset \mathbb{P}^1$. Over each of these points, the fibre is isomorphic to $\mathbb{P}(5,1)$ for a generic choice of $\alpha$. When $\alpha=0$, however, the fibre over $p_t$ degenerates and becomes a cubic curve in $\mathbb{P}(5,3,1)$. In that case, we have a hypersurface flip of type $(-5,-3,-1,2,3;-3)$ with equation $xu+z^3+ \cdots =0$.
\end{description}
\vspace{0.5cm}

Finally, there is a projection $T_3 \rightarrow \mathbb{P}^1$ whose fibres are copies of $\mathbb{P}(1,2,1,1,1)$ in both cases. Indeed this projection restricts to a del Pezzo fibration: The fibre above a general point is a del Pezzo Surface of degree 3, $dP_3$, in $\mathbb{P}^3$. Over the point $(1:0)$ we get a singular complete intersection $S_{2,3} : (\widetilde{f}(u,0,z,t,v,0,0)=\widetilde{g}(u,w,z,t,v,0,0)=0) \subset \mathbb{P}(1,2,1,1,1)$. In both cases the anticanonical divisor is linearly equivalent to a hyperplane section.
\end{proof}

We find two convenient models for $X \in \{ 123,124 \}$. After a change of variables, any quasismooth member of family 123 is given by the equations
\begin{align*}
wy+vt+z^3+f_{12}(x,z)&=0\\
wt+v^2+vg_{7}(x,y,z,t)+t^2z+tg_{9}(x,y,z)+g_{14}(x,y,z)&=0
\end{align*}
inside $\mathbb{P}(2,3,4,5,7,9)$ with homogeneous variables $x,\,y,\,z,\,t,\,v,\,w$ where $y^3 \in g_9$ or $y^4x \in g_{14}$ and we have as well $x^6 \in f_{12}$ or $x^7 \in g_{14}$. Similarly, after a change of variables, any quasismooth member of family 124 is given by the equations
\begin{align*}
wt+v^2+vf_9(x,y)+z^3+f_{18}(x,y,z)&=0\\
wv+vg_{11}(x,y,z,t)+t^2z+tg_{13}(x,y)+y^4+x^5+g_{20}(x,y,z)&=0
\end{align*}
inside $\mathbb{P}(4,5,6,7,9,11)$ with homogeneous variables $x,\,y,\,z,\,t,\,v,\,w$. 

\begin{Lem} \label{lem:123and124}
Let $X$ be quasimooth member of families 123 or 124. Then, $X$ is birational to $X'$ where $-K_{X'} \sim \mathcal{O}_{X'}(2)$ and admits a projection $\pi \colon X' \dashrightarrow \mathbb{P}^1$ whose generic fibre $S := \pi^{-1}(t)$ has isolated singularities and is uniruled.    
\end{Lem}

\begin{proof}
We look at each family separately.

\paragraph{Family 123.} By Theorem \ref{thm:cod2ExcII}, there is a Sarkisov link to a codimension 2 complete intersection Fano 3-fold $X':=Z_{6,6}$. See  Theorem \ref{thm:cod2ExcII} for the equations. Consider the projection $\mathbb{P}(1,1,2,2,3,5)~\dashrightarrow~\mathbb{P}^1$ given by $(t : y : u : z : v : w) \mapsto (t:y)$ and its restriction to $\pi \colon Z_{6,6} \dashrightarrow \mathbb{P}^1$. A generic fibre of $\pi$ is the divisor $S \colon (t- \lambda y=0) \subset Z_{6,6}$ for some $\lambda \in \mathbb{C}^*$. In other words, $S$ is the complete intersection
\begin{align*}
w y+\lambda vyu+z^3+f_{6}(u,z)&=0\\
w \lambda y+v^2+vg_{3}(u,1,y,z,\lambda y)+(\lambda y)^2zu+\lambda y g_{5}(u,1,y,z)+(\lambda y^4)u+g_{6}(u,1,y,z)&=0
\end{align*}
in $\mathbb{P}(1,2,2,3,5)$ with homogeneous variables $t,\,y,\,u,\,z,\,v,\,w$. It is easy to see that $S$ has only one singular point, $\mathbf{p_w}$ and that this is analytically isomorphic to the germ
\[
0 \in (v^2+u^2z+z^3=0) \big/\frac{1}{5}(2_u,2_z,3_v)
\]
where $\mathbb{Z}_5$ acts equivariantly on the $D_4$ singularity so that $\mathbf{p_w}$ is an isolated surface singularity. By adjunction for singular surfaces, $-K_S$ is a big divisor. Hence $S$ is uniruled. See, for instance, \cite[Lemma~2.1]{liviatiago}.


\paragraph{Family 124.} By Theorem \ref{thm:124toHyp}, there is a Sarkisov link to a Fano 3-fold hypersurface $X' := Z_{10}$. We refer to Theorem \ref{thm:124toHyp} for its equations. Consider the projection $\mathbb{P}(1,1,2,3,5)~\dashrightarrow~\mathbb{P}^1$ given by $(w : y : x : v : t) \mapsto (w:y)$ and its restriction to $\pi \colon Z_{10} \dashrightarrow \mathbb{P}^1$. A generic fibre of $\pi$ is the divisor $S \colon (w- \lambda y=0) \subset Z_{6,6}$ for some $\lambda \in \mathbb{C}^*$. In other words, $S$ is the hypersurface
\begin{multline*}
t^2+tg_{5}(x,y,v)+\lambda y v(t\lambda y+v^2+vf_3(x,y)+f_{6}(x,y,z))+\\y^4(\lambda ty+v^2+vf_3(x,y)+f_{6}(x,y,z))+x^5+g_{10}(x,y,1)=0
\end{multline*}
in $\mathbb{P}(1,2,3,5)$ with homogeneous variables $y,\,x,\,v,\,t$. In this case, the only singular point of $S$ is $\mathbf{p_v} \sim \frac{1}{3}(1,1)$ which is an isolated klt singularity. By adjunction for singular surfaces, $-K_S$ is a big divisor. Hence $S$ is uniruled. See, for instance, \cite[Lemma~2.1]{liviatiago}.
\end{proof}

\begin{Lem} \label{lem:compdp}
Let $X \in \{91,123,124,125\}$. Then $X$ is birational to a del Pezzo fibration.
\end{Lem}

\begin{proof}
By Theorem \ref{thm:cod2curve}, $X_{91}$ is birational to the singular complete intersection of two cubics $Z_{3,3} \subset \mathbb{P}(1^4,2)$. Grinenko proved in \cite{Grinenko33} that a general complete intersection of two cubics $X_{3,3}\subset \mathbb{P}(1^4,2)$ is birational to a fibration into cubic surfaces over $\mathbb{P}^1$. As expected, the same conclusion holds for $Z_{3,3}$ by blowing up its $\frac{1}{2}(1,1,1)$ singularity and then performing a number of flops. We omit the details. Notice that $Z_{3,3}$ is also birational to a quartic containing a plane.

By Proposition \ref{prop:endcod2} and Theorem \ref{thm:cod2caseI}, $X_{125}$ is not birationally rigid and is in fact birational to a (non-quasismooth) member of family 120 i.e., $Z_{8,12} \subset \mathbb{P}(1,3,4,4,5,7)$ containing a $\frac{1}{7}(1,3,4)$ singularity. Blowing up this singularity in $Z_{8,12}$ initiates a Sarkisov link to a fibration into cubic sufaces over $\mathbb{P}^1$. The computations carry exactly as in Corollaries \ref{cor:dpflip} and its conclusion in \ref{cor:dP}. See also Remark \ref{rem:120}.

Consider the models $X'$ and the projections $\pi \colon X' \dashrightarrow \mathbb{P}^1$ as in Lemma \ref{lem:123and124}. The result follows from \cite[Theorem~2.2]{liviatiago}. Compare also with \cite{hamidhyp}, where the authors use the same approach to prove non-solidity of a number of Fano 3-fold hypersurfaces of Fano index at least 2.
\end{proof}

The following theorem puts together the results of this section and establishes part of Theorem \ref{thm:main}.
\begin{Thm} \label{thm:dP}
Let $X$ be a quasismooth member of a deformation family in $I_{dP}$. Then there is a Sarkisov link 
\[
\sigma \colon X \rat Y/\mathbb{P}^1
\]
to a del Pezzo fibration of degrees $1$, $2$ or $3$. The Sarkisov link is initiated by the Kawamata blowup of a cyclic quotient singularity.
\end{Thm}
\begin{proof}
We only justify that the map $\sigma$ is indeed a Sarkisov link. Since $Y$ follows the 2-ray game on $T$, the closure of their movable cones coincide. Hence, in particular, the class of the anticanonical divisor of $Y$ is in the interior of $\overline{\Mov}(T)$ and $Y$ is a Mori Dream Space. Moreover, we have seen that each small $\mathbb{Q}$-factorial modification $Y_i \rat Y_{i+1}$ is terminal. Hence the 2-ray game on $Y$ always happens in the Mori category. By Lemma \ref{lem:linkmds}, the Kawamata blowup at $\mathbf{p_{\xi}}$ initiates a Sarkisov link.
\end{proof}

\begin{Thm} \label{thm:nS}
Let $X$ be a quasismooth member of a deformation family in $I_{nS}$. Then $X$ is non-solid.
\end{Thm}

\paragraph{A digression into special members}

There are special members of certain families that admit a structure of a del Pezzo fibration even when a general member does not seem to admit one. This is the case of families 92 and 94, for instance, as can be seen in Table \ref{tab:big} at the end. Indeed, it is expected that the number of Mori fibre space structures increases as we specialise the members further, although this relation is quite subtle. 

We show one example in this direction. Let $X$ denote a general quasismooth member of family 94. Then we have the following prediction, see also Conjecture \ref{conj:solid}.
\begin{Conj}
Let $X$ denote a general quasismooth member of family 94. Then $X$ is solid. Moreover, it admits only 3 non-isomorphic Mori fibre space structures which are given in Table \ref{tab:big}.
\end{Conj} 

Let $V \subset\mathbb{P}H^0(\mathbb{P},\mathcal{O}_{\mathbb{P}}(d_1))\cap \mathbb{P}H^0(\mathbb{P},\mathcal{O}_{\mathbb{P}}(d_2))$ be the projective space parametrising quasismooth elements of the deformation family of $X$. We denote by $X_U$ the families parametrised by $U$, whenever $U \subset V$ is a sub-variety of $V$. Let $D \subset V$ be those families for which $\frac{1}{3}(1,1,2) \in X_D$ is \emph{not} a linear cyclic quotient singularity with respect to $X$, see Subsection \ref{subsec:lin}.

\begin{Prop}
$X_D$ is non-solid. Indeed, $X_D$ is birational to a standard smooth del Pezzo fibration of degree 1.
\end{Prop}


\begin{Ex} \label{ex:pic2}

Let $X$ be a member of family 94 as in Table \ref{tab:big} with $\alpha=0$. For simplicity we can take $X$ to be
\begin{align*}
t^2y+wx+yf_{3}(y,z)+\beta z^4+x^8 &=0\\
tw+v^2+g_5(y,z)&=0
\end{align*}
where $f_3$ and $g_5$ are general polynomials, although the same conclusions are valid for general members provided $\alpha=0$. Notice that $X$ is quasismooth for any $\beta \in \mathbb{C}^*$ and in the limit $\beta \rightarrow 0$,  it has 2 isolated singular points. Using Lemma \ref{lem:lift},
\begin{align*}
t,\, v,\,w &\in H^0\Big(Y,-\frac{a_{\mu}}{2}K_Y+\frac{1}{2}E\Big), \, a_{\mu} \in \{\wt(t), \wt(v),\wt(w)\} \\
z &\in H^0(Y,-K_Y), \\ 
x,\,y &\in H^0\Big(Y,-\frac{a_{\mu}}{2}K_Y-\frac{a_{\mu}}{2}E\Big), \, a_{\mu} \in \{\wt(x),\wt(y)\}.
\end{align*}
Hence, $\widetilde{f} \in |-4K_Y|$ and $\widetilde{g} \in |-5K_Y+E|$ since $m_f = \frac{4}{3} $ and $m_g = \frac{2}{3}$. So $Y$ is 
\begin{align*}
t^2y+wx+yuf_{3}(yu,z)+\beta z^4+x^8u^4 &=0\\
tw+v^2+ug_5(yu,z)&=0
\end{align*}
in 
			\[
\begin{array}{ccccc|cccccc}
             &       & u  & t &   v & w & z & x & y & \\
\actL{T}   &  \lBr &  -2 & -1 & -1 & -1 & 0 & 1 & 2 &   \actR{.}\\
             &       & -7 & -2 &   -1 & 0 & 1 & 4 & 8 &  
\end{array}
\]
Crossing the $w$-wall on $\Mov(T)$ restricts to an antiflip $(-7,-1,1,8)$ over a point. On the other hand, crossing the $z$-wall on $\Mov(T')$, where $T'$ is defined with the same action as $T$ but has irrelevant ideal $(u,t,v,w) \cap (z,x,y)$, depends on $\beta$. If $\beta \in \mathbb{C}^*$ the $z$-wall crossing restricts to an isomorphism and otherwise, when $\beta =0$, it restricts to a blowdown of $\mathbb{P}^2$ to a point although on $T'$ this is a small contraction. This is because $Y'|_{(x=y=0)} \simeq \mathbb{P}^2$. We denote the resulting 3-fold after the blowdown by $Y_{\text{bl}}$.

The last chamber of $\Mov(T)$ is generated by $\mathbb{R}_+[(z=0)]+\mathbb{R}_+[(x=0)]$ which coincides with the last chamber of $\Eff(T)$. Hence we have a fibration to $\mathbb{P}^1(1,2)$ whose fibres are $\mathbb{P}(1,1,2,3,4)$. This restricts to a fibration on $Y_{\text{bl}}$ whose general fibre is the smooth del Pezzo surface 
\[
S \colon (t(t^2+uf_{3}(u,z)+u^4)+v^2+ug_5(u,z)=0) \subset \mathbb{P}(1_u,1_z,2_t,3_v)
\]
with $K_S^2 =1$.   If $\rank \Cl(X) = 1$, then $\rank \Cl(Y_{\text{bl}})=1$. But that is impossible since otherwise $\rank \Cl(Y_{\text{bl}}/\mathbb{P}^1)=0$ but $\rank \Cl(Y_{\text{bl}}/\mathbb{P}^1)=1$ since we are contracting an extremal ray. We conclude that $\rank \Cl(X) \geq 2$ if $\beta=0$.
\end{Ex}

\section{Sarkisov Links to complete intersection Fano 3-folds}  \label{sect:solid}

Throughout this section $X \in I$ is a quasismooth member of a deformation family in $I$. Recall that the proper transform of $X$ via the Kawamata blowup of $X$ at $\mathbf{p_{\xi}}$ is given by
\begin{equation} \label{eq:T}
\begin{array}{cccc|ccccccc}
             &       & u  & \xi &   x_0 & x_1 & x_2 & x_3 & x_4 & \\
\actL{Y \colon (\widetilde{f}=\widetilde{g}=0)\subset T }   &  \lBr &  0 & a_{\xi} &   \iota_X & a_1 & a_2 & a_4 & a_3 &   \actR{}\\
             &       & 1 & k &   \frac{k\iota_X-1}{a_{\xi}} & \frac{ka_1-\overline{ka_1}}{a_{\xi}} & \frac{ka_2-\overline{ka_2}}{a_{\xi}} & \frac{ka_3-b_3}{a_{\xi}} & \frac{ka_4-b_4}{a_{\xi}} &  
\end{array}
\end{equation}
where the rays defined by the bidegrees of each section are not necessarily ordered. Recall that $k$ and $b_i$ are defined as in the proof of Lemma \ref{lem:utbl}. For ease of notation we write the bottom row of $T$ above as
\[
1,\,k,\,\iota_X',\,a_1',\, \ldots,\, a_4'.
\]
We also notice that $\mathbf{p_{\xi}}$ is not necessarily a linear cyclic quotient singularity w.r.t. $X$, although by Lemma \ref{lem:polyform}, each $X$ has a linear cyclic quotient singularity which is a coordinate point. We prove, in particular, that every quasismooth member of a deformation family in $I_{S} \subset I$ is birational to a non-isomorphic hypersurface or codimension 2 Fano 3-fold establishing the second part of Theorem \ref{thm:main}. A byproduct of this section is that the Kawamata blowup centred at a linear cyclic quotient singularity always initiates a Sarkisov link.

\subsection{Divisorial contractions to a point} \label{sub:divtopoint}


We consider the following cases:
\[
  \begin{cases}
              \text{Case I:} \quad\,\, a_2 > a_4(\iota_X-1)>a_3(\iota_X-1)\\
              \text{Case II:}  \quad\,  a_4(\iota_X-1)>a_2 >a_3(\iota_X-1)\\
         \text{Case III:} \quad a_4(\iota_X-1)>a_3(\iota_X-1)>a_2 \\
         \text{Case IV:\,\, Exceptional cases}
            \end{cases}
\]

Let $\mathbf{p_{\xi}} \in X$ be a linear cyclic quotient singularity, see Definition \ref{def:linear}. We prove explicitly that blowing up $\mathbf{p_{\xi}} \in X$ initiates a Sarkisov link ending with a divisoral contraction to a point. 
The distinction between the three cases is made to reflect the behavior of the cone of movable divisors of $T$ which in turn explains the birational geometry of $T$. The last exceptional case consists of non-linear cyclic quotient singularities of families of Fano index 4. Each of these four cases have in common the feature that the blowup of $(X, \mathbf{p}_{\xi})$ initiates an elementary Sarkisov link ending with a hypersurface Fano 3-fold.



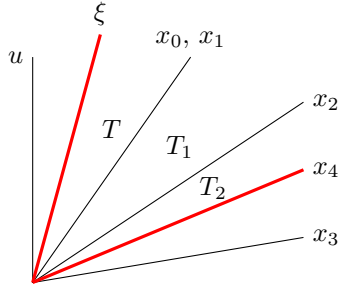
\begin{figure}%
\centering
\begin{tikzpicture}[scale=3]
 \coordinate (A) at (0, 0);
 \coordinate [label={left:$u$}] (E) at (0, 1);
 \coordinate [label={above:$x_0,\,x_1$}] (K) at (0.7, 1);
 \coordinate [label={above:$\xi$}] (5) at (0.3,1.1);
 \coordinate [label={right:$x_2$}] (2) at (1.2,0.8);
 \coordinate [label={right:$x_4$}] (4) at (1.2,0.5);
 \coordinate [label={right:$x_3$}] (3) at (1.2,0.2);

\coordinate [label={$T$}] (T) at (0.35,0.6);
\coordinate [label={right:$T_1$}] (T1) at (0.55,0.6);
\coordinate [label={$T_2$}] (T) at (0.8,0.33);

  \draw (A) -- (E);
  \draw (A) -- (K);
	\draw [very thick,color=red] (A) -- (5);
	\draw (A) -- (2);
	\draw (A) -- (3);
	\draw [very thick,color=red] (A) -- (4);
  
\end{tikzpicture}
\caption{A representation of the chamber decomposition of the cone of effective divisors of $T$. In red the subcone of movable divisors of T. The ample models $T$, $T_1$ and $T_2$ are represented in the interior of the corresponding Nef chamber.}%
\label{fig:HypCaseI}%
\end{figure}

\paragraph{Case I:} Let $X$ be a quasismooth member of a deformation family in the following table and $\mathbf{p_{\xi}} \in X$ a linear cyclic quotient singularity with respect to $X$.
\begin{center}
\resizebox{\textwidth}{!}{\begin{tabular}{lccccccccc} \toprule
    $X$ & 94 & 97  & 98 & 100 & 105 & 106 & 108 & 109 & 110  \\
    $\mathbf{p_{\xi}}$  & $\frac{1}{7}(1,1,6)$  &$\frac{1}{9}(1,1,8)$ & $\frac{1}{9}(1,2,7)$ & $\frac{1}{11}(1,2,9)$ & $\frac{1}{9}(1,1,8)$& $\frac{1}{11}(1,1,10)$& $\frac{1}{11}(1,2,9)$& $\frac{1}{13}(1,1,12)$& $\frac{1}{13}(1,2,11)$   \\ 
		\bottomrule
		\end{tabular}}
\end{center}
These are the cases for which $a_2 > a_4(\iota_X-1)>a_3(\iota_X-1)$. Notice that $\iota_X  = 2$. We can write the defining equations of $X$ as 
\begin{align*}
f\colon \xi x_3 + f' &=0\\
g\colon \xi x_4 + g' &=0
\end{align*}
where 
\begin{itemize}
	\item $X \subset \Proj \mathbb{C}[x_0,\ldots,x_4,\xi] \simeq \mathbb{P}(a_0,\ldots,a_4,a_{\xi})=:\mathbb{P}$ with no specified order on the weights.
	\item The set of all pure monomials in $(x_0,x_1)$ in $f$ and $g$ is $f(x_0,x_1)$ and $g(x_0,x_1)$, respectively. By quasismoothness of $X$,  $f(x_0,x_1) \in f$ and $g(x_0,x_1) \in g$ contain no common factors. See Lemma \ref{lem:awaybase}.
	\item By quasismoothness of $X$, we have $x_2^2 \in f$ or $x_2^2 \in g$. 
	\item By quasismoothness of $X$, $x_4^2 \in f$ if and only if $a_4-a_3=4$.
\end{itemize}
We consider the toric blowup $\Phi \colon T \rightarrow \mathbb{P}$ whose restriction to $X$ is the unique Kawamata blowup $\varphi \colon E \subset Y \rightarrow X$ centred at $\mathbf{p_{\xi}} \in X$. The assumption on the weights of $a_2,\,a_3,\,a_4$ implies that the cone of movable divisors of $T$ decomposes as in Figure \ref{fig:HypCaseI}.

We play the 2-ray game on $T$ and restrict it to $Y$. The small $\mathbb{Q}$-factorial modification $\tau \colon T \rat T_1$ restricts to an isomorphism on $Y$ by lemma \ref{lem:iso}. We continue the game:

\begin{Lem} \label{hypItoricflip}
The map $T_1 \rat T_2$ induced by crossing the $x_2$ wall replaces 
\[
\mathbb{P}\Big( a_2,a_5-\frac{a_1}{2},1,\frac{a_1}{2} \Big) \subset T_1
\]
 by 
\[
\mathbb{P}\Big( \frac{a_2-a_4}{2}, \frac{a_2-a_3}{2} \Big) \subset T_2.
\]
\end{Lem} 

\begin{proof}
Let $a_2' : = \frac{ka_2-\overline{ka_2}}{a_{\xi}}$. Then by an easy computation or direct verification, we have $\gcd(a_2,a_2') = 1$. Hence there are integers $r,\,s$ for which $ra_2+sa_2' = 1$.
Consider the matrix
\[
A=\begin{pmatrix}
r	& s \\ 
-a_2'	& a_2
\end{pmatrix} \in \SL_2(\mathbb{Z}).
\]
Then,
\[
A\cdot \begin{pmatrix}
    a_2          \\
    a_2'       
\end{pmatrix} = \begin{pmatrix}
    1         \\
    0     
\end{pmatrix}
\]

Hence, on a neighborhood of $x_2 \not = 0$, the action can be reduced to the second row of $A\cdot T$. We have the following diagram
\[
        \begin{tikzcd}[ampersand replacement=\&, column sep = 2em]
             T_1   \ar[rr, dashed ] \ar[rd,swap, "f_1" ]  \&  \&T_2 \ar[ld, "g_1" ]\\
						{} \& p_{x_4} \in \mathcal{F}_1  \& {} 
        \end{tikzcd}
    \]
where 
\[
\mathcal{F}_1 = \Proj\bigoplus_{m\geq 1} H^0(T_1,m\mathcal{O}(x_2))=\Proj \mathbb{C}[x_4, \ldots, u_iv_i, \ldots]
\]
where each $u_i$ is a monomial in $(u,\xi,x_0,x_1) $ and each $v_i$ is a monomial in $(x_4,x_3)$. The contractions $f_1$ and $g_1$ are given in coordinates by
\[
(u,\xi,x_0,x_1,x_2,x_4,x_3) \mapsto (x_4, \ldots, u_iv_i, \ldots).
\]
Thus, $f_1$ contracts the locus $(x_4=x_3=0) \subset T_1$ to $p_{x_2} = (1:0\cdots:0) \in \mathcal{F}_1$ and $g_2$ contracts the locus $(u=\xi=x_0=x_1=0) \subset T_2$ to $p_{x_2}$. Easy but tedious computations show that the contracted loci are isomorphic to
\[
\mathbb{P}\Big( a_2,a_5-\frac{a_1}{2},1,\frac{a_1}{2} \Big) \subset T_1 \quad \text{and} \quad \mathbb{P}\Big( \frac{a_2-a_4}{2}, \frac{a_2-a_3}{2} \Big) \subset T_2,
\]
respectively.  As an illustration I show $\mathbb{P}\Big( \frac{a_2-a_4}{2}, \frac{a_2-a_3}{2} \Big) \simeq (u=\xi=x_0=x_1=0) \subset T_2$. It is clear that the locus $(u=\xi=x_0=x_1=0) \subset T_2$ is a weighted projective line $ \mathbb{P}(a,b) \simeq \mathbb{P}^1$. We find $a$ and $b$. As we mentioned, the remaining action $x_4 \not = 0$ is given by the second row of $A \cdot T$. Hence we only need to compute $-a_2'a_2+a_2a_4'$ and $-a_3'a_2+a_2a_3'$ where $a_4':= \frac{ka_4-m_2}{a_{\xi}}$ and $a_3':= \frac{ka_3-m_1}{a_{\xi}}$. We have,
\begin{align*}
-a_2'a_4+a_2a_4' & = -\frac{ka_2-\overline{ka_2}}{a_{\xi}}\cdot a_4 +  a_2\cdot \frac{ka_4-m_2}{a_{\xi}} \\
&=\frac{\overline{ka_2}a_4-m_2a_2}{a_{\xi}} \\
&=\frac{(a_{\xi}\iota_X-a_1)a_4-d_2a_2}{a_{\xi}\iota_X}  && \text{Linearity of}\, \mathbf{p_{\xi}} \implies m_2 = \frac{d_2}{\iota_X},\,\, \overline{ka_2} = a_{\xi}-\frac{a_1}{\iota_X}\\
&= \frac{a_{\xi} \iota_X a_4-a_{\xi}a_2-a_4a_{\xi}}{\iota_X a_{\xi}} &&\text{Linearity of}\, \mathbf{p_{\xi}} \implies d_2 = a_{\xi} + a_4,\,\,a_{\xi}=a_1+a_2\\
&= \frac{a_4(\iota_X-1)-a_2}{\iota_X} \\
&=\frac{a_4-a_2}{2} && \iota_X=2.
\end{align*}
The computation of $b$ is completely analogous with $a_3$ replacing $a_4$. Hence, the small modification $\tau_1$ is 
\[
\bigg(-a_2,-\Big(a_{\xi}-\frac{a_1}{2}\Big),-1,-\frac{a_1}{2},\frac{a_2-a_4}{2}, \frac{a_2-a_3}{2}\bigg).
\]
\end{proof}

\begin{Lem} \label{lem:flipx2}
\begin{enumerate} The following two conditions hold:
	\item Suppose the monomial $x_2^{\alpha}x_{\mu}$ is in $f$ or $g$ where $x_{\mu} \in \{1,x_0,x_1,x_3,x_4\}$. Then it lifts to $x_2^{\alpha}x_{\mu}u^{\beta}$ where, in particular, $\beta = 1$ if $x_{\mu}=1$.
	\item Suppose $x_4^{\alpha} \in f$. Then it lifts to $x_4^2u \in \widetilde{f}$.
\end{enumerate}
\end{Lem}

\begin{proof}
\begin{enumerate}
	\item Suppose wlog, $x_2^{\alpha}x_{\mu} \in f$. By Lemma \ref{lem:lift}, we have $(\widetilde{f}=0) \in |-m_1K_Y|$ for some positive rational number $m_1$. Since $x_2^{\alpha}x_{\mu} \in H^0\big(Y,-m_1K_Y-\alpha'E \big)$ for some $\alpha' >0$, it follows that $x_2^{\alpha}x_{\mu}u^{\alpha'}\in \widetilde{f}$.

Suppose that $x_2^{\alpha} \in f$, that is, that $x_{\mu}=1$. Then, by Lemma \ref{lem:lift},  
\[
x_2^{\alpha}u^{\beta}\in H^0\bigg(Y,-m_1K_Y+ \Big(\beta-\alpha\frac{\iota_X-1}{\iota_X}\Big)E\bigg) 
\]
which in turn implies that $\beta-\alpha\frac{\iota_X-1}{\iota_X} = 0$ or, equivently, $\alpha = 2\beta$ since $\iota_X=2$ in this case. Indeed, for degree reasons and quasismoothness of $X$, $\alpha=2$, hence $\beta=1$.
	\item Suppose $x_4^{2} \in f$.  By Lemma \ref{lem:lift},
	\[
	x_4^2u^{\beta} \in H^0\bigg(Y, -\frac{d_1}{\iota_X}K_Y - \Big(1-\beta\Big)E \bigg).
	\]
	Since $(\widetilde{f}=0) \in |-m_1K_Y|$ for some positive rational number $m_1$, and $x_4^2u^{\beta} \in \widetilde{f}$ we have $\beta = 1$.
\end{enumerate}
\end{proof}

\begin{Prop}
The $x_2$-wall crossing is a terminal hypersurface flip of type 
\[
\bigg(-a_5+\frac{a_1}{2},-1,-\frac{a_1}{2}, \frac{a_2-a_4}{2}, \frac{a_2-a_3}{2};\frac{d_i}{2} \bigg)
\]
where $d_i=d_1$ if $x_2^{\alpha} \in g$ and $d_i=d_2$ if $x_2^{\alpha} \in f$.
\end{Prop}

\begin{proof}
Suppose wlog that $x_2^{\alpha} \in f$. By Lemma \ref{lem:flipx2}, it is possible to eliminate $u$ in an analytical neighborhood around $x_2\not = 0$ and indeed no other variable. The remaining equation is of the form $x_5x_4+g(x_0,x_1)+...=0$ of degree $d_2/\iota_X$. Again by Lemma \ref{lem:flipx2}, no other variable can be eliminated analytically (or globally). Hence we have a hypersurface flip

\[
\bigg(-a_{\xi}+\frac{a_1}{\iota_X},-(\iota_X-1),-\frac{a_1}{\iota_X}(\iota_X-1), \frac{a_2-a_4(\iota_X-1)}{\iota_X}, \frac{a_2-a_3(\iota_X-1)}{\iota_X};\frac{d_2}{\iota_X} \bigg)
\]
which is a type I flip $\tau_1 \colon Y \simeq Y_1 \subset T_1 \rat Y_2 \subset T_2$ according to \cite[Theorem~8]{brown}. Putting $\iota_X=2$ gives the result.
\end{proof}

\begin{Lem} \label{lem:toricendhyp}
The end of the toric 2-ray game is an extremal divisorial contraction of $(x_3=0) \subset T$ to 
\[
\mathbf{p_4} \in \mathbb{P}\Big(a_3,\frac{d_1}{2},1,\frac{a_1}{2}\frac{a_2-a_3}{2},\frac{a_4-a_3}{2}\Big)
\]
with homogeneous variables $u,\,\xi,\,x_0,\,x_1,\,x_2,\,x_4$. Moreover, the point $\mathbf{p_4}$ is an orbifold point if and only if $a_4-a_3=4$.
\end{Lem}

\begin{proof}
As we already observed, the rays 
\[
\mathbb{R}_+\Big[\left(\begin{smallmatrix} a_4 \\ a_4' \end{smallmatrix} \right)\Big] \quad \text{and} \quad \mathbb{R}_+\Big[\left(\begin{smallmatrix} a_3 \\ a_3' \end{smallmatrix} \right)\Big]
\]
are in the boundaries of $\overline{\Mov}(T)$ and $\overline{\Eff}(T)$, respectively. These are indeed distinct rays since 
\[
\det \begin{pmatrix}
a_4 & a_3 \\
a_4' & a_3' 
\end{pmatrix} = \frac{a_3-a_4}{\iota_X} < 0.
\] 
Consider the matrix
\[
A=\begin{pmatrix}
-a_3'	& a_3 \\ 
-a_4'	& a_4
\end{pmatrix} \in \GL_2(\mathbb{Z}).
\]
Then, multiplying by $A$ we get an isomorphism, see \cite[Lemma~2.4, Lemma~2.9]{hamidplia} with
\[
\begin{array}{ccccccc|cccc}
             &       & u  & \xi &   x_0 & x_1 & x_2 & x_4 & x_3 & \\
\actL{T_2}   &  \lBr &  a_3	& \frac{d_1}{\iota_X} & 1 & \frac{a_1}{\iota_X}& \frac{a_2-a_3(\iota_X-1)}{\iota_X}&\frac{a_4-a_3}{\iota_X}&0 &   \actR{.}\\
             &       & a_4	& \frac{d_2}{\iota_X} & 1 & \frac{a_1}{\iota_X}& \frac{a_2-a_4(\iota_X-1)}{\iota_X}&0&-\frac{a_4-a_3}{\iota_X}
 &  
\end{array}
\] 
There is a map 
\[
\Phi' \colon T_2 \rightarrow \mathcal{F}', \quad \mathcal{F}':= \Proj \bigoplus_{m\geq 1} H^0\bigg(T_2,\mathcal{O}\left(\begin{smallmatrix} \frac{a_4-a_3}{\iota_X} \\ 0 \end{smallmatrix} \right)\bigg) \simeq \mathbb{P}\bigg(a_3, \frac{d_1}{\iota_X}, 1, \frac{a_1}{\iota_X}, \frac{a_2-a_3(\iota_X-1)}{\iota_X},\frac{a_4-a_3}{\iota_X}\bigg).
\]
We can write this map in coordinates as 
{\small
\begin{align*}
\Phi' \colon T_2 &\longrightarrow \mathcal{F}' \\
(u,\xi,x_0,x_1,x_2,x_4,x_3) &\longmapsto \Big (x_3^{\frac{\iota_X}{a_4-a_4}\cdot a_4}u : x_3^{\frac{\iota_X}{a_4-a_4}\cdot \frac{d_2}{\iota_X}} \xi :x_3^{\frac{\iota_X}{a_4-a_4}}x_0 :x_3^{\frac{\iota_X}{a_4-a_4}\cdot \frac{a_1}{\iota_X}}x_1 :x_3^{\frac{\iota_X}{a_4-a_4}\cdot \frac{a_2-a_4(\iota_X-1)}{\iota_X}}x_2 :x_4  \Big).
\end{align*}\par}

Hence $\Phi$ contracts the divisor $(x_3=0) \subset T_2$ to $\mathbf{p_4}=(0: \ldots :0: 1) \in \mathcal{F}'$. Moreover, it is an isomorphism away from from $x_3=0$. It is clear that $\mathbf{p_4}$ is a smooth point if and only if $a_4-a_3=\iota_X$. Otherwise, $a_4-a_3 \geq 4$ since the difference $a_4-a_3$ is a positive even integer. Also notice that by looking at each family, we can see that $a_4-a_3$ is at most 4. Hence $\frac{a_4-a_3}{\iota_X} \in \{1,2\}$.
\end{proof}

\begin{Prop} \label{lem:endhyp1}
The end of the 2-ray game on $Y$ restricts to an extremal divisorial contraction $\varphi'$ to a point in a Fano 3-fold hypersurface, 
\[
X'\colon (h'\colon -x_4f'+g' = 0) \subset \mathbb{P}\Big(a_3,1,\frac{a_1}{2},\frac{a_2-a_3}{2},\frac{a_4-a_3}{2}\Big)
\]
of Fano index 1 and homogeneous variables $u,\,\xi,\,x_0,\,x_1,\,x_2,\,x_4$. The divisor is contracted to the point $\mathbf{p_4}$ which is a $\frac{1}{2}(1,1,1)$ cyclic quotient singularity if and only if $a_4-a_3=4$ (or equivalently $d_2-d_1=4$) or a cDV point otherwise, in which case $\varphi'$ has discrepancy 1. The projective line $\Gamma : (x_0=x_1=u=0)$ is in $X'$ and there is a singular point in $\partial_u h'|_{\Gamma} =0$ which is different from $\mathbf{p_4}$.
\end{Prop}

\begin{proof}
We restrict the map given in Lemma \ref{lem:toricendhyp}, which we call $\varphi'$, to the 3-fold $Y_2 \subset T_2$. Recall that
\[
Y_2 \colon (\xi x_3+f'(u,x_0,\ldots, x_4)=\xi x_4+g'(u,x_0,\ldots, x_4)=0) \subset T_2
\] 

and that $\Phi'$ contracts the divisor $x_3=0$ to the point $\mathbf{p_4} \in \mathcal{F}'$. Away from this divisor we can eliminate $\xi$ globally from the first equation which gives us
\[
X' \colon (-x_4f'(u,x_0,x_1,x_2,1,x_4)+ g'(u,x_0,x_1,x_2,1,x_4)=0 ) \subset \mathbb{P}\Big(a_3,1,\frac{a_1}{2},\frac{a_2-a_3}{2},\frac{a_4-a_3}{2}\Big)
\]
with homogeneous variables $u,\,x_0,\,x_1,\,x_2,\,x_4$. The only way for another variable to be eliminated is if $x_3^{\alpha} x_{\mu}\in g'$ for some positive integer $\alpha$. Since $g' \in |-mK_{Y_2}|$ for some $m>0$, and $x_3^{\alpha} \in H^0(Y_2,-\beta_1K_{Y_2}-\beta_2E)$ for some positive integers $\beta_i$, it follows that 
\[
x_3^{\alpha}x_{\mu} \in g' \implies \beta_2=1\,\,\text{and}\,\, x_{\mu}=u.
\]  
Moreover, this happens if and only if $\alpha=\iota_X=2$ which is not possible. Hence, $X'\colon (h'=0)$ is indeed a hypersurface of degree $\frac{d_2}{2}=\frac{a_4+a_{\xi}}{2}$ and therefore by adjunction, $-K_{X'} \sim \mathcal{O}_{X'}(1)$. We now analyse the singular point $\mathbf{p_4}$. 
\paragraph{Suppose $a_4-a_3=4$.} By Lemma \ref{lem:flipx2}, we have $x_4^2u \in f'$ and therefore $-x_4^3u \in h'$. We can then eliminate the variable $u$ locally analytically around $\mathbf{p_4}$ and we have 
\[
\mathbf{p_4} \sim \frac{1}{2}(1,1,1)
\]
with variables $x_0,\,x_1,\,x_2$. In this case, $\varphi'$ is the Kawamata blowup of $\mathbf{p_4}$. 
\paragraph{Suppose $a_4-a_3=2$.} By Lemma \ref{lem:toricendhyp}, we have a weighted blowup with weights $\wt(u,x_0,x_1,x_2)=\big(a_4,1,\frac{a_1}{2},\frac{a_2-a_4}{2} \big)$ and discrepancy $\frac{\iota_X}{a_4-a_3}=1$ by Lemma \ref{lem:discr}. In this case, $\mathbf{p_4}$ is a non-smooth 3-fold Gorenstein terminal singularity and hence a cDV point.

\phantom{.....}\newline Finally, $I_{X'} \subset (u,x_0,x_1)$: In fact, $f' \in |-m_1K_X|$ and $g' \in |-m_2K_X|$, where $m_i >0$ and therefore, monomials in $f'$ or $g'$ are either pure powers of $x_0,\,x_1$ or multiples of $u$. Hence, $ \Gamma \colon (u=x_0=x_1) \subset X'$. From this paragraph, it is also clear that $\partial_{x_{\mu}} h'$ is identically zero on $\Gamma$ if $x_{\mu} \not = u$. Since $X'$ is terminal its singularities are isolated and therefore
\[
\partial_u h'|_{\Gamma}=0
\]
is a non-trivial condition.
\end{proof}

\begin{Ex} \label{ex:97}
Consider a quasismooth member of family 97: $X_{10,14} \subset \mathbb{P}(1,2,2,5,7,9)$ with homogeneous coordinates $x,\,y,\,z,\,t,\,v,\,w$ given by 
\[
wx+t^2+f_{10}(y,z)=wt+v^2+g_{14}(y,z)+x^7v=0. 
\]
Notice that $\Gamma:=X|_{t=x=0}\colon (f_{10}(y,z)=v^2+g_{14}(y,z)=0)$ is a curve in $X$ containing the point $\mathbf{p_w}$. We blow up its only cyclic quotient singularity which is $\mathbf{p_w} \sim \frac{1}{9}(2,2,7) \sim \frac{1}{9}(1,1,8)$ with local variables $y,\,z,\,v$. This is the Kawamata blowup $\varphi \colon Y \rightarrow X$ where $Y$ is given by 
\[
	\begin{cases} 
											wx+t^2u+f_{10}(y,z)=0 &   \\
                      wt+v^2u+g_{14}(y,z)+x^7vu^4=0 & 
                 \end{cases} 
								\]
								inside
								\[
\begin{array}{cccccc|ccccc}
             &       & u  & w &  y & z & v & t & x & \\
\actL{T^- }   &  \lBr &  -2 & -1 &  0 & 0 & 1 & 1 & 1 &   \actR{}\\
             &       & -7 & -8 & -1 & -1 & 0 & 1 & 3 &  
\end{array}
\] 
The Jacobian of $Y$ is 
\[
J:=\begin{pmatrix}
t^2	& x & \partial_yf_5 & \partial_zf_5 & 0 & 2tu &w\\ 
v^2	& t & \partial_yg_7 & \partial_zg_7 & 2vu+x^7u^4 & w &7x^6vu^4
\end{pmatrix}. 
\]
The pre-image of $\Gamma$ under $\varphi$ is
\[
\varphi^{-1}\Gamma = \{u=0\}\cup C
\] 
where $C$ is the curve $\{x=t=f_5(y,z)=v^2u+g_7(y,z)=0\}\subset Y$. Notice that $J|_C$ has rank 2 whether the irrelevant ideal considered is $(u,w)\cap (y,z,v,t,x)$ or $(u,w,y,z)\cap (v,t,x)$. Hence there are no singularities in $C$. However when we flip $C$ over $\mathbf{p_v}$ we introduce a singularity. We have a flip diagram
 \[
        \begin{tikzcd}[ampersand replacement=\&, column sep = 2em]
             C^{-} \subset Y^{-}   \ar[rr, dashed ] \ar[dr, swap] \& \& Y^+  \supset C^{+} \ar[dl ] \\
             \& \mathbf{p_v} \in Z \& 
        \end{tikzcd}
    \] 
where $C^{-}\subset Y^-$ is the fibre of the non-$\mathbb{Q}-$factorial point $\mathbf{p_v} \in Z$ on the left hand side of the flip diagram. The curve $C^{-}$ is the localisation of $C$ to the chart $v=1$. On the right hand side of the diagram we extract a curve $C^{+}$ (which is isomorphic to $\mathbb{P}(1_t,3_x)$) inside
						\[
\begin{array}{ccccccc|cccc}
             &       & u  & w &  y & z & v & t & x & \\
\actL{T^+ }   &  \lBr &  -2 & -1 &  0 & 0 & 1 & 1 & 1 &   \actR{}\\
             &       & -7 & -8 & -1 & -1 & 0 & 1 & 3 &  
\end{array}
\] 
More concretely, $C^{+}$ is given by the localisation of $\{u=w=y=z=0\} \subset Y^+$ at $v=1$ and $C^+$ has a singularity precisely when $x-t^3=0$ as can be seen from the Jacobian restricted to this curve. That is, the point $\mathbf{p}:=(0:0:0:0:1:t_0:t_0^3)$ is singular.

Consider the divisor $x=0$ in $Y^+ \subset T$.  
								\[
\begin{array}{ccccccc|cccc}
             &       & u  & w &  y & z & v & t & x & \\
\actL{T}   &  \lBr &  5 & 7 &  1 & 1 & 1 & 0 & -2 &   \actR{}\\
             &       & 1 & 5 & 1 & 1 & 3 & 2 & 0 &  
\end{array}
\] 
Then $\varphi' \colon Y^+ \rightarrow Z_7$ is the Kawamata blowup of the cyclic quotient singularity $\frac{1}{2}(1,1,1) \in Z_7$ where $Z_7 \colon (h:=(-t^3+v^2)u-f_{10}(y,z)t+g_{14}(y,z)+vu^4=0) \subset \mathbb{P}(1,1,1,2,3)$ with homogeneous variables $u,\,y,\,z,\,t,\,v$. Let $\Gamma \colon (u=y=z=0)$. Clearly, $\Gamma \subset Z_7$ and it contains the image of the Kawamata contraction, the cyclic quotient singularity $\frac{1}{2}(1,1,1)$. Moreover,
\[
\partial_u h|_{\Gamma}=-t^3+v^2 \quad \partial_{x_i} h|_{\Gamma} \equiv 0. 
\]
Hence, the point $\mathbf{p'}:=(0:0:0:t:v) \in \Gamma$ satisfying $-t^3+v^2=0$ is a singular point. We now look at the the preimage of $\Gamma$ under $\varphi$ and we find that $\varphi^{-1} \Gamma = \{ux^{5/2}=yx^{1/2}=zx^{1/2}=0\}=\{x=0\} \cup \{u=w=y=z=0\}$. Above the $1/2$ point, that is the point $v=0$ in $\Gamma$, we have the divisor $x=0$ and above $\Gamma$ the curve $\{u=w=y=z=0\}$ localised at $v=1$ which is the curve $C^+$ that arises from the flip. That is 
\[
(\varphi')^{-1} \Gamma = \{x=0\}\cup C^+.
\] 
Finally we prove that $\varphi(\mathbf{p})=\mathbf{p'}$. The Kawamata blowup is given by
\[
(u:w:y:z:v:t:x) \mapsto (ux^{5/2}:wx^{7/2}:yx^{1/2}:zx^{1/2}:vx^{1/2}:t).
\]
Hence the point $\mathbf{p}$ is sent to $(0:0:0:0:t_0^{3/2}:t_0) \in \Gamma$ which satisfies $v^2-t^3=0$. This is exactly the point $\mathbf{p'}$.


\end{Ex}

\begin{Ex} \label{ex:fam100}
Let $X$ be a quasismooth member of family 100, $X_{12,14} \subset \mathbb{P}(1,2,3,4,7,11)$, given by
\[
X \colon (wx+z^4+t^3+y^6+x^{12}=wz+v^2+ty^5+x^{14}=0).
\]
Then, blowing up its unique cyclic quotient singularity $\mathbf{p_w}$, initiates a Sarkisov link which ends with a divisorial contraction to 
\[
\mathbf{p_z} \in X' \colon (z^4u^2+t^3+y^6+u^6)z+v^2u+ty^5+u^7 \subset \mathbb{P}(1,1,1,2,3)
\]
with homogeneous variables $z,\,u,\,y,\,t,\,v$. Moreover, $\mathbf{p_z}$ is a $cD_4$ point and the divisorial contraction is a weighted blowup of weights $\wt(u,y,t,v)=(3,1,2,2)$ and discrepancy 1.
\end{Ex}

We conclude the following:
\begin{Thm} \label{thm:hypIfull}
Let $X \in I_S$ be such that $a_2 > a_4(\iota_X-1)>a_3(\iota_X-1)$. Then, there is an elementary Sarkisov link to a singular Fano 3-fold orbifold hypersurface $X'$, 
 \[
        \begin{tikzcd}[ampersand replacement=\&,column sep = 2em]
             \& Y \ar[dl, "\varphi"]   \ar[r,"\simeq"] \& Y_1 \ar[rrrrrr, dashed, "{\bigg(-a_5+\frac{a_1}{2},-1,-\frac{a_1}{2}, \frac{a_2-a_4}{2}, \frac{a_2-a_3}{2};\frac{d_i}{2} \bigg)}"] \ar[rrrd, swap]\&\&\& \& \& \& Y_2 \supset (x_3=0) \ar[rd, "\varphi'"] \ar[dlll] \& \\
             \mathbf{p_{\xi}} \in X  \&  \&  \&  \& \& \mathbf{p_2} \&   \& \& \& X' \ni \mathbf{p_4} 
        \end{tikzcd}
    \]
		initiated by the Kawamata blowup  $\varphi \colon Y \rightarrow X$ of a linear cyclic quotient singularity $\mathbf{p_{\xi}}$. Moreover, $\mathbf{p_4} \sim \frac{1}{2}(1,1,1)$ if and only if $a_4-a_3=4$ (or equivalently $d_2-d_1=4$) and $\mathbf{p_4}$ is a cDV point otherwise.
\end{Thm}

\begin{proof}
The 2-ray game on $Y$ follows the game on $T$. Hence, the movable cones of $Y$ and $T$ coincide, $-K_Y \in \Int\overline{\Mov}(Y)$ and $Y$ is a Mori Dream Space. Since all the steps are in the Mori category, it follows that the composition is a Sarkisov link. See Lemma \ref{lem:linkmds}. 
\end{proof}

Notice that by \cite{chel}, the 3-fold $X'$ has to be singular as an orbifold since it is non-rigid. More than that, it also happens that the point we contract to at the end of the Sarkisov link is not the only singularity of the orbifold, in contrast with \cite{okadaI}. This happens because the flipped curve on $Y_3$ carries singularities away from the exceptional divisor $(x_3=0)$. This can be seen very explicitly in Example \ref{ex:97}.

\paragraph{Case II:} Let $X$ be a quasismooth member of a deformation family in the following table and $\mathbf{p_{\xi}} \in X$ a linear cyclic quotient singularity with respect to $X$.

\begin{center}
\begin{tabular}{lccccc} \toprule
    $X$ & 93 & 96  & 102 & 115  \\
    $\mathbf{p_{\xi}}$  & $\frac{1}{5}(1,1,4)$  &$\frac{1}{7}(1,2,5)$ & $\frac{1}{7}(1,1,6)$ & $\frac{1}{11}(1,2,9)$   \\ 
		\bottomrule
		\end{tabular}
\end{center}

The computations and conclusions are absolutely analogous to case I. The only difference is the fact that the base of the flip and the point to which we contract a divisor at the end change. This is reflected in the fact that the movable and effective cones of $T$ are the same as in the one of Figure \ref{fig:HypCaseI}, with the rays generated by $x_2$ and $x_4$ swapped. We give the main theorem and refer the reader to the computations for Case I, which are completely analogous. 

\begin{Thm} \label{thm:hypIIfull}
Let $X \in I_S$ be such that $ a_4(\iota_X-1)>a_2>a_3(\iota_X-1)$. Then, there is an elementary Sarkisov link to a singular Fano 3-fold orbifold hypersurface 
\[
X' \colon (-x_4f'+g'=0) \subset \mathbb{P}\Big(a_3,\frac{d_1}{\iota_X},1,\frac{a_1}{\iota_X},\frac{a_4-a_3}{\iota_X},\frac{a_2-a_3(\iota_X-1)}{\iota_X}\Big)
\] 
with homogeneous variables $u,\,\xi,\,x_0,\,x_1,\,x_4,\,x_2$. The link is
 \[
        \begin{tikzcd}[ampersand replacement=\&,column sep = 2em]
             \& Y \ar[dl, swap, "{\varphi}"]   \ar[r,"\simeq"] \& Y_1 \ar[rr, dashed, "{\tau}"] \ar[rd] \& \& Y_2 \supset (x_3=0) \ar[rd, "{\varphi'}"] \ar[dl] \& \\
             \mathbf{p_{\xi}} \in X    \&  \&  \&  \mathbf{p_4} \&   \&  X' \ni \mathbf{p_2} 
        \end{tikzcd}
    \]
		initiated by the Kawamata blowup of a linear cyclic quotient singularity $\mathbf{p_{\xi}}$ and ending with a divisorial contraction to a cDV point $\mathbf{p_2}$. The map $\tau$ is a flip over a point of type
		\[
		\bigg(a_4,1,\frac{a_1}{2},\frac{-a_4+a_2}{2},\frac{a_3-a_4}{2}; \frac{d_1}{2} \bigg)
		\]
		for the families  of index two, i.e., $X \in \{93,96,102\}$, and a Francia's flip
		\[
		\bigg(1,\frac{a_1}{3},\frac{-a_4+a_2}{3},\frac{a_3-a_4}{3}\bigg)= (-1,-2,1,1)
		\]
		for the index 3 family $X \in \{115\}$.
		Finally, the map $\varphi' \colon Y_2 \rightarrow X'$ is a weighted blowup of weights
		\begin{align*}
		\wt(w,x_0,x_1,x_4)&=\bigg(\overline{ka_2},\iota_X-1,\frac{a_1}{\iota_X}(\iota_X-1), \frac{a_4(\iota_X-1)-a_2}{\iota_X}\bigg)\\
		\wt(u,w,x_0,x_1,x_4)&=\bigg(a_2,\overline{ka_2},\iota_X-1,\frac{a_1}{\iota_X}(\iota_X-1), \frac{a_4(\iota_X-1)-a_2}{\iota_X}\bigg)
		\end{align*}
		for families 93, 102 and 96, 115, respectively and discrepancy
		\[
		\frac{\iota_X(\iota_X-1)}{a_2-a_3(\iota_X-1)}=\begin{cases}
1 &\text{if $X \in \{93,96,102 \}$}\\
2 &\text{if $X \in \{115\}$}
\end{cases}
		\]
\end{Thm}

\begin{Ex}
Consider the Picard rank two quasismooth 3-fold $Y$
\begin{align*}
wx+v^2+vy^2+y^4+z^3u&=0\\
wz+t^3u^2+y^5+u^5x^{15}&=0
\end{align*}
inside
	\[
\begin{array}{ccccccc|cccc}
             &       & u  & w &   y & v & z & t & x & \\
\actL{T}   &  \lBr &  5 & 9 &  2 & 4 & 1 & 0 & -1 &   \actR{.}\\
             &       & 1 & 4 &   1 & 2 & 1 & 1 & 0 &  
\end{array}
\]
The divisor $E' \colon (x=0)$ is isomorphic to
\[
(v^2+vy^2+y^4+z^3u=wz+u^2+y^5=0) \subset \mathbb{P}(5,9,2,4,1)
\]
with homogeneous variables $u,\,w,\,y,\,v,\,z$. It is contracted to the germ
\[
0 \in (u^2+v^2z+\frac{3}{4}zy^4+y^5=0) \subset \mathbb{C}^4
\]
which is a $cD_6$ point.
\end{Ex}

\paragraph{Case III:} Let $X$ be a quasismooth member of a deformation family in the following table and $\mathbf{p_{\xi}} \in X$ a linear cyclic quotient singularity with respect to $X$.


\begin{center}
\resizebox{\textwidth}{!}{\begin{tabular}{lcccccccccccccc} \toprule
    $X$ & 92 & 94 & 95 & 96 & 99& 101& 104& 105& 107& 111& 114 & 120 & 117& 125  \\
    $\mathbf{p}$  & $\frac{1}{3}(1,1,2)$  &$\frac{1}{3}(1,1,2)$ & $\frac{1}{5}(1,2,3)$ & $\frac{1}{5}(1,2,3)$ & $\frac{1}{7}(1,3,4)$& $\frac{1}{5}(1,1,4)$& $\frac{1}{9}(1,4,5)$& $\frac{1}{5}(1,1,4)$& $\frac{1}{7}(1,2,5)$ & $\frac{1}{11}(1,3,8)$& $\frac{1}{4}(1,1,3)$&$\frac{1}{5}(1,1,4)$ &$\frac{1}{8}(1,1,7)$& $\frac{1}{8}(1,1,7)$   \\ 
		\bottomrule
		\end{tabular}}
\end{center}
We can write the defining equations of $X$ as 
\begin{align*}
f\colon \xi x_3 + f' &=0\\
g\colon \xi x_4 + g' &=0
\end{align*}
where 
\begin{itemize}
	\item $X \subset \Proj \mathbb{C}[x_0,\ldots,x_4,\xi] \simeq \mathbb{P}(a_0,\ldots,a_4,a_{\xi})=:\mathbb{P}$ with no specified order on the weights.
	\item The set of all pure monomials in $(x_0,x_1)$ in $f$ and $g$ is $f(x_0,x_1)$ and $g(x_0,x_1)$, respectively. By quasismoothness of $X$,  $f(x_0,x_1) \in f$ and $g(x_0,x_1) \in g$ contain no common factors. See Lemma \ref{lem:awaybase}.
	\item By quasismoothness of $X$, either $x_4^{\alpha} \in f$ or $x_4x_2 \in f$. 
\end{itemize}

\begin{figure}%
\centering
\begin{tikzpicture}[scale=3]
 \coordinate (A) at (0, 0);
 \coordinate [label={left:$u$}] (E) at (0, 1);
 \coordinate [label={above:$x_0,\,x_1$}] (K) at (0.7, 1);
 \coordinate [label={above:$\xi$}] (5) at (0.3,1.1);
 \coordinate [label={right:$x_4$}] (2) at (1.2,0.8);
 \coordinate [label={right:$x_3$}] (4) at (1.2,0.5);
 \coordinate [label={right:$x_2$}] (3) at (1.2,0.2);

\coordinate [label={$T$}] (T) at (0.35,0.6);
\coordinate [label={right:$T_1$}] (T1) at (0.55,0.6);
\coordinate [label={$T_2$}] (T) at (0.8,0.33);

  \draw (A) -- (E);
  \draw (A) -- (K);
	\draw [very thick,color=red] (A) -- (5);
	\draw (A) -- (2);
	\draw (A) -- (3);
	\draw [very thick,color=red] (A) -- (4);
  
\end{tikzpicture}
\caption{A representation of the chamber decomposition of the cone of effective divisors of $T$. In red the subcone of movable divisors of T. The ample models $T$, $T_1$ and $T_2$ are represented in the interior of the corresponding Nef chamber.}%
\label{fig:Cod2CaseI}%
\end{figure}
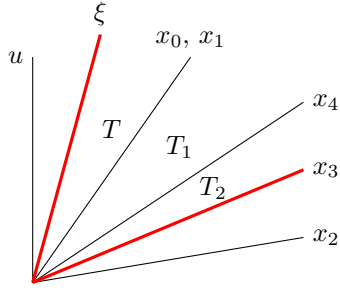

We consider the toric blowup $\Phi \colon T \rightarrow \mathbb{P}$ whose restriction to $X$ is the unique Kawamata blowup $\varphi \colon E \subset Y \rightarrow X$ centred at $\mathbf{p_{\xi}} \in X$. The assumption on the weights of $a_2,\,a_3,\,a_4$ implies that the cone of movable divisors of $T$ decomposes as in Figure \ref{fig:Cod2CaseI}.

\begin{Rem} \label{rem:fake}
If $X$ is a quasismooth member of family 117 then $\mathbf{p_{\xi}} \sim \frac{1}{7}(1,1,6) \in X$ is a linear cyclic quotient singularity satisfying the above conditions on the weights. However, in this case, there is an extra $\boldsymbol{\mu_2}$-action on $T$. See Section \ref{subsec:fake}.
\end{Rem}

We play the 2-ray game on $T$ and restrict it to $Y$. Notice that the small $\mathbb{Q}$-factorial modification $\tau \colon T \rat T_1$ restricts to an isomorphism on $Y$ due to Lemma \ref{lem:iso}. We continue the 2-ray game from $T_1$:

\begin{Lem}
The map $T_1 \rat T_2$ induced by crossing the $x_4$ wall replaces 
\[
\mathbb{P}\Big( a_4,\frac{d_2}{\iota_X},1,\frac{a_1}{\iota_X} \Big) \subset T_1
\]
 by 
\[
\mathbb{P}\Big( \frac{a_4-a_3}{\iota_X}, \frac{a_4(\iota_X-1)-a_2}{\iota_X} \Big) \subset T_2.
\]
\end{Lem} 

\begin{proof}
This is analogous to the proof of Lemma \ref{hypItoricflip}, so we are brief. On a neighborhood of $x_4 \not = 0$, the action can be reduced to the second row of $A\cdot T_1$. We have the following diagram
\[
        \begin{tikzcd}[ampersand replacement=\&, column sep = 2em]
             T_1   \ar[rr, dashed ] \ar[rd,swap, "f_1" ]  \&  \&T_2 \ar[ld, "g_1" ]\\
						{} \& p_{x_4} \in \mathcal{F}_1  \& {} 
        \end{tikzcd}
    \]
where 
\begin{align*}
\mathcal{F}_1 &= \Proj\bigoplus_{m\geq 1} H^0(T_1,m\mathcal{O}(x_2))\\
&=\Proj \mathbb{C}[x_4, \ldots, \text{monomial in}\,  (u,\xi,x_0,x_1) \times \text{monomial in}\, (x_3,x_2), \ldots]
\end{align*}
 The result follows by computing $T_1|_{(x_3=x_2=0)}$ and $T_2|_{(u=\xi=x_0=x_1=0)}$.
\end{proof}

\begin{Lem} \label{lem:flipx4}
If $x_4^{\alpha} \in f$, it lifts to $x_4^{\alpha}u$ if and only if $\alpha = \iota_X$. Otherwise, we have $x_4x_2 \in f$, which lift to $x_4x_{2}u$.
\end{Lem}

\begin{proof}
Recall that $\widetilde{f} \in |-m_1K_Y|$ for some $m_1 \in \mathbb{Q}_{>0}$ by Lemma \ref{lem:lift}. Suppose that $x_4^{\alpha} \in f$. Then, $\alpha a_4=d_1$ and 
\[ 
x_4^{\alpha}u^{\beta} \in H^0\bigg(Y,-m_1K_Y+ \Big(\beta - \alpha \frac{1}{\iota_X}\Big)E\bigg).
\]
Hence, we have $x_4^{\alpha}u \in \widetilde{f}$ if and only if $\beta=1$ and $\alpha = \iota_X$.

Suppose, on the other hand, that $x_4^{\alpha} \not \in f$. Then, by quasismoothness at $\mathbf{p_4}$, we have that $x_4^{\beta}x_{\mu} \in f$. By inspection we conclude that $\beta=1$ and $\mu=2$, that is, $x_4x_2 \in f$. This monomial lifts to 
\begin{align*}
 f^* \ni x_4x_2u^{\gamma}  &\in H^0\bigg(Y,-m_1K_Y - \Big(\frac{1}{\iota_X}+\frac{\iota_X-1}{\iota_X} - \gamma \Big)E\bigg) \\
&= H^0\big(Y,-m_1K_Y - (1-\gamma)E).  
\end{align*}
Hence, $\gamma=1$ as expected.
\end{proof}

\begin{Prop} \label{prop:cod2caseIflip}
The $x_4$-wall crossing is a toric flip of type
\[
\bigg(-1,-\frac{a_1}{\iota_X}, 1, -1+\frac{a_1}{\iota_X}\bigg)
\]
if and only if $x_4^{\iota_X} \in f$. Otherwise it is a hypersurface flip of type 
\[
\bigg(-a_4,-1,-\frac{a_1}{\iota_X}, \frac{a_1}{\iota_X}, \frac{a_4(\iota_X-1)-a_2}{\iota_X};\frac{d_1}{\iota_X}\bigg).
\]
\end{Prop}

\begin{proof}
Is is clear from Lemma \ref{lem:flipx4} that the only way to eliminate $u$ in an analytical neighborhood of $(x_4\not = 0)$ is when $x_4^{\iota_X} \in f$ and that no other variables can be eliminated using $f$. Moreover, $\xi$ can always be globally eliminated from $g$. In each case for which $x_4^{\iota_X} \in f$, we observe that the difference of the degrees of the equations defining $X$ is $\iota_X$ or, equivalently, $(a_4-a_3)/\iota_X=1$. Hence, 
\[
\frac{a_4(\iota_X-1)-a_2}{\iota_X}=\frac{a_5+a_3 -a_4-a_2}{\iota_X}=-1+\frac{a_5-a_2}{\iota_X}=-1+\frac{a_1}{\iota_X}
\]
and we conclude that this is a flip by \cite[Theorem~7]{brown}.

Suppose that $x_4^{\alpha} \not \in f$. Recall that, in that case, $x_4x_2 \in f$. Then we end up with the flip $x_{2}u+f(x_0,x_1)+\cdots=0$ which has degree
\[
\deg (x_2u)=-\frac{a_4+a_2}{\iota_X}=-\frac{a_5+a_3}{\iota_X}= -\frac{d_1}{\iota_X}.
\]
Moreover, $a_4-a_3=a_5-a_2=a_1$. Hence we have a type I flip according to \cite[Theorem~8]{brown}.
\end{proof}

\begin{Lem} \label{lem:toricendcod2}
The end of the toric 2-ray game is an extremal divisorial contraction of $(x_2=0) \subset T_2$ to 
\[
\mathbf{p_3} \in \mathbb{P}\bigg(a_2,a_{\xi}-\frac{a_1}{\iota_X},\iota_X-1,\frac{(\iota_X-1)a_1}{\iota_X},\frac{a_4(\iota_X-1)-a_2}{\iota_X},\frac{a_3(\iota_X-1)-a_2}{\iota_X}\bigg)
\]
with homogeneous variables $u,\,\xi,\,x_0,\,x_1,\,x_4,\,x_3$. Moreover, the point $\mathbf{p_3}$ is an orbifold point if and only if $a_3(\iota_X-1)-a_2 = 4$.
\end{Lem}

\begin{proof}
This is analogous to Lemma \ref{lem:toricendhyp}. The rays 
\[
\mathbb{R}_+\Big[\left(\begin{smallmatrix} a_3 \\ a_3' \end{smallmatrix} \right)\Big] \quad \text{and} \quad \mathbb{R}_+\Big[\left(\begin{smallmatrix} a_2 \\ a_2' \end{smallmatrix} \right)\Big]
\]
are in the boundaries of $\overline{\Mov}(T)$ and $\overline{\Eff}(T)$, respectively. These are indeed distinct rays since 
\[
\det \begin{pmatrix}
a_3  & a_2 \\
a_3' & a_2' 
\end{pmatrix} = \frac{a_2-a_3(\iota_X-1)}{\iota_X} < 0.
\] 
Consider the matrix
\[
A=\begin{pmatrix}
-a_2'	& a_2 \\ 
-a_3'	& a_3
\end{pmatrix} \in \GL_2(\mathbb{Z}).
\]
Then, multiplying by $A$ we get an isomorphism, see \cite[Lemma~2.4, Lemma~2.9]{hamidplia} with
\[
\begin{array}{ccccccc|ccccc}
             &       & u  & \xi &   x_0 & x_1 & x_4 & x_3 & x_2 & \\
\actL{T_2}   &  \lBr &  a_2&a_{\xi}-\frac{a_1}{\iota_X}&\iota_X-1&(\iota_X-1)\frac{a_1}{\iota_X}&\frac{a_4(\iota_X-1)-a_2}{\iota_X}&\frac{a_3(\iota_X-1)-a_2}{\iota_X}&0  & \actR{.}\\
             &       & a_3	& \frac{d_1}{\iota_X} & 1 & \frac{a_1}{\iota_X}& \frac{a_4-a_3}{\iota_X}&0&-\frac{a_3(\iota_X-1)-a_2}{\iota_X} &
 &  
\end{array}
\] 
We notice that $X$ is such that each row of the weight system of $T_2$ has no common divisor. In particular this is immediate when $\iota_X = 2$ as can be seen from the bidegree of the section $x_0$. See Remark \ref{rem:fake}. There is a map 
\[
\Phi' \colon T_2 \rightarrow \mathcal{F}'
\]
where,
\begin{align*}
 \mathcal{F}' &:= \Proj \bigoplus_{m\geq 1} H^0\bigg(T_2,\mathcal{O}\left(\begin{smallmatrix} \frac{a_3(\iota_X-1)-a_2}{\iota_X} \\ 0 \end{smallmatrix} \right)\bigg) \\
&\simeq \mathbb{P}\bigg(a_2,a_{\xi}-\frac{a_1}{\iota_X},\iota_X-1,(\iota_X-1)\frac{a_1}{\iota_X},\frac{a_4(\iota_X-1)-a_2}{\iota_X},\frac{a_3(\iota_X-1)-a_2}{\iota_X}\bigg).
\end{align*}
Let $d := \frac{\iota_X}{a_3(\iota_X-1)-a_2}$. We can write this map in coordinates as 
{\footnotesize
\begin{align*}
\Phi' \colon T_2 &\longrightarrow \mathcal{F}' \\
(u,\xi,x_0,x_1,x_2,x_4,x_3) &\longmapsto \Big (x_2^{d\cdot a_3}u : x_2^{d \cdot \frac{d_1}{\iota_X}} \xi :x_2^{d}x_0 :x_2^{d\cdot \frac{a_1}{\iota_X}}x_1 :x_2^{d\cdot \frac{a_4-a_3}{\iota_X}}x_4 :x_3  \Big).
\end{align*}\par}

Hence $\Phi'$ contracts the divisor $(x_2=0) \subset T_2$ to $\mathbf{p_3}=(0: \ldots :0: 1) \in \mathcal{F}'$. Moreover, it is an isomorphism away from from $x_2=0$. It is clear that $\mathbf{p_3}$ is a smooth point if and only if $a_3(\iota_X-1)-a_2=\iota_X$. Otherwise, $a_3(\iota_X-1)-a_2 \geq 2\iota_X$ since the difference $a_3(\iota_X-1)-a_2$ is a positive multiple of $\iota_X$. Also notice that, exactly as in Lemma \ref{lem:toricendhyp}, by looking at each family, we can see that $a_3(\iota_X-1)-a_2$ is at most 4. Hence $\frac{a_3(\iota_X-1)-a_2}{\iota_X} \in \{1,2\}$.
\end{proof}

\begin{Lem} \label{lem:auxcod2}
There is no variable $x_{\mu}$ such that $x_2^{\alpha}x_{\mu} \in \widetilde{f}$ or $\widetilde{g}$.
\end{Lem}

\begin{proof}
Recall that $ \widetilde{f}\in |-m_1K_Y|$ for some $m_1 \in \mathbb{Q}_{>0}$ by Lemma \ref{lem:lift}. We have by assumption $x_2^{\alpha} \in H^0(Y,-m_1K_Y-m_2E)$, where $m_i \in \mathbb{Q}_{>0}$. Hence if $x_2^{\alpha}x_{\mu} \in \widetilde{f}$ we have $x_{\mu} \in H^0(Y,-m_1'K_Y+m_2E)$, where $m_1' \in \mathbb{Q}_{>0}$. Hence, $x_{\mu} \in \{x_5,u \}$. Since $x_5x_3$ is the only monomial in $\widetilde{f}$ featuring $x_5$, we conclude that $x_{\mu}=u$.

If $x_2^{\alpha} \in f$ then,
\[
\widetilde{f} \ni x_2^{\alpha}u^{\gamma} \in H^0\bigg(Y,-\alpha\frac{a_2}{\iota_X}K_Y-\Big(\alpha\frac{\iota_X-1}{\iota_X}-\gamma\Big)E\bigg).
\]
Hence, $\alpha = \frac{\iota_X}{\iota_X-1}\gamma$ and if $\gamma = 1$, it follows from the fact that $\alpha$ is a (positive) integer that $\iota_X=2=\alpha$. However, by inspection, $\alpha \geq 4$.
\end{proof}

\begin{Prop} \label{prop:endcod2}
The end of the 2-ray game restricts to an extremal divisorial contraction $\varphi' \colon Y_2 \rightarrow X'$ in the Mori Category to a Fano 3-fold complete intersection
\[
X'\colon (\xi x_3+f'=\xi x_4+g'=0) \subset \mathbb{P}\bigg(a_2,a_{\xi}-\frac{a_1}{\iota_X},\iota_X-1,(\iota_X-1)\frac{a_1}{\iota_X},\frac{a_4(\iota_X-1)-a_2}{\iota_X},\frac{a_3(\iota_X-1)-a_2}{\iota_X}\bigg)
\]
of Fano index $\iota_X-1$ and homogeneous variables $u,\,\xi,\,x_0,\,x_1,\,x_4,\,x_3$. The divisor is contracted to the point $\mathbf{p_3} \in X'$ which is either
\begin{itemize}
	\item A terminal singularity of index 2 which is a $\frac{1}{2}(1,1,1)$ cyclic quotient singularity or a $cAx/2$. In this case $\varphi'$ has discrepancy $\frac{1}{2}$.
	\item A terminal singularity of index 1 (i.e., a compound du Val singularity). In this case $\varphi'$ has discrepancy $\frac{\iota_X}{a_3(\iota_X-1)-a_2}=1$.  
\end{itemize}
The rational curve $\Gamma : (u=\xi=x_0=x_1=0)$ is in $X'$ and contains the singular points of $X'$.
\end{Prop}

\begin{proof}
We restrict the map given in Lemma \ref{lem:toricendcod2}, which we call $\varphi'$, to the 3-fold $Y_3 \subset T_3$. Recall that
\[
Y_2 \colon (\xi x_3+f'(u,x_0,\ldots, x_4)=\xi x_4+g'(u,x_0,\ldots, x_4)=0) \subset T_2
\] 
and that $\varphi'$ contracts the divisor $x_2=0$ to the point $\mathbf{p_3} \in X'$ so that $\rho_{X'}=1$. By Lemma \ref{lem:auxcod2}, no variable can be eliminated away from $x_2=0$. Hence, 
\[
X'\colon (\xi x_3+f'(u,x_5,x_0,x_1,1,x_3,x_4)=\xi x_4+g'(u,x_5,x_0,x_1,1,x_3,x_4)=0) \subset \mathbb{P}'
\]
where $\mathbb{P}'$ is as in Lemma \ref{lem:toricendcod2} with homogeneous variables $u,\,\xi,\,x_0,\,x_1,\,x_4,\,x_3$. Hence $X'$ is indeed a codimension 2 3-fold given by equations of degrees $d_1\frac{\iota_X-1}{\iota_X}$ and $d_2\frac{\iota_X-1}{\iota_X}$ and by adjunction $-K_{X'} \sim \mathcal{O}_{X'}(\iota_X-1)$. That is, $-K_{X'}$ is ample since $\iota_X \geq 2$ and $X'$ is Fano.

We now analyse the singular point $\mathbf{p_3}$. By definition of $X'$, the variable $\xi$ can be eliminated in a neighborhood of $\mathbf{p_3}$. We divide it in two cases:

\begin{enumerate}
	\item Suppose $a_3(\iota_X-1)-a_2=\iota_X$. This is the case for families
	\[
	X \in \{92,95,96,99,101,107,111,114,117,125 \} \subset I.
	\]
	In this case, no other variable can be locally analytically eliminated. We conclude that $\varphi'$ is a weighted blowup of weights $\wt(u,x_0,x_1,x_4)=\big(a_3,1,\frac{a_1}{\iota_X},\frac{a_4-a_3}{\iota_X} \big)$ to a terminal Gorenstein point $\mathbf{p_3}$. The discrepancy can be readily computed by Lemma\ref{lem:discr} and it equals $\frac{\iota_X}{a_3(\iota_X-1)-a_2}$ which is $1$ by the case assumption we are in.
	\item Suppose If $a_3(\iota_X-1)-a_2=2\iota_X$. This is the case for families
	\[
	X \in \{94,104,105,120\} \subset I.
	\]
	For $X \in \{94,105,120 \}$ we have $x_3^{\iota_X}u \in g'$, proving that $\mathbf{p_3} \sim \frac{1}{2}(1,1,1)$ with variables $x_0,\,x_1,\,x_4$ and $\varphi'$ is the unique Kawamata blowup centred at $\mathbf{p_3}$. Otherwise, $X$ is a quasismooth member of family 104 and $\mathbf{p_3}$ is analytically equivalent to
	\[
	0 \in (u^2+v^2+h_{\geq 4}(y,t)=0) \subset \mathbb{C}^4_{uvyt}\bigg/\frac{1}{2}(1,0,1,1).  
	\]
	That is, $\mathbf{p_3}$ is a $cAx/2$ and $\varphi'$ is the weighted blowup $\wt(u,v,y,t)=\frac{1}{2}(5,4,1,1)$ and discrepancy $\frac{1}{2}$. See \cite[Section~8]{HayakawaI} for the classification of extremal extractions from $cAx/2$ points with discrepancy $\frac{1}{2}$.
\end{enumerate}
Finally, $I_{X'} \subset (u,\xi,x_0,x_1)$: In fact, $f' \in |-m_1K_{X'}|$ and $g' \in |-m_2K_{X'}|$ where $m_i >0$, and, therefore, monomials in $f'$ or $g'$ are either pure powers of $x_0,\,x_1$ or multiples of $u$ or $\xi$. Hence, $ \Gamma \colon (u=\xi=x_0=x_1) \subset X'$. From this paragraph, it is also clear that
\[
J(X')|_{\Gamma}=
\begin{pmatrix}
f'(0,0,0,0,1,x_3,x_4) & x_3 & 0 & 0 & 0 & 0 \\
g'(0,0,0,0,1,x_3,x_4) & x_4 &  0 & 0 & 0 & 0 
\end{pmatrix}
\]
where the columns of the matrix are the partial derivatives taken in the order $\{u,\,\xi,\,x_0,\,x_1,\,x_4,\,x_3 \}$. Since $X'$ is terminal its singularities are isolated and therefore
\[
x_4f'(0,0,0,0,1,x_3,x_4)-x_3g'(0,0,0,0,1,x_3,x_4)=0
\]
is a non-trivial condition and its finite solutions are the singular points of $X$.
\end{proof}

Hence, we conclude the following:
\begin{Thm} \label{thm:cod2caseI}
Let $X \in I_S$ together with the linear cyclic quotient singularity $\mathbf{p_{\xi}}$ such that $a_4(\iota_X-1)>a_3(\iota_X-1)>a_2$. Suppose that $(X,\mathbf{p_{\xi}})$ is not a member of family 117 with $\mathbf{p_{\xi}} \sim \frac{1}{7}(1,1,6)$. Then, there is an elementary Sarkisov link to a Fano 3-fold complete intersection $X'$, 
 \[
        \begin{tikzcd}[ampersand replacement=\&,column sep = 2em]
             \& Y_1 \ar[dl, "\varphi"]   \ar[r,"\simeq"] \& Y_2 \ar[rr, dashed, "\tau"] \ar[rd, swap] \& \& Y_3 \supset (x_2=0) \ar[rd,"\varphi'"] \ar[dl] \& \\
             \mathbf{p_{\xi}} \in X   \&    \& \& \mathbf{p_4}     \& \& X' \ni \mathbf{p_3} 
        \end{tikzcd}
    \]
		initiated by the Kawamata blowup $\varphi$ of a linear cyclic quotient singularity $\mathbf{p_{\xi}} \in X$. Moreover, the singular point $\mathbf{p_3}$ is either a compound du Val singularity or it has index 2, in which case it is the cyclic quotient singularity  $\frac{1}{2}(1,1,1)$ or a $cAx/2$ point. 
\end{Thm}

\begin{proof}
For an explanation on the exclusion of $(X, \frac{1}{7}(1,1,6))$ where $X$ is a member of family 117 see Remark \ref{rem:fake}.
The 2-ray game on $Y$ follows the game on $T$. Hence, $Y$ is a Mori Dream Space, the movable cones of $Y$ and $T$ coincide and $-K_Y \in \Int\overline{\Mov}(Y)$. Since the small modification $\tau$ is a flip, it follows that the composition is a Sarkisov link. See Lemma \ref{lem:linkmds}.
\end{proof}


\paragraph{Case IV: The exceptional cases:} The cases treated in this section are the ones in the following tables
\begin{center}
\begin{tabular}{lcccccccccccc} \toprule
    $X$ & 102 & 105  & 116 & 118 & 124 & 124   \\
    $\mathbf{p_{\xi}}$  & $\frac{1}{3}(1,1,2)$  &$\frac{1}{3}(1,1,2)$ & $\frac{1}{2}(1,1,1)$ & $\frac{1}{5}(1,2,3)$ & $\frac{1}{7}(1,3,4)$& $\frac{1}{11}(1,4,7)$  \\ 
		\bottomrule
		\end{tabular}
\end{center}
\begin{center} \resizebox{\textwidth}{!}{\begin{tabular}{lcccccccccc} \toprule
    $X$ & 95 & 98  & 107 & 108 & 110 & 114 & 121 &122 & 123 & 125\\
    $\mathbf{p}$  & $\frac{1}{3}(1,1,2)$  &$\frac{1}{3}(1,1,2)$ & $\frac{1}{5}(1,1,4)$ & $\frac{1}{5}(1,1,4)$  & $\frac{1}{7}(1,1,6)$ & $\frac{1}{2}(1,1,1)$ &$\frac{1}{3}(1,1,2)$ &$\frac{1}{5}(1,2,3)$ &$\frac{1}{9}(1,4,5)$   & $\frac{1}{2}(1,1,1)$ \\ 
		\bottomrule
		\end{tabular}}
\end{center}

These are the links initiated by blowing up a \emph{non-linear} cyclic quotient singularity or the ones where $X$ has non-prime Fano index, i.e., $\iota_X=4$. We prove that each quasismooth member of the first table is birational to a Fano 3-fold hypersurface  and each quasismooth member of the second table is birational to a Fano 3-fold complete intersection of two hypersurfaces.

\begin{figure}%
\centering
\begin{tikzpicture}[scale=3]
 \coordinate (A) at (0, 0);
 \coordinate [label={left:$u$}] (E) at (0, 1);
 \coordinate [label={above:$\xi$}] (5) at (0.3,1.1); 
 \coordinate [label={above:$v$}] (K) at (0.7, 1);
 \coordinate [label={right:$x_0,\,x_1$}] (2) at (1.2,0.8);
 \coordinate [label={right:$w$}] (4) at (1.2,0.5);
 \coordinate [label={right:$t$}] (3) at (1.2,0.2);

\coordinate [label={$T$}] (T) at (0.35,0.6);
\coordinate [label={right:$T_1$}] (T1) at (0.55,0.6);
\coordinate [label={$T_2$}] (T) at (0.8,0.33);

  \draw (A) -- (E);
  \draw (A) -- (K);
	\draw [very thick,color=red] (A) -- (5);
	\draw (A) -- (2);
	\draw (A) -- (3);
	\draw [very thick,color=red] (A) -- (4);
  
\end{tikzpicture}
\caption{A representation of the chamber decomposition of the cone of effective divisors of $T$. In red the subcone of movable divisors of T. The ample models $T$, $T_1$ and $T_2$ are represented in the interior of the corresponding Nef chamber.}%
\label{fig:HypExcpI}%
\end{figure}
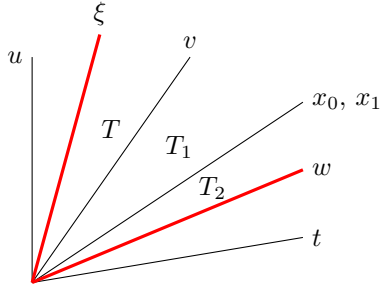

\paragraph{Families 102 and 105.} In the following, we deal with families 102 and 105. Any quasismooth member $X$ of each of these families contains a cyclic quotient singularity of index 3 and $X \subset \mathbb{P}(2,2,3,5,7,a_w)$ of homogeneous variables $x_0,\,x_1,\,\xi,\,t,\,v,\,w$, where $a_w$ is 7 or 9 if $X$ is a member of family 102 or 105, respectively. Assume that $\xi^3t \in g$, see Remark \ref{rem:102}, and for family 105 assume $tv \in f$, see Proposition \ref{prop:fake2}. It is easy to see that after a change of variables any such member can be written as
\begin{align*}
\xi w+tx_{\mu}+f_{d_1}&=0\\
\xi^3t+\xi^2g_{d_2-2a_{\xi}}+\xi g_{d_2-a_{\xi}}+v^2+wx_{\mu'}+g_{d_2}&=0
\end{align*}
where $x_{\mu}=t$ for $X_{102}$ and $x_{\mu}=v$ for $X_{105}$, and the opposite for $x_{\mu'}$. Moreover, $\xi,\,x_{\mu} \not \in f_{d_1}$. Hence, $\mathbf{p_{\xi}}\sim \frac{1}{3}(1,1,2)$ with local variables $(x_0,\,x_1,v)$. The Kawamata blowup of $X$ centred at $\mathbf{p_{\xi}}$ is an extremal weighted blowup $\varphi \colon Y\rightarrow X$ and $-K_Y \sim \varphi^*(-K_X)-\frac{1}{a_{\xi}}E$. Locally, at $\mathbf{p_{\xi}}$, the Kawamata blowup is the graph of the rational map to $\mathbb{P}(1,1,2)$ given by 
\[
\bigg(\frac{x_0^2}{\xi}:\frac{x_1^2}{\xi}:\frac{v^2}{\xi^4}\bigg).
\] 
Hence, by Corollary \ref{cor:genlift},
\[
\xi \in H^0\bigg(Y,-\frac{3}{2}K_Y+\frac{1}{3}E\bigg),\,v\in H^0\bigg(Y,-\frac{7}{2}K_Y+\frac{1}{2}E\bigg),\,x_0,\,x_1 \in H^0(Y,-K_Y).
\]
Since the section $v$ vanishes with order $\frac{2}{3}$ at $E$ and $v^2 \in g$ it follows that $m_g \leq \frac{4}{3}$. However, in each case $d_2 - 2 a_{\xi} = 8 > 2$, and therefore there are no linear terms in $g_{d_2-2a_{\xi}}$ nor in $g_{d_2-a_{\xi}}$.  Hence, $m_g = \frac{4}{3}$ and $\widetilde{g} \in |-7K_Y+E|$. We can now compute how the section $t \in H^0(X,5A)$ lifts. By Corollary, \ref{cor:genlift}, there are integers such that $t \in H^0(Y,-\frac{a_t}{2}K_Y+\frac{\alpha_t}{2}E)$. Since $\xi^3t \in g^*$ we have
\[
a_t = 7\cdot 2 - 3\cdot a_{\xi}=5, \quad \alpha_t = 1\cdot 2 - 3\cdot 1=-1.
\]
On the other hand, $f_{d_1} \in \mathbb{C}[x_0,x_1]$ and we have seen that $x_0,\,x_1$ vanish with order $\frac{1}{3}$. The monomials in $f_{d_1}$ are of the form $x_0^{\alpha_0}x_1^{\alpha_1}$ where $\alpha_i$ satisfy $\alpha_0a_1+\alpha_1a_1=2(\alpha_0+\alpha_1)=d_1$.  Hence, when pulled back by $\varphi$, the divisor $f \in |-\frac{d_1}{\iota_X}K_X|$ vanishes at $E$ with order at least 
\[
\alpha_0\cdot \frac{1}{3} + \alpha_1 \cdot \frac{1}{3} =   \frac{d_1}{6}.
\]
Similarly, $tx_{\mu}$ vanishes at $E$ with order $\frac{d_1}{6}$. Hence, $m_f=\frac{d_1}{6}$ and $\widetilde{f} \in |-\frac{d_1}{2}K_Y|$ by Corollary \ref{cor:genlift}. Finally,
\[
w \in H^0\bigg(Y,-\frac{a_w}{2}K_Y-\frac{1}{2}E\bigg).
\]
Now that we have computed how all section $x_i \in H^0(X,a_iA)$ lift we know how to construct the ambient toric variety $T$ where $Y\colon (\widetilde{f}=\widetilde{g}=0)$ lives. This is done as in Lemma \ref{lem:utbl} and its effective and movable cones are as in Figure \ref{fig:HypExcpI}.

We now use variation of GIT on the movable cone of the ambient toric variety of $Y$.
\begin{Lem} \label{lem:HypExcp1stwall}
The small $\mathbb{Q}$-factorial modification $T \rat T_1$ induced by crossing the $v$-wall restricts to an isomorphism on $Y$.
\end{Lem} 
\begin{proof}
The small modification is the diagram
 \[
        \begin{tikzcd}[ampersand replacement=\&, column sep = 2em]
             T  \ar[rr, dashed ] \ar[dr, swap, "\displaystyle{f_0}" ] \& {} \& T_1 \ar[ld, "\displaystyle{g_0}" ] \\
             {} \& \mathcal{F} \& {}
        \end{tikzcd}
    \]
		where 
		\[
		f_0 \colon T \rightarrow \mathcal{F}, \quad  \mathcal{F} = \Proj \bigoplus_{m\geq 1}H^0\big (T, m\mathcal{O}(v) \big).
		\]
That is, the map $f_0$ is given by the sections of multiples of $v$ and it contracts the locus $ \Gamma \colon (x_0=x_1=w=x=0) \subset T$ to a point in $\mathcal{F}$. 

On $Y$ the locus contracted by $f$ corresponds to the intersection of divisors $D_i \in |-m_iK_Y+n_iE|$ where $m_i>0$ and $n_i\geq 0$. In particular, $\widetilde{f}|_{\Gamma}$ vanishes identically and $\widetilde{g}|_{\Gamma} =v^2=0$. Hence, the locus on $Y$ is $(v=x_0=x_1=w=x=0) \subset Y$ which is empty.
\end{proof}

The difference to the non-exceptional cases, i.e., when the Sarkisov link is initiated by the blowup of a \emph{linear} cyclic quotient singularity, is that, in this case, we have a number of flops rather than an isomorphism:

\begin{Lem} \label{lem:flops}
In the setting above, crossing the wall generated by the characters $x_0,\,x_1$ induces a flop on $Y$ over $s$ points where $s$ is the number of components of $\widetilde{f}|_{\xi=t=v=w=0}=0$.
\end{Lem} 
\begin{proof}
Crossing the $x_0,\,x_1$ wall induces a birational map

 \[
        \begin{tikzcd}[ampersand replacement=\&, column sep = 2em]
             T_1  \ar[rr, dashed, "\displaystyle{\tau}" ] \ar[dr, swap, "\displaystyle{f_1}" ] \& {} \& T_2 \ar[ld, "\displaystyle{g_1}" ] \\
             {} \& \mathcal{F}_1 \& {}
        \end{tikzcd}
    \]
		where 
		\[
		f_1 \colon T_1 \rightarrow \mathcal{F}_1, \quad \mathcal{F}_1= \Proj\bigoplus_{m\geq 1} H^0(T_1,m\mathcal{O}_{T_1}(x_0)) 
		\]
	The sections of $\mathcal{F}_1$ are given by $x_0,\,x_1$ and monomials in the ideal $(u,\xi,v) \cap (w,t)$. Therefore, the map $f_1$ contracts the subset $\Gamma_1\colon (w=t=0) \subset T_1$ to $\mathbb{P}^1(a_0,a_1) \subset \mathcal{F}_1$ and the map $g_1$ extracts from it the subset $\Gamma_2\colon(u=\xi=v=0) \subset T_2$.
	The exceptional loci of $Y_1 \subset T_1$ and $Y_2 \subset T_2$ are $Y_1|_{\Gamma_1}$ and $Y_2|_{\Gamma_2}$, respectively. On the other hand, the base of the small modification $\tau$ restricts to
	\[
	\mathcal{B} \colon (f_{d_1}(x,y) = 0) \subset \mathbb{P}^1(a_0,a_1) \subset \mathcal{F}_1
	\]
	which is a cluster of $\frac{d_1}{a_0a_1}$ points. 
		We claim that, over each of these points, $f_1$ contracts a finite set of curves. Indeed, over a point in $\mathcal{B}$,
\[
Y_1|_{w=t=0} \colon (\xi^2g_{d_2-2a_{\xi}}(x_0,x_1)+\xi v g_{d_2-a_{\xi}}(x_0,x_1)+v^2+ug_{d_2}(x_0,x_1)=0) \subset \mathbb{P}(2_u,1_{\xi},1_v).
\]
By quasismoothness of $X$, we have that $g_{d_2}(x_0,y_0) \not = 0$. Hence, $u$ can be globally eliminated and hence, over each $(x_0 :x_1) \in \mathcal{B}$, we contract a single rational curve on each side of the diagram above.

Finally, it is clear that each contracted curve is trivial against $-K_Y$. We conclude that we have $s$ Atiyah flops.
		\end{proof}
		
		\begin{Prop} \label{prop:exchypend}
		The end of the 2-ray game restricts to an extremal divisorial contraction in the Mori Category to a Fano 3-fold hypersurface $X'\colon (h'=0)$ of Fano index 1 in a four dimensional weighted projective space. Moreover the divisor $E_t \colon (t=0)$ is contracted to the point $\mathbf{p_w}$ which is a cDV singularity for $X_{102}$ and a $cAx/2$ singularity for $X_{105}$. Also, $X'$ is quasismooth away from $\mathbf{p_w}$.
		\end{Prop}
		
\begin{proof}
Let $T':=T_2$. There is a map $\Phi' \colon T' \rightarrow \mathbb{P}(5,4,6,1,1,a_w')=:\mathbb{P}'$, where
 \[
\begin{array}{ccccccc|cccc}
             &       & u  & \xi & v & x_0 & x_1 & w & t & \\
\actL{T'}   &  \lBr &  a_w & \frac{a_w+3}{2} & \frac{a_w+7}{2} & 1 & 1 & 0 & -\frac{a_w-5}{2} &    \actR{.}\\
             &       & 5 & 4 & 6 &  1 & 1 & \frac{a_w-5}{2} & 0 &  
\end{array}
\] 
and $a_w' =\frac{a_w-5}{2}>0$ is given  by
\[
(u, z, v, x_0, x_1, w, t) \mapsto (ut^{a_w\cdot\frac{2}{a_w-5}}:\xi t^{\frac{a_w+3}{2}\cdot\frac{2}{a_w-5}}:vt^{\frac{a_w+7}{2}\cdot\frac{2}{a_w-5}}:x_0t^{\frac{2}{a_w-5}}:x_1t^{\frac{2}{a_w-5}}:w).
\]
It is clear that $\Phi'$ contracts the divisor $D_t \colon (t=0) \subset T'$ to the point $\mathbf{p_w} \in \mathbb{P}'$ and similarly when restricted to $Y'$. Hence, the 2-ray game on $Y$ follows the game on $T$ and $\Phi'$ is a extremal contraction which is $-K_Y$-negative since $-K_Y \in \Int\overline{\Mov}(Y)$.

 Now, for degree reasons $f_{d_1} \in \mathbb{C}[x_0,x_1]$ for family 102 and $f_{d_1} \in \mathbb{C}[u, x_0,x_1,t]$ for family 105. Recall that $tx_{\mu} \in f$ is $t^2$ for family 102 and $tv$ for family 105. These lift to
\[
t^2 \in H^0(Y_{102},-5K_{Y_{102}}-E), \quad \text{and} \quad tv \in H^0(Y_{105},-6K_{Y_{105}}). 
\]
Since $f^* \in |-mK_Y|$ for some $m>0$ for each of these two families, we have $t^2u \in f^*$ for family 102 and $tv \in f^*$ for family 105. Hence, away from the section $Y'|_{D_t}$, we can globally eliminate from $f^*$, the variable $u$ and $v$ for families 102 and 105, respectively. We have,
\[
u=\xi w - f_{d_1} \quad \text{and} \quad v=\xi w - f_{d_1}
\]
for families 102 and 105, respectively. No other variable can be globally eliminated. We look locally analytically at the contracted point $\mathbf{p_w}$. The variable $\xi$ can be eliminated from $f^*$ and the germ is analytically isomorphic to
\begin{align*}
&0 \in (uv+g_{14}(x_0,x_1)=0) \subset \mathbb{C}^4  \qquad &\text{for family 102}\\
&0 \in (u^2+v^2+g_{14}(x_0,x_1)=0) \subset \mathbb{C}^4/\mathbb{Z}_2(1,0,1,1)  \qquad &\text{for family 105}
\end{align*}

Finally $\varphi' \colon Y' \rightarrow X'$ is a weighted blowup of weights
\begin{align*}
&\wt(u,v,x_0,x_1)=(7,7,1,1) \qquad &\text{for family 102}\\
&\wt(u,v,x_0,x_1)=\frac{1}{2}(9,8,1,1) \qquad &\text{for family 105}
\end{align*}
 of a $cA_{13}$ and $cAx/2$ singularities with discrepancies $1$ and $\frac{1}{2}$, respectively.
\end{proof}

\begin{Thm} \label{thm:hypexcI}
Let $X \in I_S$ be as above. Then, there is an elementary Sarkisov link to a singular Fano 3-fold orbifold hypersurface $X'$, 
 \[
        \begin{tikzcd}[ampersand replacement=\&,column sep = 2em]
             \& Y \ar[dl]   \ar[r,"\simeq"] \& Y \ar[rr, dashed, "{s \times (-1,-1,1,1)}"] \ar[rd, swap]\& \&  Y' \supset (t=0) \ar[rd] \ar[dl] \& \\
             X  \&  \&  \&    \widetilde{f}|_{\Gamma_1}=0 \&   \& X' \ni \mathbf{p_w} 
        \end{tikzcd}
    \]
		initiated by the Kawamata blowup of a cyclic quotient singularity $\mathbf{p_{\xi}}\sim \frac{1}{3}(2,1,1)$ and followed by $s=\frac{d_1}{2}$ simultaneous Atiyah flops over the $s$ points $\widetilde{f}|_{\Gamma_1}=0$. The link ends with a divisorial contraction to a singular point $\mathbf{p_w} \in X'\subset \mathbb{P}'$. Moreover, $X'$ is quasismooth away from $\mathbf{p_w}$ and it has Fano index 1.
\end{Thm}

\begin{proof}
The 2-ray game on $Y$ follows the game on $T$. Hence, the movable  and pseudo-effective cones of $Y$ and $T$ coincide, $Y$ is a Mori dream space and $-K_Y \in \Int\overline{\Mov}(Y)$. The small $\mathbb{Q}$-factorial modifications of $Y$ are Atiyah flops. Hence, by Lemma \ref{lem:linkmds}, the composition $X \rat X'$ is a Sarkisov link. Moreover, notice that $Y$ and $Y'$ have the same analytic singularities by \cite[Theorem~6.15]{kollarmori} and, in particular, $Y'$ is quasismooth. The contraction to the point $\mathbf{p_w} \in X'$ is an isomorphism away from the contracted divisor and therefore $X'$ is quasismooth away from $\mathbf{p_w}$. Hence, $X'$ is smooth in codimension 2 and by adjunction, $-K_{X'} \sim \mathcal{O}_{X'}(1)$. 
\end{proof}

\begin{Rem} \label{rem:102}
By quasismoothness it follows that, if $\xi^3t \not \in g$, then $\xi^4 x_1 \in g$. Then, Corollay \ref{cor:102} shows that blowing up the $\mathbf{p_{\xi}}$ does not initiate a Sarkisov link.
\end{Rem}

\paragraph{Family 116.} Let $X$ be a \emph{general} quasismooth member of the deformation family 116. Then $X=X_{9,12} \subset \mathbb{P}(2,3,3,4,5,7):=\mathbb{P}$ with homogeneous variables $\xi,\,y,\,z,\,t,\,v,\,w$. Any such $X$ can be written as 
\begin{align*}
\xi w + vt + f_{9}(y,z) &= 0 \\
\xi ^4 t + \xi^3g_6(y,z) + \xi^2g_8(y,z,t,v)+ \xi g_{10}(y,z,t,v,w)+wv+t^3+g_{12}(y,z)  &=0.
\end{align*}
By the generality assumption, we have $v^2 \in g_{10}$. The point $\mathbf{p_{\xi}}$ is the cyclic quotient singularity $\frac{1}{2}(1,1,1)$ with orbinates $y,\,z,\,v$. Let $\Phi \colon T \rightarrow \mathbb{P}$ be the toric blowup that restricts to the Kawamata blowup $\varphi \colon E \subset Y \rightarrow X$ centred at $\mathbf{p_{\xi}}$. Then, when lifted to $Y$ via $\varphi$, the orbinates $y,\,z$ and $v$ vanish on $E$ with order $\frac{1}{2}$ each. Hence,
\[
y,\,z \in H^0(Y,-K_Y),\,\, v \in H^0\bigg(Y, -\frac{5}{3}K_Y+\frac{1}{3}E\bigg).
\]
Then, $m_f = \frac{3}{2}$ and $m_g=1$, where $m_f$ and $m_g$ are as in Lemma \ref{lem:utbl}. Hence,
\[
w \in H^0\bigg(Y,-\frac{7}{3}K_Y-\frac{1}{3}E\bigg),\,\, t \in H^0\bigg(Y,-\frac{4}{3}K_Y-\frac{1}{3}E\bigg)
\]
where these lifts are computed with Corollary \ref{cor:genlift}. Notice also that the defining equations of $Y$ are $(\widetilde{f}=\widetilde{g}=0) \subset T$ where 
\[
\widetilde{f} \in |-3K_Y|,\,\, \widetilde{g} \in |-4K_Y+E|. 
\]
Explicitly we can write them as 
\begin{align*}
\xi w + vt + f_{9}(y,z) &= 0 \\
\xi ^4 t + \xi^3g_6(y,z) + \xi^2g_8(u,y,z,t,v)+ \xi g_{10}(u,y,z,t,v,w)+wvu+t^3u^2+ug_{12}(y,z)  &=0.
\end{align*}
In particular, $\widetilde{g} \in (\xi,u)$ where $(u=0)$ is the exceptional divisor. We define the section
\begin{align*}
\eta &:=\frac{wv+t^3u+g_{12}(y,z)}{\xi} = - \frac{\xi ^3 t + \xi^2g_6(y,z) + \xi g_8(u,y,z,t,v)+ g_{10}(u,y,z,t,v,w)}{u} \\
&\in H^0\bigg(Y,-\frac{10}{3}K_Y-\frac{1}{3}E\bigg).
\end{align*}
We use this section to enlarge our ambient toric variety $T$ to $T^{\eta}$ by adjoining the section $\eta$ via the inclusion
\begin{align*}
i^{\eta} \colon T & \longrightarrow T^{\eta}\\
(u,\xi,v,y,z,w,t)&\longmapsto (u,\xi,v,y,z,\eta,w,t).
\end{align*}
Notice that $i^{\eta}$ is a morphism since the irrelevant ideal of $T$ is contained in the ideal of the indeterminacy locus of $i^{\eta}$, that is, $(u,\xi)\cap (v,y,z,w,t) \subset (u,\xi)$. Similarly, the projection away from $\mathbf{p_{\eta}}= (0,0,0,0,0,1,0,0)$ is everywhere defined since $\mathbf{p_{\eta}} \not \in T^{\eta}$. We conclude that $i^{\eta}$ is a biregular map onto its image. Define $Y^{\eta} \subset T^{\eta}$ with ideal
\begin{align*}
\xi w + vt + f_{9}(y,z) &= 0 \\
\eta \xi - wv- t^3u- g_{12}(y,z)&=0 \\
\eta u +\xi ^3 t + \xi^2g_6(y,z) + \xi g_8(u,y,z,t,v)+ g_{10}(u,y,z,t,v,w) &=0.
\end{align*}
Notice that $Y^{\eta}$ is the image of $Y$ under $i^{\eta}$. We conclude that these are biregular. Let $h_1$ and $h_2$ be the second and third equations above. Then, $h_1 \in |-4K_Y|$ and $h_2 \in |-\frac{10}{3}K_Y+\frac{2}{3}E|$. Now we play the 2-ray game on the rank two toric variety $T^{\eta}$. Since $-K_Y \sim \left(\begin{smallmatrix} 3 \\ 1 \end{smallmatrix} \right)$, we can view $T^{\eta}$ as
  \[
\begin{array}{cccc|ccccccc}
             &       & u  & \xi & v & y & z & \eta & w & t  &\\
\actL{T^{\eta} }   &  \lBr &  0 & 2 & 5 & 3 & 3 & 10 & 7 &  4 & \actR{.} \\
             &       & 1 & 1 & 2 &  1 & 1 & 3 &2 & 1  &
\end{array}
\] 
The 2-ray game on $T^{\eta}$ generates a diagram,
 \[
        \begin{tikzcd}[ampersand replacement=\&,column sep = 2em]
             \& T^{\eta} \ar[dl,swap,"\Phi^{\eta}"] \ar[dr,swap, "f"]   \ar[rr, dashed, "\displaystyle{\tau}"]\& \& T^{\eta}_1 \ar[dl,"g"]\ar[rr, dashed, "\displaystyle{\tau_1}"] \ar[rd, swap,"f_1"]\& \&  T^{\eta}_2 \ar[rd,swap,"f_2"] \ar[dl,"g_1"] \ar[rr, dashed, "\displaystyle{\tau_2}"] \& \& T^{\eta}_3 \ar[rd,"\Phi'"] \ar[dl,"g_2"] \\
             \mathbb{P}^{\eta}  \&  \& \mathcal{F} \& \&    \mathcal{F}_1 \&   \& \mathcal{F}_3 \& \& \mathcal{F}_4
        \end{tikzcd}
    \]
where the contractions are given by choosing the appropriate characters in 
\[
\overline{\Mov}(T^{\eta})= \langle\left(\begin{smallmatrix} 2 \\ 1 \end{smallmatrix} \right),\left(\begin{smallmatrix} 7 \\ 2 \end{smallmatrix} \right) \rangle. 
\]
For instance,
\[
\Phi^{\eta} \colon T^{\eta} \rightarrow \mathbb{P}^{\eta}=\Proj \bigoplus_{m\geq 1} H^0(T^{\eta},m\mathcal{O}\left(\begin{smallmatrix} 2 \\ 1 \end{smallmatrix} \right)) \simeq \mathbb{P}(2,3,3,4,5,7,10)
\]
and, recalling that $T$ is $(\eta = 0) \subset T^{\eta}$, we have that the restriction of $\Phi^{\eta}$ to $T$ is exactly $\Phi$. On the other end of the link we have 
\[
\Phi' \colon T^{\eta}_3 \rightarrow \mathcal{F}_4=\Proj \bigoplus_{m\geq 1} H^0(T^{\eta},m\mathcal{O}\left(\begin{smallmatrix} 7 \\ 2 \end{smallmatrix} \right)) \simeq \mathbb{P}(1,1,1,2,2,3,4)
\]
given in coordinates by
\[
(u,\xi,,v,y,z,\eta,w,t) \mapsto (w,yt,zt,\eta t,\xi t^3,vt^4,ut^7).
\]
In the same way, it is easy to see that the small $\mathbb{Q}$-factorial modification $\tau_i$ are 
\[
\tau = (-5,-1,1,1,5,4,3),\, \tau_1 = (-3,-1,-1,1,1,1),\, \tau_2 = (-10,-4,-5,-1,-1,1,2).
\]
We restrict this construction to $Y^{\eta}$. The map $\tau$ restricts to an isomorphism on $Y^{\eta}$ since the locus contracted by $f$ is $C \colon (y=z=\eta=w=t=0) \subset T^{\eta}$ which is not in $Y^{\eta}$ since $h_2|_C = v^2$ which implies that $v=0$. Hence, when we restrict to $Y^{\eta}$, the contracted points are all in the irrelevant ideal. Similarly, by Lemma \ref{lem:iso}, $\tau_1$ restricts to an isomorphism on $Y^{\eta}$. The map $\tau_2$ is the composition $g_2^{-1}\circ f_2$ where $f_2$ is 
\[
f_2 \colon T_2^{\eta} \rightarrow \mathcal{F}_3, \quad \mathcal{F}_3 =  \Proj \bigoplus_{m\geq 1} H^0(T^{\eta},m\mathcal{O}\left(\begin{smallmatrix} 10 \\ 3 \end{smallmatrix} \right))
\]
given in coordinates by
\[
(u,\xi,,v,y,z,\eta,w,t) \mapsto (\eta,zw,yw,z^2t,y^2t,vw^5,\ldots,ut^5). 
\]
That is, it contracts $(w=t=0) \subset T^{\eta}_2$ to the point $p=(1,0,\ldots,0)$. On an analytic neighborhood $U$ of $p \in \mathcal{F}_3$, we can eliminate the variables $\xi$ and $u$ on the fibres over $U$. Hence we get the hypersurface flip $(-5,-1,-1,1,2;3)$ with flipping equation $vt+f_3(y,z)+\cdots = 0$.

Let $Y^{\eta}_3$ be the image of $\tau_2$. The divisorial contraction $\Phi'$ restricts to a divisorial contraction 
\[
\varphi' \colon Y^{\eta}_3 \rightarrow Z_6.
\]
That is, it contracts $t=0$ to $\mathbf{p}_w$, and away from this divisor, we can globally eliminate $v$ from $f^*$ and $u$ from $h_1$. Hence we get, 
\[
Z_6 \colon (\eta(w^2\xi + wf_3(y,z) + g_4(y,z)+\eta z) + \xi^3+ \xi^2g_2(y,z) + \xi^2w^2+f_3(y,z)^2+\xi w f_3(y,z) = 0)
\]
inside $\mathbb{P}(1,1,1,2,2)$ with homogeneous variables $w,\,y,\,z,\,\eta,\xi$ and $\varphi'$ is a $(7,1,1,1)$-weighted blowup of a $cA_8$ singularity.

\paragraph{Family 118.} Let $X$ be a general quasismooth member of deformation family 118. The equations of $X$ are 
\begin{align*}
\xi x + vy + tz + f_{6}(x,y,z,t,v) &=0 \\
\xi (t+z) + v^2 + vy^2+y^4+g_8(x,y,z,t,v) &=0
\end{align*}
where we can assume $y^3 \not \in f_{6}$ and $X$ is inside $\mathbb{P}(1,2,3,3,4,5)$ with homogeneous variables $x,\,y,\,z,\,t,\,v,\,\xi$. The point $\mathbf{p_{\xi}}$ is a cyclic quotient singularity of type $\frac{1}{5}(1,2,3)$ with local coordinates $v,\,z,\,y$. We consider the toric blowup $\Phi \colon T \rightarrow \mathbb{P}$ whose restriction to $X$ is the unique Kawamata blowup $\varphi \colon E \subset Y \rightarrow X$ centred at $\mathbf{p_{\xi}}$. By Corollary \ref{cor:genlift}, we know that $v$ lifts to an anticanonical section on $Y$. Moreover,
\[
z,\,t \in H^0\bigg(Y,-\frac{3}{4}K_Y-\frac{1}{4}E\bigg), \, y \in H^0\bigg(Y,\frac{1}{2}K_Y-\frac{1}{2}E \bigg),\, x \in H^0\bigg(Y,\frac{1}{4}K_Y -\frac{3}{4}E\bigg).
\]
Also, 
\[
\xi \in H^0\bigg(Y, -\frac{5}{4}K_Y+\frac{1}{4}E\bigg).
\]
We have $\widetilde{f} \in |-\frac{3}{2}K_Y-\frac{1}{2}E|$ and $\widetilde{g} \in |-2K_Y|$. The corresponding bidegrees on $\mathbb{R}^2$ form the weight system of $T$:
  \[
\begin{array}{cccc|ccccccc}
             &       & u  & \xi & v & z & t  & y & x  &\\
\actL{T }   &  \lBr &  0 & 5 & 4 & 3 & 3 & 2 &  1 & \actR{.} \\
             &       & 1 & 4 & 3 &  2 & 2  &1 & 0  &
\end{array}
\] 
The cone of movable divisors of $T$ is 
\[
\overline{\Mov}(T) = \langle \left(\begin{smallmatrix} 5 \\ 4 \end{smallmatrix} \right), \left(\begin{smallmatrix} 2 \\ 1 \end{smallmatrix} \right) \rangle.
\]
This cone is subdivided in two chambers  separeted by the rays generated by $\left(\begin{smallmatrix} 4 \\ 3 \end{smallmatrix} \right)$ and $\left(\begin{smallmatrix} 3 \\ 2 \end{smallmatrix} \right)$. Notice that the divisors $D_y \sim D_z$ are both in the second wall as before. We play the 2-ray game on $T$. Then we have two small $\mathbb{Q}$-factorial modifications $T \rat T_1 \rat T_2$ which we call $\tau$ and $\tau_1$, respectively. We look closer at $\tau_1$. This map is the composition $g_1^{-1} \circ f_1$ where 
\[
f_1 \colon T_1 \rightarrow \Proj \bigoplus_{m\geq 1} H^0(T_1, m \mathcal{O}_{T_1}\left(\begin{smallmatrix} 3 \\ 2 \end{smallmatrix} \right)) = \Proj \bigoplus_{m\geq 1} H^0(T_2, m \mathcal{O}_{T_2}\left(\begin{smallmatrix} 3 \\ 2 \end{smallmatrix} \right)) \leftarrow T_2 \colon g_1
\]
In particular, $f_1$ contracts $(y=x=0) \subset T_1$ to $\mathbb{P}^1$ and $g_1^{-1}$ extracts from it the locus $(u=\xi=v=0) \subset T_2$. This is the small modification $(-3,-2,-1,1,2)$ with variables $u,\,\xi,\,v,\,y,\,x$ over $\mathbb{P}^1$. On the other hand, since the ray generated by $\left(\begin{smallmatrix} 2 \\ 1 \end{smallmatrix} \right)$ is in $\partial \Mov(T)$, we have a divisorial contraction,
\[
\Phi' \colon T_2 \rightarrow \Proj \bigoplus_{m\geq 1} H^0(T_2, m \mathcal{O}_{T_2}\left(\begin{smallmatrix} 2 \\ 1 \end{smallmatrix} \right)) \simeq \mathbb{P}(1,1,2,2,3,4). 
\]
This is in coordinates
\[
(u,\xi,v,z,t,y,x) \mapsto (y,ux^2,zx,tx,vx^2,\xi x^3).
\]
We now restrict the 2-ray game to $Y$. Since $v^2 \in \widetilde{g}$, the small modification $\tau$ restricts to an isomorphism. For $\tau_2$, we have that $\widetilde{g}|_{y=x=0}$ vanishes identically since $\widetilde{g} \not \in |mD_z|$ for any integer $m$. On the other hand, $\widetilde{f} \in |2D_z|$ and $\widetilde{f}$ contains the monomial $tz$ purely in $(t,z)$ and no others. Hence, the image of $f_1|_Y$ where $Y \subset T_1$ is $(tz = 0) \subset \mathbb{P}^1$ which is the two points $p_z$ and $p_t$. The fibre of $f_1$ over any of these points is analytically isomorphic to $\mathbb{P}^1(3,1)$ and the fibre of $g_1$ is $\mathbb{P}(1,2)$. 

The map $\Phi'$ restricts to a divisorial contraction. It is clear that the divisor $(x=0) \subset Y_2 \subset T_2$ is contracted to a point and away from this divisor we can set $x=1$. Hence $\xi$ can be globally eliminated using $\widetilde{f}$. We get the 3-fold $Z_6$ given by
\[
(vy+tz+f_6(u,x,y,z,t,v)) (t+z) + v^2 + vy^2u+y^4u^2+g_6(u,x,y,z,t,v) =0
\]
inside $\mathbb{P}(1,1,2,2,3)$ with homogeneous variables $y,\,u,\,z,\,t,\,v$. The map $\varphi$ can be realised as a weighted blowup of weights 
\[
\wt(u,v,z,t) = (2,3,1,1)
\]
and discrepancy 2 of the germ $0 \in (vt+u^2+h(u,z)=0)\sim cA_1 \in Z_6$.

\paragraph{Family 124.} Finally, we look at the cases of the deformation family 124. After a change of variables, any quasismooth member of family 124 is given by the equations
\begin{align*}
wt+v^2+vf_9(x,y)+z^3+f_{18}(x,y,z)&=0\\
wv+vg_{11}(x,y,z,t)+t^2z+tg_{13}(x,y,v)+y^4+x^5+g_{20}(x,y,z)&=0
\end{align*}
inside $\mathbb{P}(4,5,6,7,9,11)$ with homogeneous variables $x,\,y,\,z,\,t,\,v,\,w$. 

\begin{Thm} \label{thm:124toHyp}
Let $X$ be a quasismooth member of family 124. Then, there are two elementary Sarkisov links to singular Fano 3-fold orbifold hypersurfaces $Z_{10} \subset \mathbb{P}(1,1,2,3,5)$ and $Z_{15} \subset \mathbb{P}(1,2,3,5,7)$. 
		These are initiated by the Kawamata blowup of the cyclic quotient singularities $\frac{1}{11}(1,4,7)$ and $\frac{1}{7}(1,3,4)$, respectively and followed by a flip over a point. The link ends with a divisorial contraction to a singular point in $Z_{15}$ and $Z_{10}$, respectively.
\end{Thm}

\begin{proof} We consider two different toric blowups $T \rightarrow \mathbb{P}$ whose restriction to $X$ is, in each case, the unique Kawamata blowup of $X$ centred at some cyclic quotient singularity.  In the first case we impose that it is centred at $\mathbf{p_t}$ while in the second case we impose that it is centred at $\mathbf{p_w}$.

\paragraph{The centre of the blowup is $\mathbf{p_t}$.} We start by blowing the \emph{non-linear} cyclic quotient singularity $\mathbf{p_t} \sim \frac{1}{7}(1,3,4)$ on $X$ with local variables $x,\,y,\,v$. This is the Kawamata blowup $\varphi \colon Y \rightarrow X$ and the anticanonical divisor $-K_{Y}$ is Cartier divisor given by $x=0$ on $Y$. By Corollary \ref{cor:genlift},
\begin{align*}
v &\in H^0\bigg(Y,-\frac{9}{4}K_Y-\frac{1}{4}E\bigg),\,\, y \in H^0\bigg(Y,-\frac{5}{4}K_Y-\frac{1}{4}E\bigg) \\
w &\in H^0\bigg(Y,-\frac{11}{4}K_Y-\frac{3}{4}E\bigg), \,\, z \in H^0\bigg(Y,-\frac{3}{2}K_Y-\frac{1}{2}E\bigg).
\end{align*}
We also have $\widetilde{f} \in |-\frac{9}{2}K_Y-\frac{1}{2}E|$ and $\widetilde{g} \in |-5K_Y|$. Hence, $Y \subset T$ is defined by
\begin{align*}
tw+v^2+vf_9(x,y)+z^3u+f_{18}(x,y,z) &=0\\
t^2z+tg_{13}(x,y,v)+wvu+y^4u+x^5+g_{20}(x,y,z)&=0
\end{align*}
where $T$ is the rank two toric variety
  \[
\begin{array}{cccc|ccccccc}
             &       & u  & t & x & v & y  & w & z  &\\
\actL{T }   &  \lBr &  0 & 7 & 4 & 9 & 5 & 11 &  6 & \actR{.} \\
             &       & 1 & 2 & 1 &  2 & 1  &2 & 1  &
\end{array}
\]
where $E \colon (u=0)$ is the exceptional divisor of $\varphi$. The closure of the movable cone of $T$ is 
\[
\overline{\Mov}(T) = \langle \left(\begin{smallmatrix} 7 \\ 2 \end{smallmatrix} \right),\left(\begin{smallmatrix} 11 \\ 2 \end{smallmatrix} \right) \rangle
\]
We play the 2-ray game on $T$. Since $\overline{\Mov}(T)$ is subdivided in three chambers separated by the rays generated by the divisors $D_x$, $D_v$ and $D_y$ of $T$. Hence, there are small modifications $T \rat T_1 \rat T_2 \rat T_3$. We look closer at $\tau_2 \colon T_2 \rightarrow T_3$. This map is the composition $g_2^{-1} \circ f_2$ where 
\[
f_2 \colon T_2 \rightarrow \Proj \bigoplus_{m\geq 1} H^0(T_2, m \mathcal{O}_{T_1}\left(\begin{smallmatrix} 5 \\ 1 \end{smallmatrix} \right)) = \Proj \bigoplus_{m\geq 1} H^0(T_3, m \mathcal{O}_{T_2}\left(\begin{smallmatrix} 5 \\ 1 \end{smallmatrix} \right)) \leftarrow T_2 \colon g_2
\]
In particular, $f_2$ contracts $(w=z=0) \subset T_2$ to $\mathbb{P}^1$ and $g_2^{-1}$ extracts from it the locus $(u=t=x=v=0) \subset T_3$. This is the small modification $(-5,-3,-1,-1,1,1)$ with variables $u,\,t,\,x,\,v,\,w,\,z$ over a point. On the other hand, since the ray generated by $\left(\begin{smallmatrix} 11 \\ 2 \end{smallmatrix} \right)$ is in $\partial \Mov(T)$, we have a divisorial contraction,
\[
\Phi' \colon T_3 \rightarrow \Proj \bigoplus_{m\geq 1} H^0(T_3, m \mathcal{O}_{T_3}\left(\begin{smallmatrix} 11 \\ 2 \end{smallmatrix} \right)) \simeq \mathbb{P}(1,1,2,3,5,6). 
\]
This is in coordinates
\[
(u,t,x,v,y,w,z) \mapsto (w,yz,xz^3,vz^4,tz^8,uz^{11}).
\]
We restrict the 2-ray game to $Y$. Since $v^2 \in \widetilde{f}$ and $x^5 \in \widetilde{g}$ by quasismoothness of $X$, it follows that the maps $T \rat T_1 \rat T_2$ restrict to isomorphisms on $Y$. The fibre of $f_2$ over the point $p_y \in \Proj \bigoplus_{m\geq 1} H^0(T_2, m \mathcal{O}_{T_1}\left(\begin{smallmatrix} 5 \\ 1 \end{smallmatrix} \right))$ is locally analytically isomorphic to the curve $C \colon (wt+v^2+x^2+ \cdots = 0) \subset \mathbb{P}(3,1,1)$ with local coordinates $t,\,x,\,v$. On the other side, the fibre of $g_2$ over $p_y$ is $\mathbb{P}^1$. Hence, we have the hypersurface flip $(3,1,1,-1,-1;2)$. Let $Y_3 \subset T_3$ be the image of $Y$ under $\tau_2$. Then, $\Phi'$ restricts to a divisorial contraction $\varphi' \colon Y_3 \rightarrow Z_{10}$. The equations of $Z_{10}$ are obtained by moving away from the divisor $z=0$ on $Y_3$. That is, at $z=1$, we can globally elimiante $u$ using $\widetilde{f}$
\[
u=tw+v^2+vf_9(x,y)+f_{18}(x,y,z).
\]
We thus get,
\begin{multline*}
t^2+tg_{5}(x,y,v)+wv(tw+v^2+vf_3(x,y)+f_{6}(x,y,z))+\\y^4(tw+v^2+vf_3(x,y)+f_{6}(x,y,z))+x^5+g_{10}(x,y,1)=0
\end{multline*}
inside $\mathbb{P}(1,1,2,3,5)$ with homogeneous variables $w,\,y,\,x,\,v,\,t$. The divisorial contraction $\varphi'$ is centred at the point $\mathbf{p_w} \sim cA_4 \in Z_{10}$. Locally around it, the map $\varphi'$ is a weighted blowup with weights
\[
\wt(u,x,v,y) = (11,3,4,1)
\]
and discrepancy 3.

\paragraph{The centre of the blowup is $\mathbf{p_w}$.} The situation is similar to the previous one so we are brief. The ambient toric variety $T$ has weight system,
  \[
\begin{array}{cccc|ccccccc}
             &       & u  & w & x & y & z  & t & v  &\\
\actL{T }   &  \lBr &  0 & 11 & 4 & 9 & 5 & 6 &  7 & \actR{.} \\
             &       & 1 & 3 & 1 &  2 & 1  &1 & 1  &
\end{array}
\]
The 2-ray game on $T$ yields the diagram,
 \[
        \begin{tikzcd}[ampersand replacement=\&,column sep = 1em]
             \& T \ar[dl,swap,"\Phi"] \ar[dr,swap, "f"]   \ar[rr, dashed, "\displaystyle{\tau}"]\& \& T_1 \ar[dl,"g"]\ar[rr, dashed, "\displaystyle{\tau_1}"] \ar[rd, swap,"f_1"]\& \&  T_2 \ar[rd,swap,"f_2"] \ar[dl,"g_1"] \ar[rr, dashed, "\displaystyle{\tau_2}"] \& \& T_3 \ar[rd,"\Phi'"] \ar[dl,"g_2"] \\
             \mathbb{P}(4,5,6,7,9,11)  \&  \& \mathcal{F} \& \&    \mathcal{F}_1 \&   \& \mathcal{F}_3 \& \& \mathbb{P}(1,2,3,5,7,10)
        \end{tikzcd}
    \]
		and $Y \subset T$ is such that $\widetilde{f} \in |-\frac{9}{2}K_Y-\frac{1}{2}E|$ and $\widetilde{g} \in |-5K_Y|$. It is then easy to see that $\tau$ and $\tau_1$ restrict to isomorphisms and that $\tau_2$ is a hypersurface flip of type $(4,1,1,-1,-2;2)$ with local variables $w,\,x,\,v,\,z,\,t$ and flipping equation $wt+v^2+\cdots = 0$. We have then a divisorial contraction to a $cA_4$ point in the hypersurface $Z_{15}$ with equation
		\[
		(v^2+z^3u+\cdots)v+zu^2+y^4u+x^5+ \cdots = 0
		\]
		inside $\mathbb{P}(1,2,3,5,7)$ and homogeneous variables $z,\,y,\,x,\,v,\,u$.
		
		This contraction is a weighted blowup with weights
		\[
		\wt(w,x,v,y)=(7,2,3,1)
		\]
		and discrepancy 2.
\end{proof}

\paragraph{Families 95, 98, 107, 108 and 110.} Let $X\in I_S$ be a quasismooth member of a deformation family in \{95, 98, 107, 108, 110\}. The defining equations of each such $X$ can be written as
\begin{align*}
\xi w+f_{d_1}&=0\\
(1-\alpha)\xi^3x+\xi^2(\alpha x_1+g_{d_2-2a_{\xi}})+\xi g_{d_2-a_{\xi}}+v^2+wx_{\mu}+g_{d_2}&=0
\end{align*}

where 
\begin{itemize}
	\item $X \subset \Proj \mathbb{C}[x_0,x_1,x,\xi,v,w] \simeq \mathbb{P}(a_0,a_1,a_x,a_{\xi},a_v,a_w)$. 
	\item If $X \in \{107,108,110 \}$ we have $\xi^3x \not \in g$ for degree reasons and therefore $\alpha \not =0$ by quasismoothness of $X$. Hence we can also assume that $\xi^2g_{d_2-2a_{\xi}} \not \in g$ by a change of variables.  
	\item $f_{d_1},\,g_{d_2-a_{\xi}},\,g_{d_2-2a_{\xi}} \in \mathbb{C}[x_0,x_1,x,v,w]$
	\item $x_0 \in H^0(X,-K_X)$ and $x_1 \in H^0(X,-2K_X)$. In particular, $a_1=2a_0=2\iota_X =4$. 
	\item The monomial $x_{\mu}$ is $v$ if $a_w=a_v$ and is $x$ otherwise, where $a_w \geq a_v$ are two highest weights. 
\end{itemize}

Suppose we are in the most general case, in which $\alpha \not = 0$. Then, for families 95 and 98, we can eliminate $\xi^3x$ and $\xi^2g_{d_2-2a_{\xi}}$. In this case, the equations of $X \in \{95, 98, 107, 108, 110\}$ are
\begin{align*}
\xi w+f_{d_1}&=0\\
\xi^2x_1+\xi g_{d_2-a_{\xi}}+v^2+wx_{\mu}+g_{d_2}&=0.
\end{align*}
Hence $\mathbf{p_{\xi}} \in X$ is a cyclic quotient singularity type $\frac{1}{a_{\xi}}(\iota_X,a_v,a_x)$ with local variables $x_0,\,v,\,x$. By inspection, we have $a_v \equiv \iota_X \pmod{a_{\xi}}$. Since $\mathbf{p_{\xi}}$ is a terminal cyclic quotient singularity, there exists a unique $0 < k < a_{\xi}$ for which $k\iota_X \equiv 1 \pmod{a_{\xi}}$ and 
\[
\mathbf{p_{\xi}} \sim \frac{1}{a_{\xi}}(1,1,a_{\xi}-1)
\] 
The Kawamata blowup $\varphi \colon Y \rightarrow X$ centred at $\mathbf{p_{\xi}} \in X$ is the $(1,1,a_{\xi}-1)$-weighted blowup. Hence, locally at $\mathbf{p_{\xi}}$ it is the graph of of the rational map to $\mathbb{P}(1,1,a_{\xi}-1)$ given by
\[
\bigg(\frac{x_0^k}{\xi} : \frac{v^k}{\xi^{k_v}} : \frac{x^k}{\xi^{k_x}} \bigg)
\]
where $k_v,\,k_x$ are positive integers such that $ka_v-k_va_{\xi}=1$ and $ka_x-k_xa_{\xi} = a_{\xi}-1$, respectively. We now use Corollary \ref{cor:genlift} to lift the sections $H^0(X,a_iA)$ to $Y$. We have
\[
x_0 \in H^0(Y,-K_Y),\, v \in H^0\bigg(Y,-\frac{a_v}{\iota_X}K_Y+\frac{a_v-\iota_X}{\iota_Xa_{\xi}}E\bigg),\, x \in H^0\bigg(Y,-\frac{a_x}{\iota_X}K_Y-\frac{\iota_X-1}{\iota_X}E\bigg).
\]
Moreover, since $\varphi^*(\xi)$ does not vanish at $E$,
\[
\xi \in H^0\bigg(Y,-\frac{a_{\xi}}{\iota_X}K_Y+\frac{1}{\iota_X}E\bigg).
\]
By quasismoothness of $X$, $v^2 \in g$ (and so we know $2a_v=d_2$). Since $v$ vanishes at $E$ with order $\frac{1}{a_{\xi}}<1$ we have that $m_g = \frac{2}{a_{\xi}}$. Recall the definition of $m_g$ (and $m_f$) from the proof of Lemma \ref{lem:utbl}. Indeed, since $m_g=\frac{2}{a_{\xi}}$, then $g^*$ is a divisor in the linear system
\begin{align*}
\widetilde{g} \in \bigg|-\frac{d_2}{\iota_X}K_Y+\Big(\frac{d_2}{\iota_X a_{\xi}}-\frac{2}{ a_{\xi}} \Big)E\bigg| &= \bigg|-\frac{2a_v}{\iota_X}K_Y+\Big(\frac{2a_v}{\iota_X a_{\xi}}-\frac{2\iota_X}{\iota_X a_{\xi}} \Big)E\bigg|\\
&=\Bigg|-a_vK_Y+\frac{a_v-2}{a_{\xi}} E\bigg|\\
&=|-a_vK_Y+E|.
\end{align*}
We can now determine the lift of $x_1$ to $Y$. Again, by Corollary \ref{cor:genlift}, there are rational numbers $\alpha, \beta$ for which, $x_1 \in H^0(Y,-\alpha K_Y+\beta E)$. Hence,
\[
\xi^2x_1 \in H^0\Big(Y,-(a_{\xi}+\alpha)K_Y+ (1+\beta)E \Big).
\]
Since $\xi^2x_1 \in \widetilde{g}$, it follows that $\alpha = a_v-a_{\xi}=2$ and $\beta=\frac{a_v-2}{a_{\xi}}-1 = 0$. Hence, $x_1$ lifts to a pluri-anticanonical section in $H^0(Y,-2K_Y)$. In particular $\varphi^*(x_1)$ vanishes with order $\frac{2}{a_{\xi}}$.

By quasismoothness of $X$ there are monomials $f(x_0,x_1)$ purely in $(x_0,x_1)$ in $f$ (and in $g$). These monomials are of the form $x_0^{\alpha}x_1^{\beta}$ where $d_1 =\alpha a_0 + \beta a_1 = \alpha \iota_X + 2\beta \iota_X$. Hence, when these are pulled back by $\varphi$, they vanish at $E$ with order 
\[
\alpha \cdot \frac{1}{a_{\xi}}+\beta\cdot \frac{2}{a_{\xi}} = \frac{\alpha + 2 \beta}{a_{\xi}} = \frac{d_1}{\iota_X a_{\xi}}
\] 
It is not true that $\frac{d_1}{\iota_Xa_{\xi}}< 1$ so it is not immediate that $m_f=\frac{d_1}{\iota_Xa_{\xi}}$. However, it is easy to check that each monomial vanishes with order at least $\frac{d_1}{\iota_Xa_{\xi}}$ so we conclude that this is indeed the value of $m_f$ and therefore $\widetilde{f}=0$ is a divisor in the pluri-anticanonical linear system $|-\frac{d_1}{\iota_X}K_Y|$. Finally, we can readily see that $w$ lifts to
\[
w \in H^0\bigg(Y,-\frac{a_w}{\iota_X}K_Y-\frac{1}{\iota_X}E\bigg).
\]

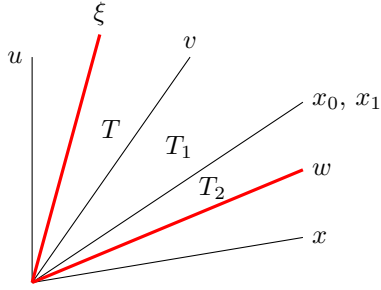
\begin{figure}%
\centering
\begin{tikzpicture}[scale=3]
 \coordinate (A) at (0, 0);
 \coordinate [label={left:$u$}] (E) at (0, 1);
\coordinate [label={above:$\xi$}] (5) at (0.3,1.1); 
\coordinate [label={above:$v$}] (K) at (0.7, 1);
 \coordinate [label={right:$x_0,\,x_1$}] (2) at (1.2,0.8);
 \coordinate [label={right:$w$}] (4) at (1.2,0.5);
 \coordinate [label={right:$x$}] (3) at (1.2,0.2);

\coordinate [label={$T$}] (T) at (0.35,0.6);
\coordinate [label={right:$T_1$}] (T1) at (0.55,0.6);
\coordinate [label={$T_2$}] (T) at (0.8,0.33);

  \draw (A) -- (E);
  \draw (A) -- (K);
	\draw [very thick,color=red] (A) -- (5);
	\draw (A) -- (2);
	\draw (A) -- (3);
	\draw [very thick,color=red] (A) -- (4);
  
\end{tikzpicture}
\caption{A representation of the chamber decomposition of the cone of effective divisors of $T$. In red the subcone of movable divisors of T. The ample models $T$, $T_1$ and $T_2$ are represented in the interior of the corresponding Nef chamber.}%
\label{fig:Cod2ExcI}%
\end{figure}

Now that we have computed how all sections $H^0(X,a_iA)$ lift, we can construct the rank 2 toric variety as in Lemma \ref{lem:utbl}. We have a representation of the cone of effective divisors of $T$ as in Figure \ref{fig:Cod2ExcI}.

We now play the 2-ray game on $T$.

\begin{Lem} \label{lem:Cod2ExcpSmallMod}
The following assertions hold:
\begin{itemize}
	\item The small $\mathbb{Q}$-factorial modification $T \rat T_1$ induced by crossing the $v$-wall restricts to an isomorphism on $Y$.
	\item The small $\mathbb{Q}$-factorial modification $T_1\rat T_2$ induced by crossing the $(x_0,x_1)$-wall restricts to $s$ simultaneous flops on $Y$ over $s$ points where $s$ is the number of solutions to $\widetilde{f}|_{w=x=0}=0$.
\end{itemize}
\end{Lem} 
\begin{proof}
The proof of Lemmas \ref{lem:HypExcp1stwall} and \ref{lem:flops} applies to each item respectively with no change to this situation.
\end{proof}

\begin{Lem} \label{lem:Cod2ExcToricEnd}
The end of the toric 2-ray game on $T$ is an extremal divisorial contraction $\Phi'$ of $(x=0) \subset T_2$ to 
\[
(0:0:0:0:0:1) \in \mathbb{P}\Big(a_x, a_{\xi}-1,    a_v-2 , 1 , 2, \frac{a_w-a_x}{\iota_X}\Big).
\]
Moreover, $\mathbf{p_w}$ is an orbifold point.
\end{Lem}

\begin{proof}
We begin to see that the rays $\mathbb{R}_+[w]$ and $\mathbb{R}_+[x]$ are actually distinct. These are generated by the vectors 
\[
\left(\begin{smallmatrix}  a_w \\ a_w' \end{smallmatrix} \right), \quad  \left(\begin{smallmatrix}  a_x \\ a_x' \end{smallmatrix} \right)
\]
where $a_w' := \frac{ka_w-m_1}{a_{\xi}}$ and $a_x' := \frac{ka_x-(a_{\xi}-1)}{a_{\xi}}$. We have
\[
\det
\begin{pmatrix}
a_w & a_x \\
a_w' & a_x' 
\end{pmatrix}
= \frac{a_x-a_w}{\iota_X}=:-\frac{1}{d} < 0
\]
and indeed $\mathbb{R}_+[w]$ and $\mathbb{R}_+[x]$ are rays in the boundaries of $\overline{\Mov}(T)$ and $\overline{\Eff}(T)$, respectively, and these are distinct. Hence, there is a extremal divisorial contraction $\Phi' \colon T_2 \rightarrow \mathbb{P}'$ which is given by the sections $\bigoplus_{m\geq 1}H^0(T_2,m\mathcal{O}(w))$. We write this map explicitly. The weight system of $T_2$ is as in \ref{eq:T}. We write it here for convenience:
		\[
\begin{array}{ccccccc|cccc}
             &       & u  & \xi &   v & x_0 & x_1 & w & x & \\
\actL{T_2 }   &  \lBr &  0 & a_{\xi} &   a_v & \iota_X & 2\iota_X & a_w & a_x &   \actR{.}\\
             &       & 1 & k &  \frac{ka_v-1}{a_{\xi}} & \frac{k\iota_X-1}{a_{\xi}} & 2\cdot \frac{k\iota_X-1}{a_{\xi}} & a_w' & a_x' &  
\end{array}
\] 
Consider the matrix 
\[
A=\begin{pmatrix}
-a_x' & a_x \\
-a_w' & a_w 
\end{pmatrix}
 \in \GL(2, \mathbb{Z}).
\]
Then multiplying the weight system $T_2$ by $A$ on the left yields the isomorphic rank 2 toric variety (see \cite[Lemma~2.4]{hamidplia})
		\[
\begin{array}{ccccccc|cccc}
             &       & u  & \xi &   v & x_0 & x_1 & w & x & \\
\actL{T_2 }   &  \lBr &  a_x & a_{\xi}-1 &   a_v-2 & \iota_X-1 & 2 (\iota_X-1) & \frac{a_w-a_x}{\iota_X} & 0 &   \actR{.}\\
             &       & a_w & \frac{d_1}{\iota_X} &  \frac{a_v+a_w}{\iota_X} & 1 & 2 & 0 & -\frac{a_w-a_x}{\iota_X}  &  
\end{array}
\] 

Since $\iota_X =2 $ the divisorial contraction $\Phi'$ is
\begin{align*}
\Phi' \colon T_2 & \longrightarrow \mathbb{P}\Big(a_x, a_{\xi}-1,    a_v-2 , 1 , 2, \frac{a_w-a_x}{\iota_X}\Big) \\
(u,\xi, v, x_0, x_1, w, x) & \longmapsto (ux^{da_w},\xi x^{d\cdot \frac{d_1}{\iota_X}},v x^{d\cdot \frac{a_v+a_{\xi}}{\iota_X}},x_0 x^{d},x_1 x^{2d},w).
\end{align*} 
Hence, $\Phi'$ contracts the divisor $(x=0) \subset T_2$ to the point $\mathbf{p_w}:=(0:0:0:0:0:1)  \in \mathbb{P}'$. Moreover, $\mathbb{P}'$ is a well-formed weighted projective space since, for instance, $a_x$ is odd. Now $\mathbf{p_w}$ is an orbifold point as a consequence of $\frac{a_w-a_x}{\iota_X}>0$.
\end{proof}

\begin{Lem} \label{lem:auxcod2Exc}
There is no variable $y$ such that $x^{\alpha}y \in \widetilde{f}$ or $\widetilde{g}$.
\end{Lem}

\begin{proof}
Recall that $\widetilde{f} \in |-\frac{d_1}{\iota_X}K_Y|$. Suppose by contradiction that there is a variable $y \in \{u,\xi,v,x_0,x_1,w,x\}$ for which $x^{\alpha}y \in \widetilde{f}$. There are rational numbers $\beta_1,\,\beta_2$ for which $\varphi^*(y) \in H^0(Y,-\beta_1K_Y+\beta_2 E)$. Hence,
\[
x^{\alpha}y \in H^0\bigg(Y,-\Big(\alpha\cdot \frac{a_x}{\iota_X}+\beta_1 \Big)K_Y-\Big(\alpha\cdot \frac{\iota_X-1}{\iota_X}-\beta_2\Big)E \bigg).
\]
In particular $\beta_2 >0$. Therefore, $y \in \{u,\xi,v \}$. However, notice that $\xi$ appears in $f$ only in the monomial $\xi w$ and so $y \in \{u,v \}$. If $y=u$, then $\beta_2=1$ and $\alpha=2$. Otherwise, $\beta_2= \frac{a_v-\iota_X}{\iota_Xa_{\xi}} = \frac{1}{\iota_X}$ and $\alpha=1$. In particular it implies that $x^2 \in f$ or $xv \in f$ which is not possible in either case for degree reasons. We conclude that there is no variable $y$ such that $x^{\alpha}y \in \widetilde{f}$. A similar argument works for $\widetilde{g} \in |-a_vK_Y+E|$.
\end{proof}

		\begin{Prop}
		The end of the 2-ray game restricts to an extremal divisorial contraction in the Mori Category to a Fano 3-fold complete intersection, $X'\colon (h_1=h_2=0) \subset \mathbb{P}'$. Moreover the divisor $E' \colon (x=0)$ is contracted to the point $\mathbf{p_w}$ which is a hyperquotient singularity. Also, $X'$ is quasismooth away from $\mathbf{p_w}$.
		\end{Prop}
\begin{proof}
We restrict the divisorial contraction of Lemma \ref{lem:Cod2ExcToricEnd} to $Y_2 \subset T_2$ and call it $\varphi'$. The map $\varphi'$ contracts the divisor $(x=0) \subset Y_2$ to a point. By Lemma \ref{lem:auxcod2Exc}, the 3-fold to which we contract has codimension 2 and is given by
\begin{align*}
\xi w+f_{d_1}&=0\\
\xi^2x_1+\xi g_{d_2-a_{\xi}}+v^2+wx_{\mu}u^{\beta}+g_{d_2}&=0
\end{align*}
inside $\mathbb{P}\Big(a_x, a_{\xi}-1,    a_v-2 , 1 , 2, \frac{a_w-a_x}{\iota_X}\Big)$ with homogeneous coordinates $u,\,\xi,\,v,\,x_0,\,x_1,\,w$, where $\beta = 1$ for family 95 and 107 and $\beta=2$ otherwise. The point to which we contract is $\mathbf{p_w} = (0:\cdots : 0 :1 )$. In an analytic neighborhood of this point, we can eliminate $\xi$ from the first equation and no other variables since $\beta > 0$. This is a singular point whose germ is
\[
0 \in (u^2+v^2+h(x_0,x_1)=0)/\bm{\mu}_{\frac{1}{d}}(a_x,a_v-2,1,2).
\]
of the hyperquotient singularity $cA/2$ for families 95 and 107 and $cAx/4$ for families 98, 108 and 110. The map $\varphi' \colon Y_2 \rightarrow X'$ can be realised as a weighted blowup with weights 
\[
\wt(u,v,x_0,x_1)=\frac{\iota_X}{a_w-a_x}\bigg(a_w,\frac{a_v+a_w}{\iota_X},1,2\bigg)
\]
and discrepancy 
\[
\frac{\iota_X}{a_w-a_x}= 
\begin{cases*}
                    \frac{1}{2} & if  $X \in \{95, 107\}$  \\
                    \frac{1}{4} & if $X \in \{98, 108,110 \}$
                 \end{cases*}
\]
by Lemma \ref{lem:discr}. These are indeed terminal extractions by \cite[Section~6 and 7]{HayakawaI}.
\end{proof}
\begin{Ex}
Let $X$ be the quasismooth Fano defined by
\begin{align*}
t^2y+wz+v^2+x^9 &= 0\\
tw+z^4+y^5+x^{10}&=0
\end{align*}
in the deformation family $X_{18,20} \subset \mathbb{P}(2,4,5,7,9,13)$ with homogenous variables $x,\,y,\,z,\,t,\,v,\,w$. The proper transform $Y$ of the Kawamata blowup centred at $\mathbf{p_t} \in X$ is 
\begin{align*}
t^2y+wzu^2+v^2+x^9u &= 0\\
tw+z^4u^2+y^5+x^{10}&=0
\end{align*}
inside 
\[
\begin{array}{ccccc|cccccc}
             &       & u  & t & v & x & y  & w & z & \\
\actL{T_1}   &  \lBr &  -2 & -1 & -1 & 0 & 0 & 1 & 1 &    \actR{.}\\
             &       & 1 & 4 & 5 &  1 & 2 & 6 & 2 &  
\end{array}
\] 
There is a small contraction $T_1 \dashrightarrow T_2$ where $T_2$ has the same Cox ring as $T_1$ but its irrelevant ideal is $(u,t,v,x,y) \cap (w,z)$. This contraction restricts to 
 \[
        \begin{tikzcd}[ampersand replacement=\&,column sep = 2em]
              Y_1 \ar[rr, dashed, "{s \times (-1,-1,1,1)}"] \ar[rd, swap, "f_1"]\& \&  Y' \ar[dl, "g_1"]  \\
         \&    \mathcal{F}_1 \&   
        \end{tikzcd}
    \]
		where $s=5$ is the number of solutions to $(x^{10}+y^5=0) \subset \mathbb{P}(1,2)$. Furthermore, the map $\Phi' \colon T' \rightarrow \mathbb{P}(1,2,4,5,6,7)$ restricts to a divisorial contraction $\varphi' \colon Y' \rightarrow X'$ to a point $\mathbf{p_w} \in X'$ where $X'$ has defining equations
		\begin{align*}
		t^2y+wu^2+v^2+x^9u &= 0\\
		tw+u^2+y^5+x^{10}&=0
		\end{align*}
		inside $\mathbb{P}(1,2,4,5,6,7)$ with homogeneous variables $x,\,y,\,w,\,u,\,t,\,v$. In an analytic neighborhood of $\mathbf{p_w}$ we can write $t=u^2+y^5+x^{10}$. So we get the germ
		\[
		(0 \in (u^2+y^5+x^{10})^2y+u^2+v^2 + x^9u = 0) \subset \mathbb{C}^4_{u,v,x,y} / \mathbb{Z}_4(1,3,1,2).
		\]
		Notice that the action of $\mathbb{Z}_4$ is equivariant on $F(u,v,x,y)=0$. The divisorial contraction $\varphi'$ is the $\frac{1}{4}(13,11,1,2)$-weighted blowup of $\mathbf{p_w} \in X'$ and discrepancy $\frac{1}{4}$. This is indeed an extraction in the Mori category. See \cite[Section~7]{HayakawaI}.  
\end{Ex}

\begin{Thm} \label{thm:HypCaseIIIExcI}
Let $X \in I_S$ be as above. Then, there is an elementary Sarkisov link to a singular Fano 3-fold orbifold hypersurface $X'$, 
 \[
        \begin{tikzcd}[ampersand replacement=\&,column sep = 2em]
             \& Y \ar[dl]   \ar[r,"\simeq"] \& Y \ar[rr, dashed, "{s \times (-1,-1,1,1)}"] \ar[rd, swap]\& \&  Y' \supset (x=0) \ar[rd] \ar[dl] \& \\
             X  \&  \&  \&    f^*|_{\Gamma_1}=0 \&   \& X' \ni \mathbf{p_w} 
        \end{tikzcd}
    \]
		initiated by the Kawamata blowup of a cyclic quotient singularity $\mathbf{p_{\xi}}$ and followed by $s>0$ simultaneous Atiyah flops over the $s$ points $\widetilde{f}|_{\Gamma_1}=0$. The link ends with a divisorial contraction to a hyperquotient singularity $\mathbf{p_w} \in X'$. Moreover, $X'$ is quasismooth away from $\mathbf{p_w}$.
\end{Thm}

\begin{proof}
The proof is the same as in Theorem \ref{thm:hypexcI}.
\end{proof}

\begin{Thm} \label{thm:cod2ExcII}
Let $X \in \{114, 117, 121, 122, 123, 125\}$ be the general member. Then, there is an elementary Sarkisov link to a singular Fano 3-fold orbifold hypersurface $X'$, 
 \[
        \begin{tikzcd}[ampersand replacement=\&,column sep = 2em]
             \& Y \ar[dl]   \ar[r,"\simeq"] \& Y \ar[rr, dashed, "\tau"] \ar[rd, swap]\& \&  Y' \supset (x=0) \ar[rd] \ar[dl] \& \\
             X  \&  \&  \&    f^*|_{\Gamma}=0 \&   \& X' \ni \mathbf{p_w} 
        \end{tikzcd}
    \]
		initiated by the Kawamata blowup of a cyclic quotient singularity $\mathbf{p_{\xi}}$ and followed by $s>0$ simultaneous Atiyah flops over the $s$ points $f^*|_{\Gamma}=0$. The link ends with a divisorial contraction to a hyperquotient singularity $\mathbf{p_w} \in X'$. Moreover, $X'$ is quasismooth away from $\mathbf{p_w}$.
\end{Thm}

\begin{proof}We do each family separately. 

\paragraph{Family 114.} Let $X$ be the general member in family 114. We have proven in Theorem \ref{thm:delPezzoMain} that any quasismooth member of this family is non-solid. W.l.o.g., we can assume that $\mathbf{p_y} \in X$ and so $\mathbf{p_y} \sim \frac{1}{2}(1,1,1)$ with variables $x,\,z,\,t$. The Kawamata blowup of $\mathbf{p_y}$ is $\varphi \colon Y \rightarrow X$ where $Y$ is 
\begin{align*}
yv+wxu+f_6(u,x,t,z)&=0 \\
y^2w+yg_{7}(u,z,t,v,w)+wvu+g_{9}(u,x,z,t,v)&=0
\end{align*}
in
\[
\begin{array}{cccccc|ccccc}
             &       & u  & y & z & t & v  & w & x & \\
\actL{T_1}   &  \lBr &  -4 & -2 & -1 & -1 & 0 & 1 & 1 &    \actR{.}\\
             &       & 1 & 1 & 1 &  1 & 1 & 1 & 0 &  
\end{array}
\] 
In this case, the movable cone of $T$ is subdivided into three Mori chambers and by Lemma \ref{lem:iso}, the first wall crossing in $\overline{\Mov}(T)$ is an isomorphism. Crossing to the last chamber in $\overline{\Mov}(T)$ amounts to contracting the curve $(w=x=0) \subset T_1$ to a point and then extract the locus $(u=y=z=t=0) \subset T_2$. The resulting variety $T_2$ has the same Cox ring as $T_1$ but its irrelevant ideal is $(u,y,z,t,v) \cap (w,x)$. By multiplying the weight system of $T_1$ by the matrix 
\[
\begin{pmatrix}
-1	& 1 \\ 
0	& 1
\end{pmatrix} \in \SL_2(\mathbb{Z})
\]
we can write $T_2$ as 
\[
\begin{array}{ccccccc|cccc}
             &       & u  & y & z & t & v  & w & x & \\
\actL{T_2}   &  \lBr &  5 & 3 & 2 & 2 & 1 & 0 & -1 &    \actR{.}\\
             &       & 1 & 1 & 1 &  1 & 1 & 1 & 0 &  
\end{array}
\] 
We have then a divisorial contraction $\Phi \colon T_2 \rightarrow \mathbb{P}^5$ to the ray $\mathbb{R}_+[\left(\begin{smallmatrix} 0 \\ 1 \end{smallmatrix} \right)]$ i.e., induced by the sections, 
\[
\bigoplus_{m\geq 1}H^0(T_2,m\mathcal{O}\left(\begin{smallmatrix} 0 \\ 1 \end{smallmatrix} \right))
\]
and is given in coordinates by
\[
(u,y,z,t,v,w,x) \mapsto (ux^5,yx^3,zx^2,tx^2,vx,w).
\] 
We restrict this construction to $Y$. Over an analytic neighborhood around $p_v$, the fibre of the contraction of $(w=x=0) \subset Y$ is $wu+g_3(z,t) + \cdots = 0$ where $\wt(u,z,t,w,x)=(-4,-1,-1,1,1)$. Hence, this is a flip of type $(-4,-1,-1,1,1;3)$. The divisorial contraction $\Phi$ restricts to the weighted blowup $\varphi' \colon Y_2 \rightarrow Z_{2,3}$ of weights
\[
\wt(y,z,t,v) = (3,2,2,1)
\]
of $cD_4 \sim \mathbf{p_w} \in Z_{2,3}$ with discrepancy 2 where $Z_{2,3}$ is given by the equations,
\begin{align*}
yv+wu+f_2(u,1,t,z)&=0 \\
y^2w+yg_{2}(u,z,t,v,w)+wvu+g_{3}(u,1,z,t,v)&=0
\end{align*}
in $\mathbb{P}^5$.

\paragraph{Family 117.} Let $X$ be the general member in family 117. W.l.o.g., we can assume that $\mathbf{p_z} \in X$ and so $\mathbf{p_z} \sim \frac{1}{4}(1,1,3)$ with local variables $x,\,y,\,t$. The Kawamata blowup of $\mathbf{p_z}$ is $\varphi \colon Y \rightarrow X$ where $Y$ is 
\begin{align*}
zw+vtu+f_{12}(x,y,z,t)&=0 \\
z^2v+zg_{11}(x,y,t,w)+wvu+t^3u+g_{15}(x,y,t,vu)&=0
\end{align*}
in
\[
\begin{array}{cccccc|ccccc}
             &       & u  & z & x & y & w  & t & v & \\
\actL{T_1}   &  \lBr &  -3 & -1 & 0 & 0 & 1 & 1 & 2 &    \actR{.}\\
             &       & -8 & -4 & -1 &  -1 & 0 & 1 & 3 &  
\end{array}
\] 
In this case, the movable cone of $T$ is subdivided into three Mori chambers and, by Lemma \ref{lem:iso}, the first wall crossing in $\overline{\Mov}(T)$ is an isomorphism. Crossing to the last chamber in $\overline{\Mov}(T)$ amounts to contracting the curve $(t=v=0) \subset T_1$ to a point and then extract the locus $(u=z=x=y=0) \subset T_2$. The resulting variety $T_2$ has the same Cox ring as $T_1$ but its irrelevant ideal is $(u,z,x,y,w) \cap (t,v)$. By multiplying the weight system of $T_1$ by the matrix 
\[
\begin{pmatrix}
3	& -2 \\ 
1	& -1
\end{pmatrix} \in \SL_2(\mathbb{Z})
\]
we can write $T_2$ as 
\[
\begin{array}{ccccccc|cccc}
             &       & u  & z & x & y & w  & t & v & \\
\actL{T_2}   &  \lBr &  7 & 5 & 2 & 2 & 3 & 1 & 0 &    \actR{.}\\
             &       & 5 & 3 & 1 &  1 & 1 & 0 & -1 &  
\end{array}
\] 
We have then a divisorial contraction $\Phi \colon T_2 \rightarrow \mathbb{P}(7,5,2,2,3,1)$ which corresponds to the contraction of the ray $\mathbb{R}_+[\left(\begin{smallmatrix} 0 \\ -1 \end{smallmatrix} \right)]$. That is induced by the sections, 
\[
\bigoplus_{m\geq 1}H^0(T_2,m\mathcal{O}\left(\begin{smallmatrix} 1 \\ 0 \end{smallmatrix} \right))
\]
and is given in coordinates by
\[
(u,z,x,y,w,t,v) \mapsto (uv^5,zv^3,xv,yv,wv,t).
\] 
We restrict this construction to $Y$. The small modification $T_1 \rat T_2$ contracts the curve $C_1 \colon (f_5(x,y)=0) \subset \mathbb{P}(1,1,8)$ with variables $x,\,y,\,w$ to a point and extracts $C_2 \simeq \mathbb{P}(1,3)$. Over an analytic neighborhood around the base of the contraction, we have then a flip of type $(8,1,1,-1,-3;5)$ with local variables $u,\,x,\,y,\,t,\,v$ and flipping equation $vu+b_5(x,y)+\cdots  = 0$.  The divisorial contraction $\Phi$ restricts to the weighted blowup $\varphi' \colon Y_2 \rightarrow Z_{8,10}$ of weights
\[
\wt(z,x,y,w) = (3,1,1,1)
\]
and discrepancy 1 to the point $cA_3 \sim \mathbf{p_t} \in Z_{8,10}$ where $Z_{8,10}$ is given by the equations,
\begin{align*}
zw+tu+f_{12}(x,y,z,t)&=0 \\
z^2+zg_{11}(x,y,t,w)+wu+t^3u+g_{15}(x,y,t,u)&=0
\end{align*}
in $\mathbb{P}(7,5,2,2,3,1)$ with homogeneous variables $u,\,z,\,x,\,y,\,w,\,t$.

\paragraph{Family 121.} Let $X$ be a general member in the family 121. We have proven in Proposition \ref{prop:121_122} that any quasismooth member of this family is non-solid. W.l.o.g., we can assume that $\mathbf{p_y} \in X$ and so $\mathbf{p_y} \sim \frac{1}{3}(1,1,2)$ with variables $x,\,z,\,t$. The Kawamata blowup of $\mathbf{p_y}$ is $\varphi \colon Y \rightarrow X$ where $Y$ is 
\begin{align*}
yw+vz+t^2+uf_9(x,z,t,v)&=0 \\
y^2v+yg_{9}(u,x,z,t,v)+z^3+v^2u+wtu+g_{13}(u,x,z,t,v)&=0
\end{align*}
in
\[
\begin{array}{cccccc|ccccc}
             &       & u  & y & z & t & v  & w & x & \\
\actL{T}   &  \lBr &  -6 & -3 & -2 & -1 & 0 & 1 & 1 &    \actR{.}\\
             &       & 1 & 1 & 1 &  1 & 1 & 1 & 0 &  
\end{array}
\] 
The curve $C^- \colon (w=x=0) \subset Y$ is isomorphic to $\mathbb{P}(3,1)$ and is contracted to a point. On the other side, a curve $C^+$ isomorphic to $\mathbb{P}^1$ is extracted. We have the following diagram: 
 \[
        \begin{tikzcd}[ampersand replacement=\&,column sep = 2em]
              Y \ar[rr, dashed, "{(-3,-1,1,1)}"] \ar[rd, swap, "f_1"]\& \&  Y' \ar[dl, "g_1"]  \\
         \&    \mathcal{F} \&   
        \end{tikzcd}
    \]
 
There is a divisorial contraction, $T^+ \rightarrow \mathbb{P}^5$ which is given by the sections 
\[
\bigoplus_{m\geq 1}H^0(T^+,m\mathcal{O}(0,1))=\mathbb{C}[ux^7,yx^4,zx^3,tx^2,vx,w]
\]
where $T^+$ is 
\[
\begin{array}{ccccccc|cccc}
             &       & u  & y & z & t & v  & w & x & \\
\actL{T^+}   &  \lBr &  7 & 4 & 3 & 2 & 1 & 0 & -1 &    \actR{.}\\
             &       & 1 & 1 & 1 &  1 & 1 & 1 & 0 &  
\end{array}
\] 
Restricted to $Y^+$, this is a divisorial contraction $\varphi' \colon Y^+ \rightarrow Z_{2,3} \subset \mathbb{P}^5$ centred at the singular point $\mathbf{p_w}$ of discrepancy 3 by Lemma \ref{lem:discr}. In fact, $\varphi'$ is the weighted blowup
\[
\wt(u,z,t,v)=(7,3,2,1)
\]
of the $cA_2$ singularity 
\[
0 \in (tu+z^3+g_{\geq 3}(z,v)=0) \subset \mathbb{C}^4_{tuzv}.
\]

In \cite[Example~6.8]{kawakitaelephants}, Kawakita gives an example of a different weighted blowup of a $cA_2$ point with discrepancy 3. See also \cite[Theorem~1.13~(ii)]{kawakitaelephants}. 

 The equations of $Z_{2,3}$ are 
\begin{align*}
yw+vz+t^2+uf_9(1,z,t,v)&=0 \\
y^2v+yg_{9}(u,1,z,t,v)+z^3+v^2u+wtu+g_{13}(u,1,z,t,v)&=0.
\end{align*}
Moreover, in \cite[Chapter~3]{iskpuk}, the birational rigidity of a general 3-fold complete intersection of a quadric and a cubic $Z$ has been obtained. This example shows that this fails when $Z$ has a $cA_2$ point. In fact, it fails quite dramatically as $Z$ is birational to del Pezzo fibration, see Proposition \ref{prop:121_122}. 

\paragraph{Family 122.} Let $X$ be a general member of family 122. After a change of variables we can write $X$ with equations
\begin{align*}
v t + wy +  f_{10}(x,y,z,t)&=0\\
(v+t) w +wg_{5}(x,y)+z^3+y^4+g_{12}(x,y,z,t) &=0
\end{align*}
 where $f_{10}$ has no pure powers of $\xi$ or $t$ and $X \subset \mathbb{P}(2,3,4,5,5,7)$ with homogeneous variables $x,\,y,\,z,\,t,\,v,w$. Then, the rational curve $(x=y=z=w=0) \subset \mathbb{P}(2,3,4,5,5,7)$ intersects $X$ in the two coordinate points $\mathbf{p_t}$ and $\mathbf{p_v}$. Each of these points is a cyclic quotient singularity of type $\frac{1}{5}(1,2,3)$. We consider the toric blowup $\Phi \colon T \rightarrow \mathbb{P}$ whose restriction to $X$ is the unique kawamata blowup centred at $\mathbf{p_v}$. The cyclic quotient singularity $\mathbf{p_v} \sim \frac{1}{5}(1,2,3)$ has local coordinates $z,\,y,\,x$. From Corollary \ref{cor:genlift}, we know how these sections lift under $\varphi$ and what the weight system of $T$ is:
\[
\begin{array}{cccc|ccccccc}
             &       & u  & v & z & w & y  & t & x & \\
\actL{T}   &  \lBr &  0 & 5 & 4 & 7 & 3 & 5 & 2 &    \actR{.}\\
             &       & 1 & 4 & 3 &  5 & 2 & 3 & 0 &  
\end{array}
\] 
We play the 2-ray on $T$. Notice that the movable cone of $T$ is
\[
\overline{\Mov}(T) = \langle \left(\begin{smallmatrix} 5 \\ 4 \end{smallmatrix} \right),\left(\begin{smallmatrix} 5 \\ 3 \end{smallmatrix} \right) \rangle
\]
and its subdivided into four GIT chambers. Choosing a character in the interior of the first chamber, that is, the chamber $\langle \left(\begin{smallmatrix} 5 \\ 4 \end{smallmatrix} \right),\left(\begin{smallmatrix} 4 \\ 3 \end{smallmatrix} \right) \rangle$, gives a variety isomorphic to $T$. Similarly, we define $T_1$, $T_2$ and $T_3$ to be the toric varieties defined by choosing characters in the interior of the second, third  and fourth chambers, respectively. These varieties are related by the small $\mathbb{Q}$-factorial modifications $T \rat T_1 \rat T_2 \rat T_3$ which are 
\[
(-4,-1,1,1,3,2),\, (-7,-3,-1,1,4,3),\, (-3,-2,-1,-1,1,1).
\]
Moreover, there is a divisorial contraction 
\[
\Phi' \colon T_3 \rightarrow \bigoplus_{m\geq 1} H^0(T_3, \mathcal{O}\left(\begin{smallmatrix} 5 \\ 3 \end{smallmatrix} \right)) \simeq \mathbb{P}(1,1,2,2,3,3).
\]
This is in coordinates
\[
(u,v,z,w,y,t,x) \mapsto (t,yx,ux^5,zx^3,vx^5,wx^4).
\]

We restrict this construction to $Y$. The defining equations of $Y$ satisfy $\widetilde{f} \in |-\frac{5}{2}K_Y-\frac{1}{2}E|$ and $\widetilde{g} \in |-3K_Y|$.  Since $z^3 \in g$ by quasismoothness of $X$, it follows that $z^3 \in \widetilde{g}$. Hence $\tau$ restricts to an isomorphism on $Y$. We have that $wy \in \widetilde{f}$ and $vw \in \widetilde{f}$ so that $y$ and $v$ can be locally eliminated around $p_w$. Hence $\tau_1 \colon T_1 \rat T_2$ restricts to the toric flip $(-7,-1,4,3)$ over a point. In the same way, $y^4 \in g$ lifts to $y^4u \in \widetilde{g}$ (and again $wy \in \widetilde{f}$). Hence, $\tau_2\colon T_2 \rat T_3$ restricts to the flip $(-2,-1,1,1)$ over a point. On the other hand $\Phi'$ restricts to a divisorial contraction $\varphi'$ to a point in a Fano 3-fold $Z_{4,6}$, where $Z_{4,6}$ is 
\begin{align*}
v t + wy +  f_{4}(u,1,y,z,t)&=0\\
(v+tu) w +wg_{3}(ux,y)+z^3+y^4u+g_{6}(u,1,y,z,t) &=0
\end{align*}
inside $\mathbb{P}(1,1,2,2,3,3)$ with homogeneous variables $t,\,y,\,u,\,z,\,v,\,w$ where $\varphi'$ can be realised as the weighted blowup 
\[
\wt(u,z,w,y) = (5,3,4,1).
\]
of the $cA_2$ singularity 
\[
0 \in (wu+z^3 + h_{\geq 4}(z,y) = 0) \subset \mathbb{C}^4_{uzwy}.
\]
Moreover, we have that the extremal divisorial contraction $\varphi' \colon Y' \rightarrow Z_{4,6}$ has discrepancy 3 by Lemma \ref{lem:discr}. Compare with the previous case of family 121.

\paragraph{Family 123.} After a change of variables, any quasismooth member of family 123 is given by the equations
\begin{align*}
\xi y+vt+z^3+f_{12}(x^2,z)&=0\\
\xi t+v^2+vg_{7}(x,y,z,t)+t^2z+tg_{9}(x,y,z)+y^4x+y^3t+g_{14}(x,y,z)&=0
\end{align*}
inside $\mathbb{P}(2,3,4,5,7,9)$ with homogeneous variables $x,\,y,\,z,\,t,\,v,\,\xi$.
We blowup the linear cyclic quotient singularity point $\frac{1}{9}(1,4,5)$ on $X_{123}$ with local variables $z,\,v,\,x$. This is the Kawamata blowup $\varphi \colon Y \rightarrow X$ and the anticanonical divisor $-K_{Y}$ is Cartier divisor given by $z=0$ on $Y_{123}$.

Moreover, by Corollary \ref{cor:genlift},
\begin{align*}
v &\in H^0\bigg(Y,-\frac{7}{4}K_Y-\frac{1}{4}E\bigg),\,\, y \in H^0\bigg(Y,-\frac{3}{4}K_Y-\frac{1}{4}E\bigg) \\
t &\in H^0\bigg(Y,-\frac{5}{4}K_Y-\frac{3}{4}E\bigg), \,\, x \in H^0\bigg(Y,-\frac{1}{2}K_Y-\frac{1}{2}E\bigg).
\end{align*}

The pseudo-effective cone of $T \supset Y$ decomposes as 
\[
\mathbb{R}_+[E]+\mathbb{R}_+[\varphi^*(-K_X)]+\mathbb{R}_+[-K_Y]+\mathbb{R}_+[v=0]+\mathbb{R}_+[y=0]+\mathbb{R}_+[t=0]+\mathbb{R}_+[x=0].
\] 

and $Y$ is given by the equations
\begin{align*}
\xi y+vtu+z^3+f_{12}(x^2u,z)&=0\\
\xi t+v^2+vg_{7}(u,x,y,z,t)+t^2zu+tg_{9}(u,x,y,z)+(y^4x+y^3t)u+g_{14}(u,x,y,z)&=0
\end{align*}
where $E \colon (u=0)$ is the exceptional divisor of $\varphi$. Since $z^3 \in \widetilde{f}$ and $v^2 \in \widetilde{g}$ by quasismoothness of $X$, it follows that the $-K_Y$-wall and $v$-wall crossings on $\Eff(T)$ are isomorphisms when restricted to $Y$. The $y$-wall crossing is a $(3,1,1,-1,-1;2)$-flip over a point and finally a divisorial contraction to a $cA_2$ point $\mathbf{p_t} \in Z_{6,6} \subset \mathbb{P}(1,1,2,2,3,5)$. This divisorial contraction can be realised by a weighted blowup with weights
\[
\wt(u,z,v,y)=(5,3,4,1)
\]
and discrepancy 3. Explicitly $Z_{6,6}$ has equations
\begin{align*}
\xi y+vtu+z^3+f_{6}(u,z)&=0\\
\xi t+v^2+vg_{3}(u,1,y,z,t)+t^2zu+tg_{5}(u,1,y,z)+(y^4+y^3t)u+g_{6}(u,1,y,z)&=0
\end{align*}
in $\mathbb{P}(1,1,2,2,3,5)$ with homogeneous variables $t,\,y,\,u,\,z,\,v,\,\xi$. Notice that we have
\[
-K_{Z_{6,6}} \sim \mathcal{O}(2).
\]

\paragraph{Family 125.} Let $X$ be the general member in family 125. W.l.o.g., we can assume that $\mathbf{p_x} \in X$ and so $\mathbf{p_x} \sim \frac{1}{2}(1,1,1)$ with local variables $y,\,z,\,t$. The Kawamata blowup of $\mathbf{p_x}$ is $\varphi \colon Y \rightarrow X$ where $Y$ is 
\begin{align*}
xw+vyu+zt&=0 \\
x^4v+x^3y^3+x^2g_{11}(y,z,t,w) +xg_{13}(y,z,t,vu,w)+wvu+y^5u+g_{15}(y,z,t,vu)&=0
\end{align*}
in
\[
\begin{array}{cccccc|ccccc}
             &       & u  & x & z & t & w  & y & v & \\
\actL{T_1}   &  \lBr &  -3 & -1 & -1 & -1 & -1 & 0 & 1 &    \actR{.}\\
             &       & -8 & -2 & -1 &  -1 & 0 & 1 & 5 &  
\end{array}
\] 
In this case, the movable cone of $T$ is subdivided into three Mori chambers and, by Lemma \ref{lem:iso}, the first wall crossing in $\overline{\Mov}(T)$ is an isomorphism. Crossing to the last chamber in $\overline{\Mov}(T)$ amounts to contracting the curve $(y=v=0) \subset T_1$ to a point and then extract the locus $(u=x=z=t=0) \subset T_2$. The resulting variety $T_2$ has the same Cox ring as $T_1$ but its irrelevant ideal is $(u,x,z,t,w) \cap (y,v)$. By multiplying the weight system of $T_1$ by the matrix 
\[
\begin{pmatrix}
-1	& 0 \\ 
-5	& 1
\end{pmatrix} \in \SL_2(\mathbb{Z})
\]
we can write $T_2$ as 
\[
\begin{array}{ccccccc|cccc}
             &       & u  & x & z & t & w  & y & v & \\
\actL{T_2}   &  \lBr &  3 & 1 & 1 & 1 & 1 & 0 & -1 &    \actR{.}\\
             &       & 7 & 3 & 4  & 4 & 5 & 1 & 0 &  
\end{array}
\] 
We have then a divisorial contraction $\Phi \colon T_2 \rightarrow \mathbb{P}(7,3,4,4,5,1)$ which corresponds to the contraction of the ray $\mathbb{R}_+[\left(\begin{smallmatrix} -1 \\ 0 \end{smallmatrix} \right)]$. That is induced by the sections, 
\[
\bigoplus_{m\geq 1}H^0(T_2,m\mathcal{O}\left(\begin{smallmatrix} 0 \\ 1 \end{smallmatrix} \right))
\]
and is given in coordinates by
\[
(u,z,x,y,w,t,v) \mapsto (uv^3,xv,zv,tv,wv,y).
\] 
We restrict this construction to $Y$. The small modification $T_1 \rat T_2$ contracts the curve $C_1 \colon (f_3(y,z)=0) \subset \mathbb{P}(1,1,8)$ with variables $y,\,z,\,u$ to a point and extracts $C_2 \simeq \mathbb{P}(1,5)$. Over an analytic neighborhood around the base of the contraction, we have then a flip of type $(8,1,1,-1,-5;3)$ with local variables $u,\,z,\,t,\,y,\,v$ and flipping equation $vu+b_3(x,y)+\cdots  = 0$.  The divisorial contraction $\Phi$ restricts to the weighted blowup $\varphi' \colon Y_2 \rightarrow Z_{8,12}$ of weights
\[
\wt(x,z,t,w) = (1,1,1,1)
\]
and discrepancy 1 to the point $cA_1 \sim \mathbf{p_y} \in Z_{8,12}$ where $Z_{8,12}$ is given by the equations,
\begin{align*}
xw+yu+zt&=0 \\
x^4+x^3y^3+x^2g_{11}(y,z,t,w) +xg_{13}(y,z,t,u,w)+wu+y^5u+g_{15}(y,z,t,u)&=0
\end{align*}
in $\mathbb{P}(7,3,4,4,5,1)$ with homogeneous variables $u,\,x,\,z,\,t,\,w,\,y$.

\end{proof}

\subsubsection{Birational Involutions} \label{sec:BI}

In the following we deal with families 117 and 125 with homogeneous coordinates $x_0,\,x_1,\,x_2,\,\xi,\,v,\,w$. These coordinates are such that $\wt(v,w)=(7,8)$ and $x_0 \in H^0(X,-K_X)$. We assume that $\xi^{\iota_X-1}x_{\mu} \not \in g$. These families can be written as
\begin{align*}
\xi w+vx_1 +f_{d_1}&=0\\
\xi^{\iota_X}x_0+\xi^{\iota_X-2}g_{d_2-a_{\xi}(\iota_X-2)}+\cdots+\xi g_{d_2-a_{\xi}}+wv+x_1^{\iota_X}+g_{d_2}&=0.
\end{align*}
where  $\wt(x_1)=a_{\xi}+1$. The point $\mathbf{p_{\xi}}$ is a cyclic quotient singularity of type $\frac{1}{4}(1,3,1)$ and $\frac{1}{2}(1,1,1)$ of families 117 and 125 respectively, with coordinates $(x_1,x_2,v)$. The Kawamata blowup of $X$ centred at $\mathbf{p_{\xi}}$ is an extremal weighted blowup $\varphi \colon Y\rightarrow X$ and $-K_Y \sim \varphi^*(-K_X)-\frac{1}{a_{\xi}}E$. Locally, at $\xi =1$, the maps are given by 
\[
\bigg(\frac{x_1^3}{\xi^2}:\frac{x_2^3}{\xi^3}:\frac{v^3}{\xi^5}\bigg)\,\, \text{and}\,\, \bigg(\frac{x_1}{\xi}:\frac{x_2}{\xi^2}:\frac{v}{\xi^3}\bigg)
\] 
for families 117 and 125, respectively. By Corollary \ref{cor:genlift},
\[
v\in H^0\bigg(Y,-\frac{7}{\iota_X}K_Y+\frac{1}{\iota_X}E\bigg), \quad x_2\in H^0(Y,-K_Y), \quad x_1\in H^0\bigg(Y,-\frac{(8-\iota_X)}{\iota_X}K_Y-\frac{1}{\iota_X}E\bigg).
\]

From the equations defining $X$ we conclude that $w \in H^0\bigg(Y,-\frac{8}{\iota_X}K_Y-\frac{1}{\iota_X}E\bigg)$ and that $x_0 \in H^0(Y,-K_Y-E)$. We now make the following observations: the monomials $\xi w,\,vx_1 \in f$ and $vw \in g$ are anticanonical in $Y$. Moreover, by quasismoothness of $X$, we have $x_2^{\alpha_1} \in f$ or $x_2^{\alpha_2} \in g$ and these lift to $Y$. 

\begin{Lem} \label{lem:isobi}
There is a sequence of small modifications $Y \rat Y'$ whose composition is a birational automorphism.
\end{Lem}

\begin{proof}
The previous considerations divide the Movable cone of $T \supset Y$ into chambers over which we apply the 2-ray game.
The map $Y\rightarrow Z$ given by the linear system generated by the multiples of $v$ contracts the curve $\Gamma_1:=\mathbb{P}(7,1) \colon (x_2=w=x_1=x_0=0) \subset Y$ to a point and extracts $\Gamma_2:=\mathbb{P}(8,1) \colon (u=\xi=0)$. 

Then $-K_Y \cdot \Gamma_1 =(D_v-D_z)\cdot \Gamma_1 =-D_z\cdot \Gamma_1  <0$ since $\Gamma_1$ intersects $D_z$ transversely and is disjoint from $D_v$. Also, $-K_Y \cdot \Gamma_2 > 0$ since $\Gamma_2$ intersects $-K_Y$ transversely. We call $Y_2$ the (anti)-flipped 3-fold.

The map $Y_2 \rightarrow Y_3$ is an isomorphism since  $x_2^{\alpha_1} \in f$ or $x_2^{\alpha_2} \in g$.

Finally, since $\xi w \in f$ and $vw \in g$ are anticanonical in $Y$, the maps given by the linear system generated by multiples of $w$ contract $\widetilde{\Gamma_1}:=\mathbb{P}(8,1) \colon (u=\xi=v=x_2=0) \subset Y$ to a point and extract $\widetilde{\Gamma_2}:=\mathbb{P}(7,1) \colon (x_1=x_0=0)$. The same argument shows that this is indeed a flip.
 \end{proof}

\begin{Lem}
There is a divisorial contraction $Y_4 \rightarrow X'$. This is a Kawamata blowup centred at a cyclic quotient singularity of index $a_{\xi}$ and $X'$ is isomorphic to $X$.
\end{Lem}
\begin{proof}
The 3-folds $X$ and $X'$ have the same defining ideal in the same ambient space.
\end{proof}

Putting together the results of this section we have the following explicit description of the birational involution induced by blowing up $\mathbf{p_{\xi}}$:

\begin{Prop} \label{prop:bi}
Let $X \in \{117, 125 \}$ be a quasismooth member of family 117 or 125 as above. Then, blowing $\mathbf{p_{\xi}}$ initiates a Sarkisov link which is the birational involution
\[
        \begin{tikzcd}[ampersand replacement=\&]
             Y   \ar[rr, dashed, "\displaystyle{(-7,-1,8,1)}" ] \ar[drrr ] \&{} \& Y_2   \ar[rr,  "\displaystyle{\simeq}" ] \& {}  \&Y_3 \ar[rr, dashed, "\displaystyle{(-8,-1,7,1)}" ]\& {} \&Y_4 \ar[dlll ]\\
						{} \&{} \& {} \&\mathbf{p_{\xi}} \in X  \& {} \& {} \& {}
        \end{tikzcd}
    \]
		\end{Prop}

\subsection{Divisorial contractions to a rational curve}

We now consider  divisorial contractions centred in a curve. 
We consider the following cases:
\[
  \begin{cases}
              \text{Case I:} \quad\,\, a_2 = a_4(\iota_X-1) > a_3(\iota_X-1)\\
              \text{Case II:}  \quad\,  a_4(\iota_X-1)= a_3(\iota_X-1) > a_2
            \end{cases}
\]

\paragraph{Case I:} This is satisfied for families 
\begin{center}
\begin{tabular}{lccccc} \toprule
    $X$ & 92 & 101    \\
    $\mathbf{p_{\xi}}$  & $\frac{1}{5}(1,1,4)$  &$\frac{1}{7}(1,1,6)$    \\ 
		\bottomrule
		\end{tabular}
\end{center}

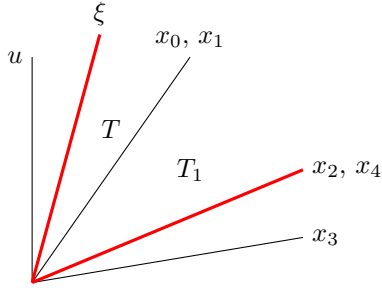
\begin{figure}%
\centering
\begin{tikzpicture}[scale=3]
 \coordinate (A) at (0, 0);
 \coordinate [label={left:$u$}] (E) at (0, 1);
 \coordinate [label={above:$x_0,\,x_1$}] (K) at (0.7, 1);
 \coordinate [label={above:$\xi$}] (5) at (0.3,1.1);
 \coordinate [label={right:$x_2,\,x_4$}] (4) at (1.2,0.5);
 \coordinate [label={right:$x_3$}] (3) at (1.2,0.2);

\coordinate [label={$T$}] (T) at (0.35,0.6);
\coordinate [label={right:$T_1$}] (T1) at (0.6,0.5);

  \draw (A) -- (E);
  \draw (A) -- (K);
	\draw [very thick,color=red] (A) -- (5);
	\draw (A) -- (3);
	\draw [very thick,color=red] (A) -- (4);
  
\end{tikzpicture}
\caption{A representation of the chamber decomposition of the cone of effective divisors of $T$. In red the subcone of movable divisors of $T$. The ample models $T$, $T_1$ are represented in the interior of the corresponding Nef chambers.}%
\label{fig:HypCaseIII}%
\end{figure}

Both families have Fano index $\iota_X=2$. The assumption is therefore equivalent to $a_2=a_4> a_3$. We consider the toric blowup $\Phi \colon T \rightarrow \mathbb{P}$ that restricts to the unique Kawamata blowup $\varphi \colon E \subset Y \rightarrow X$ centred at $\mathbf{p_{\xi}} \in X$. By Lemma \ref{lem:lift} we have in particular 
\[
x_2 \in H^0\bigg(Y,-\frac{a_2}{2}K_Y- \frac{1}{\iota_X}E \bigg) = H^0\bigg(Y,-\frac{a_4}{2}K_Y- \frac{1}{\iota_X}E \bigg) \ni x_4
 \]
and, since $a_4 > a_3$, the cone of movable divisors of $T$ is strictly contained in the interior of the cone of effective divisors of $T$, as in Figure \ref{fig:HypCaseIII}.
 From Lemma \ref{lem:iso}, the small $\mathbb{Q}$-factorial modification between the ample models $T$ and $T_1$ restricts to an isomorphism on $Y$. 

\begin{Thm} \label{thm:divtocurve}
Let $X \in I_S$ be such that $a_2= a_4 > a_3$. Then, there is an elementary Sarkisov link to a singular Fano 3-fold orbifold hypersurface $X'_{\frac{d_2}{\iota_X}} \subset \mathbb{P}(a_3,1,1,1,1)$, 
 \[
        \begin{tikzcd}[ampersand replacement=\&,column sep = 2em]
             \& Y \ar[dl, swap, "{\varphi}"] \ar[dr, "{\varphi'}"] \& \\
            \mathbf{p_{\xi}} \in X    \&  \&  X' \supset \Gamma 
        \end{tikzcd}
    \]
		initiated by the Kawamata blowup $\varphi$ of a linear cyclic quotient singularity $\mathbf{p_{\xi}}$ and ending with a divisorial contraction $\varphi'$ to a rational curve. The singular locus of $X'$ is nonempty and is contained in $\Gamma$.

\end{Thm}

\begin{proof}
We need only to analyse the divisorial contraction
\[
\Phi' \colon T_1 \rightarrow \mathcal{F}', \quad \mathcal{F}':=\Proj \bigoplus_{m\geq 1} H^0(T_1,\mathcal{O}\left(\begin{smallmatrix} 1 \\ 0 \end{smallmatrix} \right))
\]
 where 
\[
\begin{array}{cccccc|ccccc}
             &       & u  & \xi & x_0 & x_1 & x_2 & x_4 & x_3 & \\
\actL{T_1}   &  \lBr &  a_3 & \frac{d_1}{2} & 1 & \frac{a_1}{2} & 1 & 1 & 0 &    \actR{.}\\
             &       & a_4 & \frac{d_2}{2} & 1 &  \frac{a_1}{2} & 0 & 0 & -1 &  
\end{array}
\] 
We can simplify slightly by noticing by inspection that $\frac{a_1}{2}=1$. This map is given in coordinates by
\begin{align*}
\Phi' \colon T_1 & \longrightarrow \mathbb{P}\Big(a_3, \frac{d_1}{2}, 1, 1, 1, 1 \Big)\\
(u,\xi,x_0,x_1, x_2,x_4,x_3) &\longmapsto (ux_3^{a_4}:\xi x_3^{d_2/2}:x_0x_3:x_1x_3:x_2:x_4)
\end{align*}
and is a divisorial contraction to the curve $\Gamma \simeq \mathbb{P}^1 \colon (u=\xi=x_0=x_1=0) \in \mathbb{P}\Big(a_3, \frac{d_1}{2}, 1, 1, 1, 1 \Big)$.

We restrict $\Phi'$ to $Y \subset T_1$. The 3-fold $Y$ is $(\widetilde{f}=\widetilde{g}=0) \subset T_1$, where $\widetilde{f} \colon \xi x_3+f'(u,x_0,x_1,x_2,x_4,x_3) \in |-\frac{d_1}{2}K_Y|$ and $\widetilde{g} \colon \xi x_4+g'(u,x_0,x_1,x_2,x_4,x_3) \in |-\frac{d_2}{2}K_Y|$. In particular, $f'$ and $g'$ contain no monomials purely in $(x_2,x_4)$. Hence, $\Phi'$ restricts to a divisorial contraction $\varphi' \colon Y \rightarrow X'$ to $\Gamma$, where $X'$ is given by a degree $\frac{d_2}{\iota_X}$ equation:
\[
X'\colon (-x_4f'(u,x_0,x_1,x_2,x_4,1)+g'(u,x_0,x_1,x_2,x_4,1)=0) \subset \mathbb{P}(a_3,  1, 1, 1, 1 ).
\]
The singular locus of $X'$ is clearly contained in $\Gamma$. Moreover, by quasismoothness of $X$ there are monomials $f'(x_2,x_4) \in f$ and, since $f' \in |-mK_Y|$ for some $m>0$, these monomials lift to $uf'(x_2,x_4) \in \widetilde{f}$. Hence, 
\[
\partial_u (-x_4f'(u,x_0,x_1,x_2,x_4,1)+g'(u,x_0,x_1,x_2,x_4,1))|_{\Gamma}=0
\]
is a non trivial condition on $x_2,\,x_4$ while the partial derivative with any other variable vanishes identically.
\end{proof}

\paragraph{Case II:} This is satisfied for families
\begin{center}
\begin{tabular}{lccccc} \toprule
    $X$ & 91 & 113 & 123   \\
    $\mathbf{p_{\xi}}$  & $\frac{1}{3}(1,1,2)$  &$\frac{1}{4}(1,1,3)$ & $\frac{1}{5}(1,2,3)$    \\ 
		$\iota_X$ & 2 & 3 & 4 \\
 		\bottomrule
		\end{tabular}
\end{center}

We start with the families of prime Fano index.

\begin{Thm} \label{thm:cod2curve}
Let $X \in I_S$ be such that $a_4(\iota_X-1)= a_3(\iota_X-1) > a_2$. Then, there is an elementary Sarkisov link to a singular Fano 3-fold complete intersection $X'$, 
 \[
        \begin{tikzcd}[ampersand replacement=\&,column sep = 2em]
             \& Y \ar[dl, swap, "{\varphi}"] \ar[dr, "{\varphi'}"] \& \\
            \mathbf{p_{\xi}} \in X    \&  \&  X' \supset \Gamma 
        \end{tikzcd}
    \]
		initiated by the Kawamata blowup $\varphi$ of a linear cyclic quotient singularity $\mathbf{p_{\xi}}$ and ending with a divisorial contraction $\varphi'$ to a rational curve $\Gamma$. The singular locus of $X'$ is non-empty and contained in $\Gamma$.
		
\end{Thm}

\begin{proof}
This result is analogous to Theorem \ref{thm:divtocurve}. By the assumption on the weights, the movable cone of $T$ is like in Figure \ref{fig:HypCaseIII} with the rays spanned by $x_2$ and $x_3$ swapped. By Lemma \ref{lem:iso}, the small modification $T \rat T_1$ restricts to an isomorphism on $Y$. We write the weight system of $T_1$ as
\[
\begin{array}{cccccc|cccccc}
             &       & u  & \xi &   x_0 & x_1 & x_4 & x_3 & x_2 & \\
\actL{T_1}   &  \lBr &  a_2&a_{\xi}-\frac{a_1}{\iota_X}&\iota_X-1&(\iota_X-1)\frac{a_1}{\iota_X}&\frac{a_4(\iota_X-1)-a_2}{\iota_X}&\frac{a_3(\iota_X-1)-a_2}{\iota_X}&0  & \actR{.}\\
             &       & a_3	& \frac{d_1}{\iota_X} & 1 & \frac{a_1}{\iota_X}& 0 &0&-\frac{a_3(\iota_X-1)-a_2}{\iota_X} &
 &  
\end{array}
\]   
Given the small number of families treated in this result (families 91 and 113) we can simplify the the weights of $T$  noticing that $\frac{a_4(\iota_X-1)-a_2}{\iota_X}=\frac{a_3(\iota_X-1)-a_2}{\iota_X}=\frac{a_1}{\iota_X} = 1$. There is a divisorial contraction 
\[
\Phi' \colon T_1 \rightarrow \mathcal{F}', \quad \mathcal{F}': = \Proj \bigoplus_{m\geq 1}H^0(T_1,m\mathcal{O}(x_3)) \simeq \mathbb{P}\Big(a_2,a_{\xi}-1,\iota_X-1,\iota_X-1,1,1\Big).
\]
given in coordinates by
\begin{align*}
\Phi' \colon T_1 &\longrightarrow \mathbb{P}\Big(a_2,a_{\xi}-1,\iota_X-1,\iota_X-1,1,1\Big)\\
(u,\xi,x_0,x_1,x_4,x_3,x_2) & \longmapsto (ux_2^{a_3}:\xi x_2^{d_1/2}:x_0x_2:x_1x_2:x_4:x_3).
\end{align*}
This is a divisorial contraction to the curve $\Gamma \simeq \mathbb{P}^1 \colon (u= \xi=x_0=x_1=0) \subset \mathcal{F}'$.
We restrict $\Phi'$ to $Y \subset T_1$. The 3-fold $Y$ is $(\xi x_3+f'(u,x_0,x_1,x_4,x_3,x_2)=\xi x_4+g'(u,x_0,x_1,x_3,x_3,x_2)=0) \subset T$, so its image under $\varphi'$ is given by
\[
X'\colon (\xi x_3+f'(u,x_0,x_1,x_4,x_3,1)=\xi x_4+g'(u,x_0,x_1,x_4,x_3,1)=0) \subset \mathbb{P}\Big(a_2,a_{\xi}-1,\iota_X-1,\iota_X-1,1,1\Big).
\]
Since the defining equations of $Y\colon (\widetilde{g}=\widetilde{g}=0) \subset T_1$ are such that $\widetilde{f} \in |-m_1K_Y|$ and $\widetilde{f} \in |-m_2K_Y|$ for some positive integers $m_1$ and $m_2$, it follows that $\Gamma \subset X'$.
\end{proof}
\paragraph{Family 123.} Let $X$ be a general member of family 123. We claim that there are polynomials $h_7 \in \mathbb{C}[x,y,z]$ and $h_9 \in \mathbb{C}[x,y,z,v]$ for which the ideal of $X$ is 
\begin{align*}
\xi (v+h_7(x,y,z)) + y(w -h_9(x,y,z,v)) + z^3+ f_{12}(x,y,z,v) &=0\\
\xi w + v^2+y^4x+g_{14}(x,y,z,v) &=0
\end{align*}

A general member $X$ of family 123 can be written as
\begin{align*}
\xi^2x+\xi (v+f_7(x,y,z)) + wy + z^3+ f_{12}(x,y,z,v) &=0\\
\xi^2(z+x^2)+\xi (w+g_9(x,y,z,v)) + v^2+y^4x+g_{14}(x,y,z,v,w) &=0.
\end{align*}
We can assume that $g_{14} \in \mathbb{C}[x,y,z,v]$ by substituting $g$ by $xf-g$.  We make the change of variables $ v \mapsto v-f_7(x,y,z)-\xi x$ to get
\begin{align*}
\xi (v+f_7(x,y)) + wy-g_9(x,y,z,v)y + z^3+ f_{12}(x,y,z,v) &=0\\
\xi^2(z+x^2)+\xi (w+g_9(x,y,z,v)-2x(v-f_7)) + v^2+y^4x+g_{14}(x,y,z,v) &=0
\end{align*}
since $(v-f_7(x,y,z)-\xi x)^2 = v^2-2vf_7+f_7^2-2\xi x(v-f_7)+ \xi^2 x^2$. Let $g_9' := g_9-2x(v-f_7)$. We change variables again by $w \mapsto w- g_9'- \xi (z+x^2)$ to get
\begin{align*}
\xi (v+f_7(x,y)-y(z+x^2)) + y(w -(g_9+g_9')) + z^3+ f_{12}(x,y,z,v) &=0\\
\xi w + v^2+y^4x+g_{14}(x,y,z,v) &=0
\end{align*}
and the claim follows.

\begin{Thm} \label{thm:cod2curveexcept}
Let $X$ be a general quasismooth member of family 123. Then, there is an elementary Sarkisov link to a singular Fano 3-fold complete intersection $Z_{6,6}$, 
 \[
        \begin{tikzcd}[ampersand replacement=\&,column sep = 2em]
             \& Y \ar[dl, swap, "{\varphi}"] \ar[dr, "{\varphi'}"] \& \\
            \frac{1}{5}(1,2,3) \in X    \&  \&  Z_{6,6} \supset \Gamma 
        \end{tikzcd}
    \]
		initiated by the Kawamata blowup $\varphi$ of the linear cyclic quotient singularity $\frac{1}{5}(1,2,3)$ and ending with a divisorial contraction $\varphi'$ to a rational curve $\Gamma$. The singular locus of $Z_{6,6}$ is non-empty and contained in $\Gamma$.
		
\end{Thm}
\begin{proof}
The ideal defining $X$ is $(f,g)$ where $f$ and $g$ are the polynomials, 
\begin{align*}
\xi (v+h_7(x,y,z)) + y(w -h_9(x,y,z,v)) + z^3+ f_{12}(x,y,z,v) &=0\\
\xi w + v^2+y^4x+g_{14}(x,y,z,v) &=0
\end{align*}
as in the previous paragraph. $X$ is inside $\mathbb{P}:=\mathbb{P}(2,3,4,5,7,9)$ with homogeneous variables $x,\,y,\,z,\,\xi,\,v,\,w$. We consider the toric blowup $\Phi \colon T \rightarrow \mathbb{P}$ whose restriction to $X$ is the unique Kawamata blowup $\varphi \colon Y \rightarrow X$ centred at $\mathbf{p_{\xi}}\sim\frac{1}{5}(1,2,3)$ with local coordinates $z,\,y,\,x$. By Corollary \ref{cor:genlift}, we can compute how the sections of $X$ lift under $\varphi$. We have that $z$ lifts to an anticanonical section of $Y$ and
\begin{align*}
v \in H^0\bigg(Y,-\frac{7}{4}K_Y-\frac{1}{4}E\bigg),\, y \in H^0\bigg(Y,-\frac{3}{4}K_Y-\frac{1}{4}E\bigg) \\
w \in H^0\bigg(Y,-\frac{9}{4}K_Y-\frac{3}{4}E\bigg),\, x \in H^0\bigg(Y,-\frac{1}{2}K_Y-\frac{1}{2}E\bigg).
\end{align*}
Moreover $\widetilde{f} \in |-3K_Y|$ and $\widetilde{g} \in |-\frac{7}{2}K_Y-\frac{1}{2}E|$. As sections of $T$ these give the weight system,
 \[
\begin{array}{cccc|ccccccc}
             &       & u  & \xi & z & v & y & w & x & \\
\actL{T}   &  \lBr &  0 & 5 & 4 & 7 & 3 & 9 & 2 &    \actR{.}\\
             &       & 1 & 4 & 3 & 5 & 2 & 6 & 1 &  
\end{array}
\]
In particular, the cone of effective divisors is $\overline{\Eff}(T)= \langle \left(\begin{smallmatrix} 0 \\ 1 \end{smallmatrix} \right),\left(\begin{smallmatrix} 2 \\ 1 \end{smallmatrix} \right) \rangle$ and the subcone of movable divisors is $\overline{\Mov}(T)=\langle \left(\begin{smallmatrix} 5 \\ 4 \end{smallmatrix} \right),\left(\begin{smallmatrix} 3 \\ 2 \end{smallmatrix} \right) \rangle$. Moreover, the cone of movable divisors of $T$ is subdivided into three chambers separated by the rays $\mathbb{R}_+[\left(\begin{smallmatrix} 4 \\ 3 \end{smallmatrix} \right)]$ and $\mathbb{R}_+[\left(\begin{smallmatrix} 7 \\ 5 \end{smallmatrix} \right)]$.

We play the 2-ray game on $T$ by choosing characters successively in $\overline{\Mov}(T)$. We have a composition of small $\mathbb{Q}$-factorial modifications, $T \rat T_1 \rat T_2$ where $T_1$ and $T_2$ are the ample models of the different chambers in $\overline{\Mov}(T)$. For instance, $T_2$ has the same Cox ring as $T$ but its irrelevant ideal is $(u,\xi,z,v)\cap(y,w,x)$. The 2-ray game on $T$ ends with the divisorial contraction
\[
\Phi' \colon T_2 \rightarrow \Proj \bigoplus_{m\geq 1}H^0(T_2, m\mathcal{O}\left(\begin{smallmatrix} 3 \\ 2 \end{smallmatrix} \right)) = \Proj\mathbb{C}[y,ux^3,zx,w,tx^2,vx] = \mathbb{P}(1,2,2,3,3,3).
\]  
We restrict the 2-ray game on $T$ to $Y$. Since $z^3 \in |-3K_Y|$ and $v^2 \in |-\frac{7}{2}K_Y-\frac{1}{2}E|$ it follows that $z^3 \in \widetilde{f}$ and $v^2 \in \widetilde{g}$. Hence, the small modifications $T \rat T_1 \rat T_2$ all happen away from $Y$. In other words, they restrict to an isomorphism on $Y$.

The equations of $Y\subset T_2$ are 
\begin{align*}
\xi (v+h_7(ux,y,z)) + y(wu -h_9(ux,y,z,v)) + z^3+ f_{12}(u,x,y,z,v) &=0\\
\xi w + v^2+y^4xu+g_{14}(u,x,y,z,v) &=0.
\end{align*}
Then, $\Phi'$ contracts the divisor $x=0$ to the curve $\Gamma \colon (u= \xi =z = v=0) \subset \mathbb{P}(1,2,2,3,3,3)$ and, away from this divisor it is an isomorphism to $Z_{6,6}$ whose equations are
\begin{align*}
\xi (v+h_3(u,y,z)) + y(wu -h_5(u,y,z,v)) + z^3+ f_{6}(u,1,y,z,v) &=0\\
\xi w + v^2+y^4u+g_{6}(u,1,y,z,v) &=0
\end{align*}
inside $\mathbb{P}(1,2,2,3,3,3)$ with homogeneous variables $y,\,u,\,z,\,\xi,\,v,\,w$. Notice that $\Gamma \subset Z_{6,6}$. It is easy to see that the singular points of $Z_{6,6}$ as an orbifold are $\mathbf{p_w}$ and $\mathbf{p_y}$.
\end{proof}

\section{Sarkisov Links to \textit{non}-complete intersection Fano 3-folds}

In this section we construct Sarkisov links to Fano 3-folds in codimension 4 or higher. 

\subsection{Sarkisov Links to Fano 3-folds in codimenson 4}

We prove Theorem \ref{thm:cod4}. The idea of the proof is as follows: we want to obtain a Sarkisov link initiated by blowing up a non-linear cyclic quotient singularity. We do it by following the Sarkisov Program on an ambient toric variety $T \subset Y$. It turns out that the game on the 3-fold $Y$ does not follow the game played on $T$. In order to overcome this issue, we successively enlarge the ambient toric variety by means of unprojecting divisors. We do so in such a way that the unprojection happens away from the 3-fold so that these transformations can be seen as biregular changes of coordinates. 

\begin{Thm} \label{thm:cod4}
Let $X$ be a quasismooth member of the deformation family in the second column of table \ref{tab:cod4}. Then, $X$ is birational to a quasismooth Fano 3-fold $X'$ of Fano index 1 embedded as a codimension 4 as in the third column of T \ref{tab:cod4}. Moreover, there is a birational map $X \dashrightarrow X'$  which decomposes as 
\[
        \begin{tikzcd}[ampersand replacement=\&, column sep = 2em]
             Y   \ar[rr, dashed, "\displaystyle{s \times \tau}" ] \ar[d, swap, "\displaystyle{\varphi}"] \& \& Y' \ar[d,, "\displaystyle{\varphi'}" ] \\
           \mathbf{p} \in X  \& \& \mathbf{p'} \in X' 
        \end{tikzcd}
    \]
where $\varphi$ and $\varphi'$ are the Kawamata blowups of $\mathbf{p}$ and $\mathbf{p'}$, respectively and $s \times \tau$ is a number of simultaneous Atiyah flops.
\end{Thm} 

\begin{center}
\begin{table}
\begin{tabular}{lll} \toprule
    {ID} & $\mathbf{p} \in X \subset w\mathbb{P}$ &$\mathbf{p'} \in X' \subset w\mathbb{P}'$   \\  \midrule 
      99  & $2 \times \frac{1}{3}(1,1,2) \in X_{10,12} \subset \mathbb{P}(1,2,3,5,6,7)$ &  $\frac{1}{10}(1,3,7) \in X' \subset \mathbb{P}(1,1,2,3,3,4,7,10)$ \\ 
			103  & $4 \times \frac{1}{3}(1,1,2) \in X_{10,12} \subset \mathbb{P}(2,3,3,4,5,7)$ &  $\frac{1}{9}(1,2,7) \in X' \subset \mathbb{P}(1,2,3,3,4,5,7,9)$ \\ 
			104  & $\frac{1}{5}(1,1,4) \in X_{14,16} \subset \mathbb{P}(1,2,5,7,8,9)$ &  $\frac{1}{13}(1,4,9) \in X' \subset \mathbb{P}(1,1,3,4,4,5,9,13)$ \\
		106  & $2\times \frac{1}{3}(1,1,3) \in X_{18,20} \subset \mathbb{P}(2,5,6,7,9,11)$ &  $\frac{1}{12}(1,1,11) \in X' \subset \mathbb{P}(1,1,6,8,9,10,11,12)$ \\
		111  & $\frac{1}{7}(1,1,6) \in X_{18,20} \subset \mathbb{P}(2,5,6,7,9,11)$ &  $\frac{1}{14}(1,3,11) \in X' \subset \mathbb{P}(1,3,5,6,7,8,11,14)$ \\
		\bottomrule
		\end{tabular} 
		\caption{Birational Models of $X$ in codimension 4}
		\label{tab:cod4}
		\end{table}
\end{center}
\begin{proof}
Let $X$ be a quasismooth member of the deformation family in the second column of Table \ref{tab:cod4}.  Then we can write $X$ as 
\begin{align*}
\xi w + \zeta^2+\zeta f_{d_1-a_{\zeta}}(x_0,x_1,x) + f_{d_1}(x_0,x_1,x,w)&=0\\
\xi^rx_{\mu} + \xi^{r-1}g_{d_2-(r-1)a_{\xi}}+\cdots+\xi g_{d_2-a_{\xi}}(x_0,x_1,x,\zeta,w)+w\zeta +g_{d_2}(x_0,x_1,x)&=0
\end{align*}
where 
\begin{itemize}
	\item $X \subset \mathbb{P}:=\mathbb{P}(a_{x_0},a_{x_1},a_{x},a_{\xi},a_{\zeta},a_{w})$ with homogeneous variables $x_0,\,x_1,\,x,\,\xi, \, \zeta, \, w$.
	\item $x_0$ and $x_1$ are the only variables with even weight;
	\item $x_{\mu} = x_1$ for families 99 and 111 and $x_{\mu} = x$ otherwise.
	\item $\mathbf{p_{\xi}}$ is a cyclic quotient singularity of index $3,\,3,\,5,\,3,\,7$ on families 99, 103, 104, 106 and 111, respectively.
	\item $r=2$ for families 99 and 111 and $r=3$ for families 103, 104 and 106.
\end{itemize}
We divide the proof into several steps. 
\paragraph{Step 1: Lift of sections.} The point $\mathbf{p_{\xi}}$ is a cyclic quotient singularity of type $\frac{1}{a_{\xi}}(1,1,a_{\xi}-1)$. Let $\Phi \colon T \rightarrow \mathbb{P}$ be the toric blowup centred at $\mathbf{p_{\xi}}$ whose restriction to $X$ is the unique Kawamata blowup $\varphi \colon (E\subset Y) \rightarrow (\mathbf{p_t} \in X)$ centred at $\mathbf{p_{\xi}}$. We use Corollary \ref{cor:genlift} to lift the sections $\{x_0,x_1,x,\xi,\zeta,w \}$ of $X$ to $Y$ and get the effective cone of $T$ in Figure \ref{fig:codim4T1}. Moreover, $\widetilde{f} \in |-\frac{d_1}{2}K_Y+E|$ and $\widetilde{g} \in |-\frac{d_2}{2}K_Y+E|$, where $d_1$ and $d_2$ are the degrees of $f$ and $g$, respectively. Hence $Y$ is given by
\begin{align*}
\xi w + \zeta^2+u\zeta f_{d_1-a_{\zeta}}(u,x_0,x_1,x) + uf_{d_1}(u,x_0,x_1,x)&=0\\
\xi^rx_{\mu} + \xi^{r-1}g_{d_2-(r-1)a_{\xi}}+\cdots+\xi g_{d_2-a_{\xi}}+w\zeta +ug_{d_2}(u,x_0,x_1,x)&=0.
\end{align*}
We also claim that $g_{d_2-a_{\xi}} \in (u,\zeta,w)$. Indeed, $g_{d_2-a_{\xi}}$ defines a divisor in the linear system in $|-\frac{d_2}{2}K_Y+E + \frac{a_{\xi}}{2}K_Y-\frac{1}{2}E | = |-\frac{a_w+a_{\zeta}-a_{\xi}}{2}K_Y+\frac{1}{2}E|$ and in each case $a_{\zeta} > a_{\xi}$. So the ray $\mathbb{R}_+\bigg[ -\frac{a_w+a_{\zeta}-a_{\xi}}{2}K_Y+\frac{1}{2}E \bigg]$ is in between the rays generated by $-\frac{a_w}{2}K_Y+\frac{1}{2}E$ and $-m_iK_Y$. This implies that $g_{d_2-a_{\xi}} \in (u,\xi, \zeta,w) \cap (x_0,x_1,x)$ and  the conclusion follows from the assumption that $\xi \not \in g_{d_2-a_{\xi}}$.

\paragraph{Step 2: First Unprojection.} The ample models in the first two chambers of the movable cone of $T$ are related by a small modification $T \rat T'$ where $T'$ s such that $\Cox(T)=\Cox(T')$ and the irrelavant ideal of $T'$ is $(u,\xi,\zeta) \cap (w,x_0,x_1,x)$. This restricts to an isomorphism on $Y$ since $\zeta^2 \in \widetilde{f}$. On the other hand, the ample models in the second and third chambers of the movable cone of $T$ are also related by a small modification, however when restricted to $Y$ we have a divisorial extraction so the 2-ray game on $Y$ does not follow the 2-ray game on $T$.

\begin{figure}%
\centering
\begin{tikzpicture}[scale=3,font=\tiny]
  \coordinate (A) at (0, 0);
  \coordinate [label={left:$E$}] (E) at (0, 0.8);
  \coordinate [label={above:$-\frac{a_{\xi}}{2}K_Y+\frac{1}{2}E$}] (K) at (0.5,1);
	\coordinate [label={right:$-\frac{a_{\zeta}}{2}K_Y+\frac{1}{2}E$}] (5) at (1, 1);
	\coordinate [label={right:$-\frac{a_{w}}{2}K_Y+\frac{1}{2}E$}] (2) at (1.2,0.8);
	\coordinate [label={right:$-m_iK_Y$}] (4) at (1.2,0.5);
	\coordinate [label={right:$-\frac{a_{x}}{2}K_Y-\frac{1}{2}E$}] (3) at (1.2,0.2);
  \draw  (A) -- (E);
  \draw [very thick,color=red] (A) -- (4);
	\draw (A) -- (5);
	\draw (A) -- (2);
	\draw (A) -- (3);
	\draw [very thick,color=red] (A) -- (K);
\end{tikzpicture}
\caption{A representation of the chamber decomposition of the cone of effective divisors of $T$ where $-K_Y \sim \mathcal{O}(2,1)$ and $E \sim \mathcal{O}(0, 1)$. Starting from $E$ and going clockwise, we have the sections $u,\,\xi,\,\zeta, w, x_0, x_1$ and $x$. Notice that each $x_i$ is (pluri)-anticanonical.}%
\label{fig:codim4T1}%
\end{figure}
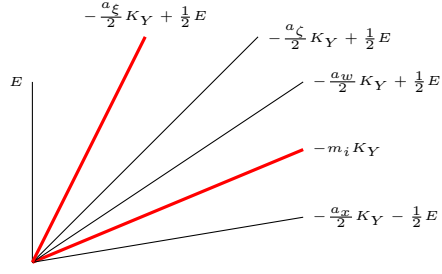

Notice that $\widetilde{f},\,\widetilde{g} \in (u,\xi,\zeta)$ and so there are $F_i, \,G_i \in \mathbb{C}[u,\xi,\zeta,w,x_0,x_1,x]$ with $1\leq i\leq 3$ such that 
\[
\begin{pmatrix}
           f^* \\
           g^* \\
         \end{pmatrix} = 
Q
\begin{pmatrix}
           u \\
           \xi \\
          \zeta \\
         \end{pmatrix},
				\quad
	Q:=	\left(
  \begin{array}{ccc}
    F_1 & F_2 & F_3 \\
    G_1 & G_2 & G_3 \\
  \end{array}
\right)		
\]
Explicitly, 
\[Q=
\left(
  \begin{array}{ccc}
    f_{d_1}+\zeta f_{d_1-a_{\zeta}} & w & \zeta \\
    g_{d_2} & \xi^{r-1}x_{\mu} + \xi^{r-2}gd_2-(r-1)a_{\xi}+\cdots+ \xi g_{d_2-2a_{\xi}}+g_{d_2-a_{\xi}} & w \\
  \end{array}
\right)
\]
where $F_1,\,G_1 \in \mathbb{C}[u,x_0,x_1,x]$ and the other $F_i,\,G_i$ are polynomials in the ideal $(u,\xi,\zeta,w)$.
Following \cite{papadakis}, for a matrix $A_{r\times (r+1)}$, we define $\bigwedge^r A $ to be a vector with $r+1$ entries whose $i^{th}$ entry is $(-1)^{i+1} |A_i|$, where $A_i$ is the submatrix of $A$ resulting from removing its $i^{th}$ column. By Cramer's Rule we have 
\[
\big(\bigwedge^2 Q\big)_i z_j-\big(\bigwedge^2 Q\big)_j z_i \in (f^*,g^*)
\]
where $z_i$ is the $i^{th}$ element of $\{u,\xi,\zeta\}$. Therefore, the ratio
\[
\eta:=\frac{(\bigwedge^2 Q)_i}{z_i} \in H^0(Y,-a_wK_Y),\quad 1\leq i\leq 3
\]
is a well-defined (pluri)-anticanonical section of $Y$. Explicitly this is
\[
\eta:=\frac{F_2G_3-G_2F_3}{u}=-\frac{F_1G_3-G_1F_3}{\xi}=\frac{F_1G_2-G_1F_2}{\zeta} \in H^0(Y,-a_wK_Y)
\]
Let $T^{\eta}= \Proj \mathbb{C}[u,\xi,\zeta,w,x_0,x_1,\eta,x]$, i.e., this is $T$ with the extra section $\eta$ and 
\[
Y^{\eta} \colon \Big(\widetilde{f}=\widetilde{g}=\eta z_i - (\bigwedge^2 Q)_i=0 \, |\, 1\leq i\leq 3\Big)
\]
inside $T^{\eta}$. Then the inclusion map
\begin{align*}
i^{\eta} \colon Y & \longrightarrow Y^{\eta}\\
(u,\xi,\zeta,w,x_0,x_1,x)&\longmapsto (u,\xi,\zeta,w,x_0,x_1,\eta,x)
\end{align*}
 is biregular since the irrelavant ideal of $Y$ is in $(u,\xi,\zeta)$.  Moreover, the ideal defining $Y^{\eta}$, $I_{Y^{\eta}}$, can be written in terms of the 5 maximal Pfaffians of the $5\times 5$ antisymmetric matrix 
\[
M=
\left(
  \begin{array}{cccc}
    u & -\xi & G_3 & F_3 \\
		 & \zeta & G_2 & F_2 \\
		 &  & G_1 & F_1 \\
		 &  &  & \eta \\
  \end{array}
\right)
\]
that we write as $I_{Y^{\eta}}=(\Pf_1,\ldots,\Pf_5)$. The unprojection equations, $\eta z_i - (\bigwedge^2 Q)_i=0$, are the first three maximal Pfaffian equations and each of these is a divisor in the linear system $|-r_iK_Y+s_iE|$ where $r_i,\,s_i >0$.

Notice that, by enlarging the ambient space in which $Y \simeq Y^{\eta}$ lives, we get in particular the extra equation, $\eta u = w^2-\zeta G_2$. Hence, crossing the $w$-wall on $T^{\eta}$ restricts to an isomorphism on $Y^{\eta}$. Furthermore, we have a divisorial contraction of $(x=0)$ to a 6-dimensional weighted projective space  $T^{\eta} \rightarrow \mathbb{P}$. The restriction to $Y^{\eta}$ is the map given by the sections of multiples of $-K_{Y^{\eta}}$. 
\[
Y^{\eta} \rightarrow \mathcal{F}, \quad \mathcal{F}:=\Proj\bigoplus_{m\geq 1} H^0(Y^{\eta},-mK_{Y^{\eta}})=\Proj \mathbb{C}\mathbb[x_0,x_1,\eta,\ldots]
\]
where the rest of the monomials are multiples of $x$. Hence $(x=0)$ is mapped to $\Proj\mathbb{C}[x_0,x_1,\eta]$ and the 2-ray game on $Y^{\eta}$ does not follow, once again, the 2-ray game on $T^{\eta}$.

\paragraph{Step 3: Second Unprojection.}

Notice that the first two rows of $M$ are polynomials in $J=(u,\xi,\zeta,w)$ and so $I_{Y^{\eta}}\subset J$ and $Y^{\eta}$ is a codimension $3$ Fano 3-fold written in Jerry$_{12}$ format. From Step 1, there are polynomials $H_i$, $1 \leq i \leq 4$ for which $G_2 = uH_1 + \xi H_2 + \zeta H_3 + w H_4$. We have,
\[
\begin{pmatrix}
           \Pf_1 \\
           \Pf_2 \\
					 \Pf_3 \\
         \end{pmatrix} = 
\left(
  \begin{array}{cccc}
    -H_1F_1 & -H_2F_1 & \eta - H_3F_1  & G_1-H_4F_1 \\
    0 & -\eta & G_1 & -F_1\\
		\eta + H_1F_3 & H_2F_3 & H_3F_3 & -G_3+H_4F_3\\
  \end{array}
\right)
\begin{pmatrix}
           u \\
           \xi \\
          \zeta \\
					w \\
         \end{pmatrix}
\]
From this matrix we produce polynomials $h_i$ which are given by
\[
 (h_1,\ldots, h_4) := \bigwedge^3 Q. 
\]

It follows from \cite[Lemma~5.11]{papadakis} that for each $1 \leq i\leq 4$ there are polynomials $K_i$ and $L_i$ such that
\[
h_i=\eta K_i + (wF_1-\zeta G_1)L_i = \eta (K_i-\xi L_i)-L_i \Pf_2
\]

Let $z_i$ be the $i^{th}$ element in $\{u,\xi,\zeta,w\}$ and $g_i=K_i-\xi L_i$ for each $1 \leq i\leq 4$. The point here is that 
\begin{align*}
I_{Y^{\eta}} &\ni h_iz_j-h_jz_i && \text{Cramer's Rule} \\
&=\eta(g_iz_j-g_jz_i)-\Pf_2(L_iz_j-L_jz_i)\\ 
&= \eta(g_iz_j-g_jz_i) \pmod{I_{Y^{\eta}}} 
\end{align*}
and so $(g_iz_j-g_jz_i) \in I_{Y^{\eta}}$ since $I_{Y^{\eta}}$ is prime. We proceed in finding $K_i$ and $L_i$. The $h_i$ can be readily computed and are
\begin{align*}                                                           
    h_1 &=\eta((F_3 H_4 - G_3) \eta  - ((F_1 H_2 + G_1 H_3) F_3 + F_1 G_3 H_3)) + (wF_1  - \zeta G_1)G_1 H_2 \\
     h_2 &=  \eta(F_1 \eta + (-F_1^2  H_3 + (F_3 H_1 - G_1 H_4) F_1 + G_1^2)) -  (wF_1 - \zeta G_1)G_1 H_1 \\
      h_3&=\eta((-F_1 H_4 + G_1) \eta + H_2 F_1^2  - F_1 G_3 H_1 + F_3 G_1 H_1)\\
         h_4&= \eta(F_1 G_1 H_2 + (F_1 H_3  - F_3 H_1 )\eta  - \eta^2)
\end{align*}

Just as before we define the ratio
\[
s:=\frac{g_i}{z_i} \in H^0\bigg(Y^u,-\frac{3a_w}{2}K_Y-\frac{1}{2}E\bigg)
\]
and let $a_s = 3a_w$. Let $T^s = \Proj \mathbb{C}[u,\xi,\zeta,w,x_0,x_1,\eta,s,x]$, i.e., this is $T^{\eta}$ with the extra section $s$ and 
\[
Y^s \colon (\Pf_j=0, sz_i-g_i=0\, | \, 1\leq j\leq 5,\, 1\leq i\leq 4)
\]
inside $T^s$. The inclusion map
\begin{align*}
i^{s} \colon Y^{\eta} & \longrightarrow Y^s\\
(u,\xi,\zeta,w,x_0,x_1,\eta,x)&\longmapsto (u,\xi,\zeta,w,x_0,x_1,\eta,s,x)
\end{align*}
is biregular. 

We now cross the wall generated by the multiples of $x_0,\,x_1,\,\eta$ on $T^{s}$. Indeed we have a diagram 
 \[
        \begin{tikzcd}[ampersand replacement=\&, column sep = 2em]
              T^s   \ar[rr, dashed ] \ar[dr, swap] \& \& {T^s}^+ \ar[dl ] \\
             \& \Proj \mathbb{C}[x_0,x_1,\eta] \& 
        \end{tikzcd}
    \] 
which restricts to the small $\mathbb{Q}$-factorial modification
\[
        \begin{tikzcd}[ampersand replacement=\&, column sep = 2em]
             C^{-} \subset Y^s   \ar[rr, dashed ] \ar[dr, swap] \& \& {Y^s}^+  \supset C^{+} \ar[dl ] \\
             \& \mathcal{F} \subset \Proj \mathbb{C}[x_0,x_1,\eta] \& 
        \end{tikzcd}
    \] 

where $\mathcal{F}$ is the restriction of $Y^s$ to $(u=\xi=\zeta=w=s=x=0)$. In each case $\mathcal{F}$ is a finite set of points, the curves contracted are isormorphic to $\mathbb{P}^1$ and
\[
K_{Y^s} \cdot C^{-}_{i} = K_{Y^s_{+}} \cdot C^{+}_{i} = 0
\]
Hence, crossing the $(x_0,x_1,\eta)$-wall induces a number of flops on $Y^s$.

\paragraph{Step 4: Final Contraction.}

The upshot of all this is that now we are finally able to carry a $K_{Y^s}$-negative divisorial contraction and terminate the link. After flopping the curves in $Y^s$, we consider the map 
\[
\varphi' \colon Y^{s} \rightarrow \mathcal{F}, \quad \mathcal{F}:=\Proj\bigoplus_{m\geq 1} H^0(Y^{s},-m\mathcal{O}(s))
\]
This map is generated by $s$ and multiples of $x$. Hence, it contracts the divisor $x=0$ in $Y^s$ to the point $\mathbf{p_s}$ and the 2-ray game on $Y^s$ follows the 2-ray game on $T^s$. In particular, we know the cone of movable divisors of $Y^s$ is
\[
\overline{\Mov}(Y^s) = \mathbb{R}_+\bigg[-\frac{a_{\xi}}{2}K_{Y^s} + \frac{1}{2}E\bigg] + \mathbb{R}_+\bigg[-\frac{a_{s}}{2}K_{Y^s} - \frac{1}{2}E\bigg] 
\]
Hence, $-K_{Y^s} \in \Int \overline{\Mov}(Y^s)$.  By construction (Kustin-Miller unprojection), the point $\mathbf{p_s} \in \mathcal{F}$ is a cyclic quotient singularity and $\varphi'$ is the Kawamata blowup centred at $\mathbf{p_s}$.
\end{proof}


\begin{Ex}
We give an example of Theorem \ref{thm:cod4} following its proof with a quasismooth member of family 111, $X_{18,20} \subset \mathbb{P}(2,5,6,7,9,11)$ with homogeneous coordinates $x,\,y,\,z,\,t,\,v,\,w$ given by 
\[
\begin{cases}
 tw+v^2+z^3 +x^9=0\\
t^2z+wv+y^4+z^2x^4=0.
\end{cases}
\]

Let $\varphi \colon (E\subset Y) \rightarrow (\mathbf{p_t} \in X)$ be the Kawamata blowup centred at $\mathbf{p_t}$. Then $E$ is an irreducible divisor in $Y$ isomorphic to $\mathbb{P}(1,1,6)$ with local variables $x,\,v,\,y$, respectively. Moreovoer, it is easy to see that $z \in H^0(Y,-3K_Y)$ and $w \in H^0\bigg(Y,\frac{-11K_Y+E}{2} \bigg)$. Hence, $Y$ is defined by
 \[
\begin{cases}
 tw+v^2+(z^3 +x^9)u=0\\
t^2z+wv+y^4u^3+z^2x^4u=0.
\end{cases}
\]
inside 
\[
\begin{array}{ccccc|cccccc}
             &       & u  & t & v & w & x & z & y & \\
\actL{T}   &  \lBr &  0 & 7 & 9 & 11 & 2 & 6 & 5 &    \actR{}\\
             &       & 1 & 4 & 5 &  6 & 1 & 3 & 2  & 
\end{array}
\] 
since crossing the $v$-wall induces an isomorphism on $Y$. Notice that $\widetilde{f},\,\widetilde{g} \in (u,t,v)$ and so there are $F_i, \,G_i \in \mathbb{C}[u,t,v,w,x,z,y]$ with $1\leq i\leq 3$ such that 
\[
\begin{pmatrix}
           \widetilde{f} \\
           \widetilde{g} \\
         \end{pmatrix} = 
\left(
  \begin{array}{ccc}
    F_1 & F_2 & F_3 \\
    G_1 & G_2 & G_3 \\
  \end{array}
\right)
\begin{pmatrix}
           u \\
           t \\
          v \\
         \end{pmatrix}.
\]
Explicitly for this example these are given by 
\[Q=
\left(
  \begin{array}{ccc}
    z^3+x^9 & w & v \\
    z^2x^4+y^4u^2 & tz & w \\
  \end{array}
\right).
\]

Let $z_i$ be the $i^{th}$ element of $\{u,t,v\}$. The ratio
\[
\eta:=\frac{(\bigwedge^2 Q)_i}{z_i} \in H^0(Y,-11K_Y),\quad 1\leq i\leq 3
\]
is a well-defined (pluri)-anticanonical section of $Y$.  Let 
\[
Y^{\eta} \colon (\widetilde{f}=\widetilde{g}=\eta z_i - (\bigwedge^2 Q)_i=0 \, |\, 1\leq i\leq 3)
\]
inside 
\[
\begin{array}{ccccc|cccccc}
             &       & u  & t & v & w & x & z & \eta & y & \\
\actL{T^{\eta}}   &  \lBr &  0 & 7 & 9 & 11 & 2 & 6 & 22 &5 &    \actR{.}\\
             &       & 1 & 4 & 5 &  6 & 1 & 3 & 11 &2  & 
\end{array}
\] 
Then the map
\begin{align*}
i^{\eta} \colon Y & \longrightarrow Y^{\eta}\\
(u,t,v,w,x,z,y)&\longmapsto (u,t,v,w,x,z,\eta,y)
\end{align*}
 is biregular.  Moreover, the ideal defining $Y^{\eta}$, $I_{Y^{\eta}}$, can be written in terms of the 5 maximal Pfaffians of the $5\times 5$ antisymmetric matrix 
\[
M=
\left(
  \begin{array}{cccc}
    u & -t & w & v \\
		 & v & tz & w \\
		 &  & z^2x^4+y^4u^2 & z^3+x^9 \\
		 &  &  & \eta \\
  \end{array}
\right)
\]
that we write as $I_{Y^{\eta}}=(\Pf_1,\ldots,\Pf_5)$. Notice that the first two rows of $M$ are polynomials in $J=(u,t,v,w)$ and so $I_{Y^{\eta}}\subset J$ and $Y^{\eta}$ is a codimension $3$ Fano 3-fold written in Jerry$_{12}$ format.  We have,
\[
\begin{pmatrix}
           \Pf_1 \\
           \Pf_2 \\
					 \Pf_3 \\
         \end{pmatrix} = 
\left(
  \begin{array}{cccc}
    0 & -z(z^3+x^9)  & \eta & z^2x^4+y^4u^2 \\
    0 & -\eta  & z^2x^4+y^4u^2 & -(z^3+x^9)\\
		\eta & zv & 0 & -w\\
  \end{array}
\right)
\begin{pmatrix}
           u \\
           t \\
          v \\
					w \\
         \end{pmatrix}.
\]

From this matrix we produce polynomials $h_i$ which are given by
\[
 (h_1,\ldots, h_4) := \bigwedge^3 Q. 
\]

For each $1 \leq i\leq 4$ there are polynomials $K_i$ and $L_i$ such that
\[
h_i=\eta K_i + (wx_2-vx_1)L_i = \eta (K_i-tL_i)-L_i \Pf_2
\]
where $x_1=z^3+x^9$ and $x_2=z^2x^4+y^4u^2$.

Let $z_i$ be the $i^{th}$ element in $\{u,t,v,w\}$ and $g_i=K_i-tL_i$ for each $1 \leq i\leq 4$. We proceed in finding $K_i$ and $L_i$. Suppose that $2 \leq i\leq 4$. Then, in our example, it is obvious that $h_i$ is a multiple of $\eta$. Hence we take $L_i$ to be the zero polynomial and $K_i = \frac{h_i}{\eta}$. Otherwise,
\[
h_1=\eta(-wx^{11}+vz^4+w\eta+vx^9z+u^2vx^2y^4+vx^6z^2-wx^2z^3)-zx_1(wx_2-vx_1).
\]
Just as before we define the ratio
\[
s:=\frac{g_i}{z_i} \in H^0\bigg(Y^u,\frac{-33K_Y-E}{2}\bigg).
\]
 Let
\[
Y^s \colon (\Pf_j=0, sz_i-g_i=0\, | \, 1\leq j\leq 5,\, 1\leq i\leq 4)
\]
inside 
\[
\begin{array}{cccccc|cccccc}
             &       & u  & t & v & w & x & z & \eta & s &y & \\
\actL{T^s}   &  \lBr &  0 & 7 & 9 & 11 & 2 & 6 & 22 &33 &5 &    \actR{.}\\
             &       & 1 & 4 & 5 &  6 & 1 & 3 & 11 & 16&2  & 
\end{array}
\] 
The map
\begin{align*}
i^s \colon Y^u & \longrightarrow Y^s\\
(u,t,v,w,x,z,y,\eta)&\longmapsto (u,x,v,y,z,w,t,\eta,s)
\end{align*}
is biregular. 

We now cross the wall generated by the multiples of $|-K_Y|$. Indeed we have a diagram
 \[
        \begin{tikzcd}[ampersand replacement=\&, column sep = 2em]
             C^{-} \subset Y^s_{-}   \ar[rr, dashed ] \ar[dr, swap] \& \& Y^s_{+}  \supset C^{+} \ar[dl ] \\
             \& \mathcal{F} \in Z \& 
        \end{tikzcd}
    \] 
		
where $\mathcal{F}$ is a cluster of 10 points. Over each of these points $p_i$ there are curves $C^{-}_{i}$ isormorphic to $\mathbb{P}^1$. Moreover,
\[
-K_{Y^s_{-}} \cdot C^{-}_{i} = -K_{Y^s_{+}} \cdot C^{+}_{i} = 0
\]
Hence, crossing the $(x,z,\eta)$-wall induces 10 flops on $Y^s$.

After flopping the curves in $Y^s$ and multiplying by a suitable matrix we are in 
\[
\begin{array}{ccccccccc|ccc}
             &       & u  & t & v & w & x & z & \eta & s &y & \\
\actL{T^s}   &  \lBr &  33 & 20 & 21 & 22 & 1 & 3 & 11 & 0 & -14 &    \actR{.}\\
             &       & 5 & 6 & 7 &  8 & 1 & 3 & 11 & 14&0  & 
\end{array}
\] 
Consider the map given by $\bigoplus_{m\geq 1} H^0(Y^s,m\mathcal{O}(s))$, that is,
\[
(uy^{33/14}:ty^{10/7}:vy^{2}:wy^{22/14}:xy^{1/14}:zy^{3/14}:\eta y^{11/14}:s).
\]
Clearly $(y=0) \subset Y^s$ is contracted to $\mathbf{p_s} \in X'$ and away from this divisor we have an isomorphism so this is a divisorial contraction. More than that, it is an extremal divisorial contraction in the Mori Category. The point $\mathbf{p_s}$ is a germ of a cyclic quotient singularity of type $\frac{1}{14}(1,3,11)$ and local variables $(x,z,\eta)$. 
\end{Ex}
\subsection{Fake Weighted Projective Spaces and higher codimension} \label{subsec:fake}


\begin{Def}[{\cite{Alfake}}]
A \textbf{fake weighted projective space} is a $\mathbb{Q}$-factorial toric variety with Picard number 1.
\end{Def}

See also \cite[Def~6.1]{fake}. By \cite[Exercise~5.1.13]{Coxbook} if $X$ is a fake weighted projective space, then there are positive integers $a_0,\ldots,a_n$ and a finite abelian group $H$ acting on $\mathbb{P}(a_0,\ldots,a_n)$ for which $X \simeq \mathbb{P}(a_0,\ldots,a_n)/H$.   

\begin{Lem} \label{lem:fakewps}
Let $T$ be a well-formed rank two toric variety of the form
\[
\begin{array}{ccccc|cccc}
             &       & v_1  &   & v_n & y & x & \\
\actL{T}   &  \lBr &  a_1 & \cdots & a_n & 0 & -1 &   \actR{}\\
             &       & b_1 & \cdots & b_n & 1 & 0 &  
\end{array}
\]
and suppose that the group of $r$th roots of unity $\boldsymbol{\mu_r}$ acts on $T$ via 
\[
(v_1,\ldots,v_n,y,x) \mapsto (v_1,\ldots,v_n,y,\epsilon x) 
\]
where $\epsilon \in \boldsymbol{\mu_r}$ is a primitive element. Then, the morphism $\Phi \colon T \rightarrow \mathbb{P}(b_1,\ldots,b_n,1)$ given by the sections of multiples of $D_y$ descends to a morphism of orbit spaces 
\[
\Phi' \colon T/\boldsymbol{\mu_r} \rightarrow \mathbb{P}(b_0,\cdots,b_n,1)/\boldsymbol{\mu_r}.
\] 
\end{Lem}

\begin{proof}
The map $\Phi$ is 
\[
\Phi \colon T \rightarrow \Proj \mathcal{F}, \quad \mathcal{F}:= \bigoplus_{m \geq 1} H^0(T, m\mathcal{O}_T(D_y)) 
\]
realised by $(v_1,\ldots,v_n,y,x) \mapsto (v_1x^{a_1}: \ldots : v_nx^{a_n} : y)$. This contracts the divisor $(x=0)$ to the point $\mathbf{p_y} \in \mathbb{P}(b_0,\ldots,b_n,1_y)$ and is an isomorphism otherwise. On the other hand, the action of $\boldsymbol{\mu_r}$ on $T$ fixes all the points in $T$ for which $x=0$ and, if $x\not = 0$, it has an orbit with $r$ elements. Hence, there is a natural quotient map $q \colon T \rightarrow T/\boldsymbol{\mu_r}$ whose fibre over each point is a $\boldsymbol{\mu_r}$-orbit. The quotient $T/\boldsymbol{\mu_r}$ is a rank 2 toric variety $T'$ that can be realised as 
\[
\begin{array}{ccccc|cccc}
             &       & v_1  &   & v_n & y & x' & \\
\actL{T'}   &  \lBr &  a_1 & \cdots & a_n & 0 & -r &   \actR{}\\
             &       & b_1 & \cdots & b_n & 1 & 0 &  
\end{array}
\]
Since $\Phi$ is an isomorphism away from $x=0$, the action is carried to  $\mathbb{P}(b_0,\ldots,b_n,1)$, namely defined by
\[
\epsilon\cdot (u_0:\ldots:u_n:y) \mapsto (\epsilon^{a_0}u_0:\ldots:\epsilon^{a_0}u_0:y).
\]
In the same way there is a quotient map $q \colon \mathbb{P}(b_0,\ldots,b_n,1)\rightarrow \mathbb{P}(b_0,\ldots,b_n,1)/\boldsymbol{\mu_r}$ and the quotient $\mathbb{P}(b_0,\ldots,b_n,1)/\boldsymbol{\mu_r}$ can be realised as 
\[
\Proj \mathbb{C}[\ldots, u_{i_0}^{\alpha_{i_0}}\cdots u_{i_k}^{\alpha_{i_k}}, \ldots, y]=\Proj \mathbb{C}[\ldots, v_{i_k}, \ldots, y]/I
\]
where, for every $i_j$, $\alpha_{i_0}a_{i_0}+\cdots+\alpha_{i_k}a_{i_k}$ is a multiple of $r$ and $I$ is the ideal of relations of the generators $u_{i}'$. On the other hand, the morphism $\overline{\varphi} \colon T/\boldsymbol{\mu_r} \rightarrow \mathbb{P}(b_0,\ldots,b_n,1)/\boldsymbol{\mu_r}$ given by the sections of multiples of the divisor $D_y$ is
\[
(v_1x^{a_1}:\ldots :v_nx^{a_n}:y) \in \mathbb{P}(b_0,\ldots,b_n,1)
\]
where $x^r$ is defined as $x'$ in $\Cox(T')$. On $x \not = 0$, we have the $\boldsymbol{\mu_r}$-action $v_i \mapsto \epsilon ^{a_i}v_i$. We conclude that the diagram 
\[
        \begin{tikzcd}[ampersand replacement=\&,column sep = 3.5em,row sep = 3.5em] \label{diag:equiv}
             T   \ar[r, "\displaystyle{q }" ] \ar[d, swap, "\displaystyle{\Phi}"] \& T'=T/\boldsymbol{\mu_2} \ar[d, "\displaystyle{\Phi'}"] \\
           \mathbb{P} \ar[r,swap, "\displaystyle{q }" ]  \&  \mathbb{P}/\boldsymbol{\mu_2}
        \end{tikzcd}
		\]
		commutes and $\varphi$ descends to a morphism of the corresponding orbit spaces.
\end{proof}

\begin{Ex} \label{ex:faketoricvar}
Consider the quasismooth rank two toric variety
\[
\begin{array}{ccccccc|cccc}
             &       & u  & v &   x & y & w & t & z & \\
\actL{T'}   &  \lBr &  5 & 4 &  1 & 1& 1 & 0 & -2 &   \actR{.}\\
             &       & 2 & 3 & 1 & 1 & 2 & 1 & 0 &  
\end{array}
\]
The map given by the sections $\Phi' \colon \bigoplus_{m \geq 1} H^0(T',m\mathcal{O}_{T'}(t=0))$ is 
\[
(u^2z^{5}:uxz^3:uyz^3:uwz^3:vz^2:x^2z:xyz:xwz:y^2z:ywz:w^2z:t)  
\]
inside $\mathbb{P}(4,3,3,4,3,2,2,3,2,3,4,1):=\mathbb{P}'$. Notice that, when $z=1$, this is precisely 
\[
\Phi' \colon T' \rightarrow \Proj \mathcal{F}^{\boldsymbol{\mu_2}}, \quad \mathcal{F}:= \mathbb{C}[u,v,x,y,w,t] 
\]
with $\wt(u,v,x,y,w,t)=(2,3,1,1,2,1)$ and $\boldsymbol{\mu_2} \times \mathcal{F} \rightarrow \mathcal{F}$ an action given by 
\[
u \mapsto \epsilon u,\, v \mapsto v,\, x \mapsto \epsilon x,\, y \mapsto \epsilon y,\, w \mapsto \epsilon w,\, t \mapsto  t.
\]
We find generators and relations of $\mathbb{P}':=\Proj \mathcal{F}^{\boldsymbol{\mu_2}}=\mathbb{P}(2,3,1,1,2,1)/\boldsymbol{\mu_2}$. From the contruction above,
\[
\mathcal{F}^{\boldsymbol{\mu_2}} = \mathbb{C}[u^2z^{5},uxz^3,uyz^3,uwz^3,x^2z,xyz,xwz,y^2z,ywz,w^2z,vz^2,t]. 
\]
Call the 12 generators $u_1,\ldots u_{12}$. These have weights 4, 3, 3, 4, 2, 2, 3, 2, 3, 4, 3, 1, respectively.
It is easy to see that the relations between the $u_i$ are generated by $\rk A = 1$, where $A$ is the symmetric $4\times 4$ matrix:
\[
\begin{pmatrix}
	u^2z^5&uxz^3&uyz^3&uwz^3\\
	&x^2z&xyz&xwz\\
	&&y^2z&ywz\\
	&&&w^2z\\
\end{pmatrix}=:
\begin{pmatrix}
	u_1&u_2&u_3&u_4\\
	&u_5&u_6&u_7\\
	&&u_8&u_9\\
  &&&u_{10}\\
\end{pmatrix}
\]

Hence
\[
\mathbb{P}(2,3,1,1,2,1)/\boldsymbol{\mu_2}=\Proj \mathcal{F}^{\boldsymbol{\mu_2}} = (\rk A = 1 ) \subset \mathbb{P}(4, 3, 3, 4, 2, 2, 3, 2, 3, 4, 3, 1).
\]
\end{Ex}

The next proposition gives examples where a similar action - that is, an action of a finite cyclic group on some toric variety - arises naturally.

\begin{Prop} \label{prop:fake1}
Let $X$ be a quasismooth member of the deformation families 116, 117 or 125. Then, $X$ is birational to a Fano 3-fold $X'$ embedded in a fake weighted projective space as in the second column of table \ref{tab:fake}. Moreover, there is a birational map $X \dashrightarrow X'$ which decomposes as 
\[
        \begin{tikzcd}[ampersand replacement=\&, column sep = 2em]
             Y   \ar[rr, dashed, "\displaystyle{ \tau}" ] \ar[d, swap, "\displaystyle{\varphi}"] \& \& Y' \ar[d,, "\displaystyle{\varphi'}" ] \\
           \mathbf{p} \in X  \& \& \mathbf{p'} \in X' 
        \end{tikzcd}
    \]
where $\varphi$ and $\varphi'$ are the Kawamata blowups of $\mathbf{p}$ and $\mathbf{p'}\sim \frac{1}{2}(1,1,1)$, respectively and $\tau$ is a flip over a point.
\end{Prop}

\begin{proof}
We can write $X$ as 
\begin{align*}
\xi x_3 +x_4x_2 + f_{d_1}(x_0,x_1)+f_{d_1}'(x_0,x_1,x_4,x_3,x_2)=0\\
\xi x_4 +x_3^{\alpha}+g_{d_2}(x_0,x_1)+g_{d_2}'(x_0,x_1,x_4,x_3,x_2)=0
\end{align*}
inside $\mathbb{P}(a_{x_0},a_{x_1},a_{x_2},a_{x_3},a_{x_4},a_{\xi})$ where $\wt(x_{\mu}) = a_{\mu}$.  Then, $\mathbf{p_{\xi}} \sim \frac{1}{a_{\xi}}(1,1,a_{\xi}-1)$ with local coordinates $(x_0,x_1,x_2)$. Let $E \colon (u=0)\cong \mathbb{P}(1,1,a_{\xi}-1)$ be the exceptional divisor under the Kawamata blowup $\varphi \colon Y \rightarrow X$ centred at $\mathbf{p_{\xi}}$. The 3-fold $Y$ can be written as
\begin{align*}
\widetilde{f} \colon \xi x_3 +x_4x_2u + f_{d_1}(x_0,x_1)+uf_{d_1}'(u,x_0,x_1,x_4,x_3,x_2)=0\\
\widetilde{f} \colon \xi x_4 +x_3^{\alpha}u+g_{d_2}(x_0,x_1)+ug_{d_2}'(u,x_0,x_1,x_4,x_3,x_2)=0
\end{align*}
where $\widetilde{f},\, \widetilde{g} \in |-\frac{d_i}{\iota_X}K_Y|$. See Lemma \ref{lem:lift}. Since $\mathbf{p_{\xi}}$ is linear and $a_4(\iota_X-1)>a_3(\iota_X-1) > a_2$ we are in Case I of Subsection \ref{sub:divtopoint}. The difference between the 3-folds treated in Section \ref{subsec:cod2} and the present ones is that in the former, the ambient toric variety in Lemma \ref{lem:toricendcod2} 
\[
T'' \colon \begin{pmatrix}
a_3	& \frac{d_1}{\iota_X} & 1 & \frac{a_1}{\iota_X}& \frac{a_4-a_3}{\iota_X}&0&-\frac{a_3(\iota_X-1)-a_2}{\iota_X}\\
a_2&a_{\xi}-\frac{a_1}{\iota_X}&\iota_X-1&(\iota_X-1)\frac{a_1}{\iota_X}&\frac{a_4(\iota_X-1)-a_2}{\iota_X}&\frac{a_3(\iota_X-1)-a_2}{\iota_X}&0
\end{pmatrix} 
\]
with irrelavant ideal $(u,\xi,x_0,x_1,x_4)\cap (x_3,x_2)$ satisfies $a_ib_j - a_jb_i \equiv 0 \bmod r$ for all $1\leq i,\,j\leq 6$, where no row consists of multiples of $r$. In this case, however, since $\iota_X \geq 3$ and $r:=\frac{a_3(\iota_X-1)-a_2}{\iota_X}=2$, the second row is a multiple of $r=2$. Define the matrix $T'$ as 
\[
T''=
\begin{pmatrix}
1&0\\ 
0	& r 
\end{pmatrix} T'. 
\] 
Then, $T''$ and $T'$ are isomorphic as quotient varieties by \cite[Lemma~2.4]{hamidplia} and $T'$ is well-formed. Let $s$ be a section such that $s^r =x_2$ and $T_s\rightarrow T'$ the corresponding double cover. We define $Y_s \subset T_s$ as 
\begin{align*}
\widetilde{f} \colon \xi x_3 +x_4s^ru + f_{d_1}(x_0,x_1)+uf_{d_1}'(u,x_0,x_1,x_4,x_3,s^2)=0\\
\widetilde{g} \colon \xi x_4 +x_3^{\alpha}u+g_{d_2}(x_0,x_1)+ug_{d_2}'(u,x_0,x_1,x_4,x_3,s^2)=0
\end{align*}
Hence $Y_s$ has an automorphism of order $r$ which is the restriction of the action of $\boldsymbol{\mu_r}$ on $T_s$ to $Y_s$. The map 
\[
\Phi \colon T_s \rightarrow \Proj \mathcal{F}=\mathbb{P}\bigg(\frac{a_2}{r},\frac{a_{\xi}-\frac{a_1}{\iota_X}}{r},\frac{\iota_X-1}{r},(\iota_X-1)\frac{a_1}{r\iota_X},\frac{a_4(\iota_X-1)-a_2}{r\iota_X},1\bigg):=\mathbb{P}.
\]
is given in coordinates by 
\[
(u,\xi,x_0,x_1,x_4,x_3,s)\mapsto(us^{a_3}:\xi s^{d_1/\iota_X}:x_0s:x_1s^{a_1/\iota_X}:x_4s^{(a_4-a_3)/\iota_X}:x_3).
\]
Notice that all powers of $s$ are positive integers. Its restriction to $Y_s$ is a 3-fold divisorial contraction which contracts $s=0$ to a point and away from this divisor, we have an isomorphism to $Z_{d_1,d_2}$
\begin{align*}
f \colon \xi x_3 +x_4u + f_{d_1}(x_0,x_1)+uf_{d_1}'(u,x_0,x_1,x_4,x_3,1)=0\\
g \colon \xi x_4 +x_3^{\alpha}u+g_{d_2}(x_0,x_1)+ug_{d_2}'(u,x_0,x_1,x_4,x_3,1)=0
\end{align*}
inside $\mathbb{P}$, where $d_1 = \frac{(a_{\xi}+a_3)(\iota_X-1)}{r\iota_X}$ and $d_2 = \frac{(a_{\xi}+a_4)(\iota_X-1)}{r\iota_X}$. 
From Lemma \ref{lem:fakewps}, $\Phi$ descends to a divisorial contraction on the corresponding orbit spaces 
\[
\Phi' \colon T/\boldsymbol{\mu_2} \rightarrow \Proj \mathcal{F}^{\boldsymbol{\mu_2}} = \mathbb{P}/\boldsymbol{\mu_2}
\]
which is $(u,\xi,x_0,x_1,x_4,x_3,s)\mapsto(us^{a_3/r}:\xi s^{d_1/r\iota_X}:x_0s^{1/r}:x_1s^{a_1/r\iota_X}:x_4s^{(a_4-a_3)/r\iota_X}:x_3)$ in coordinates. We have already noticed that $Y\subset T$ is invariant under the action of $\boldsymbol{\mu_2}$ on $T$, so it makes sense to restrict $\Phi'$ to $Y$. Let $\varphi'$ be this restriction. This is a divisorial contraction to the point $\mathbf{p_{x_3}} \in Z_{d_1,d_2} \subset \mathbb{P}/\boldsymbol{\mu_2}$ where the action of $\boldsymbol{\mu_2}$ is carried through to $\mathbb{P}$ and is given by
\[
(u,\xi,x_0,x_1,x_4,x_3)\mapsto( \epsilon^{a_3} u:\epsilon^{d_1/\iota_X}\xi:\epsilon x_0:\epsilon^{a_1/\iota_X}x_1:\epsilon^{(a_4-a_3)/\iota_X}x_4:x_3)\\
\]
Locally, around $\mathbf{p_{x_3}}$, the map $\varphi'$ is a weighted blowup of weights
\[
\wt(x_0,x_1,x_4) = \frac{1}{2}\bigg(1,\frac{a_1}{\iota_X}, \frac{a_4-a_3}{\iota_X}\bigg)=\frac{1}{2}(1,1,1)
\]
and discrepancy $\frac{1}{2}$ by Lemma \ref{lem:discr}. Hence $\varphi'$ is the Kawamata blowup of $\mathbf{p_{x_3}}$.
 \end{proof}

\begin{Ex}
Let $T'$ be as in Example \ref{ex:faketoricvar}
Let $s$ be a section such that $s^2 =z$ and $T\rightarrow T'$ the corresponding double cover. We define $Y \subset T$ as 
\begin{align*}
vt+ws^2u+a_4(x,y) &=0\\
vw+t^3u+b_5(x,y)+s^4vu&=0
\end{align*}
where $a_4,\,b_5$ are general polynomials in degrees 4 and 5, respectively. Hence, $Y$ has an automorphism of order two which is the restriction of the action of $\boldsymbol{\mu_2}$ on $T$ to $Y$. The map 
\[
\Phi \colon T \rightarrow \Proj \mathcal{F}=\mathbb{P}(2,3,1,1,2,1)
\]
is given in coordinates by $(u,v,x,y,w,t,z)\mapsto(us^5:vs^4:xs:ys:ws:t)$. Its restriction to $Y$ is a 3-fold divisorial contraction which contracts $s=0$ to a point and away from this divisor, we have an isomorphism to $Z_{4,5}$
\begin{align*}
vt+wu+a_4(x,y) &=0\\
vw+t^3u+b_5(x,y)+vu&=0
\end{align*}
inside $\mathbb{P}(2_u,3_v,1_x,1_y,2_w,1_t)$. Then, $Z_{4,5}$ has exactly two singular points which are $(w^2-t^4=0) \subset \mathbb{P}(2,1)$. From Lemma \ref{lem:fakewps}, $\varphi$ descends to a divisorial contraction on the corresponding orbit spaces 
\[
\Phi' \colon T/\boldsymbol{\mu_2} \rightarrow \Proj \mathcal{F}^{\boldsymbol{\mu_2}} = \mathbb{P}(2,3,1,1,2,1)/\boldsymbol{\mu_2}
\]
which is $(u,v,x,y,w,t,z)\mapsto(us^{5/2}:vs^2:xs^{1/2}:ys^{1/2}:ws^{1/2}:t)$ in coordinates. We have already noticed that $Y$ is invariant under this action, so it makes sense to restrict $\Phi'$ to $Y$. This is a divisorial contraction to a point in $Z_{4,5}$ inside $\mathbb{P}(2,3,1,1,2,1)/\boldsymbol{\mu_2}$ where $\boldsymbol{\mu_2}$ acts on $\mathbb{P}(2,3,1,1,2,1)$ by
\[
(u_0:v_0:x_0:y_0:w_0:t_0) \mapsto  (-u_0:v_0:-x_0:-y_0:-w_0:t_0)= (-u_0:-v_0:x_0:y_0:-w_0:-t_0).
\]
Hence, the two singular points of $Z_{4,5}$ are in the same orbit since 
\[
(0:0:0:0:t_0^2: t_0) \mapsto (0:0:0:0:-t_0^2:t_0).
\] 
\end{Ex}

The following proposition is the analogous version of Proposition \ref{prop:fake1} for families 105, 123 and 125. The Sarkisov link in this case is initiated by blowing a non-linear cyclic quotient singularity.
For the purposes of the next proposition we have the following assumptions, which we refer to Table \ref{tab:big}:
\begin{itemize}
	\item If $X$ is a member of family 105 we assume $\alpha = 0,\, (\iff tv \not \in f)$ and $y^6 \in f$.
	\item If $X$ is a member of family 123 we assume $\beta \not = 0,\, (\iff y^3t \in f,\, y^4x \not \in f)$.
	\item If $X$ is a member of family 125 we assume $\beta = \gamma = 0,\, (\iff x^5(t+z),\, x^4v \not \in g)$.
\end{itemize}

\begin{Prop} \label{prop:fake2}
Let $X$ be a quasismooth member of the deformation families 105, 123 or 125 as in the previous paragraph. Then, $X$ is birational to a Fano 3-fold $X'$ embedded in a fake weighted projective space as in the second column of Table \ref{tab:fake}. Moreover, there is a birational map $X \dashrightarrow X'$ which decomposes as 
\[
        \begin{tikzcd}[ampersand replacement=\&, column sep = 2em]
             Y   \ar[rr, dashed, "\displaystyle{ \tau}" ] \ar[d, swap, "\displaystyle{\varphi}"] \& \& Y' \ar[d,, "\displaystyle{\varphi'}" ] \\
           \mathbf{p} \in X  \& \& \mathbf{p'} \in X' 
        \end{tikzcd}
    \]
where $\varphi$ is the Kawamata blowup of $\mathbf{p}$ and $\varphi'$ is the weighted blowup of $\mathbf{p'}$ and $\tau$ is a small $\mathbb{Q}$-factorial modification over a point.
\end{Prop}

\begin{center}
\begin{table}
\resizebox{\textwidth}{!}{\begin{tabular}{llll} \toprule
    {ID} & $\mathbf{p} \in X \subset w\mathbb{P}$ &$X\subset \mathbb{P}(b_1,\ldots,b_5,b_6)/\boldsymbol{\mu_r}$ &  Action $\boldsymbol{\mu_r} \times \mathbb{P} \rightarrow \mathbb{P}$  \\
			\midrule 
      105  & $\frac{1}{5}(1,1,4) \in X_{12,14}$ &  $Z_{6,8} \subset \mathbb{P}(1,3,4,5,1,1)/\boldsymbol{\mu_2}$ & $(\epsilon v_1: v_2:\epsilon v_3: v_4:\epsilon v_5: y),\, \epsilon \in \boldsymbol{\mu_2}$ \\ 
			116  & $\frac{1}{5}(1,1,4) \in X_{9,12}$ &  $Z_{3,4} \subset \mathbb{P}(1,2,1,1,2,1)/\boldsymbol{\mu_2}$ & $(\epsilon v_1: \epsilon v_2:v_3:v_4:\epsilon v_5:\epsilon y),\, \epsilon \in \boldsymbol{\mu_2}$ \\ 
		  117  & $\frac{1}{7}(1,1,6) \in X_{12,15}$ &  $Z_{4,5} \subset \mathbb{P}(2,3,1,1,2,1)/\boldsymbol{\mu_2}$ & $(\epsilon v_1: v_2:\epsilon v_3:\epsilon v_4:\epsilon v_5: y),\, \epsilon \in \boldsymbol{\mu_2}$ \\ 
		  123  & $\frac{1}{3}(1,1,2) \in X_{12,14}$ &  $Z_{3,4} \subset \mathbb{P}(1,1,2,1,2,1)/\boldsymbol{\mu_2}$ & $( \epsilon v_1: v_2:\epsilon v_3:\epsilon v_4: \epsilon v_5:\epsilon y),\, \epsilon \in \boldsymbol{\mu_2}$ \\ 
		  125  & $\frac{1}{7}(1,1,6) \in X_{10,15}$ &  $Z_{4,6} \subset \mathbb{P}(1,3,2,2,3,1)/\boldsymbol{\mu_2}$ & $(\epsilon v_1: v_2:\epsilon v_3:\epsilon v_4:\epsilon v_5: y),\, \epsilon \in \boldsymbol{\mu_2}$ \\ 
		  125  & $\frac{1}{2}(1,1,1)\in X_{10,15}$ &  $Z_{4,6} \subset \mathbb{P}(1,1,3,2,2,3)/\boldsymbol{\mu_3}$ & $(\epsilon^8 v_1:\epsilon^2 v_2: \epsilon^3 v_3:\epsilon v_4: \epsilon v_5: y),\, \epsilon \in \boldsymbol{\mu_3}^*$ \\ 
		\bottomrule
		\end{tabular}} 
		\caption{Models in Fake Weighted Projective Spaces}
		\label{tab:fake}
		\end{table}
\end{center}

\begin{proof}
We prove each case separately for clarity since there are a few noticeable differences.
\paragraph{Family 105.} With the assumptions declared before the proposition, $X$ can be written as 
\begin{align*}
\xi^2x+\xi zf_2(x,y)+wz+z^4+f_6(x,y)+f_{12}&=0\\
\xi w+v^2+g_7(x,y)+g_{14}&=0
\end{align*}
inside $\mathbb{P}:=\mathbb{P}(2,2,3,5,7,9)$ with homogeneous variables $x,\,y,\,z,\,\xi,\,v,w$. We have $f_{12},\,g_{14} \in \mathbb{C}[x,y,z,v]$ without loss of generality. As usual, we consider the blowup $\Phi \colon T \rightarrow \mathbb{P}$ whose restriction to $X$ is the unique Kawamata blowup $\varphi \colon E \subset Y \rightarrow X$ centred at $\mathbf{p_{\xi}} \sim \frac{1}{5}(1,1,4)$ with variables $(y,v,z)$. Since $\iota_X = 2$, we have by Corollary \ref{cor:genlift}, 
\begin{align*}
\xi & \in H^0 \bigg(Y,-\frac{5}{2}K_Y+ \frac{1}{2}E \bigg), \\
v   & \in H^0 \bigg(Y,-\frac{7}{2}K_Y+ \frac{1}{2}E \bigg), \\
w   & \in H^0 \bigg(Y,-\frac{9}{2}K_Y+ \frac{1}{2}E \bigg), \\
y &\in H^0(Y,-K_Y), \\
z & \in H^0 \bigg(Y,-\frac{3}{2}K_Y- \frac{1}{2}E \bigg), \\
x & \in H^0 (Y,-K_Y- E ).
\end{align*}
We also have $\widetilde{f} \in |-6K_Y|$ and $\widetilde{g} \in |-7K_Y+E|$. These sections, from $u$ to $x$ on $T$ generate the cone $\Eff(T)$ and its Mori chamber decomposition. We also have that $\Mov(T)$ is the cone
\[
\mathbb{R}_+[\xi] + \mathbb{R}_+[z]. 
\]
We write down the contractions between toric varieties in adjacent chambers of $\Mov(T)$. The toric variety $T$ is 
\[
\begin{array}{cccc|ccccccc}
             &       & u  & \xi &  v & w & y & z & x & \\
\actL{T}   &  \lBr &  0 & 5 &  7 & 9& 2 & 3 & 2 &   \actR{.}\\
             &       & 1 & 3 & 4 & 5 & 1 & 1 & 0 &  
\end{array}
\]
and we can vary the GIT quotients by successively multiplying by convenient invertible matrices with integer entries. The contractions between the different GIT quotients are 
\[
        \begin{tikzcd}[ampersand replacement=\&,  column sep={6em,between origins}]
             T   \ar[rr, dashed, "\displaystyle{(-7,-1,1,1,5,10)}" ] \ar[d, swap, "\displaystyle{\Phi}"] \& \& T_1 \ar[rr,dashed, "\displaystyle{(-9,-2,-1,1,6,10)}" ]  \& \& T_2 \ar[rr,dashed, "\displaystyle{(-2,-1,-1,-1,1,2)}" ]\& \& T' \ar[d, "\displaystyle{\Phi'}"]\\
           \mathbf{p} \in \mathbb{P}  \& \&  \& \& \& \&  \mathbf{p'} \in \mathbb{P}'
        \end{tikzcd}
    \]
where each small contraction happens over a point. The toric variety $T'$ is 
\[
\begin{array}{ccccccc|cccc}
             &       & u  & \xi &  v & w & y & z & x & \\
\actL{T'}   &  \lBr &  3 & 4 &  5 & 6 & 1 & 0 & -2 &   \actR{.}\\
             &       & 1 & 3 & 4 & 5 & 1 & 1 & 0 &  
\end{array}
\]
We have a double cover $T_s \rightarrow T'$ given by $s^2=x$. Suppose $\boldsymbol{\mu_2}$ acts on $T_s$ by
\[
(u,\xi,v,w,y,z,s) \mapsto (u,\xi,v,w,y,z,-s).
\]
Then, $T'$ is the orbit space of $T_s$ under this action. On the other hand, the divisorial contraction
\[
\Phi_s \colon T_s \rightarrow \Proj \mathcal{F}, \quad \mathcal{F}=\bigoplus_{m\geq 1}H^0(T_s,m\mathcal{O}(z=0))
\]
is given in coordinates by $(u,\xi,v,w,y,z,s)\mapsto(us^3:\xi s^4:vs^5:ws^6:ys:z) \in \mathbb{P}(1,3,4,5,1,1)$. By Lemma \ref{lem:fakewps}, $\Phi$ descends to a map on the orbits $\Phi'\colon T'=T_s/\boldsymbol{\mu_2} \rightarrow \mathbb{P}(1,3,4,5,1,1)/\boldsymbol{\mu_2}$ where the action on $\mathbb{P}(1,3,4,5,1,1)$ is carried over from $T_s$ by $\Phi_s$. That is, the action on $\mathbb{P}(1,3,4,5,1,1)$ is $(u,\xi,v,w,y,z)\mapsto(-u:\xi:-v:w:-y:z)$.

We restrict this link to $Y \subset T$. Since $v^2 \in H^0(Y,-7K_Y+E)$ it follows that $v^2 \in \widetilde{g}$. Hence the $v$-wall crossing on $T$ restricts to an isomorphism. The $w$-wall crossing restricts to a toric anti-flip since over $\mathbf{p_w}$: This is because $wz \in \widetilde{f}$ and $tw \in \widetilde{f}$. Hence, locally around $\mathbf{p_w}$ we can eliminate $z$ and $t$. This is exactly the anti-flip $(-9,-1,1,10)$. Similarly to the $v$-wall crossing, the $y$-wall crossing restricts to an isomorphism since $y^6 \in \widetilde{f}$. 
We now define $Y_s \subset T_s$ to be $Y'|(s^2-x=0)$, where $Y' \subset T'$ is the resulting 3-fold after the small contraction over $\mathbf{p_w}$. Its equations are
\begin{align*}
\xi^2s^2+tzf_2(us^2,y)+wz+u^2z^4+f_6(us^2,y)+uf_{12}(\ldots,s^2,\ldots)&=0\\
\xi w+v^2+ug_7(us^2,y)+ug_{14}(\ldots,s^2,\ldots)&=0
\end{align*}
By construction, $Y_s$ is invariant under the action of $\boldsymbol{\mu_2}$ on $T_s$. On the other hand, the map $\varphi_s := \Phi_s|_{Y_s} \colon Y_s \rightarrow Z_{6,8}$ is a divisorial contraction to a point where $Z_{6,8}$ is
\begin{align*}
\xi^2+tzf_2(u,y)+wz+u^2z^4+f_6(u,y)+uf_{12}(\ldots,1,\ldots)&=0\\
\xi w+v^2+ug_7(u,y)+ug_{14}(\ldots,1,\ldots)&=0
\end{align*}
in $\mathbb{P}(1,3,4,5,1,1)$ with homogeneous variables $u,\,\xi,\,v,\,w,\,y,\,z$. Since $\varphi$ is an isomorphism away from from $s=0$, it follows that $Z_{6,8}$ is invariant under the action of $\boldsymbol{\mu_2}$ on $\mathbb{P}(1,3,4,5,1,1)$. Hence, we can restrict $\Phi' \colon T'=T_s/\boldsymbol{\mu_2} \rightarrow \mathbb{P}(1,3,4,5,1,1)/\boldsymbol{\mu_2}$ to $Y_s$ to get the divisorial contraction to a $cE/2$ point,
\[
\varphi_s' \colon Y' \rightarrow Z_{6,8} \subset \mathbb{P}(1,3,4,5,1,1)/\boldsymbol{\mu_2}.
\]
This is a weighted blowup with weights 
\[
\wt(u,\xi,v,y)=\frac{1}{2}(3,4,5,1)
\]
and, by Lemma \ref{lem:discr}, discrepancy $\frac{1}{2}$. See \cite[Section~10]{HayakawaI} for the classification of terminal divisorial contractions to $cE/2$ points with discrepancy smaller than 1. 

\paragraph{Family 123.} With the assumptions declared before the proposition, $X$ can be written as 
\begin{align*}
\xi w+vt+z^3+f_{12}(x,z)&=0\\
\xi^4x+\xi^2g_{8}(x,z)+\xi g_{11}(x,z,t,v) + wt+v^2+vg_{7}(x,z,t)+t^2z+tg_{9}(x,z)+g_{14}(x,z)&=0
\end{align*}
in $\mathbb{P}:=\mathbb{P}(2,3,4,5,7,9)$ with homogeneous variables $x,\,\xi,\,z,\,t,\,v,\,w$. We consider the blowup $\Phi \colon T \rightarrow \mathbb{P}$ whose restriction is the unique Kawamata blowup $\varphi \colon E \subset Y \rightarrow X$centred at $\mathbf{p_{\xi}}$. As before we use Corollary \ref{cor:genlift} and write $T$ as 
\[
\begin{array}{cccc|ccccccc}
             &       & u  & \xi &  v & z & w & t & x & \\
\actL{T}   &  \lBr &  0 & 3 &  7 & 4& 9 & 5 & 2 &   \actR{.}\\
             &       & 1 & 1 & 2 & 1 & 2 & 1 & 0 &  
\end{array}
\]
Moreover, the defining equations of $Y \colon (\widetilde{f}=\widetilde{g}=0) \subset T$ are such that $\widetilde{f} \in |-3K_Y|$ and $\widetilde{g} \in |-\frac{7}{2}K_Y+\frac{1}{2}E|$. We use variation of GIT to get the following relations between ample models in adjacent chambers of $\Mov(T)$:
\[
        \begin{tikzcd}[ampersand replacement=\&,  column sep={6em,between origins}]
             T   \ar[rr, dashed, "\displaystyle{(-7,-1,1,4,3,4)}" ] \ar[d, swap, "\displaystyle{\Phi}"] \& \& T_1 \ar[rr,dashed, "\displaystyle{(-4,-1,-1,1,1,2)}" ]  \& \& T_2 \ar[rr,dashed, "\displaystyle{(-9,-3,-4,-1,1,4)}" ]\& \& T' \ar[d, "\displaystyle{\Phi'}"]\\
           \mathbf{p} \in \mathbb{P}  \& \&  \& \& \& \&  \mathbf{p'} \in \mathbb{P}'
        \end{tikzcd}
    \]
where $T'$ is 		
\[
\begin{array}{ccccccc|cccc}
             &       & u  & \xi &  v & z & w & t & x & \\
\actL{T'}   &  \lBr &  5 & 2 &  3 & 1 & 1 & 0 & -2 &   \actR{.}\\
             &       & 1 & 1 & 2 & 1 & 2 & 1 & 0 &  
\end{array}
\]

Define $T_s$ to be the double cover of $T'$ given by $s^2=x$ and consider the action of $\boldsymbol{\mu_2}$ on $T_s$ given by
\[
(u  , \xi ,  v , z , w , t , s) \mapsto (u  , \xi ,  v , z , w , t , -s).
\]
Then $T' = T_s / \boldsymbol{\mu_2}$. On the other hand, the divisorial contraction
\[
\Phi_s \colon T_s \rightarrow \Proj \mathcal{F}, \quad \mathcal{F}=\bigoplus_{m\geq 1}H^0(T_s,m\mathcal{O}(t=0))
\]
is given in coordinates by $(u  , \xi ,  v , z , w , t , s)\mapsto(us^5:\xi s^2:vs^3:zs:ws:t) \in \mathbb{P}(1,1,2,1,2,1)$. By lemma \ref{lem:fakewps}, $\Phi$ descends to a map on the orbits $\Phi'\colon T'=T_s/\boldsymbol{\mu_2} \rightarrow \mathbb{P}(1,1,2,1,2,1)/\boldsymbol{\mu_2}$ where the action on $\mathbb{P}(1,1,2,1,2,1)$ is carried over from $T_s$ by $\Phi_s$. That is, the action on $\mathbb{P}(1,1,2,1,2,1)$ is $(u,\xi,v,w,y,z)\mapsto(-u:\xi:-v:-z:-w:z)$.

We now restrict this construction to $Y$. Since $\widetilde{g} \in |-\frac{7}{2}K_Y+\frac{1}{2}E|$ we have that $v^2 \in \widetilde{g}$. Since $\widetilde{f} \in |-3K_Y|$ we have that $yw,\, z^3 \in \widetilde{f}$. On the other hand, although $wt \in g$, it lifts to $wtu \in \widetilde{g}$ since $wt \in |-\frac{7}{2}K_Y-\frac{1}{2}E|$. For degree reasons, there is no $\alpha >1$ for which $w^{\alpha}x_{\mu} \in g$ for any $x_{\mu}$. We conclude that the map $T \rat T_1 \rat T_2 \rat T'$ restricts to an isomorphism $Y \rightarrow Y_2$ followed by a flip to $Y' \subset T'$ of type $(-9,-4,-1,1,4;-8)$. Let $Y_s \subset T_s$ be $Y'|(s^2-x=0)$. By construction, $Y_s$ is invariant under the action of $\boldsymbol{\mu_2}$ on $T_s$. On the other hand, the map $\varphi_s := \Phi_s|_{Y_s} \colon Y_s \rightarrow Z_{3,4}$ is a divisorial contraction to the point $\mathbf{p_t}\sim cA/2$ where $Z_{3,4}$ is
\begin{align*}
\xi w+vt+z^3+f_{12}(u,z)&=0\\
\xi^4+\xi^2g_{8}(1,z)+\xi g_{11}(1,z,t,v) + wtu+v^2+vg_{7}(1,z,t)+t^2zu+tg_{9}(1,z)+g_{14}(u,z)&=0.
\end{align*}
inside $\mathbb{P}(1,1,2,1,2,1)/\boldsymbol{\mu_2}$ with homogeneous variables $u,\,\xi,\,v,\,z,\,w,\,t$. Locally around $\mathbf{p_t}$, the variable $v$ can be eliminated and $\varphi_s$ is a weighted blowup of weights 
\[
\wt(u,\xi,z,w)=\frac{1}{2}(5,2,1,1)
\]
with discrepancy $\frac{1}{2}$ by lemma \ref{lem:discr}. See \cite[Section~6]{HayakawaI} for the classification of terminal divisorial contractions to $cA/2$ points with discrepancy smaller than 1. 

\paragraph{Family 125.} With the assumptions declared before the proposition, $X$ can be written as 
\begin{align*}
\xi w+yv+f_2(t,z)&=0\\
\xi^6 y+\xi^2g_{9}+\xi g_{13}+wv+vyt+g_3(t,z)+g_{15}&=0
\end{align*}
in $\mathbb{P}:=\mathbb{P}(2,3,5,5,7,8)$ with homogeneous variables $\xi,\,y,\,z,\,t,\,v,\,w$. We consider the blowup $\Phi \colon T \rightarrow \mathbb{P}$ whose restriction is the unique Kawamata blowup $\varphi \colon E \subset Y \rightarrow X$ centred at $\mathbf{p_{\xi}}\sim \frac{1}{2}(1,1,1)$ with orbinates $t,\,z,\,v$. As before we use C \ref{cor:genlift} and write $T$ as 
\[
\begin{array}{cccc|ccccccc}
             &       & u  & \xi &  v & t & z & w & y & \\
\actL{T}   &  \lBr &  0 & 2 &  7 & 5 & 5 & 8 & 3 &   \actR{}\\
             &       & 1 & 1 & 3 & 2 & 2 & 3 & 0 &  
\end{array}
\]
and $Y \colon (\widetilde{f}=\widetilde{g}=0) \subset T$, where the defining equations of $Y$ are (pluri)-anticanonical, i.e., $\widetilde{f} \in |-2K_Y|$ and $\widetilde{g} \in |-3K_Y|$. We use variation of GIT to get the following relations between ample models in adjacent chambers of $\Mov(T)$:
\[
        \begin{tikzcd}[ampersand replacement=\&,  column sep={5em,between origins}]
             T   \ar[rr, dashed, "\displaystyle{(-7,-1,1,1,3,9)}" ] \ar[d, swap, "\displaystyle{\Phi}"] \& \& T_1 \ar[rr,dashed, "\displaystyle{(-5,-1,-1,1,6)}" ]\& \& T' \ar[d, "\displaystyle{\Phi'}"]\\
           \mathbf{p} \in \mathbb{P}  \& \& \& \&  \mathbf{p'} \in \mathbb{P}'
        \end{tikzcd}
    \]
where $T_1 \rat T'$ is a small contraction over $\mathbb{P}^1$ and $T'$ is 		
\[
\begin{array}{ccccccc|cccc}
             &       & u  & \xi &  v & t & z & w & y & \\
\actL{T'}   &  \lBr &  8 & 2 &  3 & 1 & 1 & 0 & -9 &   \actR{.}\\
             &       & 1 & 1 & 3 & 2 & 2 & 3 & 0 &  
\end{array}
\]

Define $T_s$ to be the degree  $9$ cover of $T'$ given by $s^9=y$ and consider the action of $\boldsymbol{\mu_9}$ on $T_s$ given by
\[
(u  , \xi ,  v , t , z , w , s) \mapsto (u  , \xi ,  v , t , z , w ,  \epsilon s)
\]
where $\epsilon \in \boldsymbol{\mu_9}$ is a primitive third root of unity. Then $T' = T_s / \boldsymbol{\mu_9}$. The divisorial contraction
\[
\Phi_s \colon T_s \rightarrow \Proj \mathcal{F}, \quad \mathcal{F}=\bigoplus_{m\geq 1}H^0(T_s,m\mathcal{O}(t=0))
\]
is given in coordinates by 
\[
(u , \xi ,  v , t , z , w , s)\mapsto(us^{8}:\xi s^{6}:vs^3:ts:zs:w) \in \mathbb{P}(1,1,3,2,2,3):=\mathbb{P}.
\] 
In $\mathbb{P}$ we have $(u , \xi ,  v , t , z , w ) = (\lambda^i u ,\lambda^i \xi ,\lambda^{3i}  v ,\lambda^{2i} t ,\lambda^{2i} z , \lambda^{3i} w ) $ for all integers $i$ and all $\lambda \in \mathbb{C}^*$. On the other hand, the $\boldsymbol{\mu_9}$-action on $T_s$ is carried over through to  $\mathbb{P}$ via $\Phi_s$. Denote by $\epsilon \in \boldsymbol{\mu_9}$ a primitive 9th root of unity. The $\boldsymbol{\mu_9}$-action on $\mathbb{P}$ is therefore,
\[
\epsilon \cdot (u , \xi ,  v , t , z , w , s)=(\epsilon^{8}u:\epsilon^{6}\xi:\epsilon^3v: \epsilon t:\epsilon z:w).
\]
However, a generic point $\mathbf{p}$ under this action has an orbit with three elements. Indeed, $\epsilon^3 \cdot \mathbf{p} = \mathbf{p}$ and $\epsilon^3 \not =1$. Moreover, the actions of $\epsilon$ and $\epsilon^2$ are different and not trivial. Hence, we have an action of each coset of $\boldsymbol{\mu_9}/\boldsymbol{\mu_3}$ on $\mathbb{P}$ which are 
\[
K_0:=\{\epsilon^{3i}\, |\, 1 \leq i \leq 3 \},\, K_1:=\{\epsilon^{1+3i}\, |\, 1 \leq i \leq 3 \},\, K_2:=\{\epsilon^{2+3i}\, |\, 1 \leq i \leq 3 \}  
\]
given by
\begin{align*}
K_0 \cdot (u , \xi ,  v , t , z , w ) &= (u , \xi ,  v , t , z , w ) \\
K_1 \cdot (u , \xi ,  v , t , z , w ) &= (\epsilon^8 u , \epsilon^2 \xi , \epsilon^3 v , \epsilon  t ,\epsilon z , w ) \\
K_2 \cdot (u , \xi ,  v , t , z , w ) &= (\epsilon^7 u , \epsilon^4 \xi , \epsilon^6 v , \epsilon^2  t ,\epsilon^2 z , w ) \\ 
\end{align*}
Of course, $\boldsymbol{\mu_9}/\boldsymbol{\mu_3} \simeq  \boldsymbol{\mu_3}$ and we denote this group by $\boldsymbol{\mu_3}^*$ to distinguish it from $\boldsymbol{\mu_3}$, which acts trivially on $\mathbb{P}$. Hence, $\Phi$ descends to a map on the orbits $\Phi'\colon T'=T_s/\boldsymbol{\mu_9} \rightarrow \mathbb{P}(1,1,3,2,2,3)/\boldsymbol{\mu_3}^*$.

We now restrict this construction to $Y$. Since $\widetilde{g} \in |-3K_Y|$ we have that $wv \in \widetilde{g}$. On the other hand, although $yv \in f$, it lifts to $yvu \in \widetilde{f}$ since $yv \in |-2K_Y-E|$. For degree reasons there is no $\alpha >1$ for which $v^{\alpha}x_{\mu} \in f$ for any $x_{\mu}$. Hence, crossing the $v$-wall on the movable cone of $T$ restricts to an anti-flip $(-7,-1,1,1,9;2)$ where the variable $w$ can be (locally around $\mathbf{p_v}$) eliminated. The small contraction $T_1 \rat T'$ restricts to an isomorphism on $Y'$ since it happens away from $Y_1$. We call $Y'$ to $Y_1$ seen inside $T'$.

We define $Y_s \subset T_s$ to be $Y'|(s^9-x=0)$.
 By construction, $Y_s$ is invariant under the action of $\boldsymbol{\mu_9}$ on $T_s$. On the other hand, the map $\varphi_s := \Phi_s|_{Y_s} \colon Y_s \rightarrow Z_{4,6}$ is a divisorial contraction to the point $\mathbf{p_w}$ where $Z_{4,6}$ is
\begin{align*}
\xi w+vu+f_2(t,z)&=0\\
\xi^6 +\xi^2g_{9}+\xi g_{13}+wv+vtu+g_3(t,z)+g_{15}&=0
\end{align*}
inside $\mathbb{P}(1,1,3,2,2,3)/\boldsymbol{\mu_3}^*$ with homogeneous variables $u,\xi,\,v,\,t,\,z,\,w$. Locally around $\mathbf{p_w}$ we can eliminate the variables $\xi$ and $v$ and therefore $\varphi' \colon Y' \rightarrow Z_{4,6}$ is a weighted blowup with weights 
\[
\wt(u,t,z)=\frac{1}{9}(8,1,1).
\]
Moreover, $K_1 \cdot (u,t,z) = (\epsilon^8 u, \epsilon t, \epsilon z)$, hence $\mathbf{p_w}$ is a cyclic quotient singularity of type $\frac{1}{9}(8,1,1)$ and $\varphi'$ is its Kawamata blowup.
\end{proof}


\begin{Ex} \label{ex:cod6}
By Proposition \ref{prop:fake1}, a quasismooth member of family 117 is birational to $Z_{4,5} \subset \mathbb{P}(2,3,1,1,2,1)/\boldsymbol{\mu_2}$. We recall that the action of $\boldsymbol{\mu_2}$ on $\mathbb{P}(2,3,1,1,2,1)$ is given by 
\[
\epsilon \cdot (u:v:x:y:w:t) \mapsto (-u:v:-x:-y:-w:t) = (-u:-v:x:y:-w:-t).
\]
We find explicit generators and relations for the ring of polynomials invariant under this action, $\mathbb{C}[u,v,x,y,w,t]^{\mathbf{\mu}_2}$ as in the proof of Lemma \ref{lem:fakewps}. We call its generators $x_1,\ldots, x_{12}$ and they are defined as in the following matrix
\[
\left(
\begin{array}{*{16}c}
	x_1&x_2&x_3&x_4&x_5&x_6&x_7&x_8&x_9&x_{10}&x_{11}&x_{12}\\
		t^2&tw&tv&tu&w^2&wv&wu&v^2&vu&u^2&y&x\\
			2&3&4&3&4&5&4&6&5&4&1&1\\
			\end{array}
\right)
\]
where the third row is the degree of the respective generator. The relations can be easily written as $\rk A = 1$, where $A$ is the symmetric $4\times 4$ matrix:
\[
\begin{pmatrix}
	t^2&tw&tv&tu\\
	&w^2&wv&wu\\
	&&v^2&vu\\
	&&&u^2\\
\end{pmatrix}=:
\begin{pmatrix}
	x_1&x_2&x_3&x_4\\
	&x_5&x_6&x_7\\
	&&x_8&x_9\\
  &&&x_{10}\\
\end{pmatrix}
\]
Finally,
\[
\Proj \mathbb{C}[u,v,x,y,w,t]^{\mathbf{\mu}_2} \cong \Proj \mathbb{C}[x_1,\ldots,x_{12}]/(\rk A = 1).
\]
That is,
\[
\mathbb{P}(2,3,1,1,2,1)/\mathbf{\mu}_2 \cong  (\rk A = 1) \subset \mathbb{P}(2,3,4,3,4,5,4,6,5,4,1,1) =: \mathbb{P}'.
\]

The equations defining $Z_{4,5}$ are $f',\, g'$ under this contraction and can be written in the new coordinates as 
\begin{align*}
f' &\colon  x_3+x_7+a_4(x_{11},x_{12})+\ldots= 0 \\
g' &\colon  x_6+x_{10}x_{12}+b_5(x_{11},x_{12})+x_1x_4+\ldots = 0
\end{align*}
inside $\mathbb{P}'$. Hence, $x_3$ and $x_6$ can be globally eliminated. So the final 3-fold is 
\[
Z \colon (\rk A=1) \subset \mathbb{P}(2,3,3,4,4,6,5,4,1,1).
\]

This is a quasismooth member of \#5414 in the GRDB and provides an example of an explicit construction of a Fano 3-fold in weighted projective space of high codimension, in this case 6. 
\end{Ex}

\section{Exclusion}  \label{sect:excl}

Throughout this section let $X$ be a quasismooth member of a family in $I_{S}$. From the previous sections we have that $\iota_X \in \{2,3\}$. In this section we exclude certain closed subvarieties as maximal centres for $X$.


\subsection{Exclusion of curves}

Recall that the degree of a curve $\Gamma \subset X$ is the degree with respect to $A$, where $-K_X \sim \iota_X A$ and $\iota_X$ is the Fano index of $X$. That is, $\deg \Gamma = (A\cdot \Gamma)$ and that $(A^3)$ is the anticanonical degree of $X$. The following result gives bounds on the degree of $\Gamma$ provided it is a maximal centre.

\begin{Lem} \label{lem:excurves}
Suppose that the curve $\Gamma \subset X$ is a maximal centre. Then 
\[
1 \leq \deg \Gamma <\iota_X ^2 A^3.
\] 
\end{Lem}

\begin{proof}
Since $A$ is ample, $\deg \Gamma = A\cdot \Gamma > 0$ and since $\Gamma$ is in the smooth locus of $X$ by Corollary \ref{cor:kwblcor}, the intersection is an integer and therefore $\deg \Gamma \geq 1$. If $\Gamma$ is a maximal centre then, by \cite[Theorem~5.1.1]{CPR}, it follows that $(-K_X)^3 > (-K_X\cdot C)$. The conclusion follows immediately from the fact that $-K_X \sim \iota_X A$.
\end{proof}

\begin{Cor}
Let $\Gamma$ be an irreducible curve on $X \in I_S \setminus \{ 92\}$. Then $\Gamma$ is not a maximal centre. If $X$ is a quasismooth member of family 92 and $\Gamma \subset X$ is a maximal centre, then $\deg \Gamma =1$.
\end{Cor} 

\begin{proof}
All families in $I_S$ satisfy $\iota_X^2 A^3 \leq 1$ with the exception of family 92. In the case of family 92,
\[
\iota_X^2 A^3 = \frac{16}{15}.
\]
By Lemma \ref{lem:excurves} we have $1 \leq \deg \Gamma <\frac{16}{15}$ and, since $\Gamma$ is in the smooth locus of $X$, $\deg \Gamma = 1$.
\end{proof}


\begin{Ex} \label{ex:curveblow} This examples illustrates how the quasismoothness hypothesis is essential.
Consider the special singular sextic $Z_6 \subset \mathbb{P}(1,1,1,1,3)$ with homogeneous coordinates $t,\,v,\,x,\,y,\,u$ given by 
\[
u^2+ut(v^2+tx)+vf_5(x,y)+g_6(x,y)=0.
\]
Since $Z_6$ is not quasismooth we can not guarantee that if $\Gamma \subset Z_6$ is a maximal centre then it is contained in the smooth locus of $Z_6$. In fact, the singular locus of $Z_6$ is contained in $\Gamma \colon (x=y=u=0) \subset Z_6$ and the blowup of the ideal sheaf of $\mathcal{I}_{\Gamma}$ over the generic point of $\Gamma$ initiates a Sarkisov link as in Theorem \ref{thm:divtocurve}. Also notice that any smooth sextic $X \subset \mathbb{P}(1,1,1,1,3)$ is birationally super-rigid, see \cite{Iskdoublesolid}, and, in particular, no curve can be a maximal centre. 
\end{Ex}

\subsection{Exclusion of non-singular points} \label{subsec:smooth}

Throughout this subsection we assume $X \subset \mathbb{P}(a_0,\ldots,a_5)$ with $a_0 \leq \cdots \leq a_5$ and the generality assumptions in table \ref{tab:gen}.

\begin{table}
\centering
\begin{tabular}{ccc}
\toprule
  Family        & Generality Assumption I     &  Generality Assumption II          \\
\midrule 
\midrule 
  95&         $X|_{x=y=0}$ is irreducible   & $(1,4,5,5) \not \in X$    \\
\midrule
  96&  $X|_{x=y=0}$ is irreducible &         \\
\midrule
 98 &  &  $(1,4,5,7) \not \in X$   \\
\bottomrule
\end{tabular}
\caption{Generality assumptions on $X$. The notation $(c_1,c_2,c_3,c_4)$ denotes an irreducible and reduced curve defined by four homogeneous polynomials of degrees $c_1,\,c_2,\,c_3,\,c_4$, respectively. Notice that, for each 4-tuple mentioned, there is a choice of a WCI curve in $X$. See Table \ref{tab:tabisol}.}
\label{tab:gen}
\end{table}

\begin{Def}[{{\cite[Definition~5.2.4]{CPR}}}]
Let $L$ be a Weil divisor class on a variety $X$ and $p \in X$ a point. For an integer $s>0$, consider the linear system 
\[
\mathcal{L}^s_p = |\mathcal{I}^s_p(sL)| 
\]
where $\mathcal{I}_p$ is the ideal sheaf of $p$. We say that the class $L$ \textbf{isolates} $p$ or is $p$-\textbf{isolating} if $p \in \Bs \mathcal{L}^s_p$ is an isolated component for some positive integer $s$.
\end{Def}

To exclude most of the non-singular points as maximal centres we use the following criterion:

\begin{Lem}[{{\cite[Pages 210 \& 211]{CPR}}}]
Let $\mathbf{p}$ be a non-singular point of a Fano $3-$fold X. If $lA$ isolates $\mathbf{p}$ and $0<l\leq 4/\iota_X^2 A^3$, then $\mathbf{p}$ is not a maximal centre.
\end{Lem}

\begin{proof}
It follows from the fact that $-K_X \sim \iota_X A$ and \cite[Page 210 \& 211]{CPR}.
\end{proof}

We call the number $4/\iota_X^2 A^3$ the \textbf{isolating threshold} of $\mathbf{p}$. The following definition is \cite[Definition~5.6.3]{CPR}. We say that a finite set of homogeneous polynomials, $\{g_i\}$, isolates the point $\mathbf{p}$ if
\[
\bigcap_{i}(g_i=0)_X
\] 
has $\mathbf{p}$ as an isolated component.

Following \cite{okadaI}, we define the rational map $\pi_k \colon X \rat \mathbb{P}(a_0,\ldots, \widehat{a_k}, \ldots, a_5)$ to be the projection away from $\mathbf{p_k}$ and $\Exc(\pi_k) \subset X$ the locus contracted by $\pi_k$.

\begin{Lem}[{{\cite[Proposition~7.4]{okadaI}}}] \label{lem:isol2}
Let $\mathbf{p}=(\xi_0:\ldots:\xi_5)$ be a smooth point of $X$ and $j$ an index for which $\xi_j \not = 0$. Let
\[
m=\max_{0\leq l\leq 5, l\not = j} \{\lcm(a_j,a_l) \}.
\]
Then $mA$ isolates $\mathbf{p}$. Additionally, if $\mathbf{p} \not \in \Exc(\pi_k)$ for some $k\not = j$, then $m_kA$ isolates $\mathbf{p}$, where 
\[
m=\max_{0\leq l\leq 5, l\not = j,k} \{\lcm(a_j,a_l) \}.
\]
\end{Lem}


\begin{Cor}
Let $X \in I_S$ be such that 
\[
a_2a_5 \leq 4/\iota_X^2A^3.
\]
Then, no smooth point of $X$ is a maximal centre.
\end{Cor}

\begin{proof}
The families satisfying the condition are $X \in \{102,103,106,107,109,110,111 \}$ and we divide the proof in 3 cases: 
\begin{enumerate}
	\item Suppose that $a_4a_5 \leq 4/\iota_X^2 A^3$. This is satisfied for family 111. Let $m$ be as defined in Lemma \ref{lem:isol2}. Then, for all $0 \leq j,\, l \leq 5$, we have that $m \leq a_4a_5$. It follows that $mA$ isolates $\mathbf{p}$. 
	\item Suppose that $a_3a_5 \leq  4/\iota_X^2 A^3 < a_4a_5$. This is satisfied for family 107 and 110. Suppose that $\mathbf{p}=(0:0:0:0:\xi_4:\xi_5)$. Then $\mathbf{p} \in X$ if and only if $\mathbf{p}$ is $\mathbf{p_4}$ or $\mathbf{p_5}$ which is not possible since $\mathbf{p}$ is smooth. Hence, $\xi_j \not = 0$ for some $0\leq j\leq 3$. Then, $m \leq a_3a_5$ and $mA$ isolates $\mathbf{p}$. 
	\item Suppose that $a_2a_5 \leq  4/\iota_X^2 A^3 < a_3a_5$. This is satisfied for family 102, 103, 106 and 109. Suppose that $\mathbf{p}=(0:0:0:\xi_3:\xi_4:\xi_5)$. Then $\mathbf{p} \in X$ if and only if $\mathbf{p}$ is $\mathbf{p_4}$ or $\mathbf{p_5}$ which is not possible since $\mathbf{p}$ is smooth. Hence, $\xi_j \not = 0$ for some $0\leq j\leq 2$. Then, $m \leq a_2a_5$ and $mA$ isolates $\mathbf{p}$.	
	\end{enumerate}
\end{proof}

\begin{Cor}
Let $X \in I_S$ be such that 
\[
a_1a_5 \leq 4/\iota_X^2A^3 < a_2a_5 \quad \text{and} \quad a_2a_4 \leq 4/\iota_X^2A^3. 
\]
Then no smooth point of $X$ is a maximal centre.
\end{Cor}

\begin{proof}
The condition is satisfied for families  $X \in \{101, 105, 108 \}$. Let $\mathbf{p}=(\xi_0:\xi_1:\xi_2:\xi_3:\xi_4:\xi_5)$. If $\xi_0$ or $\xi_1$ is non-zero then there is an $l$ such that $lA$ is a $\mathbf{p}-$isolating class so we can assume $\xi_0=\xi_1=0$. Moreover, $\xi_2 \not =0$ since, otherwise, the point would be singular. For each of these families we consider the projection $\pi_5 \colon X \rat \mathbb{P}(a_0,a_1,a_2,a_3,a_4)$. 

For families 101 and 105, we see that $(x_0=x_1=0) \cap \Exc(\pi_5) = (x_0=x_1=x_2=x_3=0)_X$ which is equal to $\{\mathbf{p_5}\}$ - a singular point - by quasismoothness. For family 108, we have $(x_0=x_1=0) \cap \Exc(\pi_5) = (x_0=x_1=x_3=0)_X=\{\mathbf{p_5}\}$ again by quasimoothness.

By Lemma \ref{lem:isol2}, we have that $a_2a_4A$ isolates $\mathbf{p}$ and $a_2a_4 \leq 4/\iota_X^2A^3$.
\end{proof}

\begin{Lem}[{{\cite[Lemma~5.6.4]{CPR}}}] \label{lem:isol1}
Suppose that $\{g_i\}$ isolates $\mathbf{p}$. Then, $lD$ isolates $\mathbf{p}$ where $l=\max \{\deg(g_i) \}$ and $D$ is a Weil divisor in the class of a hyperplane, that is, such that $\mathcal{O}_X(D)\cong \mathcal{O}_X(1)$. 
\end{Lem}

We call $l$ the \textbf{degree} of the set $\{g_i \}$. If the positive integer $l$ does not exceed the isolating threshold of $\mathbf{p}$, then we say that $\{ g_i\}$ is a \textbf{good isolating set} for $\mathbf{p}$. 

We briefly explain how to exclude the smooth points as maximal centres for most of the remaining families $X \in I_S$. We find a curve $C \subset \mathbb{P}$ given by homogenous polynomials $\{ g_i\}$ of a certain bounded degree such that $C$ intersects $X$ in at least a smooth point of $\mathbf{p} \in X$ but $C \not \subset X$. If the degree of the polynomials defining $C$ are small enough, that is, no higher than the isolating threshold for $\mathbf{p}$, then $C$ is a good isolating set for $\mathbf{p}$ and it follows that $\mathbf{p}$ is not a maximal centre of $X$.


\begin{Cor}
Let $X \in I_S$ be such that 
\[
a_1a_5 \leq 4/\iota_X^2A^3 < a_2a_5 \quad \text{and} \quad  4/\iota_X^2A^3 < a_2a_4. 
\]
Then no smooth point of $X$ is a maximal centre.
\end{Cor}

\begin{proof}
The condition is satisfied by $X\in \{104, 117\}$. Suppose $X$ is a quasismooth member of one of these families embedded in weighted projective space with homogeneous variables $x,\,y,\,z,\,t,\,v,\,w$. Let $\mathbf{p}=(\xi_0:\xi_1:\xi_2:\xi_3:\xi_4:\xi_5) \in X$ be a smooth point. If $\xi_0$ or $\xi_1$ is non-zero then there is an $l$ such that $lA$ is a $\mathbf{p}-$isolating class so we can assume $\xi_0=\xi_1=0$. Moreover, $\xi_2,\,\xi_3,\, \xi_4,\, \xi_5 \not =0$ since, otherwise, the point would be singular. Let $I_X=(f,g) \subset \mathbb{C}[x,y,z,t,v,w]$. Then, $f,\,g \not \in (x,y)$ and so $X|_{x=y=0}$ is a curve and it contains $\mathbf{p}$. To isolate it we add further relations. 
\paragraph{Family 104:} Consider the set $S:=\{x,y,\xi_2\xi_5t^2-\xi_3^2wz, \xi_3\xi_4z^3-\xi_3^3vt\}$. Then, $X|_S \colon (x=y=\xi_2\xi_5t^2-\xi_3^2wz=wz+t^2=wt+v^2=\xi_3\xi_4z^3-\xi_3^3vt=0)$ is a zero dimensional scheme containing $\mathbf{p}$ provided that $X$ does not contain any weighted complete intersection curve of weights $(1,2,14,15)$.  In that case, $X|_{S}$ contains $\mathbf{p}$ as an isolated component  and $15A$ isolates $\mathbf{p}$ since $15 < 4/\iota_X^2A^3 = \frac{45}{2}$.

\paragraph{Family 117:} After a change of variables we can assume that $z^3 \not \in f$. Consider the relations
\[
S := (x=y=\xi_5z^2-\xi_2^2w=0) \subset \mathbb{P}(3,3,4,5,7,8).\\
\]
Then, $X|_S \colon (x=y=\xi_5z^2-\xi_2^2w=wz+tv=wv+t^3+\alpha z^2v=0)$ is a zero dimensional scheme containing $\mathbf{p}$. In fact notice that if $-t^3=v(w+\alpha z^2)$ is a multiple of $\xi_5z^2-\xi_2^2w$, then $t=0$ and so $\mathbf{p}$ is singular.  

 We conclude that $X|_{S}$ contains $\mathbf{p}$ as an isolated component  and $15A$ (resp. $8A$) isolates $\mathbf{p}$ in $X_{104}$ (resp. in $X_{117}$) since $15 < 4/\iota_X^2A^3 = \frac{45}{2}$ (resp. $8 < 4/\iota_X^2A^3 = \frac{224}{9}$).


\end{proof}

\begin{Cor} \label{cor:smthpointexcl}
Let $X \in I_S$ be such that 
\[
a_0a_5 \leq 4/\iota_X^2A^3 < a_1a_5.
\]
Then no smooth point of $X$ is a maximal centre.
\end{Cor}

\begin{proof}
The condition is satisfied for $X\in \{93,95,96,97,98,99,100\}$. Let Let $\mathbf{p}=(\xi_0:\xi_1:\xi_2:\xi_3:\xi_4:\xi_5) \in X$ be a smooth point. If $\xi_0$ is non-zero then there is an $l$ such that $lA$ is a $\mathbf{p}-$isolating class so we can assume $\xi_0=0$. 

\paragraph{Family 93:} 
By a change of coordinates we can assume $w^2,\,v^2 \not \in g_{10}$.

\begin{itemize}
	\item Assume $\xi_1=0$. We have
	\[
X|_{x=y=0} \colon (t^2+\alpha z^3=wv+zg_{8}(z,t,v)=0) \subset \mathbb{P}(2,3,5,5).
	\]
	 If $\alpha = 0$, then $\xi_3 =0$ and $z^4 \in g_{8}$ by quasismoothness of $X$. Hence, $X|_{x=y=t=0} \colon (wv+z^5=0) \subset \mathbb{P}(2,5,5)$. Notice that $\xi_4\xi_5 \not =0$.  Therefore $I_C$ is the ideal of a 1-dimensional scheme where $C:=\{x,y,t,\xi_5v-\xi_4w \}$. Moreover, it isolates $\mathbf{p}$ since $g|_C \not \in I_C$ and is indeed a good isolating set of degree 5.
	On the other hand, if $\alpha \not =0$, then $\xi_2 = 0 \iff \xi_3 = 0$ which in turn implies $\mathbf{p} \in \{\mathbf{p_v}, \mathbf{p_w} \}$. So we can assume $\xi_2,\,\xi_3 \not = 0$ and let $C:=\{x,y,\xi_2\xi_3v-\xi_4zt,\xi_2\xi_3w-\xi_5zt \}$. Then $I_C$ is the ideal of a 1-dimensional scheme and $I_X \not \subset I_C$ since $f|_C \not \in I_C$ and $C$ is indeed a good isolating set of degree 5.
		\item The case $\xi_2=0$ is completely analogous
	\item Assume $\xi_1,\,\xi_2 \not = 0$ and $\xi_3 =0$. We have
	\[
X|_{x=t=0} \colon (f_6(y,z)=wv+g_{10}(y,z)=0) \subset \mathbb{P}(2,2,5,5).
	\]
	We have $\xi_4\xi_5 \not = 0$ since, otherwise, $\xi_1=\xi_2=0$ by Lemma \ref{lem:awaybase}, contradicting our assumption. Equivalently, $f_6$ is a multiple of $\xi_1z-\xi_2y$ but $g_{10}$ is not. Hence $\{x,\xi_1z-\xi_2y, t, \xi_4w-\xi_5v \}$ is a good isolating set.
	\item Assume $\xi_1,\,\xi_2,\,\xi_3 \not = 0$. We have after a change of variables
	\[
X|_{x=v=w=0} \colon (t^2+f_6(y,z)=g_{10}(y,z)=0) \subset \mathbb{P}(2,2,3).
	\]
	By Lemma \ref{lem:awaybase}, $\xi_3 \not = 0$. Hence, $ \{x,\xi_1z-\xi_2y,v,w \} $ is a good isolating set if $\xi_4=\xi_5=0$. Otherwise let $C:=\{x, \xi_2y-\xi_1z,\xi_2\xi_3v-\xi_4zt, \xi_2\xi_3w-\xi_5zt \}$. By the assumption on the $\xi_i$, we have that $I_C$ defines a 1-dimensional scheme and $I_X \not \subset I_C$ since $f|_C \not \in I_C$. Morevoer, $C$ is indeed a good isolating set of degree 5.
\end{itemize}

\paragraph{Family 95: (With the assumption that $X|_{x=y=0}$ is irreducible ($\iff z^2t \in g$) and that $X$ contains no $(1,4,5,5)$ WCI curve)} 

After a change of variables we can assume $wv \in g$ and $v^2,\, w^2 \not \in g$.
\begin{itemize}
	\item Assume $\xi_1=0$. We have that
	\[
X|_{x=y=0} \colon ((w+v)z+t^2=wv+z^2t=0) \subset \mathbb{P}(3,4,5,5)
	\]
	is an integral curve since $z^2t \in g$. 
	 Then, $\xi_3\xi_4\xi_5  = 0 \implies \mathbf{p} \in \{\mathbf{p_z},\mathbf{p_v},\mathbf{p_w}\}$. Hence, $z(w+v)$ is not a multiple of $\xi_5 v - \xi_4 w$ and therefore $\{x,y,\xi_5 v - \xi_4 w \}$ is a good isolating set of degree 5.

	\item Assume $\xi_1\not =0,\,\xi_2=0$. We have 
	\[       
X|_{x=z=0} \colon (f_8(y,t)=wv+g_{10}(y,t)=0) \subset \mathbb{P}(2,4,5,5)
	\]
	where $y^4 \in f_8$ or $y^5 \in g_{10}$ by quasismoothness of $X$. Then, $\xi_4\xi_5 \not =0$ by Lemma \ref{lem:awaybase} or, equivalently, $f_8$ is a multiple of $\xi_1^2t-\xi_3y^2$ but $g_{10}$ is not. Then, $\{x,z,\xi_1^2t-\xi_3y^2,\xi_5v-\xi_4w \}$ is a good isolating set of degree 5.
	
	\item Suppose $\xi_1,\,\xi_2 \not = 0$. If $\xi_4,\,\xi_5  =0$, then $C:=\{x,\xi_1^2t-\xi_3y,v,w \}$ is a good isolating set for $\mathbf{p}$ of degree 5. If not both $\xi_4$ and $\xi_5$ are zero, we take $C:=\{x,\xi_3y^2-\xi_1^2t,\xi_1\xi_2v-\xi_4yz, \xi_4w-\xi_5v \}$. By the assumptions on the $\xi_i$, we have that $I_C$ is the ideal of a 1-dimensional scheme which is a $(1,4,5,5)$ WCI curve. By the assumption that $X$ contains no such curve it follows that $C$ is a good isolating set of degree 5.

\end{itemize}

\paragraph{Family 96: (With the assumption that $X|_{x=y=0}$ is irreducible ($\iff t^2 \in f$))} By a change of variables we can assume the only monomial with $v$ appearing with non-zero coefficient in $g$ is $wv$. 
\begin{itemize}
	\item Suppose $\xi_1=0$. We have,
	\[       
X|_{x=y=0} \colon (vz+ t^2=wv+z^4+tg_{8}(z,t,v)=0) \subset \mathbb{P}(3,4,5,7)
	\]
	which is an irreducible curve by assumption. We claim that $C:=\{x,y, \xi_3^2vz-\xi_4\xi_2t^2,\xi_3\xi_2w-\xi_5tz \}$ is a good isolating set of degree 7: If $\xi_2=0$, then $\xi_3 = 0$ by the assumption that $X|_{x=y=0}$ is irreducible. Then $\mathbf{p} \in \{ \mathbf{p_v},\mathbf{p_w}\}$ and can assume $\xi_2 \not  = 0$. On the other hand, if $\xi_3 = 0$, it follows that $\xi_4 = 0$, since $\xi_2 \not = 0$. This implies that $\xi_2 = 0$ from  $wv+z^4+tg_{8}(z,t,v)=0$, which is impossible. So we conclude that $\xi_2,\, \xi_3 \not = 0$. Hence, $I_C$ defines a 1-dimensional projective scheme and $I_X  \not \subset I_C$ since $g|_C \not \in I_C$.

	\item Suppose $\xi_1\not=0,\,\xi_2=0$. We have,
	\[       
X|_{x=z=0} \colon (f_{8}(y,t)=wv+g_{12}(y,t)=0) \subset \mathbb{P}(2,4,5,7).
	\]
	In this case, $C:=\{x,z,\xi_1^2t-\xi_3y^2,\xi_4\xi_1w-\xi_5vy \}$ is a good isolating set: Notice that $X|_{x=z=v=0}=\mathbf{p_w}$ by Lemma \ref{lem:awaybase} so we can assume $\xi_4 \not = 0$. It follows that $I_C$ is an ideal of a 1-dimensional scheme. Moreover,  the only way to have $I_X \subset I_C$ would be for $\xi_5=0$ which in turn implies that $f(y,t),\,g(y,t)$ are multiples of $\xi_1^2t-\xi_3y^2$. But this is impossible by Lemma \ref{lem:awaybase}.

	\item Suppose $\xi_1,\,\xi_2\not=0$ and let $C:=\{x, \xi_1^2t-\xi_3y^2,\xi_1\xi_2v-\xi_4yz,\xi_1^2\xi_2w-\xi_5y^2z \}$. By the assumption on the $\xi_i$, we have that $I_C$ is the ideal of a 1-dimensional scheme. Moreover, $C$ isolates $\mathbf{p}$ since $g|_C \not \in I_C$ and it is a good isolating set of degree 7.

 \end{itemize}

\paragraph{Family 97:}
\begin{itemize}
	\item Suppose $\xi_1 = 0$. We have,
	\[       
X|_{x=y=0} \colon (t^2+\alpha z^5=wt+v^2+zg_{12}(z,t,v)=0) \subset \mathbb{P}(2,5,7,9).
	\]
	If $\alpha = 0$, then $\xi_3=0$ and $X|_{x=y=t=0} \colon (v^2+z^7=0)$. Hence, $\xi_4 = 0 \iff \xi_2=0$ which implies, in turn, $\mathbf{p} = \mathbf{p_w}$. So we can assume $\xi_2,\,\xi_4 \not = 0$ and therefore, $C_1:=\{x,y,t,\xi_2\xi_4w-\xi_5vz \}$ is a good isolating set of degree 9.
	If, on the other hand, $\alpha \not = 0$, then $\xi_2 = 0 \iff \xi_3 =0$ which implies, similarly, $\mathbf{p} = \mathbf{p_w}$, so we can assume $\xi_3,\,\xi_2 \not = 0$. Then $C_2:=\{x,y,\xi_2\xi_3v-\xi_4tz,\xi_3\xi_2^2w-\xi_5tz^2 \}$ is a good isolating set of degree 9 as well. 
	\item The case $\xi_2 = 0$ is completely analogous.
	\item Suppose $\xi_1,\,\xi_2 \not = 0$ and $\xi_3 = 0$. We have,
	\[
	X|_{x=t=0} \colon (f_{10}(y,z)=v^2+g_{14}(y,z)) \subset \mathbb{P}(2,2,7,9).
	\]
	Then $C:=\{x,t,\xi_2y-\xi_1z,\xi_4\xi_1w-\xi_5vy \}$ is a good isolating set: Notice that $X|_{x=t=v=0}= \mathbf{p_w}$ by Lemma \ref{lem:awaybase} and therefore we can assume $\xi_4 \not = 0$. Moreover, $f_{10}$ is a multiple of $\xi_2y-\xi_1z$ and it has no common components with $g_{14}$. Therefore $I_C$ is a 1-dimensional projective scheme such that $I_X \not \subset I_C$ and $C$ is a good isolating set of degree 9.
	
	\item Suppose $\xi_1,\,\xi_2,\,\xi_3 \not = 0$ and let $C:=\{x,\xi_2y-\xi_1z,\xi_3\xi_1v-\xi_4ty,\xi_3\xi_1^2w-\xi_5ty^2 \}$. From the assumptions that $\xi_i\not = 0,\, 1\leq i \leq 3$, it follows that $I_C$ is the ideal of a 1-dimensional projective scheme. Moreover, $C$ isolates $\mathbf{p}$ since $f(0,y,z,t,v,w) =t^2+f_{10}(y,z) \not \in I_C$ and $C$ is a good isolating set of degree 9. Notice that it can happen that $\xi_4=\xi_5=0$. In this case, the isolating set $C$ reduces to $C:=\{x,\xi_2y-\xi_1z,v,w\}$.

\end{itemize}

\paragraph{Family 98: (With the assumption that $X$ contains no $(1,4,5,9)$ WCI curve)} After a change of coordinates we can assume $z^4 \not \in g$ and, indeed, $\mathbf{p_z} \in X$.
\begin{itemize}
\item Suppose $\xi_1 = 0$. We have,
	\[
	X|_{x=y=0} \colon (v^2+\alpha z^2t=wz+t^3+\beta ztv=0) \subset \mathbb{P}(3,4,5,9).
	\]
	If $\alpha = 0$, then, $v=0$ and $C:=\{x,y,v,\xi_2^3w-\xi_5z^3 \}$ defines a 1-dimensional projective scheme since $\xi_2$ and $\xi_5$ are both non-zero. Moreover $I_X \not \subset I_C$ since $wz+t^3+\beta ztv \not \in (v,\xi_2^3w-\xi_5z^3)$. 
	
	If $\alpha \not = 0$, then $\xi_3 = 0 \implies \mathbf{p} \in  \{\mathbf{p_z},\mathbf{p_w}\}$ so we can assume $\xi_3 \not =0$. Moreover, $\xi_4$ and $\xi_5$ are not both zero since otherwise $\mathbf{p}= \mathbf{p_w}$. Then, $C:=\{x,y,\xi_4\xi_2t^2-\xi_3^2vz,\xi_4\xi_3w-\xi_5vt \}$ defines a 1-dimensional projective scheme and $I_X \not \subset I_C$ since $v^2+\alpha z^2t \not \in (\xi_4\xi_2t^2-\xi_3^2vz,\xi_4\xi_3w-\xi_5vt)$. In each case ($\alpha = 0$ and $\alpha \not = 0$),  $C$ is a good isolating set (of degree 9) of $\mathbf{p}$.

		\end{itemize}
	
	\begin{itemize}
		\item 	Suppose $\xi_1 \not = 0$ and $\xi_2 = 0$. We have,
	\[
	X|_{x=z=0} \colon (v^2+f_{10}(y,t)=g_{12}(y,t)=0) \subset \mathbb{P}(2,4,5,9).
	\]
Notice that $X|_{x=z=v=0}= \mathbf{p_w}$ by Lemma \ref{lem:awaybase}, so we can assume $\xi_4 \not = 0$. In this case, $C:=\{x,z,\xi_1^2t-\xi_3y^2,\xi_4\xi_1^2w-\xi_5vy^2 \}$ defines a 1-dimensional projective scheme. Moreover, $I_X \not \subset I_C$ since $f_{10}$ and $g_{12}$ have no common components. Hence $C$ is a good isolating set of degree 9.
	\end{itemize}
\begin{itemize}
	\item Suppose $\xi_1,\xi_2 \not = 0$ and $\xi_3 = 0$. We have,
	\[
	X|_{x=t=0} \colon (v^2+f_{10}(y,z,v)=wz+g_{12}(y,z,v)=0) \subset \mathbb{P}(2,3,5,9).
	\]
	The set $C:=\{x,t,\xi_2\xi_1v-\xi_4zy,\xi_2^3w-\xi_5z^3 \}$ defines a 1-dimensional projective scheme. Moreover, $I_X \not \subset I_C$ since $y^5 \in f_{10}$ or $y^6 \in g_{12}$ from quasismoothness of $X$. Hence $C$ is a good isolating set of degree 9. 
	\end{itemize}
	\begin{itemize}
		\item Suppose $\xi_1,\,\xi_2,\,\xi_3 \not = 0$. Take $C:= \{x,\xi_1^2t-\xi_3y^2,\xi_2\xi_1v-\xi_4zy,\xi_2^3w-\xi_5z^3 \}$. The ideal $I_C$ defines a 1-dimensional projective scheme by the assumptions on the $\xi_i$. Since $X$ contains no $(1,4,5,9)$ WCI curve, it follows that $C$ isolates $\mathbf{p}$ and is, indeed, a good isolating set of degree 9.
		\end{itemize}

\paragraph{Family 99} We first assume that $\xi_1=0$. We have
\[
X|_{x=y=0} \colon (wz+t^2=wt+v^2+g_{12}(z,v)=0) \subset \mathbb{P}(3,5,6,7)
\]
where $g_{12}$ contains at least one monomial by quasismoothness. It follows that $\xi_2,\, \xi_3,\, \xi_5 \not =0$. Let $C:=\{x,y,\xi_2^2v-\xi_4z^2,\xi_3^2wz-\xi_2\xi_5t^2 \}$. We have necessarily  $wz+t^2 \in (\xi_3^2wz-\xi_2\xi_5t^2)$. On the other hand, $wt+v^2+g_{12}(z,v) \not \in (\xi_2^2v-\xi_4z^2,\xi_3^2wz-\xi_2\xi_5t^2)$. Hence $I_X \not  \subset I_C$ and $C$ is a good isolating set (of degree 10) of $\mathbf{p}$.

Suppose, on the other hand, that $\xi_1 \not =0$. The projection 
\[
\pi_5 \colon X \dashrightarrow \mathbb{P}(1,2,3,5,6)
\]
away from $\mathbf{p_w}$ has exceptional locus $\Exc(\pi_5) = (z=t=0)|_X$. Moreover,
\[
X|_{x=z=t=0} \colon (f_{10}(y,v)=g_{12}(y,v)=0) \subset \mathbb{P}(1,3,7). 
\]
By Lemma \ref{lem:awaybase}, $X|_{x=z=t=0} = \{\mathbf{p_w}\}$ and therefore $\mathbf{p} \not \in \Exc(\pi_5)$. Morevoer, 
\[
\max_{0\leq l \leq 5,\, l\not = 1,\,5}\{\lcm (a_1,a_l) \}A = 10A
\]
isolates $\mathbf{p}$ since $10 < \frac{21}{2}$.

\paragraph{Family 100}
\begin{itemize}
	\item Suppose $\xi_1 = 0$. We have,
\[
X|_{x=y=0} \colon (z^4+t^3=wz+v^2+tzg_7(z,t,v)=0) \subset \mathbb{P}(3,4,7,11). 
\]
Clearly, $\xi_2 = 0\,\, (\iff \xi_3 = 0) \implies \mathbf{p} = \mathbf{p_w}$. So we can assume $\xi_2,\,\xi_3 \not = 0$. Suppose $\xi_4 =  \xi_5 = 0$ and let  $C:=\{x,y,v,w \}$. If $g_7(z,t,0)$ is identically zero, then $C$ is clearly a good isolating set for the smooth point $(z^4+t^3=0) \subset \mathbb{P}(3,4)$ in $X$. Otherwise, $g_7(z,t,0)=t^2z^2$ and $X|_C$ is empty.

Suppose at least one of $\xi_4,\,\xi_5$ is not zero and let $C:= \{x,y,\xi_3\xi_2v-\xi_4tz,\xi_4\xi_3w-\xi_5vt \}$. Then, $C$ is a good isolating set for $\mathbf{p}$. Indeed $I_C$ is the ideal of a 1-dimensional projective scheme and $I_X \not \subset I_C$ since $z^4+t^3 \not \in (\xi_3\xi_2v-\xi_4tz,\xi_4\xi_3w-\xi_5vt)$. 
\end{itemize}
\begin{itemize}
	\item Suppose $\xi_1 \not = 0$ and $\xi_2  =0$. We have,
\[
X|_{x=z=0} \colon (f_{12}(y,t)=v^2+g_{14}(y,t)=0) \subset \mathbb{P}(2,4,7,11). 
\]
Notice that if $\xi_4 = 0$, it follows that $X|_{x=z=v=0} = \mathbf{p_w}$ by Lemma \ref{lem:awaybase}, hence we can assume $\xi_4 \not = 0$. If $\xi_3 =0$, then $y^6 \not \in f_{12}$ by the assumption that $\xi_1 \not = 0$. In particular $f_{12}(y,t)$ is a multiple of $t$. Hence, by quasismoothness of $X$, $y^7 \in g_{14}$ and $X|_{x=z=t=0} \colon (v^2+y^7=0)$. In this case, $C:= \{x,z,t,\xi_4\xi_1^2w-\xi_5vy^2 \}$ is a good isolating set. Otherwise, i.e., if $\xi_3 \not = 0$, $C:= \{x,z,\xi_1^2t-\xi_3y^2,\xi_4\xi_3w-\xi_5vt \}$ is a good isolating set. 
\end{itemize}
\begin{itemize}
	\item Suppose $\xi_1,\,\xi_2 \not = 0$ and $\xi_3  =0$. We have,
\[
X|_{x=t=0} \colon (z^4+f_{12}(y,z,v)=wz+v^2+g_{14}(y,z,v)=0) \subset \mathbb{P}(2,3,7,11). 
\]
In this case we take the good isolating set $C:=\{x,t,\xi_2\xi_1^2v-\xi_4zy^2,\xi_2\xi_1^4w-\xi_5zy^4 \}$. Notice that $I_X \not \subset I_C$ since $z^4+f_{12}(y,z,v) \not \in (\xi_2\xi_1^2v-\xi_4zy^2,\xi_2\xi_1^4w-\xi_5zy^4)$.
\end{itemize}
\begin{itemize}
	\item Suppose $\xi_1,\,\xi_2,\,\xi_3 \not = 0$ and let $C:=\{x,\xi_1^2t-\xi_3y^2,\xi_2\xi_3v-\xi_4zt,\xi_2\xi_1^4w-\xi_5zy^4 \}$. It is clear that $I_C$ defines a 1-dimensional projective scheme. Moreover, 
	\[
	z^4+z^2f_{6}(y,t,v)+zf_{9}(y,t,v)+f_{12}(y,t,v) \not \in (\xi_1^2t-\xi_3y^2,\xi_2\xi_3v-\xi_4zt,\xi_2\xi_1^4w-\xi_5zy^4).
	\]
	Hence $I_X \not \subset I_C$ and $C$ is a good isolating set of degree 11.
	\end{itemize}

\end{proof}


\begin{center}
\begin{table}
{\small
\begin{tabular}{cccc}
\toprule
  Family        & Isolating Threshold & Isolating Curve     & Isolating Divisor          \\
\midrule 
\midrule 
    93       & $5$    &  $\{x, \xi_2y-\xi_1z,\xi_2\xi_3v-\xi_4zt, \xi_2\xi_3w-\xi_5zt \}$ & $5A$   \\
\midrule
  95&         $\frac{15}{2}$   & $\{x,\xi_3y^2-\xi_1^2t,\xi_1\xi_2v-\xi_4yz, \xi_4w-\xi_5v \}$   & $5A$  \\
\midrule
  96&  $\frac{35}{4}$ &     $\{x, \xi_1^2t-\xi_3y^2,\xi_1\xi_2v-\xi_4yz,\xi_1^2\xi_2w-\xi_5y^2z \}$       & $8A$  \\
\midrule
 97 &  $9$ &  $\{x,\xi_2y-\xi_1z,\xi_1\xi_4w-\xi_5vy,\xi_1\xi_3v-\xi_4ty  \}$  & $9A$    \\
\midrule
 98 & $9$ &  $\{x,\xi_1^2t-\xi_3y^2,\xi_2\xi_1v-\xi_4zy,\xi_2^3w-\xi_5z^3 \}$  &   $9A$ \\
\midrule
 100 &  $11$ &  $\{x,\xi_1^2t-\xi_3y^2,\xi_2\xi_3v-\xi_4zt,\xi_2\xi_1^4w-\xi_5zy^4 \}$  &   $11A$\\
\bottomrule
\end{tabular}\par}
\caption{Table containing isolating curves for smooth points and the corresponding isolating divisor class for each of the families treated by Corollary \ref{cor:smthpointexcl} except family 99. The isolated sets presented assume the smooth point satisfy the conditions presented in the proof.}
\label{tab:tabisol}
\end{table} 
\end{center}


\subsection{Exclusion of singular points}
\subsubsection{Bad links}

Let $X$ be a quasimooth member of families 102, 105 or 109 with homogenous variables $x,\,y,\,z,\,t,\,v,\,w$. After a change of variables and assuming that $z^3t \not \in g$ for $X \in \{102, 105\}$ each can be written as 
\begin{align*}
             zw+t^2+f_{d_1}(x,y)&=0    \\
             z^4y+z^2g_{d_2-2a_z}(x,y)+zg_{d_2-a_z}(x,y)+wx_{\mu}+g_{d_2}(x,y,t,v,w)&=0  
\end{align*}
where $f_{d_1},\, g_{d_2-2a_z},\,g_{d_2-a_z} \in (x,y)$ and $x_{\mu}$ is $v$ for family 102 and $t$ otherwise. The point $\mathbf{p_z}$ is a cyclic quotient singularity of type $\frac{1}{3}(1,1,2),\,\frac{1}{3}(1,1,2)$ or $\frac{1}{5}(1,2,3)$. We prove that blowing up $\mathbf{p_z}$ initiates a Sarkisov link whose steps are all in the Mori Category except that the curves contracted by the last step of the program are $K_Y$-trivial, that is, we produce a \emph{bad link}.

\begin{Lem}[{\cite[Lemma.~2.18]{okadaII}}] \label{lem:badlink}
 Let $\varphi \colon Y \rightarrow X$ be an extremal divisorial contraction centred at $\mathbf{p} \in X$ and exceptional divisor $E$. Assume that there are surfaces $S$ and $T$ on $Y$ for which
\begin{enumerate}
	\item $S\sim_{\mathbb{Q}}-aK_Y-dE$ and $T\sim_{\mathbb{Q}} -bK_Y$, with $a,\,b >0$ and $d\geq 0$.
	\item The intersection $\Gamma = S \cap T $ is a 1-cycle supported on numerically equivalent irreducible and reduced curves.
	\item $T\cdot \Gamma \leq 0$
\end{enumerate}
Then $\mathbf{p}$ is not a maximal centre.
\end{Lem}

The point here is that the conditions allow us to determine the Mori cone, $\overline{\NE}(Y)$. Indeed, $\overline{\NE}(Y)=R+Q$ where $R$ is a ray spanned by curves contracted by the Kawamata blowup of $\mathbf{p}$ and $Q$ is a ray spanned by an irreducible component of $\Gamma$. See \cite[Lemma.~8.8]{okadaI} and \cite[Corollary.~5.4.6]{CPR}.

\begin{Cor} \label{cor:102}
Let $X$ be a quasismooth member of families 102, 105 or 109 as above. Then, the point $\mathbf{p_z}$ is not a maximal centre.
\end{Cor}

\begin{proof}
Let $S\sim_{\mathbb{Q}}-K_Y-E$ and $T\sim_{\mathbb{Q}}-K_Y$. Moreover, $\Gamma$ is the proper transform of $(x=y=0)_X$ which is isomorphic to 
 \[
	\begin{cases*} 
                      (zw+t^2=w^2+v^2=0) \subset \mathbb{P}(3,5,7,7) & if $X$ is member of family $102$   \\
											(zw+vt=wt+v^2=0) \subset \mathbb{P}(3,5,7,9) & if $X$ is member of family $105$   \\
                      (zw+t^2=wt+v^2=0) \subset \mathbb{P}(5,9,11,13) & if $X$ is member of family $109$
                 \end{cases*} 
								\]
In the first case, $\Gamma=\Gamma_1+\Gamma_2$ where $\Gamma_i$ are reduced and irreducible and for the second and third cases $\Gamma$ is reduced and irreducible as well. 
The Kawamata blowup centred at a cyclic quotient singularity of index $r$ has irreducible exceptional divisor $E$ with discrepancy $\frac{1}{r}$, that is, we have that $-K_Y = \varphi^*(-K_X)-\frac{1}{r}E$. Moreover, since
\[
\varphi^*(-K_X)^2 \cdot E = \varphi^*(-K_X) \cdot E^2=0
\] 
and $-K_X \sim_{\mathbb{Q}} \iota_X A$, we have
 \[
	T\cdot \Gamma = \iota_X^3A^3-\frac{r+1}{r^3}E^3=\begin{cases*} 
                      8\cdot \frac{1}{21}-4\cdot \frac{1}{3\cdot 2\cdot 1}=-\frac{2}{7} & if $X$ is member of family $102$   \\
			                8\cdot \frac{2}{45}-4\cdot \frac{1}{3\cdot 2\cdot 1}=-\frac{14}{45} & if $X$ is member of family $105$   \\
                      8\cdot \frac{1}{65}-6\cdot \frac{1}{5\cdot 3\cdot 2}=-\frac{1}{13} & if $X$ is member of family $109$
                 \end{cases*} 
								\]
\end{proof}


\begin{Ex} \label{ex:105}
We construct the Sarkisov link from $X$ a quasismooth member of family 105 for which $z^3t \not \in g$ to a non-terminal Fano variety. Let $\mathbf{p_z} \sim \frac{1}{3}(1,1,2) \in X$. Locally, over $\mathbf{p_z}$, the Kawamata blowup is given by 
\[
\bigg(\frac{y^2}{z}:\frac{t^2}{z^3}:\frac{v^2}{z^4}\bigg).
\]
This extends uniquely to a birational morphism $\varphi \colon Y \rightarrow X$ with global sections
\begin{align*}
y &\in H^0(Y,-K_Y) \\
x & \in H^0\big(Y,-\frac{1}{2}K_Y-\frac{1}{2}E\big)\\
x_{\mu} &\in H^0\big(Y,-\frac{a_{\mu}}{2}K_Y+\frac{1}{2}E\big), \quad \text{for}\,\,x_{\mu} \not = x,\,y
\end{align*}
where $\wt(x_{\mu})=a_{\mu}$. These can be seen as sections of the ambient space of $Y$, a rank 2 toric variety isomorphic to 
\[
\begin{array}{cccc|ccccccc}
             &       & u  & z &   t & v & w & y & x & \\
\actL{T_1}   &  \lBr &  0 & 3 &   5 & 7 & 9 & 2 & 2 &   \actR{.}\\
             &       & 1 & 2 &   3 & 4 & 5 & 1 & 0 &  
\end{array}
\]
Using vGIT it is easy to see that $Y$ follows the 2-ray game on $T$ and therefore $\overline{\Mov}(Y)=\langle(3,2),(2,1)\rangle$ and $-K_Y \in \partial \overline{\Mov}(Y)$. Indeed, playing the 2-ray game on $Y$ generates a a diagram
\[
        \begin{tikzcd}[ampersand replacement=\&]
             Y   \ar[rr, dashed, "\displaystyle{(-5,-1,1,6)}" ] \ar[d, swap,"\displaystyle{\varphi}" ] \&{} \& Y_2   \ar[rr,  "\displaystyle{\simeq}" ] \& {}  \&Y_3 \ar[rr, dashed, "\displaystyle{(-9,-1,1,10)}" ]\& {} \&Y_4 \ar[d, "\displaystyle{\varphi'}" ]\\
						\mathbf{p_z} \in X \&{} \& {} \&  \& {} \& {} \& \mathbf{p_y} \in X'
        \end{tikzcd}
    \] 

 The last map is a divisorial contraction  $\varphi' \colon (x=0) \subset Y_4 \rightarrow \mathbf{p_y} \in Z_{7,8} \subset \mathbb{P}(1,1,2,3,4,5)$ with homogeneous variables $y,\,u,\,z,\,t,\,v,\,w$ given by the anticanonical sections
\[
\bigoplus_{m \geq 1} H^0(Y,\mathcal{O}_Y(-mK_Y)).
\]
In particular, $\varphi'^*(-K_{Z_{7,8}})=-K_Y$ is a crepant resolution of
\[
\mathbf{p_y} \sim \frac{1}{2}(1,1,1,1;2)
\]
which is a canonical singularity. The diagram above is therefore not a Sarkisov link, since it leaves the Mori Category in the last step: Any curve in the exceptional divisor $C \subset (x=0)$ is trivial against $-K_Y$. The reason for this is the presence of anti-flips along the diagram, which worsen the singularities of $-K_Y$.
\end{Ex}

\begin{Prop} \label{prop:108}
Let $X$ be a quasismooth member of family 108 for which $\alpha = 0$ and $y^2t^2,\,y^3tx \not \in g$. Then, $\mathbf{p_y} \sim \frac{1}{3}(1,1,2) \in X$ is not a maximal centre.  
\end{Prop} 

\begin{proof}
After a change of variables we can write $X \subset \mathbb{P}(2,3,4,5,7,11)$ as 
\begin{align*}
                      yw+v^2+f_{14}(x,z,t,v,w)&=0    \\
                      y^4z+y^2g_{10}(x,z)+yg_{13}(x,z,t,v,w)+g_{16}(x,z,t,v,w)&=0  
\end{align*}
Locally, around $\mathbf{p_y}$, the Kawamata blowup is given by $(u,x,t,v) \mapsto (xu^{1/3}:tu^{1/3}:vu^{2/3})$. This extends uniquely to a birational morphism $\varphi \colon Y \rightarrow X$ where the global section $w$ vanishes on $E$ with order $\frac{4}{3}$ and, since  $y^2t^2,\,y^3tx \not \in g$, the global section $z$ vanishes on $E$ with order $\frac{5}{3}$. The proper transform of $X$ via the Kawamata blowup centred at $\mathbf{p_y}$, $Y$, lives in a rank 2 toric variety isomorphic to 
\[
\begin{array}{cccc|ccccccc}
             &       & u  & y &   t & v & w & x & z & \\
\actL{T_1}   &  \lBr &  0 & 3 &   5 & 7 & 11 & 2 & 4 &   \actR{.}\\
             &       & 1 & 2 &   3 & 4 & 6 & 1 & 1 &  
\end{array}
\]
Using vGIT it is easy to see that $Y$ follows the 2-ray game on $T$ and therefore $\overline{\Mov}(Y)=\langle(3,2),(2,1) \rangle$ and $-K_Y \in \partial \overline{\Mov}(Y)$. Indeed, playing the 2-ray game on $Y$ generates a sequence $Y_1\rat \ldots \rat Y_4$ which is given by an anti-flip of type $(-5,-1,1,1,7;2)$ followed by an isomorphism followed by another antiflip of type $(-11,-2,1,13)$. The last map is a divisorial contraction  $\varphi' \colon (z=0) \subset Y_4 \rightarrow \mathbf{p_x} \in Z_{18,20} \subset \mathbb{P}(2,4,5,7,9,13)$ with homogeneous variables $x,\,u,\,y,\,t,\,v,\,w$ given by the anticanonical sections
\[
\bigoplus_{m \geq 1} H^0(Y,\mathcal{O}_Y(-mK_Y)).
\]
In particular, $\varphi'^*(-K_{Z_{18,20}})=-K_Y$, and 
\[
\mathbf{p_x} \sim \frac{1}{2}(1,1,1,1;2)
\]
is a canonical singularity.
\end{proof}

A similar proof to Proposition \ref{prop:108} or Corollary \ref{cor:102} can be done for the case $\alpha  = 1$, $y^3t \not \in f,\, y^2t^2 \not \in g$ of family 108. However, we use a different method due to Cheltsov and Park.

\begin{Lem}[{{\cite[Lemma~3.2.8]{chel}}}] \label{lem:exclsing}
Let $\varphi \colon Y \rightarrow X$ be an extremal divisorial contraction with exceptional divisor $E$. Suppose that there are infinitely many irreducible and reduced curves $C_{\lambda} \subset Y$ such that $-K_Y \cdot C_{\lambda} \leq 0$ and $E\cdot C_{\lambda} >0$. Then $\varphi$ is not a maximal singularity.
\end{Lem}


\begin{Cor} \label{cor:excl108}
Let $X$ be a quasismooth member of family 108 where $\alpha  = 1$, $y^3t \not \in f,\, y^2t^2 \not \in g$. Then, the point $ \mathbf{p_y} \sim \frac{1}{3}(1,1,2)$ is not a maximal centre.
\end{Cor}

\begin{proof}
We have $-K_Y \sim \mathcal{O}\begin{psmallmatrix}2\\1\end{psmallmatrix}$. Let $E=\{u=0\}$ and $D_{x_i} = \{x_i=0 \}$. Then, $
-3K_Y \sim 2D_y-E$. 
Take the $1$-dimensional family of numerically equivalent curves in $Y$:
\[
C_{\lambda, \mu}= \{ \lambda x^2 + \mu z = \widetilde{f}=\widetilde{g}=v=0\}.
\]
We first compute $C_{0,1}\cdot E$: On $E$ we are forced to take $y=1$, due to the irrelevant ideal of $T_1$. The toric variety $T_1$ above is isomorphic to  
\[
\begin{array}{cccc|ccccccc}
             &       & u  & y &   t & w & x & z & v & \\
\actL{\widetilde{T_1}}   &  \lBr &  2 & 1 &   1 & 1 & 0 & 0 & -1 &   \actR{.}\\
             &       & -3 & 0 &   1 & 4 & 1 & 2 & 3 &  
\end{array}
\]
via $\big(\begin{smallmatrix}
  -1 & 2\\
  2 & -3
\end{smallmatrix}\big) \cdot T_1$. So $C_{0,1} \cdot E$ is the number of solutions of 
\[
\{\widetilde{f}(0,1,t,w,x,0,0)=\widetilde{g}(0,1,t,w,x,0,0)=0 \} \subset \mathbb{P}(1_t,4_w,1_x). 
\]

In fact, we have $\widetilde{f}=w+\alpha_1t^2x^2+\alpha_2 tx^3$ and $\widetilde{g}=tw+\beta_1wx + \beta_2x^5+\beta_3tx^4+\beta_4t^2x^3$ in $\mathbb{P}(1_t,4_w,1_x)$. We can write $w=-(\alpha_1t^2x^2+\alpha_2 tx^3)$ globally and get  
	\[
	x^2(-\alpha_1 t^3 + (-\alpha_2-\beta_1 \alpha_1 + \beta_4)t^2x+(\beta_3-\beta_1\alpha_2)tx^2+\beta_2x^3)=0.
	\]
	
	If $x^2=0$, then we have one solution with multiplicity 2 which is the coordinate point $\mathbf{p_t}$. The other solutions correspond to 
\[
	-\alpha_1 t^3 + (-\alpha_2-\beta_1 \alpha_1 + \beta_4)t^2x+(\beta_3-\beta_1\alpha_2)tx^2+\beta_2x^3=0.
	\]
	which gives us 3 points counted with multiplicity. Therefore, $C_{0,1}\cdot E = 5$.
	
	\vspace{0.5cm}
	
	I now compute $D_y\cdot C_{0,1}$: In the same way as before, $y=0$ forces $u=1$ due to the irrelevant ideal of $T_1$. From $T_1$ we can see that $D_y\cdot C_{0,1}$ is the number of solutions to 
	\[
	\{\widetilde{f}(1,0,t,w,x,0,0)=\widetilde{g}(1,0,t,w,x,0,0)=0 \} \subset \mathbb{P}(5_t,11_w,2_x).
	\]
	
	where $\widetilde{f}$ and $\widetilde{g}$ have degrees 14 and 16 respectively. Explicitly,
	\[
	\{x^2(\alpha_1t^2+\alpha_3x^5)=tw+\beta_4t^2x^3+\beta_5x^8=0 \} \subset \mathbb{P}(5_t,11_w,2_x).
	\]
	If $x^2=0$, then $t=0$ or $w=0$. These contribute with $\frac{2}{11}$ and $\frac{2}{5}$, respectively, giving $\frac{32}{55}$. The other points are the solutions to $\alpha_1t^2+\alpha_3x^5=tw+\beta_4t^2x^3+\beta_5x^8=0$ which is $\frac{10\cdot 16}{5 \cdot 11 \cdot 2}$. Summing both we conclude $D_y\cdot C_{0,1} = \frac{112}{55}$. 
		Finally,
		\[
		-3K_Y\cdot C_{0,1} = 2D_y\cdot C_{0,1}-E\cdot C_{0,1} = 2\cdot \frac{112}{55} - 5 < 0
		\]
	
	The conclusion follows by Lemma \ref{lem:exclsing}. 
	\end{proof}

\begin{Prop} \label{prop:weak115}
Let $X$ be a quasismooth member of family 115. Then $\frac{1}{2}(1,1,1) \in X$ is not a maximal centre.
\end{Prop}

\begin{proof}
Any quasismooth member $X$ of family 115 can be written as
\begin{align*}
wx+v^2+z^3+f(x,y,z,t,v) &=0\\
wz+t^3+ztv+g(x,y,z,t,v) &=0
\end{align*}
inside $\mathbb{P}:=\mathbb{P}(1,3,4,5,6,11)$ with homogenous variables $x,\,y,\,z,\,t,\,v,\,w$. $X$ has a single point $\mathbf{p} \sim \frac{1}{2}(1,1,1)$ which is the intersection of the line $\Gamma := \mathbb{P}^1 \colon (x=y=t=w=0)$ with $X$. We blowup $\Gamma \subset \mathbb{P}$ and consider the restriction to $X$. 

Hence, this restriction $\varphi \colon Y \rightarrow X$ is centred at $\mathbf{p}$ and $\varphi$ is its Kawamata blowup. In this case, $-K_Y \sim \mathcal{O}_Y(3,1)$. Let $x_{\mu}$ be $x,\,y$ or $t$ and $\wt(x_{\mu})=a_{\mu}$ the respective weight. We have,
\[
x_{\mu} \in H^0(Y,a_{\mu}\varphi^*(A)-\frac{1}{2}E)=H^0(Y,-\frac{a_{\mu}}{3}K_Y+\frac{a_{\mu}-3}{6}E).
\]
Locally, around $\mathbf{p}$, the sections $z$ and $v$ do not vanish at $E$. Let $x_{\mu'}$ be $z$ or $v$ and $a_{\mu}$ the respective weight. Then,
\[
x_{\mu'} \in H^0(Y,a_{\mu'}\varphi^*(A))=H^0(Y,-\frac{a_{\mu'}}{3}K_Y+\frac{a_{\mu'}}{6}E).
\]
Finally, since $zvt \in g$ by quasismoothness of $X$, it follows that $w$ vanishes at $E$ with order exactly $\frac{1}{2}$ and $w$ is a section in $Y$ just like $x_{\mu}$. Hence, $Y$ is 
\begin{align*}
\widetilde{f} \colon wxu+v^2+z^3+uf(u,x,y,z,t,v) &=0\\
\widetilde{g} \colon wz+t^3u+ztv+ug(u,x,y,z,t,v) &=0
\end{align*}
inside 
\[
\begin{array}{ccccc|cccccc}
             &       & u  & z &   v & w & t & y & x & \\
\actL{T}   &  \lBr &  0 & 4 &   6 & 11 & 5 & 3 & 1 &   \actR{.}\\
             &       & 1 & 2 &   3 & 5 & 2 & 1 & 0 &  
\end{array}
\]
Hence, $\widetilde{f},\, \widetilde{g} \in (u,z,v^2)$. We can therefore write
\[
\begin{pmatrix}
           f^* \\
           g^* \\
					 \end{pmatrix} = 
\left(
  \begin{array}{ccc}
    wx+f(u,x,y,z,t,v) & z^2 & 1\\
    t^3+g(u,x,y,z,t,v) & tv+w & 0\\
  \end{array}
\right)
\begin{pmatrix}
           u \\
           z \\
          v^2 \\
         \end{pmatrix}.
\]
By Cramer's Rule we have 
\[
(\bigwedge^2 Q)_i z_j-(\bigwedge^2 Q)_j z_i \in (\widetilde{f},\widetilde{g})
\]
where $z_i$ is the $i^{th}$ element of $\{u,z,v^2\}$. Therefore, the ratio
\[
\eta:=\frac{(\bigwedge^2 Q)_i}{z_i} \in H^0\bigg(Y,\frac{-11K_Y+E}{3}\bigg),\quad 1\leq i\leq 3
\]
is a well-defined section of $Y$. Explicilty this is
\begin{align*}
\eta&:=\frac{-tv-w}{u}\\
&=-\frac{t^3+g(u,x,y,z,t,v)}{z}\\
&=\frac{(wx+f(u,x,y,z,t,v))(tv+w)-z^2(t^3+g(u,x,y,z,t,v))}{v^2} \in H^0\bigg(Y,\frac{-11K_Y+E}{3}\bigg).
\end{align*}
Hence, $w$ can be globally eliminated and $Y$ is isomorphic to 
\[
Y^{\eta} \colon (\widetilde{f}=\widetilde{g}=\eta z_i - (\bigwedge^2 Q)_i=0, \,\, 2 \leq i \leq 3)
\]
inside 
\[
\begin{array}{ccccc|cccccc}
             &       & u  & z &   v  & t & \eta & y & x & \\
\actL{T^{\eta}}   &  \lBr &  0 & 4 &   6  & 5 & 11 & 3 & 1 &   \actR{.}\\
             &       & 1 & 2 &   3  & 2 & 4 & 1 & 0 &  
\end{array}
\]
Crossing the $t$-wall is an isomorphism, whereas crossing the $\eta$-wall induces a small contraction on $Y^{\eta}$; It contracts an irreducible curve $C^-$ isomorphic to
\[
C^- \colon (t^3-z=z^3+v^2=0) \subset \mathbb{P}(11_u,6_z,9_v,2_t).
\] 
The curve $C^-$ is then replaced by $C^+ \simeq \mathbb{P}(1_y,4_x)$. 

We claim that this operation is an anti-flip. Since $D_v = -2K_Y+E$ we have,
\[
-K_Y \cdot C^{-} = \frac{1}{2}(D_v \cdot C^{-} - E\cdot C^{-}) = \frac{1}{2}\bigg(\frac{9}{11}-1\bigg)=-\frac{1}{11}.
\]
On the other hand, $K_{Y^{+}} \cdot C^+=-\frac{1}{4}$.

Finally, the last contraction is $\varphi' \colon {T'}^{\eta} \rightarrow \mathbb{P}(1,1,2,2,3,4)$ with homogeneous variables $u,\,y,\,z,\,t,\,v,\,\eta$ where ${T'}^{\eta}$ is
\[
\begin{array}{ccccccc|cccc}
             &       & u  & z &   v  & t & \eta & y & x & \\
\actL{{T'}^{\eta}}   &  \lBr &  3 & 2 &   3  & 1 & 1 & 0 & -1 &   \actR{.}\\
             &       & 1 & 2 &   3  & 2 & 4 & 1 & 0 &  
\end{array}
\]
 which is given by
\[
(u:z:v:t:\eta:y:x) \mapsto (ux^3:zx^2:vx^3:tx:\eta x:y).
\]
Restricting to ${Y'}^{\eta}$ this is a map to $Z\subset \mathbb{P}(1,1,2,2,3,4)$ with equations
\begin{align*}
u^6 + u^2y^4 +  u z^2 y + u v y^2 + z^3 - u v t - u^2 \eta + v^2&=0\\
    u^6 -  u^5 y -  u y^5 - z^2 y^2 -  v y^3 + t^3 +  v t y +  u \eta y -
        z \eta&=0
\end{align*}
which contracts the divisor $x=0$ to the point $\mathbf{p_y} \in Z$. Notice that $\varphi'^*(-K_Z) \sim -K_{{Y'}^{\eta}}$ and hence $\mathbf{p_z}$ is a canonical singularity. More than that, the line $\Gamma \colon (u=z=v=t=0) \subset Z$ is a line of singularities.

We have seen that the 2-ray game on $T$ restricts to a 2-ray game on $Y$. Hence $\Mov(T)=\Mov(Y)$ and, in particular, $-K_Y \in \partial \Mov(Y)$ so that we move out of the Mori category after the last contraction is performed. This link is therefore a bad link and $\mathbf{p} \sim \frac{1}{2}(1,1,1) \in X$ is not a maximal centre.
\end{proof}

\begin{Cor} \label{cor:124}
Let $X$ be a quasismooth member of families 124. Then, the cyclic quotient singularity $\frac{1}{3}(1,1,2)$ is not a maximal centre.
\end{Cor}

\begin{proof}
The proof is identical to the one of Corollary \ref{cor:102}. Let $\varphi \colon Y \rightarrow X$ be the Kawamata blowup of $X$ centred at $\frac{1}{3}(1,1,2)$ with homogeneous variables $x,\,y,\,t$. Let $T\sim_{\mathbb{Q}}-K_Y$ and $S\sim_{\mathbb{Q}}-\frac{5}{4}K_Y-\frac{1}{4}E$ and $\Gamma = T \cap S$. Then $\Gamma$ is the proper transform of $(x=y=0)|_X$ which is isomorphic to the reduced and irreducible curve $(v^2+z^3+wt=wv+t^2z=0) \subset \mathbb{P}(6,7,9,11)$ with homogeneous variables $z,\,t,\,v,\,w$. We have 
\begin{align*}
T\cdot \Gamma &= \frac{1}{4}(-K_Y)^2 \cdot (-5K_Y-E) \\
&=\frac{1}{4}(5\cdot \iota_X^3 A^3 - 8 \cdot \frac{1}{3^3}E^3) \\
&= 2\bigg(\frac{5}{231}-\frac{1}{6}\bigg) < 0.
\end{align*}
\end{proof}

\subsubsection{The case $(-K_Y)^3 \leq 0$}

The following lemma gives a way to produce a nef divisor on $Y$ to be used with the test class method to exclude singular points satisfying $(-K_Y)^3 \leq 0$. This is adapted from \cite[Lemma.~6.6]{okadaII}. 

\begin{Lem} \label{lem:isolsing}
Let $X$ be a Fano 3-fold and $\mathbf{p} \in X$ a cyclic quotient singularity. Let $\varphi \colon E \subset Y \rightarrow X \ni \mathbf{p}$ be the Kawamata blowup centred at $\mathbf{p}$. Suppose there are prime divisors $D_1,\ldots,D_k$ in $X$ such that 
\begin{enumerate}
	\item The intersection $\cap_{i}D_i$ isolates $\mathbf{p}$.
	\item Each proper transform $\widetilde{D_i}$ of $D_i$ via $\varphi$ is $\mathbb{Q}$-linearly equivalent to $-b_iK_Y+e_iE$ for some $b_i>0$ and $e_i \geq 0$.
	\item We have $c:=\max\{e_i/b_i\}\leq 1/r$.
\end{enumerate}
Then the divisor $M := -K_Y+cE $ is nef. 
\end{Lem}

\begin{Cor} \label{cor:exc111}
Let $X$ be a quasismooth member of family 111. Then the point $\mathbf{p} \sim \frac{1}{3}(1,1,2)$ is not a maximal centre.
\end{Cor}

\begin{proof}
Consider the Kawamata blowup $\varphi \colon Y \rightarrow X$ centred at $\mathbf{p}\sim \frac{1}{3}(1,1,2)$. Locally $\mathbf{p}$ has tangent variables $\{x,y,t\}$ of weights $(1,1,2)$ respectively and therefore,
\[
x \in H^0(Y,-K_Y), \quad y \in H^0\bigg(Y,-\frac{5}{2}K_Y+\frac{1}{2}E\bigg), \quad t \in H^0\bigg(Y,-\frac{7}{2}K_Y+\frac{1}{2}E\bigg).
\]
Moreover, $w$ vanishes along $E$ with order $\frac{4}{3}$ by looking at the equations defining $Y$ and so $w \in H^0\bigg(Y,-\frac{11}{2}K_Y+\frac{1}{2}E\bigg)$.
The closed subvariety $(x=t=w=0)_X$ contains $\mathbf{p}$ but no curves passing through it, that is, $\{x,t,w\}$ is an isolating set for $\mathbf{p}$. Let $c=\max \{0,\frac{1}{7},\frac{1}{11}\}=1/7<1/3=a_E(K_X)$. Then, by Lemma \ref{lem:isolsing}, $M:=-K_Y+\frac{1}{7}E$ is a nef divisor. Moreover,
\begin{align*}
M\cdot (-K_Y)^2 &= \bigg(-K_Y+\frac{1}{7}E\bigg)\cdot (-K_Y)^2 \\
&= \bigg(\varphi^*(-K_X)-\frac{4}{7}\cdot\frac{1}{3}E\bigg)\cdot \bigg(\varphi^*(-K_X)^2-\frac{2}{3}\varphi^*(-K_X)\cdot E+\frac{1}{3^2}E^2\bigg)\\
&=\varphi^*(-K_X)^3-\frac{4}{7}\cdot\frac{1}{3^3}E^3\\
&=8\cdot\frac{2}{231}-\frac{4}{7}\cdot\frac{1}{6}\\
&=-\frac{2}{77}.
\end{align*}
By the test class method, see \cite[Lemma~5.2.1,\, Corollary~5.2.3]{CPR}, the point $\mathbf{p}$ is not a maximal centre.
\end{proof}


	
	
	
	

\newgeometry{left=1cm,right=0.5cm,top=0.5cm,bottom=1cm}

\begin{landscape}
\section{The Big Table}  \label{sect:table}

We explain how to read the Big Table. Recall that $X$ is a Fano 3-fold deformation family in $I$ and that our only assumption on $X$ is quasismoothness, unless explicitly stated. Each entry of the big table consists of information about $X \in I \setminus \{86\}$ and its birational geometry. On the first line, if the isolating threshold appears, see Section \ref{subsec:smooth}, then $X \in I_S \subset I$. Otherwise, $X \in I_{nS}$. We show typical equations for $X$ and afterwards each entry has five columns:
\begin{itemize}
	\item The first column is the reference for the results found in that row.
	\item The second column is the cyclic quotient singularity we blow up written in the form $\frac{1}{r}(1,a,r-a)$, see Theorem \ref{thm:cqs}.
	\item The third column is either empty or shows a small $\mathbb{Q}$-factorial modification. A typical entry is of the form $(a_1,a_2,a_3,-b_1,-b_2;d)$ where $a_i,\,b_i,\,d \geq 0$ and it follows the notation introduced in \cite{brown}.	If it is empty and the link exists, then it denotes an isomorphism.
	\item the fourth column is the weighted blowup corresponding to the last divisorial contraction to the new birational model. This is omitted if the centre of the contraction is a curve or a cyclic quotient singularity.
	\item The last column is the centre of the contraction together with the new model. This column is empty if and only if the cyclic quotient singularity in the second column has been excluded as a maximal centre.
\end{itemize}

{\normalsize
\begin{longtable}[l]{llcrr}
\specialrule{.2em}{.1em}{.1em} 
\multicolumn{3}{l}{\textbf{Family 87}: $X_{4,4} \subset \mathbb{P}^5(1,1,1,2,2,3), \quad \iota_X=2$}\\
    \\
		\multicolumn{3}{l}{$wx+f_4(x,y,z,t,v)=0$}\\
		\multicolumn{3}{l}{$wy+g_4(x,y,z,t,v)=0$}\\
		\\
		Prop. \ref{prop:elemcb} &$\frac{1}{3}(1,1,2)$& & & $Y/\mathbb{P}^2$\\
    \specialrule{.1em}{.05em}{.1em} 
{}&&\\
		\specialrule{.1em}{.05em}{.1em} 
		\multicolumn{3}{l}{\textbf{Family 88}: $X_{4,6} \subset \mathbb{P}^5(1,1,2,2,3,3), \quad \iota_X=2$}\\
    \\
		\multicolumn{3}{l}{$(w+v)x+vy+f_4(x,y,z,t,v)=0$} \\
		\multicolumn{3}{l}{$wv+g_6(x,y,z,t,v)=0$}\\
		\\
    Thm. \ref{thm:delPezzoMain}&$2\times\frac{1}{3}(1,1,2)$&  $(-3,-1,-1,1,1;-2)$& &$ dP_3$\\
    \specialrule{.1em}{.05em}{.1em} 
{}&&\\
		\specialrule{.1em}{.05em}{.1em} 
		\multicolumn{3}{l}{\textbf{Family 89}: $X_{6,6} \subset \mathbb{P}^5(1,1,2,2,3,5), \quad \iota_X=2$}\\
    \\
		\multicolumn{3}{l}{$wx+v^2+f_6(x,y,z,t,v)=0$}\\
		\multicolumn{3}{l}{$wy+g_6(x,y,z,t,v)=0$}\\
    {}&&\\
		Prop. \ref{prop:89}&$\frac{1}{5}(1,1,4)$& $(-4,-1-1,1,1;-3)$& & $dP_3$\\
    \specialrule{.1em}{0.05em}{.1em} 
{}&&\\
		\specialrule{.1em}{0.05em}{.1em}  
		\multicolumn{3}{l}{\textbf{Family 90}: $X_{6,8} \subset \mathbb{P}^5(1,1,2,3,4,5), \quad \iota_X=2$}\\
    \\
		\multicolumn{3}{l}{$wx+t^2+f_6(x,y,z,t,v)=0$}\\
		\multicolumn{3}{l}{$wt+g_8(x,y,z,t,v)=0$}\\
  {}&&\\
		Thm. \ref{thm:delPezzoMain}&$\frac{1}{5}(1,2,3)$& $(-1,-2,1,1)$& & $dP_2$\\
    \specialrule{.1em}{0.05em}{.1em}  
{}&&\\
		\specialrule{.1em}{0.05em}{.1em} 
		\multicolumn{3}{l}{\textbf{Family 91}: $X_{6,6} \subset \mathbb{P}^5(1,2,2,3,3,3), \quad \iota_X=2$}\\
    \\
		\multicolumn{3}{l}{$wv+f_6(x,y,z,t,v)=0$}\\
		\multicolumn{3}{l}{$wt+g_6(x,y,z,t,v)=0$}\\
{}&&\\
		Thm. \ref{thm:cod2curve} & $4\times \frac{1}{3}(1,1,2)$&  & &$\mathbb{P}^1 \subset Z_{3,3}\subset \mathbb{P}(1,1,1,1,1,2)$\\
		\specialrule{.1em}{0.05em}{.1em} 
{}&&\\
		\specialrule{.1em}{0.05em}{.1em} 
		\multicolumn{3}{l}{\textbf{Family 92}: $X_{6,8} \subset \mathbb{P}(1,2,2,3,3,5),\quad \iota_X=2, \quad \frac{4}{\iota_X^2A^3}=\frac{15}{4}$}\\
    \\
		\multicolumn{3}{l}{$wx+ vt + f_6(x,y,z,t,v)=0$}\\
		\multicolumn{3}{l}{$w(\alpha v +t)+v^2z+vg_5(x,y,z,t)+g_8(x,y,z,t,v)=0$}\\
    {}&&\\
    Thm. \ref{thm:divtocurve}&$\frac{1}{5}(1,1,4)$& & & $\mathbb{P}^1 \subset Z_4 \subset \mathbb{P}^4$\\
		Thm. \ref{thm:cod2caseI}&$2 \times \frac{1}{3}(1,1,2)$, $\alpha = 1$ & $(5,1,1,-1,-1;3)$ & $(3,1,1,1)$&$cA_3 \in Z_{3,4} \subset \mathbb{P}^5(1,1,1,1,2,2)$\\
    {}&&\\
		\multicolumn{3}{l}{In the case $\alpha = 0$, blowing up $\mathbf{p_v}$ initiates a Sarkisov Link to a del Pezzo fibration of degree 2.} \\
		\specialrule{.1em}{0.05em}{.1em} 
{}&&\\
		\specialrule{.1em}{0.05em}{.1em} 
		\multicolumn{3}{l}{\textbf{Family 93}: $X_{6,10} \subset \mathbb{P}(1,2,2,3,5,5),\quad \iota_X=2, \quad \frac{4}{\iota_X^2A^3}=5$}\\
    \\
		\multicolumn{3}{l}{$(w+v)x+t^2+f_6(x,y,z,t,v)=0$}\\
		\multicolumn{3}{l}{$wv+g_{10}(x,y,z,t,v)=0$}\\
    {}&&\\
    Thm. \ref{thm:hypIIfull}&$2\times \frac{1}{5}(1,1,4)$& $(-5,-1,-1,1,2;-3)$ & $(4,1,1,1)$ &$cA_4  \in Z_5 \subset \mathbb{P}(1,1,1,1,2)$\\
{}&&\\
		\specialrule{.1em}{0.05em}{.1em} 
		\multicolumn{3}{l}{\textbf{Family 94}: $X_{8,10} \subset \mathbb{P}(1,2,2,3,5,7),\quad \iota_X=2, \quad \frac{4}{\iota_X^2A^3}=\frac{21}{4}$}\\
    \\
		\multicolumn{3}{l}{$wx + (1-\alpha)t^2z+t(\alpha v + f_{5}(x,y,z))+ \beta z^4+f_8(x,y,z,t,v)=0$}\\
		\multicolumn{3}{l}{$wt+v^2+g_{10}(x,y,z,t,v)=0$}\\
    {}&&\\
    Thm. \ref{thm:hypIfull}&$\frac{1}{7}(1,1,6)$& $(-6,-1,-1,1,2;-4)$ & $(3,1,1,1)$ &$cA_3 \in Z_5 \subset \mathbb{P}(1,1,1,1,2)$\\
		Thm. \ref{thm:cod2caseI}&$\frac{1}{3}(1,1,2)$, $\alpha  = 1$& $(-7,-1,-1,1,3;-4)$ & & $\frac{1}{2}(1,1,1)  \in Z_{4,5} \subset \mathbb{P}(1,1,1,2,2,3)$\\
		Ex. \ref{ex:pic2}&$\frac{1}{3}(1,1,2)$, $\alpha  = 0,\, \beta = 0$& & &$\rho(X_{8,10}) \geq 2$    \\
	  Ex. \ref{ex:pic2}&$\frac{1}{3}(1,1,2)$, $\alpha  = 0,\, \beta = 1$& $(-7,-1,1,8)$ & &$dP_1$  \\
    \specialrule{.1em}{0.05em}{.1em} 
		{}&&\\
			\specialrule{.1em}{0.05em}{.1em} 
		\multicolumn{3}{l}{\textbf{Family 95}: $X_{8,10} \subset \mathbb{P}(1,2,3,4,5,5),\quad \iota_X=2, \quad \frac{4}{\iota_X^2A^3}=\frac{15}{2}$}\\
     {}&&\\
		\multicolumn{3}{l}{$ (w+v)z+f_8(x,y,z,t,v)=0$}\\
		\multicolumn{3}{l}{$wv+(1-\alpha)z^3x+z^2(\alpha t+g_{4}(x,y))+zg_{7}(x,y,t,v,w)+g_{10}(x,y,z,t,v)=0$}\\
    {}&&\\
    Thm. \ref{thm:cod2caseI}&$2\times \frac{1}{5}(1,2,3)$& $(-5,-1,-2,1,2;-4)$ & $(3,1,2,1)$ &$cA_4 \in Z_{4,5} \subset \mathbb{P}(1,1,1,2,2,3)$\\
		Thm. \ref{thm:HypCaseIIIExcI}& $\frac{1}{3}(1,2,2),\, \alpha=1$& $(-1,-1,1,1)$ & $\frac{1}{2}(5,5,1,2)$ &$cA/2 \in Z_{4,6} \subset \mathbb{P}(1,1,2,2,2,3)$\\
    {}&&\\
						\multicolumn{3}{l}{The same link is obtained if $\alpha\not = 0$.}\\
						\specialrule{.1em}{0.05em}{.1em} 
			{}&&\\
		\specialrule{.1em}{0.05em}{.1em} 
		\multicolumn{3}{l}{\textbf{Family 96}: $X_{8,12} \subset \mathbb{P}(1,2,3,4,5,7),\quad \iota_X=2, \quad \frac{4}{\iota_X^2A^3}=\frac{35}{4}$}\\
    \\
		\multicolumn{3}{l}{$wx+vz+f_8(x,y,z,t,v)=0$}\\
		\multicolumn{3}{l}{$wv+g_{12}(x,y,z,t,v)=0$}\\
    {}&&\\
   Thm. \ref{thm:hypIIfull}& $\frac{1}{7}(1,5,2)$& $(-5,-1,-2,1,2;-4)$ & $(3,5,1,2,1)$ &$ cD_4  \in Z_{6} \subset \mathbb{P}(1,1,1,2,2)$\\
		Thm. \ref{thm:cod2caseI}&$\frac{1}{5}(1,2,3)$& $(-7,-1,-2,2,3;-4)$ & $(3,1,2,2)$ &$cD_4  \in Z_{4,6} \subset \mathbb{P}(1,1,1,2,3,3)$\\
    \specialrule{.1em}{0.05em}{.1em} 
{}&&\\
\specialrule{.1em}{0.05em}{.1em} 
		\multicolumn{3}{l}{\textbf{Family 97}: $X_{10,14} \subset \mathbb{P}(1,2,2,5,7,9),\quad \iota_X=2, \quad \frac{4}{\iota_X^2A^3}=9$}\\
    \\
		\multicolumn{3}{l}{$wx+t^2+f_{10}(x,y,z,t,v)=0$}\\
		\multicolumn{3}{l}{$wt+v^2+g_{14}(x,y,z,t,v)=0$}\\
    {}&&\\
    Thm. \ref{thm:hypIfull}&$\frac{1}{9}(1,1,8)$ & $(-8,-1,-1,1,3;-5)$ & &$\frac{1}{2}(1,1,1) \in Z_{7} \subset \mathbb{P}(1,1,1,2,3)$\\
    \specialrule{.1em}{0.05em}{.1em} 
{}&&\\
		\specialrule{.1em}{0.05em}{.1em} 
		\multicolumn{3}{l}{\textbf{Family 98}: $X_{10,12} \subset \mathbb{P}(1,2,3,4,5,9),\quad \iota_X=2, \quad \frac{4}{\iota_X^2A^3}=9$}\\
    \\
		\multicolumn{3}{l}{$wx+v^2+(1-\alpha)z^3x+z^2(\alpha t+f_4(x,y))+ zf_7(x,y,t,v)+ f_{10}(x,y,z,t,v)=0$}\\
		\multicolumn{3}{l}{$wz+g_{12}(x,y,z,t,v)=0$}\\
    \\
    Thm. \ref{thm:hypIfull}&$\frac{1}{9}(1,2,7)$& $(-7,-1,-2,1,2;-6)$ & $(3,1,2,1)$ &$cA_4 \in Z_{6} \subset \mathbb{P}(1,1,1,2,2)$\\
		Thm. \ref{thm:HypCaseIIIExcI}& $\frac{1}{3}(1,1,2),\, \alpha = 0, \, v^2y \not \in g$& $3\times(1,1,-1,-1)$ &$\frac{1}{4}(9,7,1,2)$ &$cAx/4\in Z_{6,6} \subset \mathbb{P}(1,1,1,2,3,3)$\\
   \\
		\multicolumn{3}{l}{$v^2y \not \in g$ w.l.o.g.}\\
				\multicolumn{3}{l}{The same link is obtained if $\alpha\not = 0$.}\\
				\specialrule{.1em}{0.05em}{.1em} 
				{}&&\\
		\specialrule{.1em}{0.05em}{.1em} 
		\multicolumn{3}{l}{\textbf{Family 99}: $X_{10,12} \subset \mathbb{P}(1,2,3,5,6,7),\quad \iota_X=2, \quad \frac{4}{\iota_X^2A^3}=21$}\\
    \\
		\multicolumn{3}{l}{$wz+t^2+f_{10}(x,y,z,t,v)=0$}\\
		\multicolumn{3}{l}{$wt+v^2+z^2v+ g_{12}(x,y,z,t,v)=0$}\\
    \\
    Thm. \ref{thm:cod2caseI}&$\frac{1}{7}(1,3,4)$& $(-1,-3,1,2)$ & $(3,1,3,1)$ &$cA_5  \in Z_{5,6} \subset \mathbb{P}(1,1,1,2,3,4)$\\  
		Thm. \ref{thm:cod4}&$2\times \frac{1}{3}(1,1,2)$& $6\times (1,1,-1,-1)$ & &$\frac{1}{10}(1,3,7) \in Z\subset \mathbb{P}(1,1,2,3,3,4,7,10)$ \\   
	\specialrule{.1em}{0.05em}{.1em} 
{}&&\\
		\specialrule{.1em}{0.05em}{.1em} 
		\multicolumn{3}{l}{\textbf{Family 100}: $X_{12,14} \subset \mathbb{P}(1,2,3,4,7,11),\quad \iota_X=2, \quad \frac{4}{\iota_X^2A^3}=11$}\\
    \\
		\multicolumn{3}{l}{$wx+z^4+f_{12}(x,y,z,t,v)=0$}\\
		\multicolumn{3}{l}{$wz+v^2+g_{14}(x,y,z,t,v)=0$}\\
    \\
    Thm. \ref{thm:hypIfull}&$\frac{1}{11}(1,2,9)$& $(-9,-1,-3,2,3;-6)$ & $(3,1,2,2)$ &$cD_4  \in Z_{7} \subset \mathbb{P}(1,1,1,2,3)$\\
    \specialrule{.1em}{0.05em}{.1em} 
		{}&&\\
		\specialrule{.1em}{0.05em}{.1em} 
		\multicolumn{3}{l}{\textbf{Family 101}: $X_{10,12} \subset \mathbb{P}(2,2,3,5,5,7),\quad \iota_X=2, \quad \frac{4}{\iota_X^2A^3}=\frac{35}{2}$}\\
    \\
		\multicolumn{3}{l}{$wz+ vt+f_{10}(x,y,z,t,v)=0$}\\
		\multicolumn{3}{l}{$w(v+t)+z^4+z^2g_6(x,y)+zg_{9}(x,y,t,v,w)+g_{12}(x,y,z,t,v)=0$}\\
    \\
    Thm. \ref{thm:divtocurve}&$\frac{1}{7}(1,1,6)$&   && $\mathbb{P}^1 \subset Z_{6} \subset \mathbb{P}(1,1,1,1,3)$\\  
		Thm. \ref{thm:cod2caseI}&$2\times\frac{1}{5}(1,1,4)$&  $(-7,-1,-1,1,2;-5)$ & $(5,1,1,1)$ &$cA_5 \in Z_{5,6} \subset \mathbb{P}(1,1,1,2,3,4)$\\  	
		\specialrule{.1em}{0.05em}{.1em} 
{}&&\\
		\specialrule{.1em}{0.05em}{.1em} 
		\multicolumn{3}{l}{\textbf{Family 102}: $X_{10,14} \subset \mathbb{P}(2,2,3,5,7,7),\quad \iota_X=2, \quad \frac{4}{\iota_X^2A^3}=21$}\\
    \\
		\multicolumn{3}{l}{$wz+t^2+f_{10}(x,y,z,t,v)=0$}\\
		\multicolumn{3}{l}{$wv+v^2+(1-\alpha) z^4x+\alpha z^3t+z^2g_8(x,y)+zg_{11}(x,y,t,v,w)+g_{14}(x,y,z,t,v)=0$}\\
    \\
		Thm. \ref{thm:hypIIfull}&$2\times \frac{1}{7}(1,1,6)$&  $(-7,-1,-1,1,2;-5)$ & $(6,1,1,1)$ &$cA_6 \in Z_{7} \subset \mathbb{P}(1,1,1,2,3)$\\  
		Thm. \ref{thm:hypexcI}&$\frac{1}{3}(1,1,2)$, $\alpha = 1$ &  $5\times(1,1-,1,-1)$ & $(7,7,1,1)$ & $cA_{13} \in Z_{12} \subset \mathbb{P}(1,1,1,4,6)$\\
	  Cor. \ref{cor:102}&$\frac{1}{3}(1,1,2)$, 	$\alpha = 0$&   & \\
		\specialrule{.1em}{0.05em}{.1em} 
		{}&&\\
	\specialrule{.1em}{0.05em}{.1em} 
		\multicolumn{3}{l}{\textbf{Family 103}: $X_{10,12} \subset \mathbb{P}(2,3,3,4,5,7), \quad \iota_X=2, \quad \frac{4}{\iota_X^2A^3}=21$}\\
    \\
		\multicolumn{3}{l}{$wz+v^2+f_{10}(x,y,z,t,v)=0$}\\
		\multicolumn{3}{l}{$wv+y^3z+ y^2g_6(x,z,t) + yg_{9}(x,z,t,v,w) +g_{12}(x,y,z,t,v)=0$}\\
\\
	Thm. \ref{thm:delPezzoMain}&	$\frac{1}{7}(1,2,5)$& $(-1,-2,1,1)$ & &$dP_1$\\
	Thm. \ref{thm:cod4}&	$4\times \frac{1}{3}(1,1,2)$& $9\times (-1,-1,1,1)$ & &$\frac{1}{9}(1,2,7) \in Z \subset \mathbb{P}(1,2,3,3,4,5,7,9)$\\
		\specialrule{.1em}{0.05em}{.1em} 
	{}&&\\
		\specialrule{.1em}{0.05em}{.1em} 
		\multicolumn{3}{l}{\textbf{Family 104}: $X_{14,16} \subset \mathbb{P}(1,2,5,7,8,9),\quad \iota_X=2, \quad \frac{4}{\iota_X^2A^3}=\frac{45}{2}$}\\
    \\
		\multicolumn{3}{l}{$wz+t^2+f_{14}(x,y,z,t,v)=0$}\\
		\multicolumn{3}{l}{$wt+z^3x+z^2g_6(x,y)+zg_{11}(x,y,t,v)+g_{16}(x,y,z,t,v)=0$}\\
    \\
    Thm. \ref{thm:cod2caseI}&$\frac{1}{9}(1,4,5)$&  $(-1,-4,1,3)$ & $\frac{1}{2}(5,4,1,1)$ &$cAx/2  \in Z_{7,8} \subset \mathbb{P}(1,1,2,3,4,5)$\\  
		Thm. \ref{thm:cod4}&$\frac{1}{5}(1,1,4)$& $2\times (1,1,-1,-1) $ & &$\frac{1}{13}(1,4,9) \in Z \subset \mathbb{P}^7(1,1,3,4,4,5,9,13)$ \\
    \specialrule{.1em}{0.05em}{.1em} 
{}&&\\
		\specialrule{.1em}{0.05em}{.1em} 
		\multicolumn{3}{l}{\textbf{Family 105}: $X_{12,14} \subset \mathbb{P}(2,2,3,5,7,9),\quad \iota_X=2, \quad \frac{4}{\iota_X^2A^3}=\frac{45}{2}$}\\
    \\
		\multicolumn{3}{l}{$wz+(1-\alpha)t^2x+\alpha tv + f_{12}(x,y,t,v)=0$}\\
		\multicolumn{3}{l}{$wt+v^2+(1-\beta)z^4x+\beta z^3t+z^2g_{8}(x,y)+zg_{13}(t,x,y)+g_{14}(x,y,z,t,v)=0$}\\
    \\
    Thm. \ref{thm:hypIfull}&$\frac{1}{9}(1,1,8)$&  $(-8,-1,-1,1,2;-6)$ & $(5,1,1,1)$ &$cA_5  \in Z_{7} \subset \mathbb{P}(1,1,1,2,3)$\\ 
		Thm. \ref{thm:cod2caseI}&$\frac{1}{5}(1,1,4)$, $\alpha \not = 0$&  $(-9,-1,-1,1,3;-6)$ & &$\frac{1}{2}(1,1,1) \in Z_{6,7} \subset \mathbb{P}(1,1,2,3,3,4)$\\ 
		Prop. \ref{prop:fake2}&$\frac{1}{5}(1,1,4)$, $\alpha = 0, \, y^6 \in f$ & $(-9,-1,1,10)$ & &$\mathbf{p_z} \in Z_{6,8}\subset \mathbb{P}(1,1,1,3,4,5)/\boldsymbol{\mu_2}$\\
	  Thm. \ref{thm:hypexcI}&$\frac{1}{3}(1,1,2)$, $\beta \not = 0$&  $6\times (1,1,-1,-1)$ & $\frac{1}{2}(9,8,1,1)$ &$cAx/2 \in Z_{12} \subset \mathbb{P}(1,1,2,4,5)$\\ 
		Cor. \ref{cor:102} &$\frac{1}{3}(1,1,2)$, $\beta = 0$&   & \\ 
    \\
		\multicolumn{3}{l}{In the case of $\frac{1}{5}(1,1,4)$, $\alpha = 0, \, y^6 \not \in f$ a similar situation to Ex. \ref{ex:pic2} happens.}\\
		\specialrule{.1em}{0.05em}{.1em} 
{}&&\\
		\specialrule{.1em}{0.05em}{.1em} 
		\multicolumn{3}{l}{\textbf{Family 106}: $X_{14,18} \subset \mathbb{P}(2,2,3,7,9,11),\quad \iota_X=2, \quad \frac{4}{\iota_X^2A^3}=33$}\\
    \\
		\multicolumn{3}{l}{$wz+t^2+f_{14}(x,y,z,t,v)=0$}\\
		\multicolumn{3}{l}{$wt+v^2+ z^3v+z^2g_{12}(x,y)+zg_{15}(x,y,t,w)+g_{18}(x,y,z,t,v)=0$}\\
    \\
    Thm. \ref{thm:hypIfull}&$\frac{1}{11}(1,1,10)$& $(-10,-1,-1,1,3;-7)$ & &$\frac{1}{2}(1,1,1) \in Z_{9} \subset \mathbb{P}(1,1,2,3,3)$\\
		Thm. \ref{thm:cod4}&$2 \times \frac{1}{3}(1,1,2)$ &  $27\times (1,1,-1,-1)$ &&$\frac{1}{12}(1,1,11)\in Z \subset \mathbb{P}(1,1,6,8,9,10,11,12)$ \\
	\specialrule{.1em}{0.05em}{.1em} 
		{}&&\\
		\specialrule{.1em}{0.05em}{.1em} 
		\multicolumn{3}{l}{\textbf{Family 107}: $X_{12,14} \subset \mathbb{P}(2,3,4,5,7,7),\quad \iota_X=2, \quad \frac{4}{\iota_X^2A^3}=35$}\\
    \\
		\multicolumn{3}{l}{$wt+vt+f_{12}(x,y,z,t,v)=0$}\\
		\multicolumn{3}{l}{$wv+v^2+t^2z+tg_9(x,y,z,v)+ g_{14}(x,y,z,t,v)=0$}\\
    \\
    Thm. \ref{thm:cod2caseI}&$\frac{1}{7}(1,2,5)$ & $(-7,-1,-2,1,2;-6)$ & $(5,1,2,1)$&$cA_6 \in Z_{6,7} \subset \mathbb{P}(1,1,2,2,3,5)$\\
		Thm. \ref{thm:HypCaseIIIExcI}&$\frac{1}{5}(1,1,4)$ & $3\times (1,1,-1,-1)$  & $\frac{1}{2}(7,7,1,2)$& $cA/2 \in Z_{6,10} \subset \mathbb{P}(1,2,2,3,4,5)$\\
				\specialrule{.1em}{0.05em}{.1em} 
		{}&&\\
		\specialrule{.1em}{0.05em}{.1em} 
		\multicolumn{3}{l}{\textbf{Family 108}: $X_{14,16} \subset \mathbb{P}(2,3,4,5,7,11),\quad \iota_X=2, \quad \frac{4}{\iota_X^2A^3}=\frac{165}{4}$}\\
    \\
		\multicolumn{3}{l}{$wy+v^2+t^2z+tf_{9}(x,y,z,v)+f_{14}(x,y,z,t,v)=0$}\\
		\multicolumn{3}{l}{$wt+(1-\alpha)y^4z+y^3(\alpha v+tx)+y^2g_{10}(x,z,t)+yg_{13}(x,y,z,t,v)+g_{16}(x,y,z,t,v)=0$}\\
    \\
    Thm. \ref{thm:hypIfull}&$\frac{1}{11}(1,2,9)$ & $(-9,-1,-2,1,2;-8)$ & $(5,1,2,1)$ &$cA_6  \in Z_{8} \subset \mathbb{P}(1,1,2,2,3)$\\
		Thm. \ref{thm:HypCaseIIIExcI}&$\frac{1}{5}(1,1,4)$, $v^2x \not \in g$ & $4\times (1,1,-1,-1)$ & $\frac{1}{4}(11,9,1,2)$&$cAx/4  \in Z_{10,8} \subset \mathbb{P}(1,2,3,4,4,5)$\\
		Cor. \ref{cor:excl108}&$\frac{1}{3}(1,1,2)$, $\alpha  = 1$, $y^3t \not \in f,\, y^2t^2 \not \in g$ &  & \\
		Prop. \ref{prop:108}&$\frac{1}{3}(1,1,2)$, $\alpha= 0$, $y^3tx, y^2t^2 \not \in g$ &  & \\
    \\
		\multicolumn{3}{l}{$v^2x  \not \in g$ wlog. }\\
		\specialrule{.1em}{0.05em}{.1em} 
{}&&\\
		\specialrule{.1em}{0.05em}{.1em} 
		\multicolumn{3}{l}{\textbf{Family 109}: $X_{18,22} \subset \mathbb{P}(2,2,5,9,11,13),\quad \iota_X=2, \quad \frac{4}{\iota_X^2A^3}=65$}\\
    \\
		\multicolumn{3}{l}{$wz+t^2+f_{18}(x,y,z,t,v)=0$}\\
		\multicolumn{3}{l}{$wt+v^2+z^4y+z^2g_{12}(x,y)+zg_{17}(x,y,t,v)+g_{22}(x,y,t,v)=0$}\\
    \\
    Thm. \ref{thm:hypIfull}&$\frac{1}{13}(1,1,12)$ & $(-12,-1,-1,1,3;-9)$ & &$\frac{1}{2}(1,1,1) \in Z_{11} \subset \mathbb{P}(1,1,2,3,5)$\\
		Cor. \ref{cor:102}& $\frac{1}{5}(1,2,3)$   &  & \\
    \specialrule{.1em}{0.05em}{.1em} 
{}&&\\
		\specialrule{.1em}{0.05em}{.1em} 
		\multicolumn{3}{l}{\textbf{Family 110}: $X_{18,20} \subset \mathbb{P}(2,4,5,7,9,13),\quad \iota_X=2, \quad \frac{4}{\iota_X^2a^3}=91$}\\
    \\
		\multicolumn{3}{l}{$wz+v^2+t^2y+tf_{11}(x,y,z,v)+f_{18}(x,y,z,t,v)=0$}\\
		\multicolumn{3}{l}{$wt+g_{20}(x,y,z,t,v)=0$}\\
    \\
    Thm. \ref{thm:hypIfull}&$\frac{1}{13}(1,2,11)$& $(-11,-1,-2,1,2;-10)$ & $(7,1,2,1)$ &$cA_8 \in Z_{10} \subset \mathbb{P}(1,1,2,2,5)$\\
		Thm. \ref{thm:HypCaseIIIExcI}&$\frac{1}{7}(1,1,6)$, $v^2x \not \in g$ & $5\times (-1,-1,1,1)$ & $\frac{1}{4}(13,11,1,2)$ & $cAx/4 \in Z_{10,14} \subset \mathbb{P}(1,2,4,5,6,7)$\\
    \\
		\multicolumn{3}{l}{$v^2x \not \in g$ wlog.}\\
		\specialrule{.1em}{0.05em}{.1em} 
{}&&\\
		\specialrule{.1em}{0.05em}{.1em} 
		\multicolumn{3}{l}{\textbf{Family 111}: $X_{18,20} \subset \mathbb{P}(2,5,6,7,9,11),\quad \iota_X=2, \quad \frac{4}{\iota_X^2A^3} = \frac{231}{2}$}\\
    \\
		\multicolumn{3}{l}{$wt+v^2+z^3+f_{18}(x,y,z,t,v)=0$}\\
		\multicolumn{3}{l}{$wv+t^2z+tg_{13}(x,y,z,v)+\alpha v^2x+\beta z^3x+g_{20}(x,y,z,t,v)=0$}\\
    \\
    Thm. \ref{thm:cod2caseI}& $\frac{1}{11}(1,3,8)$& $(-1,-3,1,2)$ & $(7,1,3,1)$ &$cA_9  \in Z_{9,10} \subset \mathbb{P}(1,1,2,3,5,8)$\\
		Thm. \ref{thm:cod4}&$\frac{1}{7}(1,1,6)$& $10\times (-1,-1,1,1)$  && $\frac{1}{14}(1,3,11) \in Z \subset \mathbb{P}^7(1,3,5,6,7,8,11,14)$\\
		Prop. \ref{cor:exc111} &$\frac{1}{3}(1,1,2)$&  & \\
		\specialrule{.1em}{0.05em}{.1em} 
		    {}&&\\
		\specialrule{.1em}{0.05em}{.1em} 
		  \multicolumn{3}{l}{\textbf{Family 112}: $X_{6,6} \subset \mathbb{P}(1,1,2,3,3,5), \quad \iota_X=3$}\\
		\\    
    \multicolumn{3}{l}{$wx+f_{6}(x,y,z,t,v)=0$}\\
		\multicolumn{3}{l}{$wy+g_{6}(x,y,z,t,v)=0$}\\
    \\
    		Prop. \ref{prop:elemcb}&$\frac{1}{5}(1,1,4)$&  & &$Y/\mathbb{P}(1,1,2)$ \\
		\specialrule{.1em}{0.05em}{.1em} 
		{}&& \\
		\specialrule{.1em}{0.05em}{.1em} 
\multicolumn{3}{l}{\textbf{Family 113}: $X_{6,6} \subset \mathbb{P}(1,2,2,3,3,4), \quad \iota_X=3 $}\\
		\\    
    \multicolumn{3}{l}{$wy+tv+ z^3+f_{6}(x,y,z,t,v,w)=0$}\\
		\multicolumn{3}{l}{$wz+ y^2z+yg_4(x,z,t,v)+ g_{6}(x,z,t,v,w)=0$}\\
    \\
Thm. \ref{thm:cod2curve}&$\frac{1}{4}(1,1,3)$&  & &$\mathbb{P}^1 \subset Z_{4,4} \subset \mathbb{P}(1,1,1,2,2,3)$\\
Prop. \ref{prop:113_119}& $\frac{1}{2}(1,1,1)$& $(-4,-1,-1,1,2;-2)$ & &$dP_{4}$\\
		\specialrule{.1em}{0.05em}{.1em} 
		{}&&\\
		\specialrule{.1em}{0.05em}{.1em} 
\multicolumn{3}{l}{\textbf{Family 114}: $X_{6,9} \subset \mathbb{P}(1,2,3,3,4,5), \quad \iota_X=3$}\\
		\\    
    \multicolumn{3}{l}{$wx+vy+f_{6}(x,y,z,t,v)=0$}\\
		\multicolumn{3}{l}{$wv+(1-\alpha)y^4x+ (1-\alpha)y^3t+\alpha y^2w+ yg_7(x,z,t,v,w)+g_{9}(x,y,z,t,v,w)=0$}\\
    \\
    Thm. \ref{thm:delPezzoMain}&$\frac{1}{5}(1,1,4)$& $(-4,-1,-1,1,2;-2)$ & &$dP_3$ \\
		Thm. \ref{thm:cod2caseI}&$\frac{1}{4}(1,1,3)$& $(-5,-1,-1,1,3;-2)$ & $(2,1,1,1)$ &$cA_2 \in Z_{4,6} \subset \mathbb{P}(1,1,2,2,3,3)$\\
		Thm. \ref{thm:cod2ExcII}&$\frac{1}{2}(1,1,1),\, \alpha \not = 0$& $(-3,-1,1,1)$ & $(3,2,2,1)$ &$cD_4 \in Z_{2,3} \subset \mathbb{P}^5$\\
		\specialrule{.1em}{0.05em}{.1em} 
		\\
\specialrule{.1em}{0.05em}{.1em} 
		\multicolumn{3}{l}{\textbf{Family 115}: $X_{12,15} \subset \mathbb{P}(1,3,4,5,6,11),\quad \iota_X=3, \quad \frac{4}{\iota_X^2A^3}=\frac{88}{9}$}\\
		\\    
    \multicolumn{3}{l}{$wx+v^2+z^3+f_{12}(x,y,z,t,v)=0$}\\
		\multicolumn{3}{l}{$wz+t^3 +g_{15}(x,y,z,t,v)=0$}\\
    \\\
   	Thm. \ref{thm:hypIIfull}&$\frac{1}{11}(1,2,9)$& $(-1,-2,1,1)$ & $(5,9,2,4,1)$ &$cD_6  \in Z_{5} \subset \mathbb{P}(1,1,1,1,2)$\\
		Prop. \ref{prop:weak115}&$\frac{1}{2}(1,1,1)$&  & \\
		\specialrule{.1em}{0.05em}{.1em} 
		{}&&\\
		\specialrule{.1em}{0.05em}{.1em} 
		\multicolumn{3}{l}{\textbf{Family 116}: $X_{9,12} \subset \mathbb{P}(2,3,3,4,5,7), \quad \iota_X=3$}\\
		\\    
    \multicolumn{3}{l}{$wx+vt +f_{9}(x,y,z,t,v)=0$}\\
		\multicolumn{3}{l}{$wv+t^3+ x^4t+ x^3g_{6}(y,z)+x^2g_{8}(y,z,t,v)+xg_{10}(y,z,t,v,w)+g_{12}(x,y,z,t,v)=0$}\\
    \\
   	Thm. \ref{thm:delPezzoMain}&$\frac{1}{7}(1,1,6)$& $(-5,-1,-1,1,2;-3)$ & &$dP_2$\\
		Prop. \ref{prop:fake1}&$\frac{1}{5}(1,1,4)$& $(-7,-1,-1,1,4;-3)$ & &$\mathbf{p_t} \in Z_{3,4} \subset \mathbb{P}(1,1,1,2,2,3)/\boldsymbol{\mu_2}$\\
		Pp. 53 &		$\frac{1}{2}(1,1,1)$& $(-5,-1,-1,1,2;-3)$  & $(7,1,1,1)$ & $cA_8 \in Z_6 \subset \mathbb{P}(1,1,1,2,2)$ \\
		\specialrule{.1em}{0.05em}{.1em} 
	{}&&\\
		\specialrule{.1em}{0.05em}{.1em} 
			\multicolumn{3}{l}{\textbf{Family 117}: $X_{12,15} \subset \mathbb{P}(3,3,4,5,7,8),\quad \iota_X=3, \quad \frac{4}{\iota_X^2A^3}=\frac{224}{9}$}\\
		\\    
    \multicolumn{3}{l}{$wz+vt+f_{12}(x,y,z,t,v)=0$}\\
		\multicolumn{3}{l}{$wv+t^3+(1-\alpha)z^3x+\alpha z^2v+zg_{11}(x,y,t,w)+ g_{15}(x,y,z,t,v)=0$}\\
    \\
   	Thm. \ref{thm:cod2caseI}&$\frac{1}{8}(1,1,7)$& $(-7,-1,-1,1,3;-4)$ & $(4,1,1,1)$ & $cA_4 \in Z_{8,10} \subset \mathbb{P}(1,2,2,3,5,7)$\\
		Prop. \ref{prop:fake1}&$\frac{1}{7}(1,1,6)$& Sec. \ref{subsec:fake}  & &$\mathbf{p_t} \in Z_{3,4} \subset \mathbb{P}(1,1,1,2,2,3)/\boldsymbol{\mu_2}$\\
		 Thm. \ref{thm:cod2ExcII}&$\frac{1}{4}(1,1,3),\, \alpha=1$& $(-8,-1,-1,1,3;-5)$ & $(3,1,1,1)$ &$cA_3 \in Z_{8,10} \subset \mathbb{P}(1,2,2,3,5,7)$\\
		Prop. \ref{prop:bi}&$\frac{1}{4}(1,1,3),\, \alpha =0$&  & & BI \\
		\specialrule{.1em}{0.05em}{.1em}
		\\
		\specialrule{.1em}{0.05em}{.1em}
		\multicolumn{3}{l}{\textbf{Family 118}: $X_{6,8} \subset \mathbb{P}(1,2,3,3,4,5), \quad \iota_X=4$}\\
		\\    
    \multicolumn{3}{l}{$wx+tz+vy+f_{6}(x,y,z,t,v)=0$}\\
		\multicolumn{3}{l}{$w(t+\alpha z)+v^2+(1-\alpha)z^2y+g_{8}(x,y,z,t,v)=0$}\\
    \\
   Pp. 55& $\frac{1}{5}(1,2,3)$& $2 \times (-3,-1,1,2)$ & $(2,3,1,1)$ &$cA_1 \in Z_{6} \subset \mathbb{P}(1,1,2,2,3)$  \\
		Prop. \ref{prop:118_119CB}&$ \frac{1}{3}(1,1,2)$& $(-5,-1,2,3)$ & &$Y/\mathbb{P}(1,2,3)$ \\
		\\
		\multicolumn{4}{l}{Blowing up the point $\mathbf{p_v}$ also initiates a Sarkisov link to $Y/\mathbb{P}(1,2,3)$ regardless of $\alpha$.}\\
		\multicolumn{4}{l}{However, when $\alpha=0$, the toric flip $(-5,-1,2,3)$ degenerates to the hypersurface flip $(-5,-2,-1,2,3;-2)$.}\\
				\specialrule{.1em}{0.05em}{.1em}
		\\
			\specialrule{.1em}{0.05em}{.1em}
		\multicolumn{3}{l}{\textbf{Family 119}: $X_{8,10} \subset \mathbb{P}(1,2,3,4,5,7), \quad \iota_X=4$}\\
		\\    
    \multicolumn{3}{l}{$wx+zv+t^2+f_{8}(x,y,z,t,v)=0$}\\
		\multicolumn{3}{l}{$wz +v^2+g_{10}(x,y,z,t,v)=0$}\\
    \\
   	Prop. \ref{prop:118_119CB}&$\frac{1}{7}(1,3,4)$&  & &$Y/\mathbb{P}(1,2,3)$\\
		Prop. \ref{prop:113_119}&	$\frac{1}{3}(1,1,1)$& $(-7,-3,-2,1,2;-6)$ & &$dP_{4}$ \\
				\specialrule{.1em}{0.05em}{.1em}
		{}&&\\
				\specialrule{.1em}{0.05em}{.1em}
				\multicolumn{3}{l}{\textbf{Family 120}: $X_{8,12} \subset \mathbb{P}(1,3,4,4,5,7), \quad \iota_X=4$}\\
		\\    
    \multicolumn{3}{l}{$wx+vy+f_{8}(x,y,z,t,v)=0$}\\
		\multicolumn{3}{l}{$wv+g_{12}(x,y,z,t,v,w)=0$}\\
    \\
	Thm. \ref{thm:delPezzoMain}&$\frac{1}{7}(1,1,6)$&  $(-5,-1,-1,1,3;-2)$ & &$dP_3$\\
	Thm. \ref{thm:cod2caseI}&$\frac{1}{5}(1,1,4)$&  $(-7,-1,-1,1,5;-2)$ & &$\frac{1}{2}(1,1,1) \in \mathbb{Z}_{6,9} \subset \mathbb{P}(1,2,3,3,4,5)$ \\
		\specialrule{.1em}{0.05em}{.1em}
		\\
		\specialrule{.1em}{0.05em}{.1em}
				\multicolumn{3}{l}{\textbf{Family 121}: $X_{10,12} \subset \mathbb{P}(1,3,4,5,6,7), \quad \iota_X=4$}\\
		\\    
    \multicolumn{3}{l}{$wy+t^2+zv+f_{10}(x,y,z,t,v)=0$}\\
		\multicolumn{3}{l}{$wt+v^2+z^3+ y^2v+yg_{9}(x,z,t,v,w)+g_{12}(x,y,z,t,v)=0$}\\
    \\
   	Prop. \ref{prop:121_122}&$\frac{1}{7}(1,2,5)$& $ (-5,-1,1,3)$ & &$dP_3$\\
		Thm. \ref{thm:cod2ExcII}&$2\times \frac{1}{3}(1,1,2)$& $(-3,-1,1,1)$ & $(7,3,2,1)$ &$cA_2 \in \mathbb{Z}_{2,3} \subset \mathbb{P}^5$ \\
		\specialrule{.1em}{0.05em}{.1em}
		{}&&\\
		\specialrule{.1em}{0.05em}{.1em}
				\multicolumn{3}{l}{\textbf{Family 122}: $X_{10,12} \subset \mathbb{P}(2,3,4,5,5,7), \quad \iota_X=4$}\\
		\\    
    \multicolumn{3}{l}{$wy +vt+f_{10}(x,z)+\widetilde{f_{10}}(x,y,z,t,v)=0$}\\
		\multicolumn{3}{l}{$w(v+\alpha t) +(1-\alpha)t^2x+z^3+g_{12}(x,y,z,t,v)=0$}\\
    \\
   	Prop. \ref{prop:121_122}&$\frac{1}{7}(1,3,4)$& $ 2 \times (-5,-1,2,3)$ && $dP_3$  \\
		Thm. \ref{thm:cod2ExcII}&$\frac{1}{5}(1,2,3)$& $(-7,-1,4,3), \, (-2,-1,1,1)$ & $(5,3,4,1)$ &$cA_2 \in \mathbb{Z}_{4,6} \subset \mathbb{P}(1,1,2,2,3,3)$ \\
		\\
		\multicolumn{4}{l}{When $\alpha = 0$, the toric flip $(-5,-1,2,3)$ over $p_t$ degenerates to the hypersurface flip $(-5,-3,-1,2,3;-3)$.}\\
		\specialrule{.1em}{0.05em}{.1em}
		{}&&\\
		\specialrule{.1em}{0.05em}{.1em}
				\multicolumn{3}{l}{\textbf{Family 123}: $X_{12,14} \subset \mathbb{P}(2,3,4,5,7,9), \quad \iota_X=4$}\\
		\\    
    \multicolumn{3}{l}{$wy+ vt +f_{12}(x,y,z,t,v)=0$}\\
		\multicolumn{3}{l}{$wt+v^2+ y^3t + y^2g_8(x,z)+yg_{11}(x,z,t,v,w)+ g_{14}(x,y,z,t,v)=0$}\\
    \\
   	Thm. \ref{thm:cod2ExcII}&$\frac{1}{9}(1,4,5)$& $ (-3,-1,-1,1,1;-2)$ & $(5,3,4,1)$ &$cA_2 \in \mathbb{Z}_{6,6} \subset \mathbb{P}(1,1,2,2,3,5)$ \\
		Thm. \ref{thm:cod2curve}&$\frac{1}{5}(1,2,3)$&  & &$\mathbb{P}^1 \subset \mathbb{Z}_{6,6} \subset \mathbb{P}(1,2,2,3,3,3)$ \\
		Prop. \ref{prop:fake2}&$\frac{1}{3}(1,1,2)$& $(-9,-4,-1,1,4;-8)$ & & $\mathbf{p_t} \in \mathbb{Z}_{3,4} \subset \mathbb{P}(1,1,1,1,2,2)/\boldsymbol{\mu_2}$ \\
		\\
	\multicolumn{4}{l}{Blowing up the point $\frac{1}{5}(1,1,4) \in Z_{6,6}$ initiates a link to a del Pezzo fibration of degree 3. See lemma \ref{lem:compdp}.} \\
		\specialrule{.1em}{0.05em}{.1em}
		{}&&\\
		\specialrule{.1em}{0.05em}{.1em}
				\multicolumn{3}{l}{\textbf{Family 124}: $X_{18,20} \subset \mathbb{P}(4,5,6,7,9,11), \quad \iota_X=4$}\\
		\\    
    \multicolumn{3}{l}{$wt+v^2+z^3+f_{18}(x,y,z,t,v)=0$}\\
		\multicolumn{3}{l}{$wv+t^2z+tg_{13}(x,y,z,v)+g_{20}(x,y,z,t,v)=0$}\\
    \\
		Thm. \ref{thm:124toHyp}&$\frac{1}{11}(1,4,7)$& $(-4,-1,-1,1,2;-2)$ & $(7,2,3,1)$&$ cA_4 \in \mathbb{Z}_{15} \subset \mathbb{P}(1,2,3,5,7)$ \\
   	Thm. \ref{thm:124toHyp}&$\frac{1}{7}(1,3,4)$& $(-3,-1,-1,1,1;-2)$ &$(11,3,4,1)$ &$ cA_4 \in \mathbb{Z}_{10} \subset \mathbb{P}(1,1,2,3,5)$ \\
		Cor. \ref{cor:124}& $\frac{1}{3}(1,1,2)$&  &  \\
	\specialrule{.1em}{0.05em}{.1em}
    {}&&\\
		\specialrule{.1em}{0.05em}{.1em}
				\multicolumn{3}{l}{\textbf{Family 125}: $X_{10,15} \subset \mathbb{P}(2,3,5,5,7,8), \quad \iota_X=5$}\\
		\\    
    \multicolumn{3}{l}{$wx+vy+tz+f_{10}(x,y,z,t,v)=0$}\\
		\multicolumn{3}{l}{$wv+\alpha x^6y + \beta x^5 t+\gamma x^4v+x^3y^3+x^2g_{11}(y,z,t,w)+xg_{13}(y,z,t,v,w)+g_{15}(x,y,z,t,v)=0$}\\
    \\
		Thm. \ref{thm:cod2caseI}&$\frac{1}{8}(1,1,7)$& $(-7,-1,-1,1,5;-2)$ & $(2,1,1,1)$ & $cA_2 \in \mathbb{Z}_{8,12} \subset \mathbb{P}(1,3,4,4,5,7)$ \\
   	Prop. \ref{prop:fake1}&$\frac{1}{7}(1,1,6)$& $(-8,-1,-1,1,6;-2)$ & &$\mathbf{p_y} \in \mathbb{Z}_{4,6} \subset \mathbb{P}(1,1,2,2,3,3)/\boldsymbol{\mu_2}$ \\
		Thm. \ref{thm:cod2ExcII}&$\frac{1}{2}(1,1,1),\, \gamma \not = 0$& $(-8,-1,-1,1,5;-3)$ & $(1,1,1,1)$ &$cA_1 \in \mathbb{Z}_{8,12} \subset \mathbb{P}(1,3,4,4,5,7)$ \\
	Prop. \ref{prop:bi}&	 $\frac{1}{2}(1,1,1),\, \gamma = 0,\, \beta \not = 0$ &  &  &BI \\
		Prop. \ref{prop:fake2}&$\frac{1}{2}(1,1,1),\, \gamma = \beta  = 0,\, \alpha \not = 0$  & $(-7,-1,1,1,9;2)$ & &$\mathbf{p_w} \in \mathbb{Z}_{4,6} \subset \mathbb{P}(1,1,2,2,3,3)/\boldsymbol{\mu_3}^*$ \\
		\specialrule{.2em}{0.1em}{.1em}
		\label{tab:big}
		\end{longtable}
\par}

\end{landscape}

\newgeometry{footskip=1cm,margin=2.5cm}

\bibliography{aaBibliograph}{}
\bibliographystyle{abbrv}

\end{document}